\documentclass[reqno,11pt]{amsart}
\usepackage[T1]{fontenc}
\usepackage{geometry}
\makeatletter
\let\@wraptoccontribs\wraptoccontribs
\makeatother

\title[SUSY Hilbert complexes and Witten deformation on stratified pseudomanifolds]{Supersymmetric Hilbert complexes and Witten deformation on stratified pseudomanifolds} 
\date{\today}
\author{Gayana Jayasinghe} 
\email{gayana.jayasinghe@uniroma1.it}

\author{Manousos Maridakis}
\email{mmaridakis@math.auth.gr}

\author{Hadrian Quan} 
\email{hquan1@ucsc.edu}

\author{Xinran Yu}
\email{xinranyu@math.pku.edu.cn} 

\contrib[Appendix B jointly written with]{Alex R. Taylor}
\email{alexrt.math@gmail.com}

\makeatletter
\renewcommand{\tocsection}[3]{%
  \indentlabel{\@ifnotempty{#2}{\bfseries\ignorespaces#1 #2\quad}}\bfseries#3}
  
\renewcommand{\tocsubsection}[3]{%
  \indentlabel{\@ifnotempty{#2}{\ignorespaces#1 #2\quad}}#3}
\let\tocsubsubsection\tocsubsection

    \renewcommand{\tocsubsubsection}[3]{%
  \indentlabel{\@ifnotempty{#2}{\ignorespaces#1 #2\quad}}#3}

\newcommand\@dotsep{4.5}
\def\@tocline#1#2#3#4#5#6#7{\relax
  \ifnum #1>\c@tocdepth 
  \else
    \par \addpenalty\@secpenalty\addvspace{#2}%
    \begingroup \hyphenpenalty\@M
    \@ifempty{#4}{%
      \@tempdima\csname r@tocindent\number#1\endcsname\relax
    }{%
      \@tempdima#4\relax
    }%
    \parindent\z@ \leftskip#3\relax \advance\leftskip\@tempdima\relax
    \rightskip\@pnumwidth plus1em \parfillskip-\@pnumwidth
    #5\leavevmode\hskip-\@tempdima{#6}\nobreak
    \leaders\hbox{$\m@th\mkern \@dotsep mu\hbox{.}\mkern \@dotsep mu$}\hfill
    \nobreak
    \hbox to\@pnumwidth{\@tocpagenum{\ifnum#1=1\bfseries\fi#7}}\par
    \nobreak
    \endgroup
  \fi}

\AtBeginDocument{%
\expandafter\renewcommand\csname r@tocindent0\endcsname{0pt}
}
\def\l@subsection{\@tocline{2}{0pt}{2.5pc}{5pc}{}}
\def\l@subsubsection{\@tocline{2}{0pt}{4.5pc}{5pc}{}}

\makeatother

\makeatletter

\renewcommand\subsection{\@secnumfont}{\bfseries}%
\renewcommand\subsection{\@startsection{subsection}{3}
  \z@{.5\linespacing\@plus.7\linespacing}{.5em}%
  {\normalfont\bfseries}}

\renewcommand\subsubsection{\@secnumfont}{\bfseries}%
\renewcommand\subsubsection{\@startsection{subsubsection}{3}
  \z@{.5\linespacing\@plus.7\linespacing}{.5em}%
  {\normalfont\bfseries}}
  
\makeatother

\usepackage{amsmath}
\usepackage{amsfonts, mathrsfs,amssymb} 
\usepackage{hyperref}
\usepackage{slashed} 
\usepackage{mathtools}
\usepackage{tikz-cd}
\allowdisplaybreaks

\usepackage{enumerate}

\newcommand\Ker{\operatorname{Ker}}

\numberwithin{equation}{section} 

\newtheorem{theorem}{Theorem}[section]
\newtheorem{proposition}[theorem]{Proposition}
\newtheorem{lemma}[theorem]{Lemma}
\newtheorem{corollary}[theorem]{Corollary}
\newtheorem{conjecture}[theorem]{Conjecture}
\newtheorem{remark}[theorem]{Remark}
\newtheorem{example}[theorem]{Example}
\newtheorem{definition}[theorem]{Definition}
\newtheorem{question}[theorem]{Question}
\newtheorem{question/conjecture}[theorem]{Question/Conjecture}
\newtheorem{assumption}[theorem]{Assumption}

\usepackage{adjustbox}
\usepackage{graphicx}
\graphicspath{ {images/} }
\usepackage{caption}
\usepackage{subcaption}
\usepackage{float}
\usepackage{comment}
\usepackage{appendix}

\usepackage{multicol}

\date{}
\begin{document}

\maketitle

\begin{abstract}
We present a framework to study topics related to index theory and Dirac operators with $\mathcal{N}=1,2$ supersymmetry (SUSY), extending it to stratified pseudomanifolds with wedge metrics (iterated conic metrics). We use this to construct Witten deformed subcomplexes of SUSY Hilbert complexes that categorify Morse polynomials and equivariant index formulae for a large class of operators.

This allows us to prove equivariant index theorems for spin-c Dirac operators with Hamiltonian actions, and Holomorphic Morse inequalities in the K\"ahler case, in particular proving a new type of equivariant Riemann-Roch (ERR) theorem on singular spaces. We also improve on stratified Morse-Bott inequalities in the de Rham case, and review other applications including for the Kervaire semi-characteristic and the equivariant signature. We compute local cohomology groups for the Dolbeault complexes at both smooth and singular non-isolated fixed point sets which categorify the ERR numbers. 

In an Appendix joint with Alex R. Taylor, we give a broad overview of the relevance of SUSY to index theory on manifolds with boundaries and singularities, mainly focusing on the new insights and results related to $\mathcal{N}=2$ SUSY that the authors have uncovered recently.
\end{abstract}

\tableofcontents

\section{Introduction}

Supersymmetry is at the heart of Index theory. Starting from the definition of the Fredholm indices of Dirac-type operators as supertraces, to relating them to supertraces of heat kernels via spectral methods (such as the McKean singer cancellations, also known as SUSY cancellations) and more differential geometric constructions (such as Getzler's derivation of the local index formula by passing to pseudo differential calculi on supermanifolds), the role played by supersymmetry (SUSY) cannot be understated (see Appendix \ref{Appendix_SUSY_and_conjectures}). Even though the founders of index theory, notably Michael Atiyah and Isadore Singer (see, e.g., \cite{Singer1996,Yaubookfoundersindex}) embraced this perspective, it seems to remain a mystery in many communities of mathematicians who work on different problems in the area. One of the main drawbacks of this is the lack of study of the many roles played by the extra SUSY's when they are present. Thus the results that have been formulated in these instances usually involve unnecessarily complicated objects, some of which do not actually contribute to the formulas due to cancellations corresponding to these extra SUSYs. Many proofs for complexes with only one SUSY admit important technical simplifications when there are more SUSY's, while there are outright simpler proof schema that aid both understanding and generalization, as one extends index theory to singular, non-compact and foliated settings.

Two of the most important such settings is the study of the Hodge-Dirac operator (also called the Gauss-Bonnet operator) and the Dolbeault-Dirac operator on complex manifolds. Even some simple questions such as discrepancies in Riemann-Roch and equivariant Riemann-Roch formulas formulated in algebraic and analytic language that have eluded understanding in the past, for which SUSY gives enlightening answers (see Appendix \ref{Appendix_SUSY_and_conjectures}).

The work of Witten in \cite{witten1982constraints,witten1982supersymmetry} laid out the beginnings of a framework (to use semi-classical analysis in SUSY Quantum mechanical systems/ Quantum field theories), popularized as \textit{Witten deformation}, which could be used to fruitfully study and expand Morse theory and Hodge theory in contexts where it was not known earlier, including the holomorphic Morse inequalities \cite{witten1984holomorphic} and in topological quantum field theories (see, e.g., \cite{WittenTopologicalQFT88}, which discusses these perspectives in Floer theory and Donaldson theory). There has been a lot of progress in Witten deformation since then, and the techniques have been used to study Kervair\'e semi-characteristics, Hopf index theorems and Morse inequalities without non-degeneracy, geometric quantization on symplectic manifolds, analytic torsion with applications of Cheeger-M\"uller theorems and understanding Quillen metrics on complex manifolds to name some, some of which have been and we refer the reader to Subsection \ref{subsection_introduction_History} including some extensions to singular spaces.

While these results already show how SUSY plays a key role in equivariant index theory, recent studies by the authors of the present article (and their friends) have shown how they even say much more on index theory, on spaces with boundaries, singularities and non-compact ends. While some of these are implicit in \cite{jayasinghe2025supersymmetricharmonicoscillatorssingular,jayasinghe2026spectralasymmetrysupersymmetryequivariant}, we make this an article where we formalize these notions in a suitable asbtract setting, and construct generalizations of Witten subcomplexes on stratified pseudomanifolds. These subcomplexes categorify both Morse polynomials (corresponding to Morse inequalities) and Lefschetz numbers/equivariant indices, in appropriate settings.

Stratified pseudomanifolds arise in many natural contexts of study, and some common examples are algebraic varieties, quotient spaces of group actions and various moduli spaces, and there have been studies of Witten deformation for de Rham and Spin-c Dirac operators in some contexts, including by some of the authors (see Subsection \ref{subsection_introduction_History}).

Our \textit{initial goal} was to generalize the holomorphic Morse inequalities for non-isolated fixed point sets in \cite{wu1998equivariant} to some singular settings, which we have achieved and present as follows.

\begin{theorem}[\textbf{Holomorphic Morse inequalities}]
\label{Theorem_introduction_holomorphic_Morse_inequal_not_abstract}
Consider a resolution $X$ of a stratified pseudomanifold $\widehat{X}$ equipped with an asymptotically $\delta>1/2$ K\"ahler wedge metric and a K\"ahler Hamiltonian $h$, whose action lifts to a holomorphic Hermitian bundle $E$. Let $A=|dh|^2$, and denote by $Z(A)$ the set of connected components of its zero set. Assume that each $F\in Z(A)$ is admissible. Consider the twisted Dolbeault complex
$\mathcal{P}_{W}(X)=(L^2\Omega^{0,\bullet}(X;E), \mathcal{D}_{W}(P), P=\overline{\partial}_E)$, where $W$ are ideal boundary conditions for a closed extension of $P$, we have 
\begin{equation}
\label{equation_blah_440}
    \sum_{F \in Z(A)}  \sum_{q=0}^n b^q \mathrm{dim}\Big(\mathcal{H}^{q}(\mathcal{P}^{\mu}_{W,B}(\mathcal{N}_{F,\varepsilon}) \Big)  = \sum_{q=0}^n b^q \mathrm{dim}\Big(\mathcal{H}^{q}(\mathcal{P}^{\mu}_{W}(X)) \Big) + (1+b) \sum_{q=0}^{n-1} E^{\mu}_q b^q
\end{equation}
where $Q^{\mu}_q$ are non-negative integers.
Equivalently we have a power series
\begin{equation}
\label{equation_blah_441}
     \sum_{F \in Z(A)}  \sum_{q=0}^n b^q \operatorname{Tr} \Big(T_{s,\theta}|_{\mathcal{H}^{q}(\mathcal{P}_{W,B}(\mathcal{N}_{F,\varepsilon}))} \Big)  = \sum_{q=0}^n b^q \operatorname{Tr} \Big(T_{s,\theta}|_{\mathcal{H}^{q}(\mathcal{P}_{W}(X))}\Big) + (1+b) \sum_{q=0}^{n-1} E_q b^q 
\end{equation}
where the $E_q$ are power series in the variable $\lambda=se^{i\theta}$ 
converging for $|\lambda|=s<1$. The coefficients of this series are non-negative integers.

\textbf{In particular setting $b=-1$ we get an equivariant Riemann-Roch theorem}.
\end{theorem}

This result is a straightforward Corollary of the existence of Witten subcomplexes presented in Subsection \ref{Theorem_main_Witten_complex_Kahler_case}.
The objects in the statement above, such as the $\mu$ equivariant complexes $\mathcal{P}^{\mu}_W$, the normal cohomology complexes  for connected components of the critical point set of $h$, boundary conditions $B$ and cohomology groups $\mathcal{H}^q$ for local and global complexes will be introduced in a generality which holds more broadly. This generalizes the holomorphic Morse inequalities of \cite{wu1998equivariant,wu2003instanton} from the smooth setting to singular spaces with K\"ahler wedge metrics. In the case of isolated fixed points and non-isolated singularities, this was proven in \cite{jayasinghe2024holomorphic} giving Morse inequalities that generalized the equivariant Riemann-Roch numbers of Baum-Fulton-Quart \cite{baum1979lefschetz,Baumformula81} where the geometric settings overlap. \textbf{The formulation of Theorem \ref{Theorem_introduction_holomorphic_Morse_inequal_not_abstract} for non-isolated singular spaces seems new even for singular algebraic varieties, even in the Riemann-Roch situation (where $b=-1$).}

In pursuing this \textit{initial goal}, we realized that the tools we were building were applicable in greater generality. We also realized SUSY simplifies many other invariants such as $\xi_T$ and eta invariants on manifolds with boundaries and singularities (see Appendix \ref{Appendix_SUSY_and_conjectures}), building on related work by some of the authors and their friends. Thus, 
\textbf{one of the main contributions of this article consists of introducing a formal framework for Witten deformation for SUSY Hilbert complexes associated to stratified pseudomanifolds, and establishing the existence of Witten complexes for them that categorify Morse inequalities and equivariant index theorems}. We now present this.

\subsection{SUSY Hilbert complexes and Witten sub-complexes}
\label{subsection_intro_SUSY_Hilbert_and_Witten}

The following formalize a setting introduced by Witten, for instance in \cite{witten1982constraints,witten1982supersymmetry}, when discussing Witten deformation.

\begin{definition}[\textbf{$\mathcal{N}=N$ SUSY Hilbert complex}]
\label{definition_N_Hilbert_complex}
    An abstract $\mathcal{N}=N$ \textbf{SUSY Hilbert complex} is a tuple $\mathcal{Q}= \left(H, \{(\mathcal{D}(Q_j), Q_{j})\}_{j=1}^N,\gamma \right)$ where $H$ is a Hilbert space equipped with an involution $\gamma$ (called the chirality/grading operator) that induces a $\mathbb{Z}_2$ grading on $H$. Here each $Q_j$, called a \textbf{SUSY operator} is a densely defined self adjoint operator with domain $\mathcal{D}(Q_j)$, which satisfies the conditions
    \begin{equation}
    \label{equation_SUSY_complex_condition_1}
        \{Q_i,Q_j\}=0 \text{   if $i \neq j$, and    } \{Q_i,\gamma\}=0  \text{    for all i=1,...,N}
    \end{equation}
    and 
    \begin{equation}
    \label{equation_SUSY_complex_condition_2}
        Q_1^2=Q_2^2=...=Q_N^2=\Delta
    \end{equation}
    where $\Delta$ is called the \textbf{Hamiltonian operator}, with the induced domain
    \begin{equation}
        \mathcal{D}(\Delta)=\{ u \in \mathcal{D}(Q_j)| Q_ju \in \mathcal{D}(Q_j) \}
    \end{equation}
    for any choice of $j$. 
    We define the \textbf{cohomology group} of $\mathcal{Q}$ to be 
    \begin{equation}
    \mathcal{H}(\mathcal{Q}):= \{u \in \mathcal{D}(\Delta_{\varepsilon})| \Delta_{\varepsilon}u=0\}
    \end{equation}
    and we say the complex is \textbf{Fredholm} if $\mathcal{H}(\mathcal{Q})$ is finite dimensional, and that the complex has \textbf{discrete spectrum} if $\mathcal{D}(\Delta_{\varepsilon})$ has discrete spectrum.
    
    If in addition the Hilbert space is the space of sections of a resolved stratified pseudomanifold $X$, we say the \textbf{$\mathcal{N}=N$ SUSY Hilbert complex is associated to $X$}, and denote it $\mathcal{Q}(X)$. 
    If the chirality operator $\gamma$ can be refined to one inducing a $\mathbb{Z}$ grading of $H$, we call the Hilbert complex \textbf{$\mathbb{Z}$ graded}.
\end{definition}

In \cite{witten1982supersymmetry}, Witten explains that such \textit{SUSY algebras} must be generalized further when studying relativistic quantum field theory, since the Lorentz transformations generate spatial translations (corresponding to momentum operators) that are not compatible with the operator identities in \eqref{equation_SUSY_complex_condition_1} and \eqref{equation_SUSY_complex_condition_2}. We do not study these here, and study the aspects related to Hilbert complexes on singular spaces which are more nuanced on stratified pseudomanifolds.

\begin{remark}[Resolved stratified pseudomanifolds and wedge structures]
\label{Remark_resolved_stratified_pseudomanifolds_and_wedge_metrics}    

Given a stratified pseudomanifold $\widehat{X}$, we can associate to them manifolds with corners with iterated fibration structures $X$ (also known as manifolds with fibered boundaries, and as fibered manifolds) which resolve them in a suitable way for analysis. There are functorial equivalences between two categories which contain these as objects and we shall study this in more detail in Section \ref{section_stratified_and_wedge}. For brevity, we shall call such a resolved space as \textbf{resolutions $X$ of a stratified space}.

Since the open dense regular part of a stratified pseudomanifold is dense in its resolved space, and is of full measure with respect to the volume measure corresponding to a wedge metric, we will study the $L^2$ cohomology groups of such complexes on the stratified spaces in terms of the groups on their resolutions.   

Part of the philosophy in geometric microlocal analysis propogated by Melrose is to treat the structures underlying singular metrics. Similar to one studies smooth structures by studying invariants defined using smooth Riemannian metrics, one can consider a rescalings of tangent bundles (called adapted tangent bundles in \cite{jayasinghe2025supersymmetricharmonicoscillatorssingular}) adapted to singular metrics, such that said singular metrics are non-degenerate bundle metrics on these rescaled bundles. In the case of incomplete iterated conic (edge) metrics, also called wedge metrics, these are the wedge tangent bundles. Then all other geometric structures such as symplectic forms, almost complex structures and Clifford modules get their wedge counterparts.   

\textbf{In this article, unless mentioned otherwise, we shall always consider resolved stratified pseudomanifolds equipped with wedge structures.}
\end{remark}

We refer the reader to Section \ref{section_stratified_and_wedge} for a more comprehensive explanation of the resolutions, wedge metrics and operators, and to \cite{jayasinghe2025supersymmetricharmonicoscillatorssingular} for some insights on the general philosophy of studying structures related to rescaled tangent bundles. Smooth manifolds with smooth Riemmanian metrics are wedge metrics and this may aid the reader not used to the subtleties of pseudomanifolds to digest the content here in a first reading. Much of this formalism is built to generalize easily for other singular metrics.

We shall drop the symbol $\gamma$ from the notation $\mathcal{Q}$ where it is clear by context. 
For an $\mathcal{N}=N$ SUSY Hilbert complex associated to a stratified pseudomanifold $X$, the operators $Q_j$ acting on the smooth sections on the interior smooth strata of the space do not admit unique self-adjoint extensions, and the choice of extension $\mathcal{D}(Q_j)$ induces a choice of extension on the square $Q_j$ for each $j$. Then the equality in \eqref{equation_SUSY_complex_condition_2} implies that the domains of each operator $Q_j^2$ coincides as well. We shall study this in more detail in Section \ref{sec:domains}.

Since $\Delta$ is a square of at least one self-adjoint operator, it is clearly non-negative. 
When studying non-negative second order operators, it is standard to see if there is a metric $g$ on the space such that it is a Laplace-type operator, that is, $\sigma_2(\Delta)=g$. This forces any self-adjoint first order operators that square to it to be Dirac-type, $\sigma_1(Q)=cl$ for an appropriate Clifford action on the cotangent bundle. Then the geometry aids the study of the operator. This motivates certain parts of Definition \ref{definition_family_of_N_2_SUSY_complexes} below.

\begin{definition}[\textbf{Family of $\mathcal{N}=N$ SUSY Hilbert complexes}]
\label{definition_family_of_N_2_SUSY_complexes}
    A family of $\mathcal{N}=N$ SUSY complexes is a family of tuples $\mathcal{Q} = \left( H, \{ ( \mathcal{D}(Q_{j,\varepsilon}), Q_{j,\varepsilon} ) \}_{j=1}^N ,\gamma \right)$ where for each fixed $\varepsilon \in \mathbb{R}^+_0$, the tuple defines a $\mathcal{N}=N$ SUSY Hilbert complex.
    We refer to $\Delta_{\varepsilon}:=Q_{j,\varepsilon}^2$ as \textbf{the Hamiltonian} of the family.
    
    Families of $\mathcal{N}=N$ SUSY Hilbert complexes associated to a resolved stratified pseudomanifold $X$ can be defined as in Definition \ref{definition_N_Hilbert_complex}, as well as $\mathbb{Z}$ graded families. Then we demand that each $Q_{j,\varepsilon}$ is of Dirac type (in particular $\Delta_{\varepsilon}$ is of Laplace-type), and that $H=L^2(X;\mathcal{E})$ for a Clifford module $\mathcal{E}$.
    In this case if the Hamiltonian can be expressed as
    \begin{equation}
    \label{equation_introduction_deformed_Hamiltonian}
        \Delta_{\varepsilon}=\Delta+\varepsilon K+ \varepsilon^2 A
    \end{equation}
    where $\Delta$, $K$ and $A$ are independent of $\varepsilon$, then we call $A$ the \textbf{potential} associated to the deformation. We refer to the connected components of \textbf{the vanishing locus of the potential} defined by $\{x \in X| A(x)=0 \}$, as the \textbf{vanishing set} and denote it by $Z(A)$. We assume that $Z(A)$ is finite (the number of connected components of the vanishing locus is finite).
    
    If the operator $K$ is $0$-th order, we say that \textbf{the Hamiltonian is of $0$-th order type}. 
    If $K=K'+iL_V$, where $L_V$ is a 1st order operator while $K'$ is 0-th order, then we say that \textbf{the Hamiltonian is of $1$-st order type}.    
\end{definition}

\begin{remark}
    In some contexts (see, e.g., \cite{maridakis2017spinor,maridakis2021localization,parker2022concentrating}) Hamiltonians corresponding to deformed families of $\mathcal{N}=1$ SUSY operators are studied, where the deformed operator is $D_{s}=D+s \mathcal{A}$ for some operator $\mathcal{A}$, the ``singular set'' of which, denoted by $Z_\mathcal{A}= \{ x \in \widehat{X}| \mathrm{ker}\mathcal{A}(x) \neq 0\}$ where $\mathcal{A}$ fails to be injective. There those studies have been conducted with conditions that yield Hamiltonians of $0$-th order type, similar to the case of Witten deformation of the de Rham operator (see Remark \ref{Remark_de_Rham_conjugation_N_2}).   
    This fits in our framework under assumptions that reduce the term $\mathcal{A}^*\mathcal{A}$ into a a scalar valued potential, since when squarframeworking the operator $D_s$ yields $\mathcal{A}^*\mathcal{A}=A$.
\end{remark}

In Witten deformation, families of SUSY Hilbert complexes arise from deformation of the operators $Q$ associated to some \textbf{\textit{dynamical system, or a classical mechanical system}} on an underlying space $X$, where the Hilbert complex should be considered as its \textbf{\textit{quantization}} (see Section 3.2 of \cite{jayasinghe2025supersymmetricharmonicoscillatorssingular}).
The limit as $\varepsilon \rightarrow \infty$ is a semi-classical limit. The idea then is that the quantum ground states at zero energy (zero-eigenvalue eigensections) of the Hamiltonians concentrate in $L^2$ mass near the classical ground states, which are at the vanishing sets of the potential.
The vanishing point set would match the turning points in the WKB formalism, of which Witten deformation is a SUSY version. 

The families of $\mathcal{N}=1$ complexes that we study in this article are of the following form.

\begin{definition}[\textbf{Witten deformation by vector fields, and geometric endomorphisms}]\label{def:witten-def-cmpx-vfs}
Given a family of $\mathcal{N}=1$ complexes $\mathcal{Q}_{\varepsilon}(X)=(H=L^2(X;\mathcal{E}), \mathcal{D}(D_{\varepsilon}),D_{\varepsilon},\gamma)$ associated to a resolved stratified pseudomanifold $X$ where $D_{\varepsilon}=D+i\varepsilon cl(V)$ for a Dirac type operator $D$ with principal symbol $\sigma(D) =  cl(-)$ and a \textit{compatible} Killing vector field $V$ generating an $S^1$ action, we say that it is a \textbf{deformation by a vector field.}

Given such a Killing vector field $V$ inducing an $S^1$ action, there is an endomorphism $T_f$, (also denoted $T_g$) on the space of sections induced by pulling back by the self map $f_g: X^{reg} \rightarrow X^{reg}$ for $g \in S^1$, that is
\begin{equation}
    T_g:=f^*_g: \mathcal{D}(D_{\varepsilon}) \rightarrow \mathcal{D}(D_{\varepsilon})
\end{equation}
which is extended to the domain by density of the smooth sections with compact support in the interior of the Domain, and we require that $[T_g,D]=0$ on the domain and this is the compatibility mentioned above.
We refer to such endomorphisms as \textbf{geometric endomorphisms} on the SUSY Hilbert complexes. 
\end{definition}

\begin{remark}
\label{Remark_nuances_for_general_geometric_endomorphisms}
    The definition of geometric endomorphisms given above works only for self maps associated to isometric $S^1$ actions. 
    For instance note that $\overline{z}^{-1}=z$ iff $z$ takes values in the unit circle, and for general $\mathbb{C}^*$ actions one needs to compensate for this when pulling back sections on the cohomology of the Dolbeault complexes (this issue can be checked by applying the pullback and taking traces instead of taking geometric endomorphisms in equation 7.58 of \cite{jayasinghe2023l2}.

    A more general construction that works more broadly is given in Section 4 of \cite{jayasinghe2023l2}, and Section 3.3 of \cite{jayasinghe2024holomorphic}.    
\end{remark}

\begin{definition}[\textbf{Circle actions and equivariant subcomplexes}]
\label{Definition_circle_actions_and_equivariant_subcomplexes}
Given a family of $\mathcal{N}=N$ SUSY complexes $\mathcal{Q}_{\varepsilon}(X) = \left( H, \{ ( \mathcal{D}(Q_{j,\varepsilon}), Q_{j,\varepsilon} ) \}_{j=1}^N ,\gamma \right)$ where each $(H, \{ ( \mathcal{D}(Q_{j,\varepsilon}), Q_{j,\varepsilon} ) \} ,\gamma)$ (for any fixed $j$ is a deformation by a vector field with a Hamiltonian of $1$-st order type, we will assume that the operator $iL_V$ in Definition \ref{definition_family_of_N_2_SUSY_complexes} is the generator of the infinitesimal action on the sections of the Clifford module, corresponding to the circle action generated by $V$.

Then for each $\mu$, we define the \textbf{equivariant Hilbert (sub)complexes}
\begin{equation}
    \mathcal{Q}_{\varepsilon}^{\mu}(X)= \Big( H,\{ (\mathcal{D}^{\mu}(Q_{j,\varepsilon}), Q_{j,\varepsilon})\}_{j=1}^N,\gamma \Big)
\end{equation}
where 
\begin{equation}
    \mathcal{D}^{\mu}(Q_{j,\varepsilon}) =\{ u \in \mathcal{D}^{\mu}(Q_{j,\varepsilon})\;|\; u \in \emph{Ker}(iL_V-\mu) \}
\end{equation}
and the Hamiltonian of $\mathcal{Q}_{\varepsilon}(X)$ restricts to one on $\mathcal{Q}_{\varepsilon}^{\mu}(X)$ that is of $0$-th order type.

We show that these are well defined as Hilbert complexes in Proposition \ref{Proposition_equivariant_complexes}.  
\end{definition}

We use Definition \ref{Definition_circle_actions_and_equivariant_subcomplexes} to reduce the study of families with Hamiltonians of 1-st order type to families of Hamiltonians with $0$-th order type. This allows us to give equivalent formulations for Morse inequalities, as in Theorem \ref{Theorem_introduction_holomorphic_Morse_inequal_not_abstract}, where those for a given Hilbert complex with a Hamiltonian of 1-st order type corresponds to countably many equivariant Hilbert complexes with Hamiltonians of 0-th order type after restricting to eigenspaces of $iL_V$.

\begin{remark}[More general actions of Lie groups]
Whenever we study group actions in this article, \textbf{unless otherwise stated we shall consider those of compact connected Lie groups $G$}. In particular these will have a maximal torus, and if non-abelian will carry an action of a Weyl group on the action chambers. Since elements that generate circle actions are dense in $G$, we shall build our main analytic results for circle actions and use density arguments to prove more general results, such as those for toric actions in Subsection \ref{Subsection_Toric_actions}.    

In the case of more general actions than those of compact connected Lie groups (such as algebraic toric actions, see Remark \ref{remark_local_extention_c_star}), we remind the reader that geometric endomorphisms are more nuanced and refer to Remark \ref{Remark_nuances_for_general_geometric_endomorphisms}.
\end{remark}

\begin{definition}[\textbf{SUSY refinements}]
    Given a $\mathcal{N}=N$ SUSY Hilbert complex $\mathcal{Q}$, if we can find another Dirac-type operator that makes it a $\mathcal{N}=N+1$ SUSY Hilbert complex, we say the that \textbf{$\mathcal{Q}$ admits a $\mathcal{N}=N+1$ SUSY refinement}.
\end{definition}
With this definition in hand, it suffices to prove most of the technical results in this article for $\mathcal{N}=1$ SUSY complexes and then take into account the improvements in the results when there are $\mathcal{N}=2$ refinements. When there are such refinements, there are simplifications due to second SUSY cancellations, coming via either boundary conditions or other mechanisms according to the proof (see Appendix \ref{Appendix_SUSY_and_conjectures}).

Since Witten deformation and equivariant index theory attempts to understand global invariants in terms of local ones, we require suitable assumptions for the vanishing loci of potentials that we study. We summarize what we study in this article.

\begin{assumption}[\textbf{localizing neighbourhoods have K\"ahler Hamiltonian actions}]
    For the families of SUSY complexes arising from deformation by vector fields that we study in this article, for each $F \in Z(A=|V|^2)$ (fixed point sets of $V$) we assume that there are \textbf{localizing neighbourhoods} $U_{t}(F)$ which admit locally K\"ahler Hamiltonian actions. This adds more structure on the critical point sets of stratified Morse-Bott type studied in \cite{jayasinghe2024witten} and we refer the reader to Subsection \ref{subsection_fixed_point_sets} for more precise details.    

    In the case of de Rham complexes, the ``fundamental neighbourhoods'' studied for stratified Morse-Bott functions in \cite{jayasinghe2024witten} are used as localizing neighbourhoods and our main technical results hold for those as well.    
\end{assumption}

In the case where $F$ and $X$ are both smooth, $U_{\infty}(F)$ is the normal bundle of $F$ in $X$, and $U_{t}(F)$ for finite $t$ are the disc bundles of the normal bundle with radius $t$. In the singular setting, roughly this can be thought of as a fiber bundle over the set $F$, where the fibers are the normal cones to the set at each point, in other words a suitable generalization of a normal bundle. 

Given such a localizing neighbourhood, we can decompose the model Dirac operators $D_{\varepsilon}$ (with respect to a model metric on $F$) into horizontal and vertical operators $D^{\mathcal{H}}=D_F, D_{\varepsilon}^{\mathcal{V}}=D_{C(Z),\varepsilon}$ that act on horizontal and vertical Clifford modules $\mathcal{E}_{\mathcal{H}}, \mathcal{E}_{\mathcal{V}}$, where $Z$ is the resolved link of the conic fibers over $F$. We define fiber-wise complexes for the vertical Clifford module and take their cohomology groups and see that they form \textbf{normal cohomology bundles}, which we denote by $\mathcal{N}_{F,\varepsilon}$ using which we define normal cohomology complexes. \textit{The following definition should be treated as an introductory statement for the benefit of the reader, and it will be made precise in Subsections \ref{subsec:localizing-domains} and \ref{subsec:cohomology-truncated-neighbourhoods}.}

\begin{definition}[\textbf{Localized and normal cohomology complexes}]
\label{definition_deformed_local_complexes_intro}
Given a family of $\mathcal{N}=N$ SUSY complexes $\mathcal{Q}_{\varepsilon}(X)$ with a potential $A$, for any $F \in Z(A)$, we define the \textbf{family of localized complexes of $\mathcal{Q}_{\varepsilon}(X)$ associated to $F$} to be a family 
\begin{equation}
    \mathcal{Q}_{\mathcal{B},\varepsilon}(U_{t}(F)) = \Big( H|_{U_{t}(F)}, \{ ( \mathcal{D}_{\mathcal{B}}(Q_{j,\varepsilon}), Q_{j,\varepsilon} ) \}_{j=1}^N , \gamma \Big)
\end{equation}
where $H|_{U_{t}(F)}=L^2(U_{t}(F);\mathcal{E})$ is the Hilbert space given by the restriction of the sections of $H$ with support on $U_{t}$. Here the domains $\mathcal{D}_{\mathcal{B}}(Q_{j,\varepsilon})$ are uniquely determined by the domains for the operators $Q_j$ on $X$, and an $L^2$ boundary condition at infinity on $U_{\infty}(F)$ (recorded by $\mathcal{B}$) corresponding to a growth condition for the potential of the Hamiltonian. For finite $t$, $\mathcal{B}$ stands for boundary conditions that are explained in Subsection \ref{subsec:cohomology-truncated-neighbourhoods}.

A \textbf{family of normal cohomology complexes over $F$ associated to $\mathcal{Q}_{\varepsilon}(X)$} is a family of complexes
\begin{equation}
    \mathcal{Q}_{\mathcal{B},\varepsilon}(\mathcal{N}_{F,\varepsilon}) = \Big( L^2(F; \mathcal{E_H} \otimes \mathcal{N}_{F,\varepsilon}), \{ ( \mathcal{D}_{\mathcal{B}}(Q_{F,j,\varepsilon}), Q_{F,j,\varepsilon} ) \}_{j=1}^N , \gamma \Big)
\end{equation}
where $Q_{F,j,\varepsilon}$ are the horizontal Dirac-type operators on $F^{reg}$ acting on the horizontal Clifford module twisted by the normal cohomology complex.
\end{definition}

Some of the main results in Section \ref{sec:domains} are concerned with defining the normal cohomology complexes appropriately and giving alternate characterizations of their cohomology groups in various cases, which in turn give alternate characterizations of the Witten complexes, Morse inequalities and Lefschetz fixed point theorems that we present below.
One can think of the normal cohomology complexes as complexes defined for the Horizontal Clifford module on $F^{reg}$ twisted by a local system (the cohomology of the fibers) where the precise definitions are involved, and we defer the definitions to Subsections \ref{subsec:cohomology-truncated-neighbourhoods} and \ref{subsec:localizing-domains}.

\begin{remark}[Domains and notational decorations for complexes]
Since Dirac-type and Laplace-type operators are not essentially self-adjoint on stratified pseudomanifolds, we shall study choices of self-adjoint extensions for them in Section \ref{sec:domains}, and we shall decorate the domains as $\mathcal{D}_{W}(Q_j)$, and $\mathcal{D}_{W}(\Delta)$ for choices of ideal boundary conditions at the strata that correspond to choices of self-adjoint extensions at the level of cohomology (an equivalence class of self-adjoint extensions that give the same cohomology).

When we study localizing neighbourhoods with truncated fibers, and fiberwise complexes, and normal cohomology complexes corresponding to such complexes, we denote the boundary conditions we use by $B$. Then the domains for the Dirac-type and Laplace-type operators are denoted by 
$\mathcal{B}=(W,B)$. When studying only a single component $F \in Z(A)$, we use $\mathcal{B}_F=(W,B_F)$.

Then we notate the corresponding complexes by similar decorations, for instance $\mathcal{Q}_{W,\varepsilon}(X),\mathcal{Q}_{\mathcal{B},\varepsilon}$.
\end{remark}

In this article we mainly focus on the cases where $\mathcal{N}=2$, and $\mathcal{N}=1$.
In the case where $\mathcal{N}=2$, Witten defines the nilpotent operators $Q_{\pm}=(Q_1 \pm \sqrt{-1} Q_2)/\sqrt{2}$ (see, e.g., \cite[\S 3]{witten1982constraints}), which can be easily checked to satisfy the nilpotency condition $Q_{\pm}^2=0$, and defines the cohomology groups $\Ker(Q_{\pm})/\text{Im}(Q_{\pm})$. While these are similar to co-chain complexes, a more precise notion that keeps track of the Hilbert spaces and domains are \textbf{Hilbert complexes} in the sense of Br\"uning and Lesch, introduced in \cite{bru1992hilbert}.

Since the operators $Q_{i}$ are self-adjoint, one can also check that $Q_{\pm}$ are adjoints of each other with respect to the inner product on the Hilbert space.
\textbf{We will denote $Q_+$ by $P$ and it is easy to check that $Q_-$ is $P^*$}.

The families of $\mathcal{N}=2$ SUSY Hilbert complexes that we study in this article arise as \textbf{families arising from conjugation}, and as such are equipped with quasi-isomorphisms of Hilbert complexes (see Proposition \ref{Proposition_cohomology_level_equivalences}). The notion of Hilbert complex we use is that introduced in \cite{bru1992hilbert}.

\begin{definition}[Hilbert complexes and families arising from conjugation]
\label{Definition_Hilbert_complexes_and_conjugated_families}
Given a $\mathcal{N}=2$ SUSY complex $\mathcal{Q}=(H,\{\mathcal{D}(Q_j),Q_j\}_{j=1}^2,\gamma)$ with $\mathbb{Z}$ grading, it determines a  
it determines a \textbf{Hilbert complex} $\mathcal{P}=(H,\mathcal{D}(P),P)$ in the sense of \cite{bru1992hilbert}, where $P:=(Q_1 + \sqrt{-1} Q_2)/\sqrt{2}$ with the domain given by taking the graph closure of $P$ on $\mathcal{D}(Q_1) \cap \mathcal{D}(Q_2)$. One can check that $\mathcal{D}(Q_j)=\mathcal{D}(P) \cap (\mathcal{D}(P))^*$ and that $Q_1=P+P^*, Q_2=i(P-P^*)$.

For a Fredholm complex, the cohomology groups for $P$ and $P^*$ are the cohomology groups for the Hilbert complex $\mathcal{P}$ ($ker(P)/im(P)$), and its adjoint Hilbert complex $\mathcal{P}^*=(H,\mathcal{D}(P^*),P^*)$. 

Given such a Hilbert complex corresponding to a $\mathcal{N}=2$ SUSY Hilbert complex, $\mathcal{P}=(H,\mathcal{D}(P),P,\gamma)$, we say that there is a \textbf{family of complexes arising from conjugation} $\mathcal{P}_{\varepsilon}=(H,\mathcal{D}(P_{\varepsilon}),P_{\varepsilon},\gamma)$ where $P_\varepsilon=M_\varepsilon^{-1}PM_\varepsilon$ equipped with the induced domain, if there exists a family of bounded operators $M_{\varepsilon}$ on $H$ with bounded inverses.
\end{definition}

The two main examples for such families arise from de Rham and Dolbeault complexes (twisted by flat and holomorphic bundles), where $M_{\varepsilon}=e^{\varepsilon h}$ where $h$ is a suitable smooth function in the de Rham case (see \cite{Kirwan_Morse_21} for functions with degeneracy) and $h$ is a K\"ahler Hamiltonian function in the Dolbeault case (cf. Remarks \ref{Remark_on_conjugation} \ref{Remark_de_Rham_conjugation_N_2}).

The $\mathcal{N}=2$ complexes arising from conjugation for the twisted Dolbeault complexes are also families arising from deformation using vector fields, and in fact holds more generally for spaces with symplectic forms and Hamiltonian vector fields where the spin-c Dirac operator can be deformed as a $\mathcal{N}=1$ complex. However this is not the case for the de Rham complexes (see Remark \ref{Remark_de_Rham_conjugation_N_2}).

\begin{definition}[\textbf{Witten complexes, Morse inequalities and Lefschetz fixed point theorems}]
\label{definition_Witten_complexes_and_more_intro}
Given a family of $\mathcal{N}=N$ SUSY Hilbert complexes $\mathcal{Q}_{W,\varepsilon}(X)$ associated to a stratified pseudomanifold $X$, with a potential $A$, we say it \textbf{admits a Witten subcomplex} $\mathcal{W}(\mathcal{Q}_{W,\varepsilon}(X))$ if there is a subcomplex 
\begin{equation}
    \mathcal{W}(\mathcal{Q}_{W,\varepsilon}(X)) = \Big( \mathcal{V}_{\varepsilon}^k, \{\mathcal{D}_W(Q_{j,\varepsilon})|_{\mathcal{V}_{\varepsilon}}, Q_{j,\varepsilon} \}_{j=1}^N,  \gamma \Big)
\end{equation}
where

\begin{equation}
     \mathcal{V}_{\varepsilon}^k:=\bigoplus_{F\in Z(A)}\mathcal{H}^k(\mathcal{Q}_{\mathcal{B},\varepsilon}(\mathcal{N}_{F,\varepsilon}) ).
\end{equation} 
\noindent Given a geometric endomorphism $T_{f}$ for $\mathcal{Q}_{W,\varepsilon}(X)$, we define the \textbf{Morse polynomial} to be 
\begin{equation}
    M(\mathcal{Q}_{W,\varepsilon}(X),T_f, b):=\sum_{k}b^k \operatorname{Tr} T_f|_{\mathcal{V}_{\varepsilon}^k}
\end{equation}
for both $\mathbb{Z}$ graded and $\mathbb{Z}_2$ graded complexes.
The \textbf{Poincar\'e polynomial} is defined to be 
\begin{equation}
    P(\mathcal{Q}_{W,\varepsilon}(X), T_f,b)=\sum_{k}b^k \operatorname{Tr} T_f|_{\mathcal{H}^k(\mathcal{Q}_{W,\varepsilon}(X))}.
\end{equation}
For both polynomials, one needs conditions for the trace to exist, for instance the finiteness of the vector spaces on which the traces are taken.

\noindent The \textbf{error polynomial} $E(b)$ is defined by the equality 
\begin{equation}
\label{definition_introduction_abstract_Morse_inequality}
    M(\mathcal{Q}_{W,\varepsilon}(X),b)=P(\mathcal{Q}_{W,\varepsilon},b)+(1+b)E(b)
\end{equation}
where $E(b)$ is a polynomial in $b$, the coefficients of which satisfy a \textbf{non-negativity condition} corresponding to $T_f$, the precise statement of which is equivalent to the content of a \textbf{Morse inequality}. 
For instance if $T_f=Id$, and the traces in the two polynomials are finite (this is the case if the Hamiltonian is of $0$-th order type), then the coefficients are non-negative integers.

A \textbf{Lefschetz fixed point theorem} in this context is equivalent to the content of equation \eqref{definition_introduction_abstract_Morse_inequality} for $b=-1$, whence the error polynomial does not contribute since the factor of $(1+b)$ appearing with it vanishes. Notice however this does not necessitate the existence of a Witten subcomplex for $\varepsilon=0$ just that there are local cohomology groups for the localizing neighbourhoods.
\end{definition}

\begin{remark}
    Given other geometric endomorphisms on the Witten complexes, there are more general Lefschetz-Morse inequalities (see \cite[\S 6.5.3]{jayasinghe2023l2}) which can be easily formulated in this generality.

    The $(1+b)$ factor arising in front of the error polynomial allows for the \textbf{Lacunary principle} of Morse theory to be used to better understand the Poincar\'e polynomial. For instance if the Morse polynomial has coefficients $0$ for the terms $b^{k-1}$ and $b^{k+1}$, then the $b^k$ term of the Morse polynomial and the Poincar\'e polynomial are equal. 
\end{remark}

Definition \ref{definition_Witten_complexes_and_more_intro} shows that finding good characterizations of the local cohomology groups allow for a better understanding of the Witten complexes, via which one can draw results on spaces admitting them, and even for explicit computations.
Such characterizations in Section \ref{sec:domains} connects Theorem \ref{Theorem_introduction_holomorphic_Morse_inequal_not_abstract} with Theorem \ref{Theorem_Witten_complex_N_2_abstract_intro} below (see Theorem \ref{Theorem_main_Witten_complex_Kahler_case}).

\begin{assumption}[Metrics and domains]   
\label{assumption_on_metrics_and_domains}
    We study \textbf{asymptotically $\delta$ wedge metrics} in this article, where in the symplectic setting we assume there are triples $(g_{w}, J_w, \omega_w)$ where $\omega_w$ is a wedge symplectic form and $J_w$ is a compatible wedge almost complex structure that we assume is integrable in the K\"ahler setting (see Subsection \ref{subsection_asymptotics_metric} for more details).
    \textbf{The two main results below are proven assuming $\delta>1/2$.}

    \textbf{We assume that the complexes are Fredholm and have discrete spectrum}, as in Definition \ref{definition_N_Hilbert_complex}.
    
    These assumptions are key to the proofs of the technical result in Proposition \ref{prop:block-estimates}, which are used to prove the main theorems that we state after the next Remark
\end{assumption}

\begin{remark}[Discrete spectrum and Fredholmness in geometric settings]
In this paper we do not attempt to exhaustively categorize all self-adjoint domains of $D_\varepsilon$ which satisfy these spectral property, leaving such analysis for future work. At present we content ourselves with observing that there are important examples of domains where this known already e.g. the minimal domain (see equation \ref{equation_minimal_domain_definition_1}) for the Dolbeault complex spaces with wedge Kähler metrics (see Proposition
3.6 of \cite{AnalyticToddBeiPiazza}) and Dirac-type operators equipped with the vertical APS domain (see Definition \ref{definition_VAPS_domain}) living on stratified pseudomanifolds (see \cite{Albin_2017_index}). See also sections \S \ref{subsubection_Fredholm_spectral} and \S \ref{subsection-spectral-background} for further discussion of these functional analytic considerations.

On the other hand, we can weaken this \emph{analytic} hypothesis if we are willing to strengthen our \emph{geometric} assumptions. Namely:
\begin{enumerate}
    \item If $Z(A)$ is a discrete set of isolated points with smooth links (singularities of depth 1), then for all possible (self-adjoint) extensions then have Theorem \ref{Theorem_spectral_properties_of_isolated_vanishing} yielding discrete spectrum and Fredholmness of the complex. In particular we can argue as in \cite[Prop 5.2]{jayasinghe2024witten} to conclude a stronger version of the block estimates proven here in Proposition \ref{prop:block-estimates}.
    \item If $Z(A)$ has non-isolated connected components, and the localizing neighborhood $U(F)$ satisfies that the induced product-type metric $g_{w,pt}$ gives $U(F)$ the structure of a flat (fiber) bundle, then again for all possible (self-adjoint) extensions we can prove \ref{Theorem_spectral_properties_of_isolated_vanishing} for such a SUSY Hilbert complex. Thus we can again argue as in \cite[Prop 5.2]{jayasinghe2024witten} to conclude stronger versions of the block estimates proven here in Proposition \ref{prop:block-estimates}.
\end{enumerate}
Having shown a stronger version of Proposition \ref{prop:block-estimates}, we can thus conclude all previously stated results up through Theorem \ref{Theorem_Witten_complex_N_2_abstract_intro}. The conclusions of Proposition \ref{prop:block-estimates} should be compared with the previous results of \cite{bismut1991complex,wu1998equivariant}, and note that the global discreteness of $\text{Spec}(\Delta_\varepsilon)$ are automatic in their geometric contexts; essentially being a consequence of the usual elliptic estimates on closed manifolds.
\end{remark}

The following two theorems are versions of Theorem \ref{Theorem_Witten_complex_N_2_abstract} and Theorem \ref{Theorem_N_1_Witten_complex}, and we refer to the proofs there.

\begin{theorem}[\textbf{Witten subcomplexes for $\mathcal{N}=2$ SUSY families}]
\label{Theorem_Witten_complex_N_2_abstract_intro}
Given a family $\mathcal{P}_{\varepsilon}(X)$ of Hilbert complexes associated to a family of $\mathcal{N}=2$ SUSY complexes arising from conjugation, with a Hamiltonian of zero-th order type where each element $F \in Z(A) \neq \emptyset$ admits a localizing neighbourhood $U_{t}(F)$ We assume that the complex and the localized complexes are Fredholm. Then there exists some large $\varepsilon_0 >0$ such that for each $\varepsilon > \varepsilon_0$, there exist finite dimensional Witten subcomplexes.

These Witten subcomplexes are isomorphic to each other each $\varepsilon > \varepsilon_0$, 
and are quasi-isomorphic to the complex $\mathcal{P}_{\varepsilon}(X)$.

Moreover for each $\varepsilon \geq 0$ there exists subcomplex $\mathcal{P}_{\varepsilon}(X)$ which is quasi-isomorphic to the Witten subcomplex for $\varepsilon >0$, that we also refer to as the Witten subcomplex by abuse of nomenclature (since all such Witten subcomplexes are isomorphic as complexes).

In particular, there are Morse polynomials and Morse inequalities.
\end{theorem}

\begin{theorem}[\textbf{Equivariant index theorem for $\mathcal{N}=1$ SUSY families}]
\label{Theorem_N_1_Witten_complex_intro}
Consider a $\mathcal{N}=1$ SUSY complex $\mathcal{Q}_{\varepsilon}=(L^2(X;\mathcal{E}),\mathcal{D}(D_{\varepsilon}),D_{\varepsilon},\gamma)$ where $\gamma$ is a $\mathbb{Z}_2$ grading, where the Hamiltonian is of $0$-th order type with a potential $A$, and where $Z(A) \neq \emptyset$.
Then there exists $\varepsilon_0>0$ such that for each $\varepsilon>\varepsilon_0$, there is a Witten subcomplex.

The equivariant index of $D_{\varepsilon}$ restricted to the Witten subcomplex is equal to the index of the operator $D_{\varepsilon}$ for any $\varepsilon >\varepsilon_0$ and for any $\varepsilon' \geq 0$, and
\begin{equation}
    \operatorname{Ind} _g(D):= \emph{Str}\, T_g|_{\mathcal{H}^k(\mathcal{Q}_{0}(X))}=\sum_{F \in Z(A)} \emph{Str}\,T_g|_{\mathcal{H}^k(\mathcal{Q}_{\varepsilon'}(\mathcal{N}_{F,\varepsilon'}))}.
\end{equation}
\end{theorem}

These two theorems apply to many cases of complexes.
As stated earlier, our initial motivation in this study was to prove the Holomorphic Morse inequalities in Theorem \ref{Theorem_introduction_holomorphic_Morse_inequal_not_abstract}, and this follows from Theorem \ref{Theorem_main_Witten_complex_Kahler_case} that we state next.

\begin{theorem}[\textbf{K\"ahler Hamiltonian actions give $\mathcal{N}=2$ families}]
\label{Theorem_main_Witten_complex_Kahler_case}
Consider a resolution $X$ of a stratified pseudomanifold $\widehat{X}$ equipped with a K\"ahler wedge metric and a K\"ahler Hamiltonian action $V$ that lifts to a holomorphic Hermitian bundle $E$.
Then there is a family of $\mathcal{N}=2$ SUSY Hilbert complexes arising from conjugation, equivalently a family of Hilbert complexes
\begin{equation}
    \mathcal{P}_{\varepsilon}(X)=(L^2(X;\mathcal{E}), \mathcal{D}_W(P_\varepsilon),P_\varepsilon), \quad P_\varepsilon=e^{-\varepsilon h} P e^{+\varepsilon h},
\end{equation}
where $W$ is any choice of ideal boundary condition that leads to a Fredholm complex with discrete spectrum. Here $h$ is the Hamiltonian function.
In particular Theorem \ref{Theorem_Witten_complex_N_2_abstract_intro} holds for each equivariant subcomplex.

For the complexes without restricting to the eigenvalues of the operator $\sqrt{-1}L_V$, we have subcomplexes
\begin{equation}
    \mathcal{W}(\mathcal{Q}_{\varepsilon}(X)) = \Big( \mathcal{V}_{\varepsilon}, \{\mathcal{D}(Q_{j,\varepsilon})\big|_{\mathcal{V}_{\varepsilon}}, Q_{j,\varepsilon} \}_{j=1}^N,  \gamma \Big), \quad \mathcal{V}_{\varepsilon}=\bigoplus_{\mu \in I_V}\mathcal{V}_{\varepsilon}^\mu
\end{equation}
where $\mathcal{V}_{\varepsilon}^\mu$ are the finite dimensional Hilbert spaces for each fixed eigenvalue $\mu$ of $\sqrt{-1}L_V$ that has non-trivial local cohomology complexes (indexed by $ I_V$)

Then given the geometric endomorphism $T_{\theta}$ associated to the group action, we define the equivariant Morse polynomials
\begin{equation}
    \sum_{F \in Z(A)}  \sum_{q=0}^n b^q \operatorname{Tr} \Big(T_{s,\theta}|_{\mathcal{H}^{q}(\mathcal{P}_{W,B}(\mathcal{N}_{F,\varepsilon}))} \Big) :=\sum_{F \in Z(A)}  \sum_{q=0}^n \sum_{\mu \in I_V} b^q e^{-\mu(s+i\theta)}\mathrm{dim}\Big(\mathcal{H}^{q}(\mathcal{P}^{\mu}_{W,B}(\mathcal{N}_{F,\varepsilon}) \Big) 
\end{equation}
and the Poincar\'e polynomials
\begin{equation}
    \sum_{F \in Z(A)}  \sum_{q=0}^n b^q \operatorname{Tr} \Big(T_{s,\theta}|_{\mathcal{H}^{q}(\mathcal{P}_{W}(X))}\Big) := \sum_{F \in Z(A)}  \sum_{q=0}^n  \sum_{\mu \in I_V} b^q e^{-\mu(s+i\theta)}\mathrm{dim}\Big(\mathcal{H}^{q}(\mathcal{P}^{\mu}_{W}(X)) \Big).
\end{equation}
Theorem \ref{Theorem_introduction_holomorphic_Morse_inequal_not_abstract} is a straightforward corollary.
\end{theorem}

\begin{proof}
    That there is such a family of Hilbert complexes arising from conjugation is a result of the discussion in Subsection \ref{subsection_Witten_deformation_Dirac}, in particular Proposition \ref{Proposition_operator_identities_Gollum} shows that the Hamiltonian is of the form we assume in Definition \ref{definition_family_of_N_2_SUSY_complexes}. Then the results for the associated $0$-th order complexes follow from Theorem \ref{Theorem_Witten_complex_N_2_abstract_intro}.

    The complexes for the $1$-st order complexes are straightforward constructions taking direct sums over the countable indexing set for $\mu$ when taking direct sums of stratified spaces. 
    
    Then the Morse and Poincar\'e polynomials are defined formally. A detailed proof in the case of isolated critical points is given in \cite{jayasinghe2024holomorphic} and the same proof goes through, \textit{mutatis mutandis} in this case, and we refer to that for the convergence of the power series in the polynomials, as well as for the statements on the error polynomial in Theorem \ref{Theorem_introduction_holomorphic_Morse_inequal_not_abstract}.   
\end{proof}

\begin{remark}[local cohomology groups]
In the case of the Dolbeault-Dirac operators, the cohomology of the normal cohomology complex over $F$ corresponds to the cohomology of the Dolbeault complex over $F$ with values twisted by the cohomology of the fibers. Here the cohomology of the fibers correspond to the tensor product of the cohomology groups (see \eqref{eq_polarized_normal_complex}) of a factorization of the conic fibers into attracting and expanding cones (based on the metric gradient flow of the hamiltonians) see Subsection \ref{subsec:cohomology-truncated-neighbourhoods} for more details.

These cohomology groups are much easier to compute on algebraic varieties and are analytic analogues of local rings, and techniques for computation developed in \cite{jayasinghe2023l2,jayasinghe2024holomorphic} can be used to compute them on normal cones.    
\end{remark}

\begin{remark}
\label{Remark_Hamiltonian_functions_exist}
It is well known that K\"ahler circle actions with fixed points on smooth manifolds are Hamiltonian by Frankel's theorem, and this extends to wedge setting (see, e.g., \cite[\S 3.2]{Mazzeo_2015}). Therefore, if there is a K\"ahler circle action which has fixed points, it is a Hamiltonian action. This is not true in general for non-K\"ahler actions which have non-trivial first Betti number (see for instance \cite{tolman2017non}).
\end{remark}

\begin{theorem}[\textbf{Hamiltonian actions on wedge symplectic spaces yielding $\mathcal{N}=1$ families}]
\label{Theorem_Hamiltonian_non_Kahler_intro}
Consider a resolution $X$ of a stratified pseudomanifold $\widehat{X}$ equipped with a triple $(g,J,\omega)$ where $g$ is an asymptotically $\delta>1/2$ wedge metric, $\omega$ a wedge symplectic form and $J$ a compatible wedge almost complex structure. Assume there exists a Hamiltonian action $V$ that lifts to a flat Hermitian bundle $E$ and that the twisted spin-c Dirac operator sastisfies the spectral Witt condition (see Definition \ref{definition_witt_spectral_witt}), i.e. is essentially self-adjoint).
Then there is a family of $\mathcal{N}=1$ SUSY Hilbert complexes arising deformation by vector fields, 
\begin{equation}
    \mathcal{Q}_{\varepsilon}(X)=(L^2(X;\mathcal{E}), \mathcal{D}_{VAPS}(Q_\varepsilon),Q_\varepsilon=D+i\varepsilon cl(V)),
\end{equation}
Here $D$ is the spin-c Dirac operator and $cl$ is the corresponding Clifford action, $W=VAPS$ is the VAPS ideal boundary condition in Definition \ref{definition_VAPS_domain}.

In particular Theorem \ref{Theorem_N_1_Witten_complex_intro} holds for each equivariant subcomplex, and we have an equivariant index theorem for the spin-c Dirac operator. 
\end{theorem}

\begin{proof}
That there is such a family of Hilbert complexes follows from the content in Subsection \ref{subsection_Dirac_Operators} (see Remark \ref{Remark_non_Kahler_Hamiltonian_2}). The form of the Hamiltonian in the smooth setting has been worked out (see, e.g., Theorem 1.4 of \cite{tian1998sympletic}) where it can be easily checked that it fits the form we assume. Then the discussion in Subsection \ref{subsection_Dirac_Operators} explain how it extends to the singular setting we study for any choices of ideal boundary conditions which give rise to self-adjoint extensions for $D_{\varepsilon}$.
Since it is not straightforward to see if such families have compatible domains, we restrict to the spectral Witt setting where the operators are essentially self adjoint.

Then the equivariant index theorem follows from applying Theorem \ref{Theorem_N_1_Witten_complex_intro}.
\end{proof}

There is work related to extending local cohomology groups for local neighbourhoods with only almost complex structures in      \cite{jayasinghe2023l2} (see Sections 5.2.4 and 7.2.5 of \textit{loc. cit.}), but here we only prove the above result for the local neighbourhoods we study, where we assume that the almost complex structures are integrable.

This deformed family has been usd to study geometric quantization, and questions of quantization commuting with reduction in the smooth setting, and we discuss relevant literature in Subsection \ref{subsection_introduction_History}.
The symplectic reduction (in particular K\"ahler reduction) of a smooth compact symplectic manifold will always yield a stratified pseudomanifold (see Theorem 2.1 of \cite{sjamaar1991stratified}), and we study actions in a generality which includes many actions coming from reductions on smooth spaces. In particular any iterated affine cone constructions will carry wedge K\"ahler metrics. Things are of course more flexible in the general symplectic setting. 

Notions such as wedge Hamiltonian vector fields will be defined in Subsection \ref{subsection_fixed_point_sets}.
In \cite{jayasinghe2024witten}, some of the authors of this article proved the existence of a Witten subcomplex for the de Rham complex, for \textbf{\textit{flat} stratified Morse-Bott functions} $h$ with non-isolated critical point sets (see Theorem 1.1 of \cite{jayasinghe2024witten}). This flatness, roughly, is that the vertical conic bundles at every non-isolated critical point set could be taken to be flat, which in turn implies that the fibration is locally a metric product over any chart of the regular part of the critical point set.
The proof techniques in this article allows one to drop the flatness condition and prove that for \textit{stratified Morse-Bott functions}, the Witten deformed de Rham complexes in that article (with various choices of domains) admit a Witten complex, which is finite dimensional (see the discussion around Theorem \ref{prop:specgap-flatness}. We formalize this as follows for completion, referring to \cite{jayasinghe2024witten} for more details.

\begin{theorem}[\textbf{Witten complex for stratified Morse-Bott functions for the de Rham complex}]
\label{theorem_de_Rham_Morse_Bott}
Let $X$ be a resolution of a stratified pseudomanifold equipped with a model wedge metric and a stratified Morse-Bott function $h$. Let E be a flat vector bundle on X, and let 
$\mathcal{P}_{W,\varepsilon}(X)=(L^2\Omega^k(X;E),\, \mathcal{D}_{W}(P_{\varepsilon}),\, P_{\varepsilon})$ be a twisted de Rham complex on $X$ with a choice of mezzo-perversity $W$, where the operator $P_{\varepsilon}=e^{-\varepsilon h}P e^{\varepsilon h}$, where $P=d_E$.
Then this family of $\mathcal{N}=2$ SUSY complexes admits a Witten subcomplex which categorify the stratified Morse-Bott inequalities.
\end{theorem}

The Morse inequalities given in Theorem 1.3 and Corollary 1.6 of \cite{jayasinghe2024witten}, as well as Lefschetz-Morse inequalities (see Theorem 6.35 of \cite{jayasinghe2023l2}) are easy Corollaries of the de Rham instanton complex presented above, and now hold without the flatness condition for the stratified Morse-Bott functions assumed in \cite{jayasinghe2024witten}.

\begin{proof}[sketch of proof]
The proof we give for the Dolbeault complex in Theorem \ref{Theorem_main_Witten_complex_Kahler_case} can be easily extended to the de Rham case. We have worked out the local spectral theory without the flatness condition now, which in \cite{jayasinghe2024witten} was done only with the flatness condition assumed, since Theorem \ref{Theorem_Witten_complex_N_2_abstract_intro} holds in this generality.

We note that in the de Rham setting, since the cohomology groups are independent of the choice of wedge metric (within the class of wedge metrics), the metrics on the fibers of the fundamental neighbourhoods can be assumed to be rigid, which reduces the operator on the fibers to the model operator, and one can avoid the polynomial estimates we used to prove the Born-Oppenheimer approximations in Section \ref{section_Witten_complexes}.
\end{proof}

\begin{remark}
    Note that Theorem \ref{theorem_de_Rham_Morse_Bott} holds for any wedge structure, and for any choice of ideal boundary condition (corresponding to mezzo-perversities as introduced in \cite{Albin_hodge_theory_cheeger_spaces,albin2025smoothatlasstratifiedspaces}, and used in \cite{jayasinghe2024witten}). We emphasize that this is because given any wedge structure one can always find a metric that is quasi-isometric to one which has model metrics of the type utilized in Section \ref{section_Witten_complexes}, and so we can take it to be asymptotically $\delta$ for arbitrarily large $\delta$.
    We do not need the spaces and the singular strata to be even dimensional or the metrics to be K\"ahler for the results in Section \ref{section_Witten_complexes}. In fact the choices of mezzo-perversities correspond to the cohomology groups of the links in the degree corresponding to $l/2$ where $l$ is the dimension of the link, and this is non-trivial only when there are singular strata of odd codimension.   
    
    The de Rham cohomology groups are invariants of the wedge structures and the mezzo-perversities, which follows for instance from their correspondence with intersection homology groups (see \cite{BanaglLclassnonWitt06,RefinedintersectionAlbinBanagl2015}).
\end{remark}

Beyond these cases, we present
\textbf{the following list of results that have been studied with Witten deformation on smooth manifolds, and in some cases in singular spaces, that the results and techniques developed in this article can help generalize.}

\begin{enumerate}
    \item The Poincar\'e index theorem for de Rham complexes (see \cite{Zhanglectures,zhang2000n}).
    \item The Kirwan-Novikov inequalities (see \cite{braverman2006kirwan} the techniques of which have also been used to construct Witten subcomplexes for non-Morse functions in \cite{Kirwan_Morse_21}).
    \item Atiyah's vanishing theorem for the Kervaire Semi-characteristic and generalizations of other equivariant mod 2 index theorems (see Section 7 of \cite{Zhanglectures}, and \cite{zhangkervaire,zhangkervaire99}). We study a narrow version of this briefly in Subsection \ref{subsection_Kervaire}.
    \item The case of Spin-c Dirac operators with Hamiltonian circle actions on spaces with symplectic structures. Theorem 1.4 of \cite{tian1998sympletic} can be used to see that this gives rise to a $\mathcal{N}=1$ complex with a $\mathbb{Z}_2$ grading with a Hamiltonian of 1-st order type as we study here. 
    \item Equivariant signature formulas \cite[\S 3]{witten1982supersymmetry} given a Killing vector field with isolated zeros. This includes localization with polyforms.      
    On spaces with K\"ahler Hamiltonian actions this follows from Poincar\'e-Hodge polynomials (see Subsection \ref{subsubsection_Poincare_Hodge_polynomial}). 
    \item Rarita-Schwinger and Spin Dirac, and other Dirac-type operators with compatible actions (see, e.g., \cite{witten1983fermion,Witten_elliptic_genera87}). Our results can be applied directly to generalize this to stratified pseudomanifolds. In the setting where the actions are K\"ahler Hamiltonian, we have Morse inequalities (for studies in the case of isolated fixed points, see \cite{jayasinghe2024holomorphic} where various notions of Rarita-Schwinger operators appearing in the literature are explained).
    
\end{enumerate}

Beyond these, Witten deformation has been used to understand Morse theory on infinite dimensional spaces corresponding to Topological Quantum Field Theories (TQFTs). For instance, Floer studied the Witten deformation corresponding to the conjugation of the analog of the de Rham operator on the space of all connections, by the functional $e^{\varepsilon h}$ where $h$ is the Chern-Simons functional (see \cite[\S 2]{WittenTopologicalQFT88}). Obviously our framework only works for QM systems, which are $1$-dimensional QFTs.

\subsection{Overview of the construction of the Witten complexes}

We give a high level overview on Witten deformation, and explain the main ideas in the proof that yield the Witten subcomplexes. The proof goes through SUSY quantum mechanics (QM), and the key idea is to quantize a dynamical system after formulating it as a classical mechanical system, and to then take a semi-classical limit to understand the quantum ground states in terms of geometric objects near the classical ground states. That is, given a flow (gradient flow, Hamiltonian flow, other flow generated by some vector field), one considers a potential $A$ associated to the flow and studies a semi-classical Schr\"odinger operator $\Delta_{\varepsilon}$. The classical ground states correspond to the minima of the potential, which if $A$ is non-negative are at $0$ where the classical dynamics come to a halt. Then the associated QM system, corresponding to the Hamiltonian $\Delta_{\varepsilon}$ (acting on a suitable Hilbert space, and with a suitable self-adjoint extension) will have quantum ground states corresponding to the minimal eigenvalue eigensections. If $\Delta_{\varepsilon}$ is non-negative but not positive, then these are the harmonic sections, and with sufficient assumptions on the spectrum (a spectral gap), one can show that there is a a subset of small eigenvalue eigensections which should contain the harmonic sections.
In the semi-classical limit as $\varepsilon \rightarrow \infty$, one can find suitable \textit{quasi-modes} for such small eigenvalue sections using semi-classical methods such as using a WKB ansatz, and for non-isolated vanishing sets of $A$, using Born-Oppenheimer approximations.
Under certain assumptions for the potential, it can be shown that the quasi-modes can be constructed by extending to the entire space, the local harmonic sections of model Hamiltonians near the vanishing set of the potential.

\begin{figure}[h]
    \centering
    \includegraphics[scale=0.60]{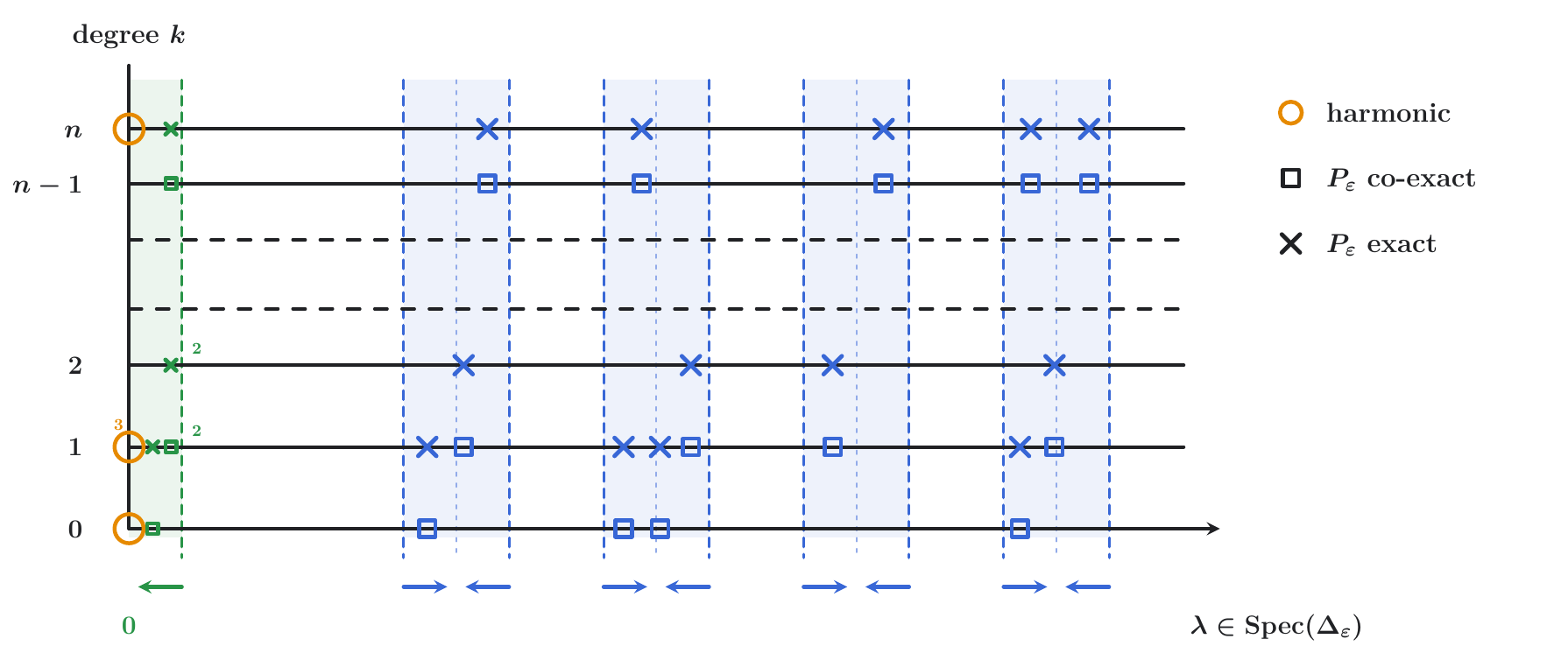}
    \caption{A spectral abacus for a $\mathbb{Z}$ graded  $\mathcal{N}=2$ supersymmetric complex corresponding to non-isolated critical points. The bands are around the eigenvalues of the Laplace operators on the horizontal Clifford modules of the vanishing sets $F$. The green band corresponds to the extensions to the space $X$ of the cohomology of the normal cohomology complexes, and includes the global cohomology.}
    \label{figure_spectral_abacus_non_isolated}
\end{figure}

Now when one imposes that there is SUSY, and takes a semi-classical QM system that respects the SUSY, then these small eigenvalue sections can to give rise to quasi-modes that correspond to a \textit{small eigenvalue subcomplex}. Then the local harmonic sections of the model Hamiltonians near the vanishing sets correspond to the cohomology of the normal cohomology complexes.

First it is shown that the global harmonic sections should localize near the vanishing set in the semi-classical limit, by a straightforward semi-classical estimate (see Proposition \ref{Propostion_growth_estimate_witten_deformed}).

Then the operator is shown to be well approximated by a model operator on any localizing neighbourhood (see Propositions \ref{proposition_rescaled_estimates-0} and \ref{prop:poly-est-delta-geq-1/2}), and then semi-classical estimates are proven in Subsection \ref{subsec:resolvent-ests} that show that the cohomology of the normal cohomology complex form a sub-complex including the global harmonic sections in the semi-classical limit. 
These latter estimates are essentially a singular Born-Oppenheimer approximation. The ``polyhomogenous approximation'' of the deformed Dirac-type operators using model operators generalizes work by the second author (see \cite{maridakis2017spinor,maridakis2021localization}), that circumvents some of the difficulties in extending the estimates in \cite{bismut1991complex} (estimating along geodesic flows in normal directions and the study of Jacobi fields is technically complicated when the geodesics may hit singular strata on singular links in finite time). 

The illustration in Figure \ref{figure_spectral_abacus_non_isolated} is inspired by Bott (see \cite{Bott1985,BottMorsetheoryindomitable}), analogous to figures 1 and 2 of  \cite{witten1983fermion}, and other such illustrations appearing in studies of SUSY. 
Since they encode the information of the Witten complexes used for spectral counting problems, we call such a diagram a \textbf{spectral abacus}.
For complexes with $\mathbb{Z}_2$ gradings, there are only two horizontal rails corresponding to the even and odd degrees, and the SUSY operator is $D_{\varepsilon}$.

When there are $\mathcal{N}=2$ refinements and the deformation arises from a conjugation then the conjugation gives a quasi-isomorphism of the entire SUSY Hilbert complexes, which gives isomorphisms for the local and global cohomology groups (which correspond to the harmonic sections). This allows one to take the small eigenvalue subcomplex for large values of $\varepsilon$, and get a quasi-isometric subcomplex at $\varepsilon=0$. Doing this for the normal cohomology complexes at a vanishing set, one obtains normal cohomology complexes with truncated normal fibers with boundary conditions that are dictated by the dynamical system (see Subsection \ref{subsec:cohomology-truncated-neighbourhoods}), which correspond to the generalized Neumann and generalized Dirichlet boundary conditions $(N$ and $D$), that unify the absolute/relative boundary conditions with the $\overline{\partial}$-Neumann condition and its adjoint condition for the Laplace-type operators (for the Hodge Laplacian and its Dolbeault counterparts respectively) with suitable twists. The conditions were studied an expounded on in detail in \cite{jayasinghe2023l2,jayasinghe2024holomorphic} for the case of cones using the Cheeger ansatz, and for other metrics in \cite{jayasinghe2025supersymmetricharmonicoscillatorssingular}. 
We shall say more about the spectral properties of such complexes in Subsection \ref{subsubection_Fredholm_spectral}.

\begin{remark}
\label{Remark_on_conjugation}
    Conjugation of $\mathcal{N}=2$ complexes for general QFT was already studied by Witten in \cite{witten1982constraints}, and in the case of QM (which is 1-dimensional QFT) in \cite{witten1982supersymmetry}, in which case arises naturally from the WKB formalism, when switching from the Schr\"odinger picture to the Heisenberg one, preserving SUSY.

    Thus for $M_{\varepsilon}:=e^{-\varepsilon h}$ one can think of $\varepsilon=\tau/\hbar$, where $\tau=it$ can be thought of one real imaginary time after a Wick rotation (also related to perspectives in \cite{gibbons1979classification}). One can think of the limit as $\varepsilon$ goes to $\infty$ as taking the semi-classical limit after fixing imaginary time (see \cite[\S 3.2]{jayasinghe2025supersymmetricharmonicoscillatorssingular}). Separately, the dualities associated to taking the inverse flow (the negative Hamiltonian function) in Subsection \ref{Subsection_Toric_actions} corresponds to a time reversal symmetry (see also \cite{jayasinghe2024holomorphic,jayasinghe2024witten}, and Remark \ref{Remark_21_on_dualities}).
\end{remark}

The relations between the cohomology groups with boundary conditions and certain \textit{Dirac cohomology groups} are explained in Appendix \ref{Appendix_SUSY_and_conjectures}, and will be explored in greater depth in the upcoming revision  of \cite{jayasinghe2026spectralasymmetrysupersymmetryequivariant}.
In the case of the localizing neighbourhoods, the second SUSY operator yields a decomposition of the spectrum for the Dirac operator on a cone into a direct sum over the cohomologies of the $P$ and $P^*$ complexes, and a complement which has a ``second SUSY". This is the mechanism behind the reduction of the $\xi_T$ invariants that are defined more generally when there is only 1 SUSY, to the case of supertraces over local cohomology groups, and we refer to Appendix \ref{Appendix_SUSY_and_conjectures} for more details where we present many more spectral abacuses which give a clearer picture of the SUSYs.

Lagrangians for the classical mechanical systems corresponding to the gradient flow of Morse functions (that quantizes to the Witten deformed de Rham complex, so to speak) were presented in \cite{witten1982supersymmetry} (see equation (23) of \textit{loc. cit.}), and those corresponding to the action after rewriting in an Euclidean metric and discarding fermions to simplify to a model case as well, where the instantons were studied. Correspondences with the dynamical Morse-Smale complexes were also studied there, which have since been developed in some cases on singular spaces (see \cite{ludwig2024morse,ludwig2017comparison}).

The relevant semi-classical analysis here is similar to the study of tunneling when there are potential wells, and there the particle that is tunneling \textit{spends very little time on (an instant on) the highest probability path} (path of steepest descent for a gradient flow) away from neighbourhoods of the vanishing sets, and such paths are therefore called instantons (see \cite{BottMorsetheoryindomitable}), and related complexes are sometimes called instanton complexes.

For more general QFT's we refer to \cite{witten1981dynamical,witten1982constraints} and \cite{WittenTopologicalQFT88} for TQFT's.
We also refer the interested reader to \cite[\S 3.2]{jayasinghe2025supersymmetricharmonicoscillatorssingular} for an attempt to formalize the SUSY QM systems involved, and for a quick summary of what it means for SUSY to be spontaneously broken.

\subsubsection{Domains and their Fredholm and spectral results}
\label{subsubection_Fredholm_spectral}

A complete theory for self-adjoint extensions of symmetric operators on stratified pseudomanifolds is yet to emerge and this has been the focus of much prior literature.
We refer to the many references in \cite{Albin_hodge_theory_cheeger_spaces,Albin_hodge_theory_on_stratified_spaces} for an overview.
Our motivation here was to set a framework which would work for arbitrary self-adjoint extensions of Dirac-type operators, with compatibility conditions for $\mathcal{N}=2$ SUSY.

\textbf{One of the most important points of view that we present is that it suffices to pick domains at the level of cohomology for index theoretic and Morse theoretic problems, and that \textit{extra SUSY cancellations} play a role in these choices of domains for $\mathcal{N}=2$ SUSY Hilbert complexes.}
\textbf{We shall comment on how the procedure we use to pick domains at the level of cohomology implicitly uses a second supersymmetric cancellation in \ref{subsection_appendix_second_susy_cancellation}.}

We follow a method influenced by the works of \cite{bera2025diracoperatorstwistedramified,booss1998maslov,jayasinghe2026spectralasymmetrysupersymmetryequivariant,Albin_hodge_theory_cheeger_spaces,Albin_hodge_theory_on_stratified_spaces,Albin_2017_index}), and use Gelfand-Robbin quotients to characterize choices of self-adjoint extensions for the SUSY operators that we study (see Subsection \ref{subsubsec:gelfand-robbin-quotient}). In fact, we only characterize these domains at the level of cohomology.
We then give a procedure for identifying domains, which seems to work in practice for de Rham and Dolbeault complexes (with suitable twists of coefficients) in Subsection \ref{subsubsec:ideal-bc-stratified}. This works especially well with analytic metrics such as Fubini-Study and affine metrics on complex varieties. In the de Rham setting, this procedure is different from the iterated procedure for picking choices of domains given in \cite{Albin_hodge_theory_cheeger_spaces,Albin_hodge_theory_on_stratified_spaces}.
\textbf{and we make no claim that we give such concrete characterizations as in \textit{loc. cit.} of all possible choices of domains for the operators we study}.

We now explain the contributions of this article to advancing spectral theory and Fredholm results.

In \cite{jayasinghe2024holomorphic}, the first author proved the holomorphic Morse inequalities for isolated critical point sets on stratified pseudomanifolds with wedge K\"ahler metrics, generalizing some techniques of Zhang which assumed that there is discrete spectrum for the associated Hamiltonians on compact smooth manifolds.
While writing this article, the authors realized that in fact the proof there goes through without assuming any spectral properties for the global complex, and in fact that it can be used to show that for any choice of ideal boundary condition for a family of $\mathcal{N}=2$ SUSY complexes arising from conjugation, that the $0$-th order complexes are Fredholm and have a spectral gap. This is the content of Theorem \ref{Theorem_spectral_properties_of_isolated_vanishing}. This can be quite useful in different ways. Given any stratified pseudomanifold with wedge and even more general types of metrics (see, e.g., \cite{jayasinghe2025supersymmetricharmonicoscillatorssingular,Jesus2018Wittens,Jesus2018Wittensgeneral}) there are stratified Morse functions that can be used to show that the de Rham complexes with choices of self-adjoint extensions are Fredholm. There are more choices of domains for the Hodge Dirac operator for general metrics than those characterized in \cite{Albin_hodge_theory_cheeger_spaces,albin2025smoothatlasstratifiedspaces} and in light of Question \ref{question_cohomological_domains_dictated_at_critical_points} that a theory could be developed more broadly. For Dolbeault-Dirac complexes, this requires that the space has a K\"ahler Hamiltonian action, and one could pursue similar questions where this is available, for instance on toric varieties with Hamiltonian actions.

The Witten complexes in \cite{jayasinghe2024witten} assumed that the normal fibers to the critical point sets formed flat fiber bundles to prove the spectral theory necessary, and we have summarized the main result for localized complexes there in Proposition \ref{prop:specgap-flatness}. Some of the main technical results of this work was to prove estimates that avoided such assumptions.
This has been accomplished, but as we phrased our results in Subsection \ref{subsection_intro_SUSY_Hilbert_and_Witten} that the choice of ideal boundary conditions $W$ that we pick for the complexes needs to give rise to a Fredholm Hamiltonian (on the $0$-th order complexes) and have discrete spectrum. The main reason is that we follow \cite{bismut1991complex,wu1998equivariant} and prove Proposition \ref{prop:block-estimates} in order to get estimates for quasi-modes. However not all the sections of the normal cohomology complexes we introduce are quasi-modes, and any non-zero eigenvalue eigensections of the Dirac-type operators on the vanishing sets $F$ play a role. This forces one to assume that the Dirac-type operators on the normal cohomology complexes have discrete spectrum. Given this condition, using the spectral results for Schr\"odinger operators in Theorem \ref{Theorem_spectral_discreteness_Schrodinger} our proof goes through. 

For instance the minimal domain for the Dolbeault complex on spaces with wedge K\"ahler metrics (see Proposition 3.6 of \cite{AnalyticToddBeiPiazza}), and the VAPS domain of \cite{Albin_2017_index} (which for the Dolbeault-Dirac operator on K\"ahler manifolds is a rolled up extension by results in \cite{jayasinghe2024holomorphic}) give Fredholm domains.
The VAPS domain can be defined for any choice of Dirac-type operator with satisfying a Witt condition, which holds for the Dolbeault complexes. If in addition the spectral Witt condition is met (which is the case if the links at all strata have small volume), then the complexes are essentially self-adjoint.

Theorem \ref{Theorem_spectral_discreteness_Schrodinger} is an extension of a classical argument to stratified spaces, based on the idea that \textit{non-negative potentials that grow as one approaches infinity are confining, and quantize the spectrum}. We believe that this can be used to study non-compact stratified pseudomanifolds with similar potential growth in a fairly straightforward manner using the results of this article. 

Of course when there are actions with isolated critical points available, we can prove the results for the case of isolated critical points for any choice $W$. One interesting case of study for the initial development of the theory is toric actions which we go over in detail in Subsection \ref{Subsection_Toric_actions}, and there the non-isolated fixed point sets (vanishing sets) occur at walls of the action chamber, and the methods in Section \ref{section_examples} can be used to compute the Morse inequalities and Equivariant index formulae at such fixed points using the formulae for the generic action.

We have mainly referred to prior work by the first author for the techniques to prove discreteness of spectrum and Fredholm properties (for the $0$-th order complexes), and we give an overview on the different approaches to study spectral theory on metric cones corresponding to $\mathcal{N}=2$ complexes.
In Cheeger's article \cite{cheeger1983spectral}, section 3 (title of \textit{functional calculus on cones}) develops an ansatz for eigensections of the Hodge Laplacian (using SUSY, even though it is not explicitly expressed as such), and proceeds to study the essential spectrum on infinite cones and construct expressions for heat kernels on them, eventually using it to obtain expressions for the eta invariant (of the signature operator). In \cite[\S 5.2.4]{jayasinghe2023l2} (title of \textit{functional calculus on truncated cones with boundary conditions}), this ansatz was studied using SUSY and extended to wedge metrics on metric cones (over singular links with appropriate ideal boundary conditions) for both the de Rham and Dolbeault complexes, where the boundary conditions described above were imposed to obtain a description of the spectrum, and the cohomology groups with the boundary conditions (where those for the $N$ and $D$ conditions are dual to each other, by a duality of Hilbert complexes) were understood. In fact it was noticed that it works with Witten deformation by conjugation as well, and this was used to construct Witten complexes for de Rham complexes, and in \cite{jayasinghe2024holomorphic} in the Dolbeault case (cf Section 4 \cite{ludwig2017index} for the de Rham case on infinite cones). A thorough study of these in a more general metric setting was done in expository style in \cite{jayasinghe2025supersymmetricharmonicoscillatorssingular} and we refer the reader to this text in lieu of repeating the details in this article. The case of infinite cones with $L^2$ boundary conditions is similar as well since the growth of the potential at infinity leads to discreteness of the spectrum and is studied in \textit{loc. cit.}.
We believe that Cheeger's ansatz in tandem with the second SUSY cancellations explained in Appendix \ref{Appendix_SUSY_and_conjectures} could be used to approach these problems as well, though we do not pursue this line of approach in the present article. 

\subsection{Historical context and motivations}
\label{subsection_introduction_History}

This subsection partly presents some historical context and motivations for this work, while we discuss motivations from different angles in Appendix \ref{Appendix_SUSY_and_conjectures}.

In certain investigations of SUSY breaking in \cite{witten1981dynamical,witten1982constraints}, Witten discovered that one could prove the classical Morse inequalities using semi-classical analysis, and categorified the Morse polynomials using what is now called Witten deformation in \cite{witten1982supersymmetry} (cf. \cite{Bott1985,BottMorsetheoryindomitable}), in which article he also explained how the perspective generalized more broadly.
In \cite{witten1983fermion,witten1984holomorphic} he used these techniques to study equivariant index formulae, and formulated the holomorphic Morse inequalities in the setting of K\"ahler manifolds. 
Rigorous proofs of these results and a wider investigation of their consequences were carried out in for the de Rham complex in \cite{helffer1985puits}. For the Dolbeault complex this was done in \cite{mathai1997equivariant,wu1996equivariant,wu1998equivariant,wu2003instanton} more than a decade after Witten's papers, which had by then been inspired more techniques in equivariant cohomology which had already unlocked a lot more formulas which are still widely used today.

Witten deformation was also used to great effect to study analytic torsion and Quillen metrics (see \cite{bismut1991complex,BismutZhangCheegerMuller1992}). Notions of $\mathbb{Z}_2$ graded torsion has also been studied in \cite{mathai2010twisted,mathai2012topology,mathai2011analytic}.

It has also been used to study equivariant index formulas and geometric quantization for Hamiltonian actions on symplectic manifolds, and certain mod 2 invariants in the presence of circle actions including the Kervaire semi-characteristic (see \cite{zhang1998holomorphic,zhangkervaire,zhangkervaire99,ShubinNovikov,tian1998sympletic,wu2000note,tian1998holomorphic,braverman1997symplectic,braverman1999holomorphic,braverman2000index}).

The study of linear operators on singular spaces is an ongoing endeavor where a wide variety of tools are being used to build different approaches of study, and different phenomena are being investigated from various points of view. 
In addition to revealing what happens in singular settings, these investigations bolster our understanding of phenomena and methods on smooth spaces.

Index theory on smooth spaces is widely used in investigations of various mathematical and physical theories, and situates itself at the core of many landmark results in representation theory, symplectic geometry, gauge theory, mathematical physics. Part of the appeal in taking index theory as the motivation to further the study of analysis on singular manifolds is because the SUSY aids the analysis, yielding many cancellations of quantities that may otherwise be ill behaved. \textbf{Having multiple SUSY makes things even easier.}

Equivariant index theory on singular spaces is still in relative infancy, and various subtleties are still being probed with different methods. In the de Rham case, Witten deformation was extended to singular settings in \cite{Jesus2018Wittens,Jesus2018Wittensgeneral,ludwig2013analytic,ludwig2017comparison,ludwig2017index,ludwig2020extension}, proving not only the Morse inequalities but also extending Cheeger-M\"uller type results.
In \cite{jayasinghe2024holomorphic}, the first named author generalized the holomorphic Morse inequalities for K\"ahler Hamiltonian actions with isolated fixed points on stratified pseudomanifolds with K\"ahler wedge metrics, a first on any class of singular space, in particular constructing the instanton complexes. Moreover its potential was demonstrated by generalizing various results of interest to mathematicians and physicists to these spaces, identifying new objects that appear in such results in the singular setting and the subtleties that accompany them. This work built on the holomorphic Lefschetz fixed point theorem in \cite{jayasinghe2023l2}, and it lays the foundations for investigations of complex analytic torsion, generalizing work in \cite{bismut1991complex,bismut1988analytic} with Bismut.
In \cite{jayasinghe2024witten} some of the authors generalized the Witten complex for singular spaces with wedge metrics and functions with \textit{stratified} Morse-Bott type singularities. 

Methods developed by the second author in the proof of the localization theorem in \cite{maridakis2017spinor,maridakis2021localization} and which used different ideas from that in \cite{bismut1991complex} have been generalized to the singular setting in this article.

The work in \cite{banagl2025gsignaturetheoremwittspaces} generalizes the equivariant signature theorem to fixed point sets that are normally non singular in the sense introduced by Goresky and MacPherson. MacPherson developed Riemann-Roch theorems using the method of deformation to the normal cone, together with Fulton and Baum in \cite{baumfultonmacKtheoryRiemannRoch,baumfultonmacpheresonRiemannRoch}, equivariant formulas for which were proven in \cite{baum1979lefschetz,Baumformula81}. We discuss some extensions of these results to the singular setting in our article in Subsection \ref{Subsection_comparison_Lefschetz_smooth} and Subsection \ref{subsection_appendix_index_theory}.

In \cite{Kirwan_Morse_21}, the authors prove Morse inequalities for \emph{degenerate} Morse functions, and do Witten deformation on a modified space that simplifies the analysis, following the ideas of \cite{braverman2006kirwan}. Here the choice of metric really is like a choice of gauge since it is the de Rham complex being studied, but even for the de Rham complex, in the singular world the choices of ideal boundary conditions and corresponding cohomology theories are only invariant upto quasi-isometry classes of metrics and more care must be taken. It is natural to study degenerate functions on stratified spaces as well, and for certain types of degeneracies variations of the notions of potentials studied in \cite{jayasinghe2025supersymmetricharmonicoscillatorssingular} together with the methods introduced in this article already yield such results.
It is easy to see that the holomorphic Witten complexes exist in general for more classes of metrics than the wedge ones, when they are K\"ahler and admit isometric K\"ahler Hamiltonian group actions (see also \cite{jayasinghe2025supersymmetricharmonicoscillatorssingular}, and Subsection \ref{subsubection_Fredholm_spectral}).

We also discuss the holomorphic Morse inequalities and compare with the smooth and singular formulae in Subsection \ref{Subsection_comparison_Lefschetz_smooth}, and study further conjectures in Appendix \ref{Appendix_SUSY_and_conjectures}.
In Section 6.2 of \cite{jayasinghe2024holomorphic}, the underlying philosophy that was proposed by Witten, of categorifying a Lefschetz fixed point formula in a manner that admits Morse inequalities is presented. Conjectural Morse inequalities for the case of the Lefschetz-Riemann-Roch formulas of \cite{baum1979lefschetz,Baumformula81} are presented and a similar description goes over to the case of non-isolated critical point sets. In particular, on spaces where the local and global cohomology groups corresponding to $L^2$ cohomology and the cohomology groups studied by Baum-Fulton-Quart are isomorphic, the holomorphic Morse inequalities we study here realize such conjectural Morse inequalities. In Appendix \ref{Appendix_SUSY_and_conjectures} we conjecture this holds for non-equivariant Riemann-Roch numbers as well, with proofs in certain cases in an upcoming revision of \cite{jayasinghe2026spectralasymmetrysupersymmetryequivariant}.

In a different way of looking at things, we replace the computation of equivairant index formula contributions from superdeterminant formulas on tangent spaces of fixed points as in \cite{AtiyahBott1,AtiyahBott2} (and normal bundles to fixed point sets) to a more categorical framework. It seems plausible that Guillemin-Sternberg and Atiyah-Bott-Guillemin formulas may also yield such categorifications that enable generalization. This is more plausible if one considers the case where the foliation of a contact manifold corresponds to that given by Hamiltonian flows on its symplectizing cone.

\subsection{Overview of Sections}

We briefly explain the content of each section to aid the reader navigate the content, based on their preferences, background, and reading styles.

Section \ref{section_stratified_and_wedge} gives a review of stratified spaces and their resolutions, and proceeds to explain wedge metrics and wedge versions of geometric structures including complex structures, symplectic forms and K\"ahler structures. Then the local forms of vector fields corresponding to the deformations we study are introduced together with the structure of the localizing neighbourhoods near fixed point sets.

In Section \ref{section_Witten_deformation}, we frame our results for general SUSY Hilbert complexes by describing the \emph{core mathematical object} we use to formalize the notion of supercharges $Q_j$, that of Dirac-type operators; afterwards we specialize to the case of spin-c Dirac and Dolbeault operators, and highlight how these arise as various examples of $\mathcal{N}=1$ and $2$ SUSY. After we recall the construction of Witten deformations of Dirac operators; the case of deformation of Dolbeault-Dirac operators motivates the study of equivariant families of Hilbert complexes and we present this alongside the Bochner-Schr\"odinger formulas and arguments for localization.

In Section \ref{sec:domains}, we begin the study of domains for SUSY Hilbert complexes, in order to later establish analytic properties of our operators. In order to describe closed extensions of operators we recall the theory of ideal boundary operators, before rephrasing this concept via the Gelfand-Robbin quotient; this construction formalizes a concept very familiar to physicists, viewing the boundary term in an integration-by-parts formula as a symplectic form on a certain Hilbert space, and realizing self-adjoint extensions as Lagrangian subspaces. After observing that this construction extends to the case of Hilbert complexes, and use this to obtain Proposition \ref{prop:cohom-splitting}, a notable tool for explicit computation of cohomology of local and global complexes (as it is used in \S \ref{subsection_non_isolated_singularity}). After recalling the familiar concepts of the vertical APS domain of a Dirac operator, and the (geometric) Witt conditions, we go over how to pick ideal boundary conditions iteratively in the case of stratified spaces with a stratification of depth $\geq 1$; we provide an explicit and archetypal example of such a space in \S \ref{Example_singular_toric_variety_domains}, a singular complex projective variety with a maximal stratum of depth 2. Finally we introduce the essential new definition of localizing domains, and use this notion to define localized complexes, fiberwise complexes, and normal cohomology complexes; these definitions are summarized in Definition \ref{Definition_all_localized_domains}. This final concept of normal cohomology complexes should be thought of as the analogue of the complexes along critical sets/fixed point sets implicitly studied in  the much simpler smooth setting in \cite{bismut1991complex, wu1998equivariant} which can themselves be thought of as formalizations of Witten's observation \cite{witten1984holomorphic} that the spectrum of the Witten-deformed (Dolbeault-)Laplacian behaves in the limit as if it is well-approximated by the \emph{Born-Oppenheimer approximation}. In this analogy, the sections on the vanishing loci play the role of ``slow modes'' and the general fiberwise sections are ``fast modes.'' We conclude the section with a brief discussion various of isomorphisms between these different complexes. 

In section \ref{section_Witten_complexes} we complete the remaining technical work for concluding the main results, by constructing the Witten complexes needed to prove Theorems \ref{Theorem_Witten_complex_N_2_abstract} and \ref{Theorem_N_1_Witten_complex}. We begin with \S \ref{subsection_polyhomogenous_estimates} which considers the hypothesis of our metrics being asymptotically $\delta$-wedge metrics, and from this derives polyhomogeneous estimates of the difference between the true Witten-deformed $D_{X,\varepsilon}$ and a model operator. This analysis should be viewed as a substitute for normal coordinates in the singular setting; if $X$ was smooth we could simply obtain such estimates by considering the normal bundle of a vanishing set $F\in Z(A)$ (as is done in e.g. \cite{bismut1986witten,bismut1991complex}). However, for a stratified space, given the lack of a well-defined normal bundle we must instead work much more carefully, our construction ultimately relying again on the notion of localizing neighborhoods which we have introduced. Finally we apply these polyhomogeneous estimates, and the normal cohomology complex of \S \ref{sec:domains} to obtain block operator estimates of the projection of $D_{X,\varepsilon}$ onto/off a natural subspace of $L^2(X;\mathcal{E})$. This permits a direct comparison of the spectrum of $D_{X,\varepsilon}$ and an induced $D_F$ in the limit $\varepsilon\to \infty$. Having done so we produce Witten subcomplexes for families of $\mathcal{N}=1,2$ SUSY Hilbert complexes, and highlight how the number of supersymmetries influences the kind of theorem one can prove (strong Morse inequalities vs. equivariant index formulae). 

In \S \ref{section_Morse_inequalities} we consider applications of our main results to stratified $X$ including K\"ahler Hamiltonian actions of compact connected Lie groups, and the de Rham operator on stratified wedge pseudomanifolds, and relate these  to local and global manifold invariants; this relates to \S \ref{Subsection_comparison_Lefschetz_smooth} which offers many insights on the equivariant Riemann-Roch formulae. For K\"ahler Hamiltonian actions we present many results, including Theorem \ref{Theorem_strong_Morse_Dolbeault_version_groups}, which gives a strengthening of the Holomorphic morse inequalities and generalizes \cite[Prop 2.6]{wu1998equivariant}. In this section we also prove several rigidity and vanishing theorems; new computations of Bolt charges for singular examples and discussion of their relations to local contributions to equivariant signatures; and new results generalizing formulae of \cite{zhangkervaire,zhang2000n} which compute Kervaire semi-characteristics and which reformulate them in terms of local cohomology groups.

In \S \ref{section_examples} we present several examples and offer a guide for computations. This is done compute several of the invariants of study including the holomorphic Morse and holomorphic Poincar\'e polynomials; these  examples of the local cohomology groups and study toric actions on both smooth and stratified spaces. We strongly encourage any readers interested in computing any of these (or related) invariants in singular settings to read this section carefully. 

In \S \ref{Appendix_spectral_theory} we collect several technical functional analytic results and discussions of previous work on the spectral theory of Dirac operators acting on stratified spaces. 

Finally in \S \ref{Appendix_SUSY_and_conjectures} we aim to present the role SUSY has previously and could further play in Index theory. In particular, we gather several examples of significant phenomena in which the definitions and computations various invariants built from Dirac operators can be fruitfully understood and simplified in the case of extra SUSY's. We conclude with several conjectures in Index theory and questions coming from physics. This appendix arguably encapsulates the motivation for this work, and we welcome further questions or discussion regarding it.

\subsection{Acknowledgements and AI statement}

\textbf{AI statement}: During earlier stages of this work, GJ attempted to use various LLMs to try and generate proofs of several technical results (mostly constrained to Section \ref{section_Witten_complexes}), in the hopes that those arguments were by now routine enough in the literature to be ``automated''. The results were mostly unsatisfactory, and other than his use in some revision and rewriting, all of the authors thought up and wrote down their own proofs. The technical core of all arguments, and all key insights in this work were those of the authors themselves, unassisted by any LLM-aided suggestions or discussions, and had been developed \emph{substantially before} any recent developments in the capabilities of various models. In preparing this manuscript, ChatGPT Astra was used in the calculation of specific numerical computations in Sections \ref{subsection_bolt_charges} and \ref{section_examples}, and were used to save time in producing the various figures.

\textbf{Acknowledgements}: GJ would like to thank Alex Taylor, Sebastian Goette, Siye Wu, Pierre Albin, Paolo Piazza, Rafe Mazzeo, and Lorenzo Foscolo for many interesting conversations that influenced this work, and thanks The Aristotle University of Thessaloniki for its hospitality and support.

MM and GJ thank Iakovos Androulidakis for many interesting conversations that helped propel this work.

HQ would like to thank Pierre Albin, Fran\c{c}ois Monard, and Marika Taylor for many useful discussions that furthered his understanding of many aspects of this work, and especially thanks Fran\c{c}ois for his unceasing support and mentoring during HQ's first year as a parent. 

XY is grateful to Pierre Albin and Alex Taylor for many insightful discussions and helpful comments.

GJ was supported by the Ministero dell’Università e della Ricerca of Italy – Avviso FIS2, grant FIS-2023-0395, CUP B53C25001070001.
XY was supported by Boya Postdoctoral Fellowship at Peking University.

\subsection{Notation}

\textbf{Section 1}
\begin{multicols}{2}
\begin{enumerate}
    \item $\widehat{X}$ stratified pseudomanifold, $X$ resolution
    \item $h$ Kähler Hamiltonian, $A = |dh|^2$ potential, $Z(A)$ set of connected components of the zero set
    \item $E$ Hermitian bundle
    \item $F$ connected component of critical point set, $U(F)$ localizing neighbourhood
    \item $\mathcal{P}_W(X)$ global Hilbert complex, $\mathcal{P}_W(X) = (L^2 \Omega^{0,\bullet}(X;E), \mathcal{D}_W(P), \overline{\partial}_E)$ twisted Dolbeault complex, $\mathcal{P}^\mu_W(X)$ 
    \item $\mathcal{P}_{\mathcal{B}}(U(F))$ local Hilbert complex
    \item $\mathcal{B} = (W,B)$, $W$ ideal boundary condition, $B$ local boundary condition
    \item $\mathcal{H}^q(\mathcal{P})$ cohomology group
    \item $\mu$ eigenvalue/weight of the infinitesimal $S^1$-action
    \item $M(\mathcal{Q}_\varepsilon(X), T_f,b)$ Morse polynomial
    \item $P(\mathcal{Q}_\varepsilon(X), T_f,b)$ Poincare polynomial
    \item $E(B)$ error polynomial  $Q_q^\mu$ coefficients of the error polynomial (nonnegative integer) 
    \item $T_{s,\theta}$ equivariant operator induced by the $S^1$-action
    \item $T_f$ geometric endomorphisms, also denoted by $T_g$, $f = f_g$ self map, $g \in S^1$ 
    \item $\mathcal{Q} = \Big( H, \{(\mathcal{D}(Q_j), Q_j)\}_{j=1}^N, \gamma \Big)$ SUSY Hilbert complex, $H$ $\mathbb{Z}_2$-graded Hilbert space, $\mathcal{D}(Q_j)$ domain, $Q_j$ SUSY operator, $\gamma$ chirality operator
    \item $\Delta = Q_j^2$ Hamiltonian, $\Delta_\epsilon = \Delta + \epsilon K + \epsilon^2 A$, $A$ potential associated to the deformation
    \item $Q_\pm = \frac{1}{\sqrt{2}}(Q_1\pm\sqrt{-1}Q_2) = P,P^*$
    \item $\mathcal{P} = (H,\mathcal{D}(P),P,\gamma)$
    \item $P_\varepsilon = M_\varepsilon^{-1} P M_\varepsilon$, $M_\varepsilon$ is a family of bounded, invertible operators on $H$
    \item $\mathcal{W}(\mathcal{Q}_\varepsilon(X)) = \Big(\mathcal{V}, \{(\mathcal{D}(Q_{j,\epsilon}|_\mathcal{V}), Q_{j,\epsilon})\}_{j=1}^N, \gamma \Big)$ Witten subcomplex, \\$\mathcal{V}^k = \bigoplus_{F\in Z(A)} \mathcal{H}^k_\varepsilon \Big(\mathcal{Q}_\varepsilon(U(F)) \Big)$
\end{enumerate}
\end{multicols}

\textbf{Section 2}
\begin{multicols}{2}
\begin{enumerate}
    \item $X^{reg}$ regular part, sometimes denoted as $\widehat{X}^{reg}$
    \item $\mathcal{S}(\widehat{S})$ union of strata, $Y_\alpha$ closed stratum, $\mathcal{T}_Y, \pi_Y, \rho_Y$ control data (tubular neighbourhood of $Y$, continuous retraction, radial function)
    \item $\mathcal{M}_1(X)$ set of boundary hypersurfaces
    \item $\mathfrak{B}_Y$ collective boundary hypersurface
    \item $\beta: X \to \widehat{X}$ blow down map
    \item $depth_X(Y)$ depth of the strata
    \item $H$ boundary hypersurface
    \item $\rho_Y$ collective boundary defining function, $\rho_X$ total boundary defining function
    \item $\mathcal{C}_\Phi^\infty(X)$ smooth function on $X$
    \item $g_w$ asymptotically $\delta$ wedge metric, $g_w - g_{w,p.t.} = x^{2\delta}E$ where $E$ is the error term
    \item $\mathcal{V}_w^*$ wedge one forms, $\prescript{w}{}{T^*X}$ wedge cotangent bundle
    \item $\mathcal{V}_e$ edge vector fields, $\prescript{e}{}{TX}$ edge tangent bundle
    \item $\operatorname{Diff}_w^*(X;E,F) = \rho_X^{-k} \operatorname{Diff}_e^*(X;E,F)$ wedge differential operators,\\ $\operatorname{Diff}_e^*(X;E,F)$ set of edge differential operators 
    \item $\sigma_k = \rho_{RC}^k \overline{\sigma_k}$ edge symbol, $\overline{\sigma_k}$ symbol map, $RC(\prescript{e}{}{TX})$ radial compactification of $\prescript{e}{}{TX}$
    \item $(g,\omega,J)$ wedge Kähler structure, $\alpha$ contact form
    \item $\mathcal{F}_\xi$ integrable foliation generated by the flow of the vector field $\xi$ on $Z$
    \item $R$ Reeb vector field 
    \item $h'\in C^0(\widehat{X})$, $h$ lift of $h'$, $\operatorname{crit}(h)$ critical point set
    \item $V_{F,k}$ Killing vector fields restrcted to each fiber $C(Z_j)$, $\rho_F^{2\delta_E} V_{F,E}$ error term
    \item $\mathcal{T}_{F^{reg}}$ 
    \item $U_t(\widehat{F}), U_t(F)$ singular localizing neighbourhood of $\widehat{F}, F$, $t \in (0,\infty]$ 
    \item $\partial U$ metric boundary
    \item $g_F$ model Kähler wedge metric on $F$
    \item $f_{\theta}(r_j,z_j)$ self map on tangent cone generated by the flow at time $\theta$
    \item $G$ compact connected Lie group, $\mathfrak{g}$ Lie algebra, $g\in G$
    \item $V$ Hamiltonian vector field 
\end{enumerate}
\end{multicols}

\textbf{Section 3}
\begin{multicols}{2}
\begin{enumerate}
    \item $cl$ Clifford multiplication
    \item $\mathcal{E} = \bigoplus_{q=0}^N \Lambda^q(\prescript{w}{}{TX}^{1,0} \otimes E)$ spinor bundle
\end{enumerate}
\end{multicols}

\textbf{Section 4}
\begin{multicols}{2}
\begin{enumerate}
    \item $\mathcal{D}_{\min}(D), \mathcal{D}_{\max}(D)$ minimal and maximal domain
    \item $\mathcal{GR}(D)$ Gelfand-Robbin quotient, \\$q_D: \mathcal{D}_{\max}(D) \to \mathcal{GR}(D)$ quotient map 
    \item $\mathcal{D}'_W(D)$ predomian 
    \item $\omega_D(u,v)$ Green-Stokes pairing
    \item $W^\perp$ annihilator w.r.t. $\omega_D$, $\mathcal{D}_{W^\perp}(D^*)$ domain of the adjoint operator
    \item $\mathcal{GR}(\mathcal{P})$ Gelfand-Robbin complex
    \item $\mathcal{H}(Q)$ cohomological quotient for a symmetric operator, $H^k(\mathcal{GR}(\mathcal{P}))$ cohomology of the Gefand-Robbin complex, $H_Y(P)$ cohomological boundary data along $Y$, $S$ section of the bundle $H_Y(Y)$
    \item $\mathcal{D}_{\max,W_1}(P)$ maximal domain of $P$ subject to the domain condition $W_1$ along the singular stratum $Y^1$, $\mathcal{D}_{\max,W_1}(P)_{C_x(Z^2)}$ maximal domain on the cone
    \item $U_{a,y;t} = \phi^{-1}(F(y))$ normal fiber over the point, where $\phi_F: U_t(F) \to F^{reg}$, $\partial U_{a,y;t}$ metric boundary
    \item $\mathcal{E}_{\mathcal{H}}, \mathcal{E}_{\mathcal{V}}$ horizontal and vertical Clifford modules
    \item $\Delta_\varepsilon^V$ vertical Laplacian
    \item $\mathcal{N}_{F,\varepsilon}$ normal cohomology bundle over $F^{reg}$
    \item $U_{a,y}^s \times U_{a,y}^u$ attracting and expanding directions of the flow of $h_a$
    \item $\mathcal{P}_1 \widehat{\otimes} \mathcal{P}_2$ tensor product of two Hilbert complexes
\end{enumerate}
\end{multicols}

\textbf{Section 5}
\begin{multicols}{2}
\begin{enumerate}
    \item $\mathscr{C}(\mathfrak{B}_Y)$ collar neighbourhood of $\mathfrak{B}_Y$
    \item $\nabla^M$ model connection w.r.t. $g_{w,p.t.}$
    \item $\{w_A=(z_i,y_\alpha)\}_{A=(i,\alpha)}$ for the submersion $\phi_Y:\mathfrak{B}_Y\to Y$ with $\{z_i \}_i$ the fiber coordinates on the model fiber $(Z^{reg}, g_Z)$ and $\{y_\alpha\}_\alpha$ the coordinates of $Y^{reg}$
    \item $x^2g_{\mathfrak{B}_Y/Y} +\phi^*_Y g_Y$ on the slice $\{x\}\times \mathfrak{B}_Y \xrightarrow{\phi_Y} Y$ horizontal and link wedge Fermi fields
    \item $p^{\mathcal V}, p^{\mathcal H}$ $g_{w,p.t.}$-orthogonal projections onto vertical and horizontal wedge-Fermi fields
    \item $O_w(x^k\partial^{\mathcal H}), O_w(x^k\partial^{\mathcal V})$ expressions of horizontal and link model wedge Fermi fields accordingly with coefficient components vanishing up to order $x^k$ when $x\to 0^+$
    \item $g_N, g_F$ normal metric and metric along $F$
    \item $g_0 = g_{w,p.t.}$, $g_1 = g_w$, $\nabla^0, \nabla^1$ connections 
    \item $[e_A^0] = [e_j^0, e_\alpha^0]_{j,\alpha}$, $e_1^0 = \partial_x$ $g_0$-orthonor-mal wedge frame, $[f_j^0 = xe_j^0]_{j\geq 2}$ $g_{\mathfrak{B}_Y/Y}$-orthonormal frame, $[e_\alpha^0]_\alpha$ horizontal lift of $g_Y$-orthonormal frame, $T$ transition matrix, similar notations for $g_1$
    \item $\prescript{w}{}{\mathcal{O}}, \prescript{w,pt}{}{\mathcal{O}}$ $SO(2n,\mathbb{R})$-subbundles of orthonormal wedge frames, $\tau$ bundle isomorphism between them
    \item $H: \prescript{w}{}{\mathcal{O}} \to \mathrm{GL}^+(\prescript{w}{}{T\mathscr{C} (\mathfrak{B}_Y)})$ principal bundle morphism
    \item $\prescript{w}{}{\mathcal{P}}, \prescript{w,pt}{}{\mathcal{P}}$: the principal $\operatorname{Spin}(2n)$-bundle lifting $\prescript{w}{}{\mathcal{O}}, \prescript{w,pt}{}{\mathcal{O}}$.
    \item $\mathcal{E}, \prescript{w,pt}{}{\mathcal{E}}$ general Clifford modules, $cl_w, cl_{w,p.t.}$ Clifford actions, $\tau: \mathcal{E} \to \prescript{w,pt}{}{\mathcal{E}}$
    \item $\nabla^{\mathcal{E}}, \nabla^{\prescript{w,pt}{}{\mathcal{E}}}$ spin$^c$-connections from $\nabla, \nabla^M$
    \item $D_X, D_M = D_{\mathcal{V}} + D_{\mathcal{H}}$ spin Dirac operators or Dirac-type operators
    \item $\mathcal{G}^+(\prescript{w}{}{T\mathscr{C}\mathfrak{B}_Y})$ principal bundle of positively oriented frames
    \item $H_{\mathcal{B}}^k(X;\mathcal{E}) = \mathcal{D}_{\mathcal{B}}(D^k)$ iteratively defined domain of $k$-th power of $D_X$.
    \item $\gamma_{t_0}^F$ cutoff function
    \item $R_\varepsilon$ first order differential operator on $\mathcal{N}_{F,\varepsilon}$
    \item $\Pi_\varepsilon R_\varepsilon$ projection onto the vertical harmonic subspace
    \item $\sigma_F$ section of $H^1_{\mathcal{B},\varepsilon}$
    \item $E = E_+ \oplus E_-$
    \item $A, B$ bounded operators in $\mathcal{L}(H_\mathcal{B}^0(X;E))$, $\mathcal{L}(H_\mathcal{B}^0(\mathcal{N}(F);E))$, $d(A,B)$ disctance
\end{enumerate}
\end{multicols}

\textbf{Section 6}
\begin{multicols}{2}
\begin{enumerate}
    \item $T$ maximal torus, $\mathfrak{t}$ Lie algebra, $\mathcal{L}$ integral lattice in $\mathfrak{t}$
    \item $\mathbb{Z}[\mathcal{L}^*]$  the ring of formal characters of $T$
    \item $\mathbb{Z}_E[b][\mathcal{L}^*]$ rings of formal power series with a shifted integral lattice, $\mathbb{Z}_E[b][\mathcal{L}^* \otimes \mathbb{C}]$ complexified integral lattice
    \item $\bigcap_{k=1}^l M_k$ intersection of polynomials
    \item $\mathcal{W} = N_G(T)/T$ Weyl group
    \item $X^T$ generically fixed points
    \item $\alpha_{a,j}$ linear functionals incorporate the constant rescaling
    \item $\mathcal{M}$ moment map
    \item $M_a(g)(b)$ local equivariant Morse polynomial
    \item $C$ action chamber 
    \item $P(g)(b)$ equivariant Poincaré polynomial 
    \item $M_\cap(b)$ combined upper bound 
    \item $R(b)$ reminder term
    \item $K_{W,\varepsilon,k}^{[0,c]}$ vector space generated by the eigenspace of $\Delta_\varepsilon$ with eigenvalues in $[0,c]$, $\mathcal{P}_{W, inst}$ associated subcomplex
    \item $\operatorname{Ind}$ index of an operator
    \item $\mathcal{P}^d_{\min}(X)$ the minimal de Rham complex,  $\prescript{p}{}{\mathcal{P}}_W(X)$ Dolbeault complex of $(p,\bullet)$-forms, $\prescript{\bullet}{}{\mathcal{P}}_{\min}(X)$ direct sum over all $p$ 
    \item $\mathcal{H}^{p,q}(\prescript{*}{}{\mathcal{P}}_W) = \mathcal{H}^{q}(\prescript{p}{}{\mathcal{P}}_W)$
    \item $\chi_{y,b}(\mathcal{P},T_g)$ Poincaré Hodge polynomial
    \item $\mathcal{HD}(\mathcal{P}^d_{\min})$ Hodge deviation of the complex
    \item $\tau_F(\theta)$ local contribution to the equivariant signature, $N_F$ charge, $T_F$ constant signature contribution
    \item $k_W(X)$ Kervaire semi-characteristic
\end{enumerate}
\end{multicols}

\textbf{Appendices}
\begin{multicols}{2}
\begin{enumerate}
    \item $\Sigma(\mathfrak{B}F)$ topological suspension of the boundary fibration
    \item $G_{\varepsilon,\mu}$ quasimodes
    \item $R(G)$ representation ring of the group $G$
\end{enumerate}
\end{multicols}
 \section{Stratified pseudomanifolds and wedge geometry}
\label{section_stratified_and_wedge}

In this section, we review stratified spaces and their resolutions, as well as wedge K\"ahler structures and the spin$^{\mathbb{C}}$ Dirac operator, before discussing local structures in neighbourhoods of fixed points of the circle actions that we study. In subsection \ref{subsection_stratified_spaces} we will go over the definition of stratified pseudomanifolds $\widehat{X}$ we work with, and the fibered manifolds $X$ which resolve them. In subsection \ref{subsubsection_wedge_metrics_related_structures} we introduce wedge geometry. A wedge structure is to a wedge metric, what a smooth structure is to a smooth Riemannian metric, and we identify it with a wedge tangent bundle on a fibered manifold. 

Wedge metrics are smooth Riemannian metrics on the regular part of $\widehat{X}$, which we also denote by $X^{reg}$, that extends to a non-degenerate bundle metric on the wedge tangent bundle. We will study such metrics, where the expansions of the metric in the boundary defining function may deviate from model metrics by $\mathcal{O}(\rho^{2\delta})$, made more precise in Subsubsection \ref{subsection_asymptotics_metric}.

We review the notions of wedge K\"ahler structures, wedge clifford structures here, before explaining precise asymptotics near the singularities for the wedge metrics that we consider. 
Most of this content is essentially the same as that in \cite{jayasinghe2024holomorphic} and is restated for the sake of completion. 

In Subsection \ref{subsection_fixed_point_sets} we describe the localizing neighbourhoods for fixed point sets of the group actions that we study, where we constrain ourselves in this article to those with circle actions that are locally K\"ahler Hamiltonian.
The model forms of the group actions are also described here, while their toric extensions are discussed in 
Subsection \ref{Subsection_Toric_actions}.

\subsection{Stratified pseudomanifolds and their resolutions.}
\label{subsection_stratified_spaces}

We review some of the background for stratified spaces and their resolutions, referring the reader to \cite{Albin_2017_index,jayasinghe2023l2,jayasinghe2024holomorphic} for more details. There are other notions of stratified spaces, and we note that the smooth atlas stratified spaces in \cite{albin2025smoothatlasstratifiedspaces} are also equivalent to the spaces studied here.

\subsubsection{Stratified pseudomanifolds}

All topological spaces we consider will be Hausdorff, locally compact topological spaces with a countable basis for their topology. Recall that a subset $W$ of a topological space $V$ is locally closed if every point $a \in W$ has a neighborhood $\mathcal{U}$ in $V$ such that $W \cap \mathcal{U}$ is closed in $\mathcal{U}$. A collection ${S}$ of subsets of $V$ is locally finite if every point $v \in V$ has a neighborhood that intersects only finitely many sets in ${S}$.

\begin{definition}
    \label{Thom-Mather stratified space}    
A Thom-Mather stratified space $\widehat{X}$ is a metrizable, locally compact, second countable space which admits a locally finite decomposition into a union of locally closed strata $\mathcal{S}(\widehat{X})=\left\{Y_{\alpha}\right\}$, where each $Y_{\alpha}$ is a smooth manifold, with dimension depending on the index $\alpha$. We require the following:
\begin{enumerate}
    \item If $Y_{\alpha}, Y_{\beta} \in \mathcal{S}(X)$ and $Y_{\alpha} \cap \overline{Y_{\beta}} \neq \emptyset$, then $Y_{\alpha} \subset \overline{Y_{\beta}}$.

    \item Each stratum $Y$ is endowed with a set of `control data' $\mathcal{T}_Y, \pi_{Y}$ and $\rho_{Y}$; here $\mathcal{T}_Y$ is a neighbourhood of $Y$ in $X$ which retracts onto $Y, \pi_{Y}: \mathcal{T}_Y \longrightarrow Y$ is a fixed continuous retraction and $\rho_{Y}: \mathcal{T}_Y \rightarrow[0,2)$ is a `radial function' in this tubular neighbourhood such that $\rho_{Y}^{-1}(0)=Y$. Furthermore, we require that if $\widetilde{Y} \in \mathcal{S}(X)$ and $\widetilde{Y} \cap \mathcal{T}_Y \neq \emptyset$, then
\begin{equation}    
    \left(\pi_{Y}, \rho_{Y}\right): \mathcal{T}_Y \cap \widetilde{Y}  \setminus Y \longrightarrow Y \times (0,2)
\end{equation}
is a differentiable submersion.

\item If $W, Y, \widetilde{Y} \in \mathcal{S}(X)$, and if $p \in \mathcal{T}_Y \cap \mathcal{T}_{\widetilde{Y}} \cap W$ and $\pi_{\widetilde{Y}}(p) \in \mathcal{T}_Y \cap \widetilde{Y}$, then $\pi_{Y}\left(\pi_{\widetilde{Y}}(p)\right)=\pi_{Y}(p)$ and $\rho_{Y}\left(\pi_{\widetilde{Y}}(p)\right)=\rho_{Y}(p)$.
\end{enumerate}
\end{definition}

From this definition, and due to Thom's first isotopy lemma, we have structure of neighborhoods as locally trivial fibrations over strata $Y$, and the following are some direct consequences.

\begin{enumerate}
\item If $Y, \widetilde{Y} \in \mathcal{S}(X)$, then
$$
\begin{aligned}
Y \cap \overline{\widetilde{Y}} \neq \emptyset & \Leftrightarrow \mathcal{T}_Y \cap \widetilde{Y} \neq \emptyset, \\
\mathcal{T}_Y \cap \mathcal{T}_{\widetilde{Y}} \neq \emptyset & \Leftrightarrow Y \subset \overline{\widetilde{Y}}, Y=\widetilde{Y} \quad \text {or } \widetilde{Y} \subset \overline{Y} .
\end{aligned}
$$

\item For each $Y \in \mathcal{S}(X)$, the restriction $\pi_{Y}: \mathcal{T}_Y \rightarrow Y$ is a locally trivial fibration with fibre the cone $C\left(Z_{Y}\right)$ over some other stratified space $Z_{Y}$ (called the \textbf{\textit{link}} over $Y)$, with atlas $\mathcal{U}_{Y}=\{(\phi, \mathcal{U})\}$ where each $\phi$ is a trivialization 
\begin{equation}
\label{equation_chartlike_map}
\pi_{Y}^{-1}(\mathcal{U}) \rightarrow \mathcal{U} \times C\left(Z_{Y}\right),    
\end{equation}
and the transition functions are stratified isomorphisms of $C\left(Z_{Y}\right)$ which preserve the rays of each conic fibre as well as the radial variable $\rho_{Y}$ itself, hence are suspensions of isomorphisms of each link $Z_{Y}$ which vary smoothly with the variable $y \in \mathcal{U}$.
\end{enumerate}

This structure shows that there is a filtration of $\widehat{X}$ of the form
\begin{equation}
    \widehat{X}_n \supseteq  \widehat{X}_{n-1} \supseteq  \widehat{X}_{n-2}... \supseteq \widehat{X}_{0}
\end{equation}
where $\widehat{X}_{j}$ is the union of strata $Y_{\alpha} \subset S(\widehat{X})$ of dimension less than or equal to $j$. It is clear that $\widehat{X}_j \backslash \widehat{X}_{j-1}$ is a smooth manifold of dimension $j$. We define the \textit{\textbf{regular part}} of $\widehat{X}$ to be $X^{reg}:=\widehat{X}_n \setminus \widehat{X}_{n-1}$. 

If in addition if $\widehat{X}_{n-1}=\widehat{X}_{n-2}$ and $\widehat{X} \backslash \widehat{X}_{n-2}$ is dense in $\widehat{X}$, 
then we say that $\widehat{X}$ is a \textbf{stratified pseudomanifold}.

There is a functorial equivalence between Thom-Mather stratified spaces and manifolds with corners and iterated fibration structures (see Proposition 2.5 of \cite{Albin_signature}, Theorem 6.3 of \cite{Albin_hodge_theory_cheeger_spaces}).

\subsubsection{Manifolds with corners with iterated fibration structures}

In section 1 of \cite{Albin_2017_index}, there is a detailed description of iterated fibration structures on manifolds with corners. We explain some of these structures from said article which we will use and refer the reader to the source for more details. In \cite{kottke2022products}, these are referred to as manifolds with fibered corners that are interior maximal.

An $n$-dimensional manifold with corners $X$ is an $n$-dimensional topological manifold with boundary, with a smooth atlas modeled on $(\mathbb{R}^+)^n$ whose boundary hypersurfaces are embedded. We denote the set of boundary hypersurfaces of $X$ by $\mathcal{M}_1(X)$. A collective boundary hypersurface refers to a finite union of non-intersecting boundary hypersurfaces.

\begin{definition}
\label{iterated_fibration_structure}
An \textit{\textbf{iterated fibration structure}} on a manifold with corners $X$ consists of a collection of fiber bundles
\begin{center}
    $Z_Y -\mathfrak{B}_Y \xrightarrow{\phi_Y} Y$
\end{center}
where $\mathfrak{B}_Y$ is a collective boundary hypersurface of $X$ with base and fiber manifolds with corners such that:

\begin{enumerate}
    \item Each boundary hypersurface of $X$ occurs in exactly one collective boundary hypersurface $\mathfrak{B}_Y$.
    \item If $\mathfrak{B}_Y$ and $\mathfrak{B}_{\widetilde{Y}}$ intersect, then dim $Y \neq$ dim $\widetilde{Y}$, and we write $Y < \widetilde{Y}$ if dim $Y < \text{dim} \widetilde{Y}$.
    \item If $Y < \widetilde{Y}$, then $\widetilde{Y}$ has a collective boundary hypersurface $\mathfrak{B}_{Y\widetilde{Y}}$ participating in a fiber bundle $\phi_{Y\widetilde{Y}} : \mathfrak{B}_{Y\widetilde{Y}} \rightarrow Y$ such that the diagram 
    
\[\begin{tikzcd}
	{\mathfrak{B}_Y \cap \mathfrak{B}_{\widetilde{Y}}} && {\mathfrak{B}_{Y\widetilde{Y}} \subseteq \widetilde{Y}} \\
	& Y
	\arrow["{\phi_{\widetilde{Y}}}", from=1-1, to=1-3]
	\arrow["{\phi_Y}"', from=1-1, to=2-2]
	\arrow["{\phi_{Y\widetilde{Y}}}", from=1-3, to=2-2]
\end{tikzcd}\]
commutes.
\end{enumerate}
\end{definition}

The base can be assumed to be connected but the fibers are in general disconnected.
As we mentioned above, there is an equivalence between Thom-Mather stratified spaces and manifolds with corners with iterated fibration structures. 

If we view a cone over a link $Z$ as the quotient space of $[0,1]_x \times Z$ under the identification where the points of the link at $\{x=0\} \times Z$ are identified, then the quotient map is a blow-down map.
More generally, there is an inductive desingularization procedure which replaces the Thom-Mather stratified space with a manifold with corners with iterated fibration structures.
This corresponds to a \textit{\textbf{blow-down}} map $\beta : X \rightarrow \widehat{X}$, which satisfies properties given in Proposition 2.5 of \cite{Albin_signature}.
See Remark 3.3 of \cite{Albin_2017_index} for an instructive toy example.

Under this equivalence, the bases of the boundary fibrations correspond to the different strata, which we shall denote by
\begin{equation}
    \mathcal{S}(X) = \{Y: Y \textit{ is the base of a boundary fibration of X} \}.
\end{equation}
The bases and fibers of the boundary fiber bundles are manifolds with corners with iterated fibration structures (see for instance Lemma 3.4 of \cite{albin2010resolution}). The condition dim $Z_Y > 0$ for all $Y$ corresponds to the category of pseudomanifolds within the larger category of stratified spaces. The partial order on $\mathcal{S}(X)$ gives us a notion of depth
\begin{center}
$ depth_X(Y) = max \{ n \in \mathbb{N}_0 : \exists Y_i \in \mathcal{S}(X)$ s.t. $Y=Y_0 <Y_1 < ... <Y_n \}$.
\end{center}
The depth of $X$ is then the maximum of the integers $depth_X(Y)$ over $Y \in \mathcal{S}(X)$. 

We now introduce some auxiliary structures associated to manifolds with corners with iterated fibration structures.
If $H$ is a boundary hypersurface of $X$, then because it is assumed to be embedded, there is a smooth non-negative function $\rho_H$ such that 
$\rho_{H}^{-1}(0)=H$ and $d\rho_H$  does not vanish at any point on $H$. We call any such function a \textbf{\textit{boundary defining function for $H$}}.
For each $Y \in \mathcal{S}(X)$, we denote a \textbf{\textit{collective boundary defining function}} by
\begin{center}
    $\rho_Y = \prod_{H \in \mathfrak{B}_Y} \rho_H$,
    and by
    $\rho_X = \prod_{H \in \mathcal{M}_1(X)} \rho_H$
\end{center}
 a \textbf{\textit{total boundary defining function}}, where $\mathcal{M}_1(X)$ denotes the set of boundary hypersurfaces of $X$.

When describing the natural analogues of objects in differential geometry on singular spaces, the iterated fibration structure comes into play. For example,
\begin{equation}
\label{equation_smooth_functions_on_stratified_spaces}
\mathcal{C}_{\Phi}^{\infty}(X)= \{ f \in \mathcal{C}^{\infty}(X) : f \big|_{\mathfrak{B}_Y} \in 
\phi_{Y}^*\mathcal{C}^{\infty}(Y) \text{ for all } Y \in \mathcal{S}(X) \}
\end{equation}
corresponds to the smooth functions on $X$ that are continuous on the underlying
stratified space. 

\subsection{Wedge metrics and wedge geometry}
\label{subsubsection_wedge_metrics_related_structures}

On Thom-Mather stratified pseudomanifolds we can define metrics which are \textit{locally conic}. For instance, consider the model space $\widehat{X} = \mathbb{R}^k \times {Z}^+$, where ${Z}^+$ is a cone over a smooth manifold $Z$. The resolved manifold with corners that corresponds to the blowup of this space is $\mathbb{R}^k \times [0,\infty)_x \times Z$. We have the model wedge metric
\begin{equation}
\label{equation_definition_homogeneous_metric_cone}
    g_{w} = g_{\mathbb{R}^k} + dx^2 + x^2 g_Z,
\end{equation}
which is a product metric on the product space and which we call a \textbf{\textit{rigid/ product type}} wedge metric, degenerating as it approaches the stratum at $x=0$. These metrics are degenerate as bundle metrics on the tangent bundle but we can introduce a rescaled bundle on which they are non-degenerate. 

Formally, we proceed as follows. Let $X$ be a manifold with corners and iterated fibration structure. Consider the ‘wedge one-forms’
\begin{center}
    $\mathcal{V}^{*}_{w}= \{ \omega \in \mathcal{C}^{\infty}(X; T^*X) : \text{ for each } Y \in \mathcal{S}(X), \text{  }i^*_{\mathfrak{B}_Y} \omega (V) =0 \text{ for all } V \in \text{ker } D\phi_Y \}$.
\end{center}
We can identify $\mathcal{V}^{*}_{w}$ with the space of sections of a vector bundle which we call $\prescript{{w}}{}{T}^*X$, the wedge cotangent bundle of $X$, together with a map
\begin{equation}
    i_{w} : \prescript{{w}}{}{T}^*X \rightarrow T^*X
\end{equation}
that is an isomorphism over $X^{reg}$ such that,
\begin{center}
    $(i_{w})_*\mathcal{C}^{\infty}(X;\prescript{{w}}{}{T}^*X) =  \mathcal{V}^{*}_{w} \subseteq \mathcal{C}^{\infty}(X; T^*X)$.
\end{center}
In local coordinates near the collective boundary hypersurface, the wedge cotangent bundle is spanned by
\begin{center}
    $ dx$, $xdz$, $dy$
\end{center}
where $x$ is a boundary defining function for $\mathfrak{B}_Y$, $dz$ represents covectors along the fibers and $dy$ covectors along the base.
The dual bundle to the wedge cotangent bundle is known as the wedge tangent bundle, $\prescript{{w}}{}{T}X$. It is locally spanned by
\begin{center}
    $\partial_x$, $\frac{1}{x}\partial_z$, $\partial_y$
\end{center}
A \textit{\textbf{wedge metric}} is simply a bundle metric on the wedge tangent bundle.

The notion of a \textit{\textbf{wedge differential operator}} 
$P$ of order $k$ acting on sections of a vector bundle $E$, taking them to sections of a vector bundle $F$ is described on page 11 of \cite{Albin_2017_index}. We first define the \textit{\textbf{edge vector fields on X}} by
\begin{center}
    $\mathcal{V}_e = \{ V \in  \mathcal{C}^{\infty}(X;TX): V \big|_{\mathfrak{B}_Y}$ is tangent to the fibers of $\phi_Y$ for all $Y \in \mathcal{S}(X) \}$.
\end{center}
There is a rescaled vector bundle that is called the \textit{\textbf{edge tangent bundle}} $\prescript{e}{}TX$ together with a natural vector bundle map $i_e : \prescript{e}{}TX \rightarrow TX$ that is an isomorphism over the interior and satisfies
\begin{center}
    $(i_e)_* \mathcal{C}^{\infty}(X;\prescript{e}{}TX)=\mathcal{V}_e$.
\end{center}
In local coordinates near a point in $\mathfrak{B}_Y, (x, y, z)$,
a local frame for $\prescript{e}{}TX$ is given by 
\begin{center}
    $ x\partial_x, x\partial_y, \partial_z$
\end{center}
Note that the vector fields $x\partial_x$ and $x\partial_y$ are degenerate as sections of $TX$, but not as sections of $\prescript{e}{}TX$. 

The universal enveloping algebra of $\mathcal{V}_e$ is the ring $\text{Diff}^*_e(X)$ of edge differential operators. That is, these are the differential operators on $X$ that can be expressed locally as finite sums of products of elements of $\mathcal{V}_e$. They have the usual notion of degree and extension to sections of vector bundles, as well as an edge symbol map defined on the edge cotangent bundle (see \cite{Mazzeo_Edge_Elliptic_1,Albin_signature,Albin_hodge_theory_cheeger_spaces}).
Similarly, $\text{Diff}^*_e(X; E, F)$ denotes the edge differential operators acting on sections of a vector bundle $E$ and taking them to a sections on a vector bundle $F$.
The edge symbol (which can be defined on the space of edge pseudo-differential operators in general) is used to define the notion of ellipticity used in this article. We follow \cite{Albin_2017_index} and define the map
\begin{equation}
    \overline{\sigma_k} : \text{Diff}^k_e(X;E,F) \rightarrow \rho_{RC}^{-k}\mathcal{C}^{\infty}(RC(\prescript{e}{}T^*X), \pi^*hom(E,F))
\end{equation}
which is the usual symbol map where $RC(\prescript{e}{}T^*X)$ denotes the radial compactification of the edge cotangent bundle and $\pi: \prescript{e}{}T^*X \rightarrow X$ is the projection map. We denote by $\rho_{RC}$ a boundary defining function for the
boundary at radial infinity. Multiplying $\overline{\sigma_k}$ by $\rho_{RC}^k$ 
the resulting map $\sigma_k$ is called the \textbf{\textit{edge symbol}} (see Section 3.3 of \cite{Albin_2017_index}).
This fits into the following short exact sequence
\begin{equation}
\label{edge symbol sequence}
    0 \rightarrow \text{Diff}^{k-1}_e(X;E,F) \rightarrow \text{Diff}^{k}_e(X;E,F) \xrightarrow{\sigma} \mathcal{C}^{\infty}(^e\mathbb{S}^*X,\pi^*hom(E,F)) \rightarrow 0.
\end{equation}
If the edge symbol of a differential operator is invertible away from the zero section we call such an operator  \textbf{\textit{edge elliptic}}.
We define \textit{\textbf{wedge differential operators}} by 
\begin{equation*}
    \text{Diff}^k_w(X;E,F) = \rho^{-k}_X \text{Diff}^k_e(X;E,F)
\end{equation*}
following, e.g., \cite{Albin_2017_index}, where $\rho_X$ is a total boundary defining function for $X$. 
If a wedge differential operator $D_k$ of order $k$ can be written as $\rho_X^{-k}A_k$ where $A_k$ is an edge elliptic differential operator of order $k$, then $D_k$ is said to be a  \textbf{\textit{wedge elliptic}} operator.

\subsubsection{Asymptotically wedge metrics and K\"ahler structures}
\label{subsection_asymptotics_metric}

An \textit{\textbf{asymptotically $\delta$ wedge metric}} is defined inductively. On depth zero spaces, which are just smooth manifolds, a wedge metric is a Riemannian metric. Assuming we have defined asymptotically $\delta$ wedge metrics on spaces of depth less than $k$, let us assume $X$ has depth $k$.
A wedge metric $g_w$ on $X$ is an 
asymptotically $\delta$ wedge metric if, for every $Y \in \mathcal{S}(X)$ of depth $k$ there is a collar neighbourhood $\mathscr{C}(\mathfrak{B}_Y) \cong [0,1)_x \times \mathfrak{B}_Y$ of $\mathfrak{B}_Y$ in $X$, a metric $g_{w, pt}$ of the form
\begin{equation}
\label{equation_wedge_metric_pt}
    g_{w, pt} = dx^2 +x^2 g_{\mathfrak{B}_Y/Y} +\phi^*_Yg_Y
\end{equation}
where $g_Y$ is an asymptotically $\delta$ wedge metric on $Y$, $g_{\mathfrak{B}_Y/Y}+\phi^*_Yg_Y$ is a submersion metric for $\mathfrak{B}_Y \xrightarrow{\phi_Y} Y$ and $g_{\mathfrak{B}_Y/Y}$ restricts to each fiber of $\phi_Y$ to be an asymptotically $\delta$ wedge metric on $Z_Y$ and
\begin{equation}
\label{equation_structure of the metric}
    g_{w} - g_{w, pt} = x^{2\delta}E\in x^{2\delta} \mathcal{C}^{\infty} \big(\mathscr{C}(\mathfrak{B}_Y); S^2(^{w}T^*X)\big).
\end{equation}
When $\delta=1/2$, we further assume that the error term $E\in \mathcal{C}^{\infty} \big(\mathscr{C}(\mathfrak{B}_Y); S^2(^{w}T^*X)\big)$ satisfies
\begin{equation}
\label{wedge_metric_assumption}
    E(\partial_x, \partial_x)|_{x=0}= E(\partial_x, V_1)|_{x=0}= E(\partial_x, V_2)|_{x=0}=0, 
\end{equation}
for every $V_1\in \mathcal{C}^{\infty} \big(\mathscr{C}(\mathfrak{B}_Y); \ker d\phi_Y\big)$ and every $V_2\in\mathcal{C}^{\infty} \big(\mathscr{C}(\mathfrak{B}_Y); \phi_Y^*(^{w} TY)\big)$. Assumptions \eqref{wedge_metric_assumption} hold true when the stratum $Y$ and the wedge metric are both smooth, as explained in Remark~\ref{rem:smooth_case_analogy}. 

If at each step  $g_{w} = g_{w, pt}$, we say $g_w$ is a \textit{\textbf{rigid}} or \textit{\textbf{product-type wedge metric}}. Off of these collar neighborhoods of the stratum of depth $k$, the form of the metric is fixed by the induction as a $\delta$ wedge metric of lower depth. 

\begin{remark}
\label{remark_totally_geodesic_wedge}
In \cite{Albin_2017_index} the case of $2\delta=1$ (called \textit{\textbf{exact wedge metrics}}) and $2\delta=2$ (called totally geodesic wedge metrics) were studied, the distinction only relevant in Getzler rescaling arguments at the strata.
In this article we required $\delta >1/2$ for the estimates in Subsection \ref{subsection_polyhomogenous_estimates}, hence we make Assumption \ref{assumption_on_metrics_and_domains}.
\end{remark}

Given a stratified pseudomanifold with a wedge metric, a complex structure on the wedge tangent bundle is called a \textit{\textbf{wedge complex structure}}.

\begin{definition}[wedge K\"ahler structures]
\label{definition_wedge_kahler_structure}
    If there is a symplectic form $\omega$ on $X^{reg}$ extending to a non-degenerate closed wedge two form such that it is tamed by a wedge complex structure $J$, yielding a wedge metric $g=\omega(J \cdot, \cdot)$ on $X$, then we say that $(g,\omega,J)$ is a \textit{\textbf{wedge K\"ahler structure}} on $X$.
\end{definition}

As in the smooth case, any two of the three constituents of a triple $(g,\omega,J)$ determines the third.
Projective algebraic varieties with wedge metrics such as conifolds and even the cusp curve, as we studied in \cite[\S 7.3]{jayasinghe2023l2} are examples of such spaces.

We briefly review some background on conic K\"ahler metrics, referring to \cite[\S 7.1]{jayasinghe2023l2,boyer&galicki,dragomir2007differential,blair2010riemannian} for more details on complex structures on cones and their CR structures, as well as Sasaki structures of K\"ahler cones.

It is well known that a complex structure on a cone $(C(Z), dx^2 +x^2 g_Z)$ induces a CR structure on the link $(Z,g_Z)$, which in particular yields an almost contact metric structure on the link if it is smooth (see for instance Section 1.1.4 of \cite{dragomir2007differential}). In the case where the CR structure is pseudoconvex, the almost complex structure is a contact structure. 
A Sasaki structure on $M$
is equivalent to a K\"ahler structure on the metric cone over $M$.
In this case, the CR structure is pseudoconvex and has a contact structure. Taking the quotient by the action of the Reeb vector field $\xi$, one gets (in general) an orbifold $\Sigma = Z / \xi$ which has a K\"ahler structure, usually known as the transversal K\"ahler structure on the Sasaki manifold. The Reeb foliation $\mathcal{F}_\xi$ on $Z$ happens to be a taut Riemannian foliation, and the contact distribution has a splitting into holomorphic and anti holomorphic components. 
The K\"ahler form can be written as 
\begin{equation}
    2dx \wedge x\alpha + x^2 d\alpha
\end{equation}
where $\alpha$ is a contact form. Here $d\alpha$ is a transversal K\"ahler form that is non-degenerate on the contact distribution $TZ/T\mathcal{F}_{\xi}$.
The K\"ahler form of a disc $2dr \wedge rd\theta$ corresponds to a \textit{trivial} Sasaki structure on the circle.

\begin{definition}[wedge Sasaki structure]
    Given a stratified pseudomanifold $\widehat{Z}$ with a wedge metric, if the resolution of the cone $C_x(\widehat{Z})$ is a wedge K\"ahler metric, we say that the resolution $Z$ of $\widehat{Z}$ has a wedge Sasaki structure, and $V:=xJ(\partial_x)$ is called the Reeb vector field.    
\end{definition}
These definitions generalize the usual notion of wedge Sasaki spaces to our setting in the obvious manner.

A wedge symplectic form is a closed wedge 2 form, that is a smooth section $\omega$ of $\Lambda^2(^{w}T^*X)$ which satisfies $d\omega=0$.
On a collar neighbourhood 
\begin{equation}
\label{equation_collar_neighbourhood_asymptotics}
    \mathscr{C}(\mathfrak{B}_Y) \cong [0,1)_x \times \mathfrak{B}_Y
\end{equation}
of $\mathfrak{B}_Y$ in $X$, a closed wedge 2 form of the form 
\begin{equation}
    \omega_{w,pt}=  2dx \wedge x\alpha +x^2 d\alpha +\phi^* \omega_Y
\end{equation}
is called a \textbf{\textit{product type wedge symplectic structure}}, where $\alpha$ is a contact form for the link $Z$ and $\omega_Y$ is a 2 form on $Y$.
In this article we consider \textbf{\textit{asymptotically $\delta$ wedge symplectic structures}} $\omega$ which satisfy
\begin{equation}
\label{equation_structure_of_the_symplectic_form}
    \omega - \omega_{w, pt} \in x^{2\delta} \mathcal{C}^{\infty} \big(\mathscr{C}(\mathfrak{B}_Y); \Lambda^2(^{w}T^*X)\big).
\end{equation}
for some $\delta>0$ in a collar neighbourhood as described in equation \eqref{equation_collar_neighbourhood_asymptotics}.

An integrable almost complex structure on the wedge tangent bundle is called a wedge complex structure.
Given a triple $(g_{w,pt},J_{w,pt},\omega_{w,pt})$ where $g_{w,pt}$ is a product wedge metric and $\omega_{w,pt}$ is a wedge symplectic structure, and $J_{w,pt}$ is a wedge complex structure such that $g_{w,pt}(J_{w,pt} \cdot, \cdot) = \omega_{w,pt} (\cdot, \cdot)$, we say that $J_{w,pt}$ is a \textit{\textbf{product type wedge complex structure}}. 

Using a connection in a collar neighbourhood as described in equation \eqref{equation_collar_neighbourhood_asymptotics} where we can extend the fibration $\phi$ and fix a splitting, we can write this in coordinates as  
\begin{equation}
    J_{w,pt}=\partial_x \otimes x\alpha - \frac{1}{x} \xi \otimes dx + J_{TZ/T\mathcal{F}_{\xi}} + J_Y.
\end{equation}
Here $x\alpha$ is a unit wedge covector field on $C(Z)$, $\frac{1}{x} \xi$ is the dual wedge vector field with respect to a product type metric $g_{w,pt}$ on the collar neighbourhood, $J_{TZ/T\mathcal{F}_{\xi}}$ is an almost complex structure on the distribution $TZ/T\mathcal{F}_{\xi}$ where $\mathcal{F}_{\xi}$ is the integrable foliation generated by the flow of the vector field $\xi$ on $Z$, and where $J_Y$ is a complex structure on $TY$ restricted to the stratum.
Given a metric product of $\mathbb{C}^k \times C(Z)$ where $Z$ has a Sasaki structure, the corresponding K\"ahler structure is a product type wedge complex structure, where $\alpha$ is a contact form and $\xi$ is the Reeb vector field for the Sasaki structure and $J_{TZ/T\mathcal{F}_{\xi}}$ is an almost complex structure on the contact distribution. In the next subsection we shall study the case of K\"ahler Hamiltonian actions on such spaces.

In this article we consider wedge complex structures $J$ where
\begin{equation}
\label{equation_structure_of_the_complex_structure}
    J - J_{w, pt} \in x^{2\delta} \mathcal{C}^{\infty} \big(\mathscr{C}(\mathfrak{B}_Y); End(^{w}TX)\big)
\end{equation}
for some $\delta>0$, and we call them \textbf{\textit{asymptotically $\delta$ wedge complex structures}}.

\subsection{Fixed point sets and localizing neighbourhoods}
\label{subsection_fixed_point_sets}

In \cite{jayasinghe2024witten}, the \textbf{fundamental neighbourhoods} of critical point sets of stratified Morse-Bott functions were studied, and in this article we choose to rename them as \textbf{localizing neighbourhoods}.
We study the localizing neighbourhoods of Killing vector fields, and consider specializations to the Hamiltonian and K\"ahler Hamiltonian cases. 

Localizing neighbourhoods for other types of vector fields, such as non-vanishing vector fields and those transverse to such non-vanishing vector fields in studies of mod 2 invariants (see Subsection \ref{subsection_Kervaire}) deserve a more detailed study, and we do not attempt that in this article. However, the localizing neighbourhoods in this article and in \cite{jayasinghe2024witten} give many examples of localizing neighbourhoods for these other complexes as well, while showing that the principle of localization and Witten deformation continues to work in those settings.

In the smooth setting, it is known that K\"ahler Hamiltonians are Morse-Bott in the sense of Kirwan \cite{holm2016morse} admitting a slightly different local description than that given by the usual Morse-Bott lemma in local coordinates. However we will define our localizing neighbourhoods in terms of the group actions, instead of the Hamiltonian functions, referring to results in the literature that show that such functions exist (see Remark \ref{Remark_Hamiltonian_functions_exist})

We briefly study the form of the Hamiltonian functions to motivate the Definitions we present and the assumptions we make.

\begin{definition}[Hamiltonian functions and vector fields]
\label{Definition_wedge_Kahler_Hamiltonian_function}
If $X$ is a resolved stratified pseudomanifold equipped with a wedge symplectic form $\omega$ and $V$ is a wedge vector field, we say it is the Hamiltonian vector field if there exists a function $h \in C^{\infty}_{\Phi}(X)$ such that $\iota_{V}\omega=-dh$, in which case we say $h$ is the \textbf{Hamiltonian} function.

If there is a wedge K\"ahler structure $(g,J,\omega)$, if in addition $L_VJ=0$, we say call them K\"ahler Hamiltonian (the vector field and the function).    
\end{definition}

Consider a cone $C_x(Z)$ where $Z$ is a resolution of a link, determined by a Reeb vector field $\xi$ on $Z$. We can find a contact form $\alpha$ restricted to $Z^{reg}$ satisfying the usual properties 
\begin{equation}
    \iota_\xi \alpha=1, \quad \iota_\xi d\alpha=0
\end{equation}
and which extends to a wedge contact form on $Z$. 

We refer the reader to Section 2 of \cite{Martellisparksyau08} for a discussion on Hamiltonian Killing vector fields on K\"ahler cones. Since the metric is preserved, the radial distances are preserved and such vector fields are necessarily Killing vector fields on the link of every cone. Since they also preserve the complex structure $J$ and $x\partial_x$, they preserve the Reeb vector field. Then it is easy to see they preserve the transverse K\"ahler structure. 

This shows that a Hamiltonian (even non-K\"ahler) on a cone can be written as $h=(1/2)x^2f \alpha$, for some function $f$ and the contact form $\alpha$ on the link, and $dh=dx \wedge f\alpha +(x^2/2)(df \wedge \alpha+ f d\alpha)$. 

Thus 
\begin{equation}
\label{equation_Hamiltonian_relevant_equation_1}
    dh=xf dx+(x^2/2)df=-\iota_X (dx \wedge x \alpha +(x^2/2)d\alpha)
\end{equation}
which shows that $f=-\alpha(X)$ and $df=-(d\alpha)(X)$.
In the K\"ahler case, there are extra conditions on $f$ since the transverse K\"ahler structure will be preserved.

\subsubsection{Vanishing sets for potentials}

We begin with a quick review of the notions of stratified critical point sets defined in \cite{jayasinghe2024witten}, refering to \textit{loc. cit.} for more details and the philosophy behind the definitions. Here we rename them as vanishing sets, as we generalize them further.

\textbf{In the following three definitions}, we shall consider  a stratified pseudomanifold $\widehat{X}$ with a wedge metric $g_w$ and that there is a Killing vector field $V$, which determines a potential $A=|V|_{g_w}^2$. We assume that the vanishing set $Z(A)$ has only has finitely many connected components.

These definitions are related to the stratified Morse-Bott critical point sets of types Smooth, I and II in \cite{jayasinghe2024witten} in the following way. When the Killing vector field $V$ is \textbf{wedge Hamiltonian}, then
there is a continuous function $h' \in C^0(\widehat{X})$ which lifts to a map $h \in C_{\Phi}^{\infty}(X)$ (see \eqref{equation_smooth_functions_on_stratified_spaces}) on the resolved manifold with corners, i.e., $h=h' \circ \beta$. We write $\text{crit}(h):=\beta(|dh|_{g_w}^{-1}\{0\})$.
Then $A:=|dh|^2$ for a Morse-Bott type function $h$, satisfying $dh=-\iota_V\omega$ for the wedge symplectic form $\omega$. 

We note that despite the appearance of the metric $g_w$ in the following definitions, the condition of being critical, i.e. $\widehat{F}\subset \beta(|dh|_{g_w}^{-1}(0))$, is independent of the choice of wedge metric $g_w$. In fact we only need that the differential of $h$ vanishes along $\overline{\beta^{-1}\left(\widehat{F}\right)}$.

\begin{definition}[\textbf{Smooth Morse-Bott type vanishing set}]
    If $F \subset \widehat{X}$ is a connected component of $Z(A)$, which is wholly a subset of $X^{reg}$, then we say it is a smooth non-degenerate vanishing set if $F$ is a smooth submanifold of $X^{reg}$ and the $V$ is a Killing vector field on the normal fibers.

    In the Hamiltonian case, the Hessian of $h$ is non-degenerate in the directions of the normal bundle of $F$ in $X^{reg}$.
\end{definition}

\begin{definition}[\textbf{Type I stratified vanishing set}]
    If $F\subset \widehat{X}$ is a connected component of $Z(A)$ which is disjoint from $\widehat{X}^{reg}$, then we say that it is a Type I stratified vanishing set if $F$ is a smooth submanifold of $Y_\alpha\in \mathcal{S}(X)$ where $F \hookrightarrow Y_\alpha$ is a smooth embedding. 
        
    Since $\pi_\alpha: \mathcal{T}_{Y_\alpha}\rightarrow Y_\alpha$ is a locally trivial fibration with fiber $C(\widehat{Z})$, (the cone over a Thom-Mather stratified space $\widehat{Z}$), we conclude that $F$ has an open neighborhood $U(\widehat{F})$ in $\widehat{X}$ which arises as a fiber bundle 
   \begin{equation}
       \pi_{F}: U(\widehat{{F}}) \to F,
   \end{equation}
    with fibers $C(\mathbb{S}^{\ell})\times C(\widehat{Z}_a)$, where $\ell+1>0$ is the codimension of the smooth embedding $F \hookrightarrow Y_\alpha$.
    
    On each fiber we define $\rho_N$ to be the distance from the set $\{0\} \times C(\widehat{Z})$, where ($\{0\}$ is the origin of the smooth disc $C(\mathbb{S}^{\ell})$) to a point on the fiber.
    We equip this vanishing set with the \textbf{localizing control data} $(\mathcal{T}_{F}, \pi_{F}, \rho_{F})$
    where $(\mathcal{T}_{F}, \pi_{F}):=(U(\widehat{{F}}) , \pi_{F})$ and $\rho_{F}:=\rho_Y^2+\rho_N^2$ where $\rho_Y$ is the boundary defining function of $Y$ corresponding to the radial function given by the control data.
\end{definition}

Since $F$ is embedded in $Y_{\alpha}$, the normal bundle of the submanifold within $Y_{\alpha}$ has $\mathbb{R}^l$ fibers, of which we take a ball neighbourhood $C(\mathbb{S}^{\ell})$ and take the product of these fibers at each point with the fibers $C(\widehat{Z})$ of the fibration $\pi_\alpha$ to form the fibers of the new fibration $\pi_{F}$.
The case of normally non-singular fixed point sets in \cite{banagl2025gsignaturetheoremwittspaces}.

In the following we will assume that the components $Y_{\alpha}$ of the strata $\mathcal{S}(\widehat{X})$ are connected, as can be arranged by refining the stratification by simply separating the connected components.

\begin{definition}[\textbf{Type II stratified vanishing set}]
Consider $\widehat{F}\subset \widehat{X}$, a connected component of $Z(A)$ (not necessarily disjoint from $\widehat{X}^{reg}$).
We say $\widehat{F}$ is a Type II stratified vanishing set if $\widehat{F}$ is a Thom-Mather stratified space which is equal to some $\widehat{Y_{\alpha}}$, by which we denote the metric closure (with the wedge metric obtained by restriction) of a connected stratum $Y_{\alpha}$ of $\widehat{X}$.
 
We denote by $F^{reg}=\widehat{F}^{reg}:= \widehat{F}_m \setminus \widehat{F}_{m-1}$ the regular part of $\widehat{F}$.
In this case we define the \textbf{localizing control data} $(\mathcal{T}_{F^{reg}}, \pi_{F}, \rho_{F})$ for $F^{reg}$ to be the restriction of the control data $(\mathcal{T}_Y,\pi_Y, \rho_Y)$ of $Y$. 
\end{definition}

We note that the Type II stratified vanishing sets are a generalization of the first two types, where the key is $\mathbb{R}^k$ is a cone over a sphere with the round metric, and thus the fibers of the normal bundle of a smooth embedding are also conic fibers. Together with the fact that a product of cones is the cone over the joing of the links of the factors, it is easy to see that the preceding two definitions fall into this third type. Nevertheless we make the distinction to help the reader match the smooth and normally non-singular cases to our work.

In what follows, we will always define and study objects related to $F^{reg}$ (in the case of Type I stratified vanishing sets $F=F^{reg}$ since in that case $F$ is a smooth manifold (a pseudomanifold of depth $0$).

\begin{definition}[\textbf{Non-degeneracy of vanishing sets}]
\label{definition_non_degeneracy_critical_point_sets_Kahler}
If $F \in Z(A)$ is a stratified vanishing set which is of Type I or Type II, we say it is nondegenerate if there is the following compatibility between the vector field $V$ and the localizing control data that we defined in each case.

Namely, given the vanishing control data $(\mathcal{T}_{F^{reg}}, \pi_{F}, \rho_{F})$, consider a decomposition of the locally trivial fibration as a Whitney sum of fiber bundles, that is 
\begin{equation}
\label{equation_fibration_critical_components}
    \pi_{F}:\mathcal{T}_{F^{reg}}\to F^{reg}, \quad \mathcal{T}_{F^{reg}} = \bigoplus_{j} E_j,
\end{equation}
where $E_{j}:=C(\widehat{Z_{j}})$, and with respect to this splitting each fiber is a product of cones 
\begin{equation}
\pi_{F}^{-1}(q) :=  \prod_j C(\widehat{Z_{j}}), \quad C(\widehat{Z}) =  \widehat{Z} \times [0,2)_{\rho_{F,j}} /_{\sim}
\end{equation}
where $\sim$ is the usual identification at $\{\rho_{F,j}=0\}$ over other Thom-Mather stratified spaces $\widehat{Z_{j}}$. 
 
The link $Z$ of the normal conic fiber is the join of the links of the conic factors, with the induced metric.
Denoting the restriction of $\rho_{F}$ to each subbundle by $\rho_{F}\big|_{E_j} =: \rho_{F,j},$
we require their compatibility with the `radial functions' $\rho_{F}$ defined for the stratified vanishing sets, via
\begin{equation}
\label{eq:distance_function_components}
\rho^2_{F}\big|_{\mathcal{T}_{F^{reg}}} = \sum_j \rho_{F,j}^2,
\end{equation}
We demand that the vector field $V$ satisfies
\begin{equation}
\label{Morse_Bott_Kirwan_modify_hamiltonian_vf}
V \big|_{\mathcal{T}_{F^{reg}}} = V_F+\rho_{F}^{2\delta_E}V_{F,E}, \quad V_F:=\sum_k \gamma_{F,k} V_{F,k}
\end{equation}
for $\delta_E>1/2$ and for some constants $\gamma_{F,k} \neq 0$, and where all $V_{F,k}$ are Killing vector fields restricted to each fiber $C(Z_j)$, and where the \textbf{error vector field} $V_{F,E}$ is a wedge vector field on the resolved localizing neighbourhood. Here $V_F$ is called the \textbf{model vector field at $F$}.

In the Hamiltonian case, we have that there exists a function $\widehat{h_F}=h_F+ \mathcal{O}(\rho^{2\delta_E})$ on $\mathcal{T}_{F^{reg}}$, where $h_F:=(1/2) \rho_F^2 f$ for a function $f \in C^{\infty}_{\Phi}(X)$ on the link $Z$ (see equation \ref{equation_Hamiltonian_relevant_equation_1}). 

We will \textbf{assume} $h_F=(1/2) \rho_F^2 f=(1/2) \rho_1^2 f_1 + (1/2) \rho_2^2f_2$ where $f_1$ and $-f_2$ are both strictly positive functions on resolved links $Z_1, Z_2$ respectively, where $C(Z)=C(Z_1)\times C(Z_2)$ where $Z_1$ and $Z_2$ have induced contact structures.
\end{definition}

Note that the wedge vector field $V_{F,E}$ restricted to $F^{reg}$ is a section of the wedge tangent bundle.
\begin{remark}[Notation and nomenclature, vanishing control data]
\label{Remark_notation_and_nomenclature_vanishing_control_data}
    The nomenclature and notation for the vanishing control data for a vanishing set $F \in Z(A)$ matches those for the control data of a stratum $F$ in Definition \ref{Thom-Mather stratified space}, and in particular are the same when $F$ is a singular stratum.
\end{remark}

We make a few remarks to clarify Definition \ref{definition_non_degeneracy_critical_point_sets_Kahler} further.

\begin{remark}[Remark 1 on Definition \ref{definition_non_degeneracy_critical_point_sets_Kahler}]
If $\omega- \omega_{w,pt}= x^{2\delta} \omega_E$ where we denote by the right hand side the error term in equation \eqref{equation_structure_of_the_symplectic_form}, then 
\begin{align*}
    dh = - \iota_V \omega &= -\iota_{V_F}\omega_{w,pt} - x^{2\delta_E}\iota_{V_{F,E}} \omega_{w,pt} - x^{2\delta} \iota_{V_F} \omega_E \\
    &= dh_F + O(x^{2\delta'}) , \qquad \text{where $\delta'=\min\{\delta_E, \delta\}$}.
\end{align*}   
If $\delta'>1/2$ then the error terms are $O(x^{2\delta'})$. Thus we need both $\delta_E>1/2$ and $\iota_{V_{F}}\omega_E =O(x)$, as $x\to 0^+$. The compatibility of the almost complex structure will impose a further condition on $J$ in the triple $(g,J,\omega)$.
\end{remark}

\begin{remark}[Remark 2 on Definition \ref{definition_non_degeneracy_critical_point_sets_Kahler}]
The decomposition in equation \ref{equation_fibration_critical_components} is not assumed to be unique. For instance, one can take the trivial decomposition, where the non-degeneracy of the vector field simply reduces to the content of equation \ref{Morse_Bott_Kirwan_modify_hamiltonian_vf}.  

The generality of this definition is partly motivated to tackle more general settings, as in singular warped product metrics (see \cite{jayasinghe2024witten}), and to deal with issues addressed in Remark \ref{Remark_clarify_2}.
\end{remark}

\begin{remark}[Remark 3 on Definition \ref{definition_non_degeneracy_critical_point_sets_Kahler}]
\label{Remark_clarify_2}

In the Hamiltonian situation where $F \in Z(|dh|^2)$ is neither a local minimum nor a local maximum, then Definition \ref{definition_non_degeneracy_critical_point_sets_Kahler} forces at least a two fold decomposition, where locally $h_F=(1/2) \rho_F^2 f=(1/2) \rho_1^2 f_1 + (1/2) \rho_2^2f_2$ where $f_1$ and $-f_2$ are both strictly positive functions on the resolved links $Z_1, Z_2$ respectively.

\textbf{We do not claim that there exists such a decomposition in general.} 
There is such a product decomposition in the smooth setting, where the normal fibers $\mathbb{C}^l=\Pi_{k=1}^l \mathbb{C}$ have actions $\sum_{k=1}^l \gamma_j \partial_{\theta_k}$ where $\theta_k$ are the angular variables on the factors. There are also spaces such as the conifold (see, for instance, Subsection 7.3.4 of \cite{jayasinghe2023l2}) where there are isometric K\"ahler Hamiltonian $\mathbb{T}^3$ actions, but where it is not a product of three cones over circles with the trivial Sasaki structures (note that the K\"ahler metric on the conifold is not flat).     

This \textbf{product assumption} allows one to use separation of variables to characterize the local cohomology in see (Subsection \ref{subsec:cohomology-truncated-neighbourhoods}). Other geometric perspectives on the structure of normal fibers seem better suited for certain types of degenerations in examples, but the analysis becomes significantly complicated.
\end{remark}

In the above definitions for vanishing sets $F$ and their degeneracies, we have not specified localizing neighbourhoods on the stratified spaces as being neighbourhoods fibering over all of $F$, only over $F^{reg}$. To commence with laying analytic foundations for the local cohomology groups appearing in Witten complexes, we now define several notions of localizing neighbourhoods, which carry minor technical decorations that become relevant in the proof, particularly when understanding ideal boundary conditions.

\begin{definition}[singular localizing neighborhood]
Let $\widehat{F}$ be a stratified nondegenerate vanishing set of a potential $A=|V|^2$ on a stratified pseudomanifold $\widehat{X}$ as above. On the tubular neighborhood $\mathcal{T}_{F^{reg}} = \{ \zeta\in \widehat{X} : \rho_{F}(\zeta)< 4t \}$ 
we can define the closed subset 
\[ \left\{ \zeta\in \mathcal{T}_{F^{reg}}: \max\{\rho_{F,u}(\zeta),\rho_{F,s}(\zeta)\} \leq t \right\} , \]
which is a fiber bundle whose fibers $\widehat{U}$ are a product of truncated cones 
\[ \widehat{U} :=  \prod_j \widehat{U_{j}} \]
where each factor $\widehat{U_{j}}$ is homeomorphic to the one point compactification at $\rho_{j}=0$ of  $\widehat{Z_{a,j}}\times [0,t]_{\rho_{F,j}}$, where $\rho_{F,j}$ is the restriction of 
$\rho_{{F},j}$ to each fiber.
Taking the metric completion in $\widehat{X}$ of this fiber bundle over $F^{reg}$, we denote this completion by $U_t(\widehat{F})$ and refer to it as the \textbf{\textit{singular localizing neighbourhood of }$\widehat{F}$}.
The spaces $U_t(\widehat{F})$ can be extended for any $t \in (0,\infty]$ once one is found with a small enough $t$, which always exists by the frontier condition of stratified space.
\end{definition}

\begin{remark}
\label{Remark_not_a_neighbourhood}
We emphasize that we do not demand any fiber bundle structure on $U_t(\widehat{F})$, since any such notion would carry with it additional complications since the candidate base is stratified. See Remark 2.11 of \cite{jayasinghe2024witten} for more details.
The smooth atlas stratified space structures in \cite{albin2025smoothatlasstratifiedspaces} can be used to describe $F$ and to define a fibration over the entire set $F$, but we do not pursue this route here.
\end{remark}

\begin{definition}[resolved localizing neighbourhood]
\label{definition_resolved_fundamental_neighbourhood}
    Given a critical point set $\widehat{F}$ with a singular localizing neighbourhood $U_t(\widehat{F})$, we define the \textbf{(resolved) localizing neighbourhood} $U_t(F)$ as the fiber bundle 
\begin{equation}
\label{equation_fibration_critical_components_resolved}
    {U} - {U_t(F)} \xrightarrow{{\phi_{F}}} F^{reg}
\end{equation}    
    where the fibers $U$ are the products of resolutions  of the fibers $C(\widehat{Z_{j}})$ i.e. 
\begin{equation}
    U = \prod_j U_{j} := \prod_j \left( [0,1]_{r_{F,j}} \times Z_{j} \right),
\end{equation}   
    where $r_{F,j}=\rho_{F,j}$ are new variables on the resolution. We denote the resolution relating these localizing neighborhoods by $\beta_F:U_t(F)\to U_t(\widehat{F})$, where we use that $\widehat{F^{\text{reg}}}=F^{\text{reg}}$ so these fiber bundles live over the same space. Finally we equip $U(F)$ with the restriction of the wedge metric $g_w$.
    We will refer to these resolved localizing neighbourhoods as \textbf{localizing neighbourhoods} when it is clear by context and notation.
\end{definition}

\begin{remark}
    We note that the neighbourhood $U_t(\widehat{F})$ arises as the metric completion in $\widehat{X}$ of a fiber bundle, and thus arises as a subset of $\widehat{X}$. On the other hand, while $U_t(F)$ is still a bundle over $\widehat{F}^{reg}$, we do not consider this resolution as a subset of $X$ (the resolution of $\widehat{X}$), and indeed it may not embed into $X$ in a natural way. There are certain ``alternate resolutions" that can be associated to $\widehat{X}$, in which the resolved neighbourhoods are embedded in those spaces, but we do not require them for our study here.
    
    For our later analysis, it will suffice to consider the space $U_t(F)$ separately from $X$, where the association is via the spaces of sections defined on each space.
\end{remark}

\begin{remark}
\label{remark_different_length_neighbourhoods}
To keep the notation concise, we shall denote by localizing neighbourhoods that are small enough by $U_t(\widehat{F})$, $U({F})$ in the singular and resolved cases, other than in places where it becomes relevant. The subscript $t$ plays no major role, especially after passing to $0$-th order equivariant Hilbert complexes.
\end{remark}

The following is an adaptation of the notion of the \textbf{metric boundary} of the \textit{resolved fundamental neighborhood} $U(F)$, presented in Definition 2.12 of \cite{jayasinghe2024witten}. We refer to \cite[\S 2.1]{Jesus2018Wittensgeneral} for a discussion of the product of cone metrics being a cone metric (c.f., \cite[\S 3.2.1]{Jesus2018Wittens}) where a refinement of the original stratification of $\widehat{U}$ and a refinement of the corresponding manifold with fibered boundary structure to the resolved space $U_p$ is presented for products of cones.

\begin{definition}
\label{definition_smooth_boundary_and_singularities}
Given a localizing neighbourhood $U({F})$, we define the fibers of the metric boundary to be 
to be 
\begin{equation}
    \partial U = \bigcup_k \bigg[{ \prod_{j \neq k} C_{r_{a,j}}({Z_{a,j}})} \times \{r_{a,k}=1\} \times Z_{a,k} \bigg] 
\end{equation}
which are the fibers of the fibration 
\begin{equation}
    {\partial U} - {\partial U(F)} \xrightarrow{{\partial \phi_{F}}} F^{reg}
\end{equation}
where this is the fibration obtained by restricting the fibration ${\phi_{F}}$ in \eqref{equation_fibration_critical_components_resolved} to the metric boundaries of the fibers.
\end{definition}

It is clear that under the blow-down map, the image of the metric boundary consists of the boundary components of co-dimension 1 (as opposed to the singular components of higher codimension) of the singular localizing neighbourhood. This is where we shall apply boundary conditions (supplementing the ideal boundary conditions chosen for the sections on the entire space $X$) to fix a domain for the local Hilbert complexes.

In \cite{jayasinghe2024witten}, Proposition 2.16 served to show that there exists a model metric in every fundamental neighbourhood of a critical point set where the metrics restricted to the fibers were isometric. Here we only use the fact that there is a model metric obtained by freezing coefficients of the expansions of the metrics near the critical point sets, yielding a metric where the fiberwise metrics are product type/rigid.
For isometric actions, we can use the $G$ action to integrate the metrics and obtain model metrics which are equivariant under the action of the model vector field $V_F$ in Definition \ref{definition_non_degeneracy_critical_point_sets_Kahler}.
However, just freezing the leading order coefficients of the metric, the vector fields, Hamiltonian functions and other objects in in typical microlocal geometric fashion gives the following formal statement. 

\begin{proposition}
\label{Proposition_model_metric_existence_Kahler}
Given a resolved fundamental neighbourhood of a critical point set $F$ of a K\"ahler Hamiltonian function, there exists a model K\"ahler wedge metric $g_F$ on $F$ such that 
\begin{equation}
\label{equation_equivariant_model_metric}
    g=g_F+\mathcal{O}(\rho_{F}^{{\delta_E}})
\end{equation}
where the model vector field $V_F=\sum_j \gamma_{F,j} V_{F,j}$ generates a K\"ahler-Hamiltonian isometry of $g_F$.
\end{proposition}

In fact, we can define an alternate stratification near the critical point set, compatible with the group action (similar to those for Bialynicki-Birula decompositions) such that the form of the metric in equation \eqref{equation_equivariant_model_metric} will be just that given in equation \eqref{equation_structure of the metric} for the compatible stratification.

\begin{remark}[Local extension of circle actions to $\mathbb{C}^*$ action]
\label{remark_local_extention_c_star}
Given a self map $f_{\theta}(r_j,z_j)$ on the tangent cone generated by the flow at time $\theta$ of a K\"ahler Hamiltonian Killing vector field, the map $f_{\theta}(s r_j,z_j)$ for $s \in \mathbb{R}$ gives an extension of the $S^1$ action to a $\mathbb{C}^*$ action, where $\lambda=se^{i\theta} \in \mathbb{C}^*$ acts on the tangent cone. Such extensions can be done on localizing neighbourhoods of critical point sets for K\"ahler Hamiltonian actions that we study. 
\end{remark}

\begin{definition}[K\"ahler moment map]
    In the same setting as in Definition \ref{Definition_wedge_Kahler_Hamiltonian_function}, but for $G$ a compact connected Lie group, we define a wedge \textbf{\textit{K\"ahler moment map}} as a moment map $M$ in the usual sense on $X^{reg}$ which extends to a continuous function on $X$ which we denote 
    \[M: X \to \mathfrak{g}\]
    where $\mathfrak{g}$ is the Lie algebra of $G$, where we demand that if $g \in G$ generates circle actions on $X^{reg}$ corresponding to K\"ahler Hamiltonian functions on $X^{reg}$ they extend to K\"ahler Hamiltonian functions as in Definition \ref{Definition_wedge_Kahler_Hamiltonian_function}. 
    Then we denote a geometric endomorphism corresponding to an element $g \in S^1 < G$ as $T_g$.
    This can be extended to a complexified moment map locally near any fixed point set (see Remark \ref{remark_local_extention_c_star}), with image the complexified Lie algebra $\mathfrak{g}_{\mathbb{C}}$.
\end{definition}

\begin{remark}
We use the convention for the moment map so that in the case of circle actions, it matches the Hamiltonians. 
In terms of \textit{algebraic moment maps} studied in the literature, ours would correspond to the norm squared of such maps. 
The Duistermaat-Heckman measures are given by the pushforward of the symplectic measure under the moment map in the sense we consider.
\end{remark}

\begin{remark}[Sign convention for Hamiltonian actions]
\label{Remark_sign_convention_Hamiltonians}
Here we follow the sign convention in \cite{witten1984holomorphic} (different from that in \cite{wu1998equivariant}) and take the Hamiltonian vector field to be $V=-J \nabla h$ (where we denote the gradient of $h$ by $\nabla h$), equivalently $dh=-\iota_V \omega$.

The difference in this sign convention is the source of the difference between the local deformed cohomology groups on page 326 of \cite{witten1984holomorphic}, and those in Proposition 3.2 of \cite{wu1998equivariant}.
\end{remark}

\section{Dirac operators and their Witten deformation}
\label{section_Witten_deformation}

In this section we describe wedge Spin-c Dirac operators, and study some important identities for their Witten deformed versions, in particular showing that the Hamiltonians are of admissible type.
In subsection \ref{subsection_spin_c_Dirac} we review Dirac operators on wedge structures and in Subsection \ref{subsection_Witten_deformation_Dirac}, their Witten deformation. Here we also present some key technical arguments that are at the heart of the formalism which we built for localization.

\subsection{Dirac-type operators}
\label{subsection_Dirac_Operators}

We briefly review Dirac-type and spin-c Dirac operators in the wedge context, referring to \cite{jayasinghe2024holomorphic,jayasinghe2023l2,Albin_2017_index} for more details.
We follow the conventions in \cite{wu1998equivariant}, referring to chapters 3 and 5 of \cite{duistermaat2013heat} for a more detailed construction of the spin$^\mathbb{C}$ Dirac operator where in the K\"ahler case this is also equal to the Dolbeault-Dirac operator. Differences between these to in the non-K\"ahler setting are  discussed in \cite{duistermaat2013heat,epstein2006subelliptic,tian1998sympletic}.

Given a wedge metric, consider a Clifford action
\begin{equation}
    cl: \prescript{w}{}{TX} \rightarrow \text{End}(\mathcal{E}), \quad cl(u)cl(v)+cl(u)cl(v)=-2g_w(u,v)\operatorname{Id}_\mathcal{E}
\end{equation}
which makes $\mathcal{E}$ into a wedge Clifford module (see Section 1.3 of \cite{Albin_2017_index} where the vertical wedge Clifford modules over a point restrict to what we study here).
Consider a connection $\nabla^\mathcal{E}$ compatible with the Clifford action in the sense that
\begin{equation}
    \nabla^\mathcal{E}\big(cl(v)s\big)=
cl\big(\nabla^{{LC}}v\big)s+cl(v)\nabla^\mathcal{E}s,
\end{equation}
where $\nabla^{{LC}}$ is the Levi-Civita connection.
Given a (local) orthonormal frame of the wedge tangent bundle $\{e_1, e_2, \ldots, e_{2n} \}$, we have the Dirac operator
\begin{equation}
    D=\sum_{i=1}^{2n} cl(e_{i}) \nabla^\mathcal{E}_{e_{i}}.
\end{equation}
Intrinsically, $D$ is the composition
\begin{equation}
    \Gamma(\mathcal{E})
\xrightarrow{\nabla^\mathcal{E}}
\Gamma(T^*M\otimes \mathcal{E})
\xrightarrow{cl}
\Gamma(\mathcal{E}),
\end{equation}
where the second map uses the metric identification
$^wT^*M \cong ^wTM$ and is defined by $c(\alpha\otimes s):=c(\alpha^\sharp)s$.

If for some $Y\in \mathcal{S}(X)$, we restrict to a collar neighborhood of $\mathfrak{B}_Y$ as in Subsection \ref{subsection_asymptotics_metric} where we have an asymptotically $\delta$ wedge metric, the Dirac operator takes the form 
\begin{equation}
\label{temp_lah_di_dah}
cl(\partial_x)\nabla^\mathcal{E}_{\partial_x} + \sum_i {cl}\left(\frac{1}{x}\partial_{z_i}\right)\nabla^\mathcal{E}_{\frac{1}{x}\partial_{z_i}} + \sum_j {cl}(\partial_{y_j})\nabla^\mathcal{E}_{\partial_{y_j}}
\end{equation}
up to a differential operator in $x^{\delta}\text{Diff}^1_w(X;\mathcal{E})$.
Here $x$ is a boundary defining function for $\mathfrak{B}_Y$, and we recognize \eqref{temp_lah_di_dah} as a wedge differential operator of order one. 

We shall study settings where the Clifford modules that admit splittings corresponding to fibrations that respect the module structure.
Consider a fibration $\pi:N \rightarrow B$ over a smooth base $B$ with fibers $U=C_x(Z)$, where $Z$ is a resolved link and the fibers are $U=(0,t)_x\times Z$ for some $t>0$. Then the wedge tangent bundle admits a splitting as a direct sum of the vertical wedge tangent bundle $V(\prescript{w}{}{TX}):= \text{ker} d\pi$ and the horizontal wedge tangent bundle $H=V^{\perp}$.

We shall assume there is a corresponding factorization of the Clifford module 
\begin{equation}\label{eq:clifford-mod-factorization}
    \mathcal{E}=\mathcal{E_H} \otimes \mathcal{E_V}
\end{equation}
where $\mathcal{E_H}, \mathcal{E_V}$ are closed Clifford submodules of $\mathcal{E}$ where the restrictions of the horizontal and vertical wedge tangent bundles act. The Dirac operator can be written as the sum of two operators $D^{\mathcal{H}}, D^{\mathcal{V}}$ that act on these bundles. We shall study these in detail in Subsection \ref{subsection_discreteness_of_spectrum}.

For instance in the de Rham setting, any submersion of a smooth manifold in another gives such a splitting of the tangent space on the base, where the exterior algebra of the cotangent bundle $\mathcal{E}$ on $M$ can be constructed using those on the base and the fibers which are the horizontal and vertical submodules.

\subsubsection{Spin$^\mathbb{C}$ Dirac and Dolbeault operators}
\label{subsection_spin_c_Dirac}

Given a wedge metric $g_w$ and an almost complex structure on the wedge tangent bundle (which we call a wedge almost complex structure), the complexified wedge tangent bundle ${}^{w}T_{\mathbb{C}}X$ has an orthogonal splitting into the holomorphic and anti-holomorphic tangent bundles as ${}^wTX^{1,0} \oplus {}^wTX^{0,1}$ induced by the wedge almost complex structure.
Given $v \in \Gamma({}^{w}T_{\mathbb{C}}X)$ where $^{w} T_{\mathbb{C}}X=\prescript{w}{}{TX} \otimes_{\mathbb{R}} \mathbb{C}$,
we can write it as $v=v^{1,0}+v^{0,1}$ where $v^{1,0} \in \Gamma(\prescript{w}{}{TX}^{1,0})$ and $v^{0,1} \in \Gamma(\prescript{w}{}{TX}^{0,1})$. Then a Clifford action of the complexified wedge Clifford algebra is determined by
\begin{equation*}
cl(v^{1,0})=\sqrt{2}\, (v^{1,0})^{*} \wedge, \quad cl(v^{0,1})=-\sqrt{2}\, \iota_{v^{0,1}},
\end{equation*}
where $(v^{1,0})^{*}=\xi^{0,1} \in \Gamma(^{w} T^*X^{0,1})$ corresponds to $v^{1,0}$ via the wedge metric.

Given a Hermitian bundle $E$ on $X$ equipped with a compatible connection, this extends to a wedge Clifford action on $\text{End}(\mathcal{E})$ where $\mathcal{E}= \bigoplus_{q=0}^n \Lambda^{q} (\prescript{w}{}{TX}^{1,0}) \otimes E$ by $cl \otimes Id$, which we denote by $cl$ with abuse of notation. 

Given wedge vector fields $v, e \in \Gamma(^{w} T_{\mathbb{C}}X)$, we have the anti-commutation relation $\{cl(v), cl(e)\}=-2 g(v, e)$. The holomorphic Hermitian connection on $\prescript{w}{}{TX} \otimes \mathbb{C}$ corresponds to the Levi-Civita connection on $\prescript{w}{}{TX}$, and we choose a compatible connection on $\mathcal{E}=\Lambda^{q} (\prescript{w}{}{TX}^{1,0}) \otimes E$ which we denote by $\nabla^\mathcal{E}$. Composing $\nabla^\mathcal{E}$ with the Clifford action we obtain the spin$^\mathbb{C}$ Dirac operator $D$ acting on sections supported on $\mathring{X}$.

It is well known that in the K\"ahler case we can write the operator as $D=\sqrt{2}(\overline{\partial}+\overline{\partial}^*)$ where
\begin{equation}
\label{equation_d_bar_spin_c_correspondence_1}
\bar{\partial}=\frac{1}{2 \sqrt{2}} \sum_{i=1}^{2 n} cl(e_i-\sqrt{-1} J e_i) \nabla^\mathcal{E}_{e_i}, \quad \bar{\partial}^{*}=\frac{1}{2 \sqrt{2}} \sum_{i=1}^{2 n} cl(e_i+\sqrt{-1} J e_i) \nabla^\mathcal{E}_{e_i}.
\end{equation}
Here $J$ is an (integrable) complex structure that tames the wedge K\"ahler form.

\subsection{Witten deformation for Dirac operators}
\label{subsection_Witten_deformation_Dirac}

Given a Dirac operator and a Killing vector field $V$ that respects the action of the operator $D$ in the sense that $[D,L_V] = 0$ (where $L_V$ is the infinitesimal generator of the action generated by the flow of $V$ on the sections of the Clifford module), Witten defined a deformed operator in \cite{witten1983fermion}
\begin{equation}
D_\varepsilon := D + \sqrt{-1}\,\varepsilon\,cl(V)    
\end{equation}
such that $[D_{\varepsilon},L_V] = 0$, for which the Hamiltonian is 
\begin{equation}
    D_{\varepsilon}^2:=D^2+i\varepsilon\{ cl(V), D\}+\varepsilon^2 |V|^2
\end{equation}
where the operator $i\varepsilon\{cl(V),D\}$ has a first order part given by $iL_V$, where $L_V$ is the infinitesimal generator of the action on the sections of the bundle.

In the case of the spin-c Dirac complex in the K\"ahler Hamiltonian setting, a simple computation shows that
\begin{equation}
    D_\varepsilon = P_{\varepsilon}+P_{\varepsilon^*}=\sqrt{2}\,(\overline{\partial}_\varepsilon + \overline{\partial}_\varepsilon^*) = D + \sqrt{-1}\,\varepsilon\,cl(V)
\end{equation}
where the deformed $\mathcal{N}=2$ SUSY operators are
\begin{equation}
    \overline{\partial}_\varepsilon = \overline{\partial} + \sqrt{-2}\,\varepsilon\,(V^{1,0})^*\wedge,
\qquad
\overline{\partial}_\varepsilon^* = \overline{\partial}^*-\sqrt{-2}\,\varepsilon\,\iota_{V^{0,1}},
\end{equation}
which can also be shown to arise by conjugation, as
\begin{equation} \label{equation_witten_def_dbar}
    \overline{\partial}_\varepsilon := e^{-\varepsilon h}\,\overline{\partial}\,e^{\varepsilon h},
    \qquad
    \overline{\partial}_\varepsilon^* := e^{\varepsilon h}\,\overline{\partial}^*\,e^{-\varepsilon h}
\end{equation}
where $h$ is the Hamiltonian function.
In more abstract notation for $\mathcal{N}=2$ SUSY complexes 
\begin{equation} 
    P_\varepsilon := e^{-\varepsilon h}  Pe^{\varepsilon h},
    \qquad
    P_\varepsilon^* := e^{\varepsilon h}P^*e^{-\varepsilon h}
\end{equation}

\begin{remark}[The de Rham case]
\label{Remark_de_Rham_conjugation_N_2}
In the de Rham case, the conjugation by $e^{\varepsilon h}$ for functions $h$ (such as in \cite{jayasinghe2024witten} for Stratified Morse-Bott functions) leads to the Dirac operator $D_{\varepsilon}=D+i\varepsilon \widehat{cl}(dh^{\sharp})$, where $\widehat{cl}$ is the principal symbol of the operator $Q_2=i(d-d^*)$ (see equation 4.12 of \cite{Zhanglectures}). The anti-commutation relation $\{Q_1,Q_2\}=0$ can be used to see that the principal symbols of the operators also anti-commute, and that the operator $D_{\varepsilon}^2=D^2+\varepsilon\{Q_1,\widehat{cl}(dh^{\sharp})\}+\varepsilon^2|dh|^2$ is of zero-th order type.    
\end{remark}

\subsubsection{Operator identities}

We review the Bochner identities for the Witten-deformed Dolbeault-Dirac operator. The formulas in the smooth setting are proven in \cite[\S3]{wu1998equivariant}, and extensions for the stratified pseudomanifolds are proven in \cite[\S4]{jayasinghe2024holomorphic}. 

Let $\overline{L_V}=\{\iota_V,\nabla^E\}$ denote the Lie derivative on $E$-valued anti-holomorphic forms induced by the lifted $S^1$-action. We denote by $\sqrt{-1}L_V$ the infinitesimal generator of the $S^1$-action on $L^2(X;\mathcal{E}) = L^2\Omega^{0,\bullet}(X;E)$, where $\mathcal{E} = \Lambda^{0,\bullet}({}^{w}T^*X)\otimes E$. 

\begin{proposition}[Operator identities]
\label{Proposition_operator_identities_Gollum}
In the setting above we have the identities
\begin{equation}
\label{equation_operator_identities_1}
    D=\sqrt{2}\left(\overline{\partial}+\overline{\partial}^{*}\right),
\end{equation}
\begin{equation}
\label{equation_operator_identities_2}
    D_{\varepsilon}:= D + \sqrt{-1} \varepsilon cl(V) =\sqrt{2}(\overline{\partial}_{\varepsilon}+\overline{\partial}_{\varepsilon}^{*}),
\end{equation}
and
\begin{equation}
\label{equation_operator_identities_3}
    \Delta_{\varepsilon} := D_\varepsilon^2 = D^{2}+\varepsilon K+\varepsilon^2|V|^{2} ,
\end{equation}
where 
\begin{equation}
\label{equation_operator_identities_4}
    K:=\frac{1}{2} \sqrt{-1} \sum_{i=1}^{2 n} cl\left(e_{i}\right) cl\left(\nabla_{e_{i}} V\right)+\left.\sqrt{-1} \operatorname{tr} \nabla \cdot V\right|_{^wT^{0,1} X} -2\sqrt{-1} (\overline{L_V})
\end{equation}
where $\nabla( \cdot )$ 
is the divergence operator on vector fields on $X^{reg}$ which extends uniquely to an operator on wedge vector fields.
\end{proposition}
\begin{proof}
    The proof is the same as that of \cite[Proposition 4.1]{jayasinghe2024holomorphic}. The computation on $X^\mathrm{reg}$ uses \cite[Lemma 3.1]{wu1996equivariant} and extends directly to $X$, as every term appearing in the formula is a wedge differential operator (where bundle endomorphisms are 0-th order operators). We only highlight the identities needed in the argument.
    
    The first identity is the standard expression for the Dolbeault-Dirac operator $D$, following from equation \eqref{equation_d_bar_spin_c_correspondence_1}. The second identity on $D_\varepsilon$ follows from the decomposition of $\overline{\partial}_{\varepsilon}$ and $\overline{\partial}_{\varepsilon}^*$:
    \begin{equation}
        \overline{\partial}_{\varepsilon}=\overline{\partial}+\frac{1}{2 \sqrt{2}} \sum_{i=1}^{2 n} cl\left(e_{i}-\sqrt{-1} J e_{i}\right) (-\nabla^E_{e_i}h), \quad \overline{\partial}_{\varepsilon}^{*}=\overline{\partial}^{*}-\frac{1}{2 \sqrt{2}} \sum_{i=1}^{2 n} cl\left(e_{i}+\sqrt{-1} J e_{i}\right) (-\nabla^E_{e_i}h).
    \end{equation}
    Substituting the expression for $D_\varepsilon$ and using
    \begin{equation}
        \label{Baila_gahapang_2}
        \left\{D, cl(V)\right\}=\sum_{i=1}^{2 n}\left(\left\{cl\left(e_{i}\right), cl(V)\right\} \nabla^E_{e_{i}}+cl\left(e_{i}\right) cl\left(\nabla^E_{e_{i}} V\right)\right)=-2 \nabla^E_{V}+\sum_{i=1}^{2 n} cl\left(e_{i}\right) cl\left(\nabla^E_{e_{i}} V\right) 
    \end{equation}
    gives the last two equations.
\end{proof}

In this case, the operator $K$ has both a first order part and a zero-th order part.
The key difference from the de Rham case for Morse functions is the appearance of the Lie derivative $L_V$, reflected in the equation 
\begin{equation}
\label{equation_de_Rham_Dolbeault_difference}
    D^2_\varepsilon= \frac{1}{2}\Delta_{\varepsilon,dR} - \sqrt{-1} \varepsilon L_V
\end{equation}
first observed in \cite{witten1982constraints} (see also, \cite[Equation 15]{witten1984holomorphic} and the introduction of \cite{mathai1997equivariant}). 
In the wedge setting we have $|V|_{^wTX}^2=|dh|_{^wT^*X}^2$, similar to the smooth case. Restricting to a $\mu$-eigenspace of $\sqrt{-1}L_V$, the operator is bounded. 
We see that for the family of deformed Dolbeault complexes, the Hamiltonian is of the form in \eqref{equation_introduction_deformed_Hamiltonian} where the potential is $|V|^2$, the vanishing set of which is the zero set of the vector field $V$, equivalently the critical points of the hamiltonian $h$.

\begin{remark}[Bochner-Schr\"odinger identities]
\label{Remark_Lichnerowicz_Schrodinger}
The identity for the expansion of the Hamiltonian $\Delta_{\varepsilon}$ in equation \eqref{equation_operator_identities_3} admits two important conceptual descriptions, well known by now to experts.
\noindent First, as a Schr\"odinger operator where the potential $\varepsilon^2|V|^2$ plays a key role in the localization,

\noindent The second as a Bochner formula, where the Witten deformation corresponds to a perturbation of the metric of the bundle $E$, which leads to a compatible perturbation of the connection on the bundle, which leads to the corresponding perturbation of the Dirac operator $D_{\varepsilon}$.
This is clearest in the case of the de Rham complex, Witten deformed by conjugation $e^{\varepsilon h}$ (as in \cite{witten1982supersymmetry}) for any function $h$, the deformed operator $Q^+_{\varepsilon}=d+\varepsilon dh \wedge$, where $\varepsilon dh$ is the connection one form.
This perspective of changing the bundle metric to localize the topological invariant quantities ties into the description of Witten deformation as a topological quantum field theory.

Such identities are referred to by many names, of mathematicians and physicists, including Lichnerowicz, Weitzenbock and others, and we refer the interested reader to \cite{kay2020editorial}.   
The question of realizing potentials of Schr\"odinger operators as curvatures of geometric data has a rich history, notably in the Nirenberg problem and other prescribed curvature problems, and carries natural positive answers in these settings of Witten deformation.
\end{remark}

\begin{remark}[The case of Hamiltonian actions that are non-K\"ahler]
\label{Remark_non_Kahler_Hamiltonian_2}
    In the case of a symplectic manifold with a Hamiltonian vector field $V$, the deformation of the spin-c Dirac operator $D$ to $D+i\varepsilon cl(V)$ was used to understand the harmonic spinors in \cite{tian1998sympletic}, where Theorem 1.4 gives the Bochner-Schr\"odinger operator $D_{\varepsilon}^2$ in that setting, and it is easy to check that it is of first order type with potential $|V|^2$.

    However in this case, the complex only has $\mathcal{N}=1$ SUSY in general.
\end{remark}

\subsubsection{Localization}

In this case we use the following proposition to get associated Hilbert complexes for which the Hamiltonian has a term $K$ with a 0-th order operator.

\begin{proposition}[1-st order to 0-th order]
\label{Proposition_equivariant_complexes}
Consider a family of $\mathcal{N}=N$ Hilbert complexes \begin{equation*}
    \mathcal{Q}_{\varepsilon}(X) = \Big( H=L^2(X;\mathcal{E}),\{ (\mathcal{D}(Q_{j,\varepsilon}), Q_{j,\varepsilon})\}_{j=1}^N,\gamma \Big)
\end{equation*} associated to a fibered manifold $X$, where the Hamiltonian is of the form $\Delta_{\varepsilon}=\Delta+\varepsilon (K'+iL_V)+\varepsilon^2 A$, where $K$ and $A$ are bounded 0-th order operators, and where $iL_V$ is the infinitesimal generator of the action on the sections, for an stratum preserving isometric action generated by the flow of a vector field $V$.

Then there are equivariant Hilbert complexes 
$\mathcal{Q}_{\varepsilon}^{\mu}(X)= \Big( H,\{ (\mathcal{D}^{\mu}(Q_{j,\varepsilon}), Q_{j,\varepsilon})\}_{j=1}^N,\gamma \Big)$
where 
\begin{equation}
    \mathcal{D}^{\mu}(Q_{j,\varepsilon}) =\{ u \in \mathcal{D}^{\mu}(Q_{j,\varepsilon})\;|\; u \in \emph{Ker}(iL_V-\mu) \}.
\end{equation}

In particular, the Hamiltonians on these complexes are of 0-th order type.
\end{proposition}

\begin{proof}
We check that the Kernel of $iL_V-\mu$ is closed for any constant $\mu$. Since $V$ is assumed to generate isometries, and preserves strata, in any collar neighbourhood of a boundary hypersurface with boundary defining function $\rho$, it will be an edge vector field. 

The Lie derivative $L_V$ for an edge vector field is a bounded operator on any self adjoint extension of a Dirac operator since $\mathcal{D}_{\max}(D) \subseteq H^1_e(X;\mathcal{E})$ (see \eqref{equation_inclusions_of_function_spaces})
Then it is standard that the kernel of a bounded linear operator is closed.
\end{proof}

This proposition allows us to associate complexes with Hamiltonians of $0$-th order type to complexes with Hamiltonians of $1$-st order type. The cohomology groups for the spaces we study are defined only on their domains, and the estimates we shall prove to construct the Witten complexes too are only on the domains.

\begin{proposition}[Localization]
\label{Propostion_growth_estimate_witten_deformed}
Consider a family of $\mathcal{N}=N$ complexes \begin{equation}
    \mathcal{Q}_{\varepsilon}(X)=\Big( H=L^2(X;\mathcal{E}),\{ (\mathcal{D}(Q_{j,\varepsilon}), Q_{j,\varepsilon}) \}_{j=1}^N, \gamma \Big),
\end{equation} associated to a fibered manifold $X$, where the Hamiltonian $\Delta_{\varepsilon}$ is of type 0 and with potential $W$.
Then for any section $s \in \mathcal{D}(Q_{j,\varepsilon})$ for any  $j \in \{1, 2,..., N\}$ with 
\begin{equation}
\label{equation_support_away_from_ZW}
    \emph{supp}(s) \cap \{x \in X \; |\; A(x) \leq c \}=\emptyset
\end{equation}
for any small $c>0$ there exist constants $C>0$, $\varepsilon_0>0$ such that for all $\varepsilon \geq \varepsilon_0$, one has 
\begin{equation}
    ||D_{\varepsilon}s||_{L^2(X;\mathcal{E})} \geq C \sqrt{\varepsilon}||s||_{L^2(X;\mathcal{E})}.
\end{equation}
\end{proposition}

We note that the condition in \eqref{equation_support_away_from_ZW} is satisfied if $\text{supp}(s) \subset \big(X \setminus \bigcup_{F \in Z(A)} U(F) \big)$ for neighbourhoods $U(F_a)$ of spaces $F_a$.

\begin{proof}
Using the identity \ref{equation_introduction_deformed_Hamiltonian}, we see that there exists a finite constant $K_1$ such that, for all $s \in \mathcal{D}(\Delta_\varepsilon)$,
\begin{equation}
||Q_{j,\varepsilon} s||_{L^2(X;\mathcal{E})}^2  = \langle Q_{j,\varepsilon} s, Q_{j,\varepsilon} s \rangle_{L^2(X;\mathcal{E})} \geq (\varepsilon^2 c^2-\varepsilon |K_1|)||s||^2_{L^2(X;\mathcal{E})}
\end{equation}
where we use the fact that $\Delta=Q^2_{j,0}$ is a non-negative operator.
For any $C>0$ we may, by taking $\varepsilon$ sufficiently large, bound the right hand side by $\varepsilon^2 C^2||s||^2_{L^2(X;\mathcal{E})}$. Finally, as the inequality holds for all $s \in \mathcal{D}(\Delta_\varepsilon)$ it holds by density on the closure of $\mathcal{D}(\Delta_\varepsilon)$ in the graph norm of $\mathcal{D}(Q_{j,\varepsilon})$.
\end{proof}

A simple consequence of this proposition is that as one considers larger and larger values of $\varepsilon$, any harmonic section of the Hamiltonian $\Delta_{\varepsilon}$ will have $L^2$-mass concentrated in smaller and smaller neighbourhoods of $Z(A)$. In terms of semi-classical analysis, the $0$-energy eigenmodes (also known as the ground states) of the SUSY quantum mechanical system will converge in $L^2$-mass to the classical ground states that occur at $Z(A)$.

Hence we can already arrive at the heuristic that the global harmonic sections must be spanned by a linear combination of the local harmonic sections in the semi-classical limit, precise versions of which statement will be made in Section \ref{section_Witten_complexes}.

\begin{remark}
\label{Remark_compatible_domains_1}
It is easy to check that for any choice of self-adjoint extension for a Dirac operator $D_{\varepsilon}$ as above, the identity \ref{equation_operator_identities_3} extends to an identity of operators acting on the sections of the domain. For the interested reader, detailed computations at the boundary hypersurfaces resolving the strata were worked out in detail for the de Rham complex in Proposition 6.22 of \cite{jayasinghe2023l2}, Proposition 4.5 of \cite{jayasinghe2024witten}, and for the Dolbeault complex in \cite{jayasinghe2024holomorphic}.
We refer to Remark 6.23 of \cite{jayasinghe2023l2} for remarks on the Morse inequalities in \cite{ProkhorenkovShubin2009Morseboundary} where the differences in domains for the deformed and undeformed operators on a manifold with boundary leads to extra terms at boundaries in the formulas, eventually contributing to Morse polynomials.
In this article we shall pick domains for the deformed and undeformed Dirac operators that are compatible.
\end{remark}

\section{Domains for operators} \label{sec:domains}

In this section, we introduce the framework of SUSY Hilbert complexes, focusing on the case $\mathcal{N}=2$. We also discuss global and local domains, together with the choice of ideal boundary conditions.

\subsection{Abstract $\mathcal{N}=N$ Hilbert complexes} \label{subsec:abstract-complexes}

We begin by discussing the domains of the operators defining an $\mathcal{N}=N$ SUSY Hilbert complex, following Witten's formulation of supersymmetric quantum mechanics\cite{witten1982constraints, witten1982supersymmetry}. Recall that the supercharges $Q_1,\ldots,Q_N$ are odd self-adjoint operators on a $\mathbb Z_2$-graded Hilbert space $\mathcal{H}$, with respective domains $\mathcal{D}(Q_j)$, satisfying the SUSY relations
\begin{equation}
    Q_iQ_j+Q_jQ_i=2\delta_{ij}\Delta,
\end{equation}
where $\Delta$ is a nonnegative self-adjoint operator. Although these relations may hold formally at the level of differential expressions, their extension as identities of unbounded operators depends on the choice of domains. In particular, the requirement that $Q_j^2=\Delta$ for every $j$ imposes compatibility conditions on the domains of the supercharges.

For each $j$, we set
\begin{equation}
    \mathcal{D}(\Delta) = \{ u\in\mathcal{D}(Q_j): Q_ju\in\mathcal{D}(Q_j) \}.
\end{equation}
Since the squared operators are required to coincide with $\Delta$, these domains must agree:
\begin{equation}
    \mathcal{D}(\Delta) = \mathcal{D}(Q_j^2), \qquad j=1,\ldots,N.
\end{equation}
We refer to a collection of self-adjoint extensions of the supercharges satisfying the SUSY relations, including the equality of their squares as self-adjoint operators, as a compatible choice of domains for the $\mathcal{N}=N$ SUSY Hilbert complex.

\begin{lemma}[Domains for SUSY operators] 
\label{lem:susy-domains} 
    Let $\mathcal{Q} = \Big(H, \{(\mathcal{D}(Q_j), Q_j)\}_{j=1}^N, \gamma \Big)$ be an $\mathcal{N}=N$ SUSY Hilbert complex with Hamiltonian $\Delta$, whose spectrum is discrete in $[0,C)$ for some $C>0$ (i.e. note that the operators $Q_j$ are self adjoint and have real spectrum).   
    Then for each $j$ there exists an orthonormal basis of $H$ consisting of simultaneous eigenvectors of $Q_j$ and $\Delta$ for the spectrum below $C$. 
    
    \noindent In particular for each even/odd parity eigenvector $\psi_{\lambda}$ of $\Delta$ with respect to $\gamma$, $Q_j\psi_{\lambda}$ is an odd/even degree eigenvector of $\Delta$ with the same eigenvalue called the \textbf{superpartner of $\psi_{\lambda}$}.

    If $\Delta$ has discrete spectrum (on $[0,\infty)$), then 
    \begin{equation}
    \mathcal{D}(Q_j)=\mathcal{D}(\Delta^{1/2}), \qquad j=1,\ldots,N, 
    \end{equation} 
    and the corresponding graph norms coincide. Moreover the finite linear span of any such basis is a common core for all the supercharges.
\end{lemma}

\begin{remark}[SUSY, McKean-Singer cancellations and superpartners]
    Lemma \ref{lem:susy-domains} is key to what are known as SUSY cancellations in index theory, also called McKean Singer cancellations. When there are multiple SUSY operators, there are more cancellations in invariants arising in the study of index theory, and for $\eta$ invariants this is discussed in \cite{jayasinghe2026spectralasymmetrysupersymmetryequivariant}.
    The pairing of eigenvalues in the abacuses (see Figures \ref{figure_spectral_abacus_isolated}, \ref{figure_spectral_abacus_non_isolated}) is also  a consequence, where we use the refinement for $\mathcal{N}=2$ complexes given by Lemma \ref{lem:exact-coexact-eigenvectors}.
\end{remark}

\begin{proof}[Proof of Lemma \ref{lem:susy-domains}]
    Since $Q_j$ is self-adjoint and $Q_j^2=\Delta$, the spectral theorem gives 
    \begin{equation}
        |Q_j|=(Q_j^2)^{1/2}=\Delta^{1/2}. 
    \end{equation} 
    Hence 
    \begin{equation}
        \mathcal{D}(Q_j)=\mathcal{D}(\Delta^{1/2}), \qquad     \|Q_ju\|=\|\Delta^{1/2}u\|. 
    \end{equation} 
    Thus the graph norms of the supercharges coincide. Suppose now that $\Delta$ has discrete spectrum. Since $Q_j$ strongly commutes with $\Delta=Q_j^2$, each eigenspace of $\Delta$ is invariant under $Q_j$, and there are superpartners in even and odd degree. The restriction of $Q_j$ to each eigenspace is self-adjoint and hence admits an orthonormal eigenbasis. 
    
    Combining these bases yields an orthonormal basis of simultaneous eigenvectors of $Q_j$ and $\Delta$. Finally, the spectral theorem implies that the finite linear span of any orthonormal eigenbasis of $\Delta$ is dense in $\mathcal{D}(\Delta^{1/2})$ with respect to its graph norm. Since this graph norm coincides with that of every $Q_l$, the assertion follows.     
\end{proof}

\subsubsection{The case of $\mathcal{N}=2$} \label{subsubsec:n2-case}

As explained in the introduction, a $\mathcal{N}=2$ SUSY Hilbert complex with self-adjoint supercharges $Q_1$ and $Q_2$, the complex supercharges
\begin{equation}
    P:=Q_+=(Q_1+iQ_2)/\sqrt{2}, \qquad P^*=Q_-=(Q_1-iQ_2)/\sqrt{2}
\end{equation}
define a closed, densely defined differential $P$ satisfying $P^2=0$, together with its adjoint. Thus, a $\mathcal{N}=2$ SUSY Hilbert complex defines a Hilbert complex $\mathcal{P} = (H,\mathcal{D}(P),P)$ in the sense of \cite{bru1992hilbert}, which also determines its adjoint complex $\mathcal{P}^* = (\mathcal{H},\mathcal{D}(P^*), P^*)$. The compatibility of the domains of operators translates to the fact that the domains for $Q_1, Q_2$ are the \textbf{rolled-up domains}
\begin{equation}
\label{equation_rolled_up_extensions}
    \mathcal{D}(Q_j) = \mathcal{D}(P)\cap\mathcal{D}(P^*)
\end{equation}
for both $j=1,2$.
The Laplacian is $\Delta=Q_1^2=Q_2^2=PP^*+P^*P$ with the corresponding domain.

In this setting we say that a section $u$ is (respectively) \textbf{exact, co-exact, closed, co-closed} if
\begin{equation}
    u=Pv, u=P^*v, Pu=0, P^*u=0
\end{equation}
respectively, for some sections $v$ in the Hilbert complex (a straight forward generalization of the nomenclature for forms and the de Rham operator).

We next record a consequence for the spectral decomposition of a $\mathbb Z$-graded Hilbert complex. This will be useful in relating the exact and co-exact eigensections of the Laplacian.

\begin{lemma}\label{lem:exact-coexact-eigenvectors}
    Let $\mathcal{P} = (H,\mathcal{D}(P),P)$ be a Hilbert complex. Then in addition to the spectral results in Lemma \ref{lem:susy-domains}, 
    the operator $|\lambda|^{-1/2}P$ gives a unitary isomorphism from co-closed eigenspaces of $\Delta$ of degree $k$ with any fixed positive eigenvalue $\lambda$, to the exact eigenspaces of $\Delta$ of degree $k+1$ with eigenvalue $\lambda$, and its inverse is given by $|\lambda|^{-1/2}P^*$.     
\end{lemma}

\begin{proof}
    Consider a co-closed unit norm eigenvector $v_{\lambda}$, where co-closeness is that $P^*v_{\lambda} = 0$. Consequently, $P^*Pv_{\lambda} = \Delta v_{\lambda} = \lambda v_{\lambda}$. Thus
    \begin{equation}
        \left\langle     \frac{Pv_{\lambda}}{\sqrt{\lambda}},     \frac{Pv_{\lambda}}{\sqrt{\lambda}}     \right\rangle  = \frac{1}{{\lambda}}     \left\langle     v_{\lambda},P^*Pv_{\lambda} \right\rangle = \frac{1}{{\lambda}}     \left\langle     v_{\lambda}, \lambda v_{\lambda} \right\rangle=1.
    \end{equation}
    Moreover, the Hilbert complex identities including the nilpotencies of $P,P^*$ can be used to see that 
    \begin{equation}
        \Delta (P v_{\lambda})=(PP^*+P^*P)Pv_{\lambda}=PP^*Pv_{\lambda}= P(PP^*+P^*P)v_{\lambda}=P\Delta v_{\lambda} = \lambda (P v_{\lambda}),
    \end{equation}
    that is, the $Pv_{\lambda}$ is an exact eigenvector of eigenvalue $\lambda$.
\end{proof}

We now discuss Witten deformations defined by conjugation. Let $M_\varepsilon$ be a bounded invertible operator on $H$ with bounded inverse, and set
\begin{equation}
    P_\varepsilon=M_\varepsilon^{-1}PM_\varepsilon, \qquad \mathcal{D}(P_\varepsilon)
    =M_\varepsilon^{-1}\mathcal{D}(P).
\end{equation}
Then $P_\varepsilon$ is again a differential satisfying $P_\varepsilon^2=0$, and multiplication by $M_\varepsilon$ induces an isomorphism between the deformed and undeformed Hilbert complexes, hence also between their cohomology groups. This holds both for local and global complexes.
The adjoint domain $\mathcal{D}(P_\varepsilon^*)$ is determined by Hilbert space adjunction. Consequently, the deformed rolled-up extensions as in \ref{equation_rolled_up_extensions},
$\mathcal{D}(Q_{j,\varepsilon}) =\mathcal{D}(P_\varepsilon) \cap\mathcal{D}(P_\varepsilon^*)$ for both $j=1,2$, provides the common domain of the deformed supercharges, preserving the $\mathcal{N}=2$ domain compatibility conditions.

\subsection{Ideal boundary conditions at singular strata} \label{subsec:ideal-boundary-conditions}

We describe \textit{ideal boundary conditions} for differential operators on the regular strata of stratified pseudomanifolds $X$, acting between sections of Hermitian vector bundles. The theory of closed extensions in this setting builds on the extensive literature on wedge operators. We take a point of view inspired by \cite{bera2025diracoperatorstwistedramified,Albin_hodge_theory_cheeger_spaces} and developed to some extent in \cite{jayasinghe2026spectralasymmetrysupersymmetryequivariant}.

\subsubsection{Gelfand-Robbin quotients, complexes and closed extensions} \label{subsubsec:gelfand-robbin-quotient}

Let $X$ be a stratified pseudomanifold equipped with a wedge structure, and let $E, F$ be Hermitian vector bundles over $X^{\mathrm{reg}}$. Consider \textit{any} 1st order differential operator $D$ defined between vector bundles $E$ and $F$ There are two canonical extensions of the operator $D$, whose domains are 
\begin{equation}
\label{equation_minimal_domain_definition_1}
    \mathcal{D}_{\min}(D) = \text{ graph closure of } C_c^\infty(X^{\mathrm{reg}}; E), 
\end{equation}
and
\begin{equation}
    \mathcal{D}_{\max}(D) = \{u \in L^2(X^{\mathrm{reg}};E): Du \in L^2(X^{\mathrm{reg}};E) \text{ distributionally} \}.
\end{equation}
Both extensions are closed and satisfy $\mathcal{D}_{\min}(D)\subset \mathcal{D}_{\max}(D)$. Every closed extension of $D$ contained in $D_{\max}$ has a domain lying between these two canonical domains.

The freedom in choosing an intermediate extension is encoded by the quotient of the maximal domain by the minimal domain \cite{bru1992hilbert,jayasinghe2024witten, jayasinghe2026spectralasymmetrysupersymmetryequivariant}. We define the \textbf{Gelfand–Robbin quotient} to be
\begin{equation}
    \mathcal{GR}(D) := \mathcal{D}_{\max}(D)/\mathcal{D}_{\min}(D),
\end{equation}
which is a Hilbert space with the quotient graph norm. We denote the quotient map by
\begin{equation}
    q_D:\mathcal{D}_{\max}(D)\longrightarrow\mathcal{GR}(D).
\end{equation}
Every closed intermediate extension of $D$ determines a closed subspace $W\subset \mathcal{GR}(D)$. Conversely, every closed subspace $W\subset \mathcal{GR}(D)$ determines a closed extension of $D$ as follows. We define the  \textbf{predomain} 
\begin{equation}
    \mathcal{D}'_W(D)=q_D^{-1}(W),
\end{equation}
where $q_D^{-1}$ denotes the set theoretic preimage, and the domain $\mathcal{D}_W(D) = (\mathcal{D}'_W(D))^{**}$.  We refer to such a choice of $W$ as an  \textbf{ideal boundary condition} following Cheeger.

Let $D^*$ denote the formal adjoint of $D$. For $u\in\mathcal{D}_{\max}(D)$ and $v\in\mathcal{D}_{\max}(D^*)$, there is a natural \textit{Green-Stokes pairing} (sometimes called the boundary pairing) defined by
\begin{equation}
    \omega_D(u,v) = \langle D_{\max}u,v\rangle_{L^2} - \langle u,(D^*)_{\max}v\rangle_{L^2}.
\end{equation}
This pairing vanishes whenever either argument belongs to its corresponding minimal domain, and hence descends to a continuous pairing
\begin{equation}
    \omega_D: \mathcal{GR}(D)\times\mathcal{GR}(D^*) \longrightarrow\mathbb C.
\end{equation}
The induced pairing is nondegenerate and identifies the (ideal) boundary data of $D$ with the dual (ideal) boundary data of its formal adjoint.
Given a closed subspace $W\subset \mathcal{GR}(D)$, define its annihilator with repsect to $\omega_D$ by
\begin{equation}
    W^{\perp} := \left\{ [v]\in\mathcal{GR}(D^*): \omega_D([u],[v])=0 \quad\forall [u]\in W \right\}.
\end{equation}
Then the Hilbert-space adjoint of the extension of $D$ with domain $\mathcal{D}_W(D)$ is the closed extension associated to $W^\perp$, whose domain is given by 
\begin{equation}
    \mathcal{D}_{W^\perp}(D^*) = (\mathcal{D}'_{W^\perp}(D^*))^{**} = (q_{D^*}^{-1}(W^\perp))^{**}.
\end{equation}
Given a first order symmetric operator $Q$, the Green-Stokes pairing defines a symplectic form. In this context, Lemma 3.3 of \cite{booss1998maslov} (c.f. Proposition 2.14 of \cite{bera2025diracoperatorstwistedramified}) shows that if $D$ is a \textit{Dirac-type} operator, any closed extension arises from picking a closed subspace of the Gelfand-Robbin quotient, and if the closed subspace is a Lagrangian with respect to the symplectic form $\omega_Q$, we have that the extension is self-adjoint.

If we are interested in these choices at the level of the harmonic sections, then it suffices to study the \textbf{cohomological quotient} 
\begin{equation}
    \mathcal{H}(Q):= q_Q(\text{ker}(Q_{\max})) \cap\mathcal{GR}(Q)
\end{equation}
which we define thus only for \textit{symmetric operators} $Q$. Since the kernel of an operator is a closed subspace, we see that $\mathcal{H}(Q)$ is closed and inherits a symplectic pairing from $\omega_Q$, which we continue to denote by that by abuse of notation.
Thus a closed Lagrangian subspace of $\mathcal{H}(Q)$ gives a closed subspace of the closed subspace $\mathcal{GR}(Q)$, but the closed extension corresponding to this domain may not necessarily be self-adjoint.

We now consider the corresponding construction for Hilbert complexes in the sense of \cite{bru1992hilbert}.

Let $\mathcal{P} = (H, P)$ be a pre-Hilbert complex (which is the data of a Hilbert complex without choices of domains for the operator $P$, which can be thought of as acting distributionally) arising from a differential complex on $X^{\mathrm{reg}}$. The maximal extension satisfies,
\begin{equation} 
    P_k\bigl(\mathcal{D}_{\max}(P_k)\bigr) \subset \mathcal{D}_{\max}(P_{k+1}), \qquad P_{k+1}P_k=0.
\end{equation} 
and the minimal extensions satisfies $\mathcal{D}_{\min}(P_k) \subset \mathcal{D}_{\max}(P_k)$. Applying the construction described above in each degree, we obtain the Gelfand-Robbin quotients 
\begin{equation} 
    \mathcal{GR}^k(P) = \mathcal{D}_{\max}(P_k)/\mathcal{D}_{\min}(P_k). 
\end{equation} 
Since the minimal complex is a subcomplex of the maximal complex, the differentials induce well-defined maps 
\begin{equation} 
    \widetilde P_k: \mathcal{GR}^k(P) \longrightarrow \mathcal{GR}^{k+1}(P), \qquad \widetilde P_k[u]=[P_k u]. 
\end{equation} 
Moreover, the induced maps are bounded with respect to the quotient graph norms and satisfy $\widetilde P_{k+1}\widetilde P_k=0$. Thus we obtain the \emph{Gelfand-Robbin complex} $\mathcal{GR}(\mathcal{P})$: 
\begin{equation} 
    \cdots\longrightarrow \mathcal{GR}^{k-1}(P) \xrightarrow{\widetilde P_{k-1}} \mathcal{GR}^{k}(P) \xrightarrow{\widetilde P_k} \mathcal{GR}^{k+1}(P) \longrightarrow\cdots. 
\end{equation}

We next introduce a class of ideal boundary conditions determined by the cohomology of the Gelfand-Robbin complex. For each degree $k$, consider the space 
\begin{equation}
    \mathcal H^k(P) := \left( \ker(P_k)_{\max} \cap \ker(P^*_{k-1})_{\min} \right) \cap \mathcal{D}_{\min}(P_k)^{\perp_{L^2}},
\end{equation}
where $(P^*_{k-1})_{\min}$ is the minimal domain for the formal adjoint of $P_{k-1}$, which we again refer to as the \textbf{cohomology quotient}.
The quotient map identifies $\mathcal H^k(P)$ with
\begin{equation}
\label{equation_GR_cohomology}
    H^k(\mathcal{GR}(\mathcal P)) = \frac{\ker\widetilde P_k} {\operatorname{im}\widetilde P_{k-1}}.
\end{equation}
The Green-Stokes pairing is compatible with the differentials of the Gelfand-Robbin complexes of $P$ and its formal adjoint $P^*$, and descends to a pairing on the cohomology of the Gelfand-Robbin complexes. 

A \textbf{cohomological predomain} is a graded closed linear subspace $W \subset\mathcal H (P)$ (which inherits the grading), associated with a degreewise operator predomain 
\begin{equation} 
    \mathcal{D}'_W(P) := \mathcal{D}_{\min}(P) \oplus W, 
\end{equation} 
and the corresponding closed domain $\mathcal{D}_W(P) :=  (\mathcal{D}'_W(P))^{**}$. Since $P_ku \subset \ker(P_{k,\max})$ for $u \in W$, we have 
\begin{equation} 
    P_k\bigl(\mathcal{D}'_W(P_k)\bigr) \subset \mathcal{D}_{\min}(P_{k+1}) \subset \mathcal{D}'_W(P_{k+1}). 
\end{equation} 
Passing to the graph closures, it follows that 
\begin{equation} 
    P_k\bigl(\mathcal{D}_W(P_k)\bigr) \subset \mathcal{D}_W(P_{k+1}). 
\end{equation}

Thus a cohomological predomain determines a closed intermediate Hilbert complex 
\begin{equation} 
    \mathcal P_W = (H,\mathcal{D}_W(P),P). 
\end{equation} 
A choice of ideal boundary conditions $W$ for $P$ also determines an extension of the Supersymmetry operators $Q=Q_{1}= P+P^*$ and $Q_2=i(P-P^*)$, with their rolled up extensions that we denote by $\mathcal{D}_W(Q_k):=\mathcal{D}_W(P)\cap \mathcal{D}_{W^\perp}(P^*)$ for both $k=1,2$, where $W^\perp$ is the annihilator of the Green-Stokes pairing.

It is easy to follow the definitions and see that 
\begin{equation}
    \mathcal{H}(P) \oplus\mathcal{H}(P^*) \subseteq \mathcal{H}(Q)
\end{equation}
is an inclusion of a closed subspace, where the restriction of the symplectic form $\omega_Q$ admits a natural splitting into maximally isotropic subspaces $\mathcal{H}^k(P)$ and $\mathcal{H}^k(P^*)$.
Rolled up extensions correspond to picking subspaces $W$ within the maximal isotropic space $\mathcal{H}(P)$, and taking $W^{\perp}+W$ gives the maximal Lagrangian extension for a subspace $W$ within $\mathcal{H}(D)$.
Since $W \subseteq \mathcal{H}(P) \subseteq \mathcal{GR}(P)$ is an inclusion of closed subspaces, and so is $W^{\perp} \subseteq \mathcal{H}(P^*) \subseteq \mathcal{GR}(P^*)$, we have that the Lagrangian $W\oplus W^{\perp}$ is a closed subspace of the $\mathcal{H}(Q)$, and in particular $\mathcal{GR}(Q)$.

Recall that a core subspace for a domain of any first order operator is a linear subspace which is dense in the graph norm, and a \textbf{core subcomplex} of a Hilbert complex is one which is dense in the complex in the graph norm. In particular the cohomology of a complex is isomorphic to the cohomology of a core subcomplex (see Theorem 2.12 of \cite{bru1992hilbert}).

Thus we see that the cohomological predomain $\mathcal{D}'_W(P)$ is a dense core and can be used to compute the cohomology of $\mathcal{D}_W(P)$. We formalize what we have proven here as follows.

\begin{proposition}\label{prop:cohom-splitting}
    Given a $\mathcal{N}=2$ SUSY pre-Hilbert complex associated to a resolved stratified pseudomanifold $X$, determined by the pre-Hilbert complex $\mathcal{P} = (H=L^2(X;
    \mathcal{E}),P)$, a choice of ideal boundary condition $W \subseteq \mathcal{H}(P)$ determines Hilbert complex $\mathcal{P}_{W}$ with domain $\mathcal{D}_W(P)$, for which the complex restricted to the cohomological pre-domain $\mathcal{D}'_W(P)=W \oplus \mathcal{D_{\min}}(P)$ yields a core subcomplex $\mathcal{P}'_{W}$ of the Hilbert complex $\mathcal{P}_{W}$, and 
    $\mathcal{H}^k(\mathcal{P}_W)=\mathcal{H}^k(\mathcal{P}'_W)$.

    In particular, 
    \begin{equation}       
    \mathcal{H}^k(\mathcal{P}_W)=\mathcal{H}^k(\mathcal{P}'_W)=\mathcal{H}^k(\mathcal{P}_0)\oplus W^k
    \end{equation}
    where $\mathcal{P}_0=\mathcal{P}_{\min}$ is the Hilbert complex with the minimal domain for $P$.
\end{proposition}
Here we use the $\mathcal{P}_0=\mathcal{P}_{\min}$ to indicate the complex with the minimal extension for $P$, which arises if one picks $W=\{0\}$.

This Proposition gives a handle on the cohomology groups, both for local and global complexes, and its utility is clear in the computations in Subsection \ref{subsection_non_isolated_singularity}.

\subsubsection{Iterated ideal boundary conditions at resolved strata} \label{subsubsec:ideal-bc-stratified}

For Dirac-type operators on stratified spaces with wedge metrics, the behavior of the domains near a singular stratum is closely related to the Witt and spectral Witt conditions, and we refer to \cite{Albin_2017_index} for more details.

On spaces $X$ with wedge metrics we have the inclusion
\begin{equation}
    \label{equation_inclusions_of_function_spaces}
    \rho_X H^1_{{e}}(X;E)
\subseteq \mathcal{D}_{\min}(D)
\subseteq \mathcal{D}_{\max}(D)
\subseteq H^1_{{e}}(X;E).
\end{equation}
where $\rho_X$ is a collective boundary defining function and 
\begin{equation}
H^1_{{e}}(X;E)= \{u \in L^2_{{w}}(X;E) :
Vu \in L^2_{{w}}(X;E)
\text{ for all } V \in C^\infty(X;{}^{{e}}TX)\}.
\end{equation}
is the edge Sobolev space.

This suggests the following natural domain corresponding to \textit{cutting the spectrum down the middle}, and was introduced in \cite{Albin_2017_index} (see Definition 2.3 of \textit{loc. cit.}.

\begin{definition}[VAPS domain]
\label{definition_VAPS_domain}
The \textbf{vertical APS domain} of $D$ is the graph
closure of $\rho_X^{1/2}H^1_{\mathrm{e}}(X;E)
\cap \mathcal{D}_{\max}(D_X)$, that is,
\begin{multline}
\mathcal{D}_{VAPS}(D)= \{u \in L^2_{w}(X;E) : \exists (u_n) \subseteq \rho_X^{1/2}H^1_{e}(X;E) \cap \mathcal{D}_{\max}(D) \\ \textit{s.t. } u_n \to u \textit{ and } (Du_n)
\textit{ is } L^2_{w}\textit{-Cauchy} \}.
\end{multline}
\end{definition}

\begin{definition}[Witt and spectral Witt conditions]
\label{definition_witt_spectral_witt}
The operator $(\mathcal{D}_{VAPS}(D),D)$
satisfies the \textbf{geometric Witt condition} 
\begin{equation}
    Y \in \mathcal{S}(X),\, y \in Y
\implies
\text{Spec}(D_{Z_y}) \cap (-\tfrac{1}{2},\tfrac{1}{2})= \emptyset.
\end{equation}

If, instead, we only require
\begin{equation}
    Y \in \mathcal{S}(X),\, y \in Y \implies \text{Spec}(D_{Z_y}) \cap \{0\} = \emptyset,
\end{equation}
then we say that $(\mathcal{D}_{VAPS}(D),D)$ satisfies the \textbf{Witt condition}.
Tacitly, we take the domain of the links to be the VAPS domains.
\end{definition}

This is Definition 2.4 of \cite{Albin_2017_index}. In particular, it shows that as one scales the links to be smaller and there are compatible Dirac-type operators, then they are essentially self-adjoint.

We now give a more systematic way to understand the choices of ideal boundary conditions for complexes with $\mathcal{N}=2$ SUSY.

Let $X$ be a resolution of a stratified pseudomanifold.
Suppose first that $X$ has depth one, and let $Y$ be a boundary fibration corresponding to a singular stratum 
\begin{equation}
    Z\longrightarrow H\xrightarrow{\phi}Y
\end{equation}
on the resolution $X$. The Gelfand–Robbin complex may be understood in terms of the asymptotic data of sections near $Y$, akin to Cauchy data of a boundary value problem. The admissible choices of boundary data depend on the geometry of the stratum and the asymptotic behavior of the differential operator. We refer the reader to \cite{Albin_hodge_theory_cheeger_spaces,Albin_hodge_theory_on_stratified_spaces} for similar descriptions for the de Rham and signature operators, which were used in \cite{jayasinghe2024witten}. 

We are primarily interested in the twisted Dolbeault complex $P=\bar\partial_E$, but our discussion extends to give a different perspective on the work in \cite{Albin_hodge_theory_cheeger_spaces,Albin_hodge_theory_on_stratified_spaces} with mezzo-perversities. Our choices are at the level of equivalence classes of domains that have the same cohomology.

In this setting, the Gelfand–Robbin complex carries additional geometric structure arising from the complex geometry of the strata and their links. Here we explain a procedure to pick certain domains at the level of cohomology, which in practice works well with analytic metrics such as on the singular projective variety we study in subsection \ref{Example_singular_toric_variety_domains}.

Suppose that the cohomological boundary data along $Y$ form a vector bundle
\begin{equation}
    H_Y(P)\longrightarrow Y
\end{equation}
where
\begin{equation}
        H_Y(P) = \bigsqcup_{y\in Y^{\mathrm{reg}}} \mathcal{H}(P)_{C_x(Z_y)}
\end{equation}
where $C_x(Z_y)$ is the normal fiber over $y$ equipped with the model metric at $x=0$ (i.e. the product type part of the metric in \eqref{equation_wedge_metric_pt}). Here $\mathcal{H}(P)_{C_x(Z_y)}$ corresponds to the cohomology quotient for the operator $P$ induced on the cone (also called the normal operators in geometric microlocal analysis, see \cite{Albin_hodge_theory_cheeger_spaces,Albin_hodge_theory_on_stratified_spaces}). The vertical Dirac operator acts on a vertical Clifford module, as in Subsection \ref{subsection_Dirac_Operators} (see also Subsection \ref{subsection_discreteness_of_spectrum}).

Since we assume there are local trivializations on any chart of $\widehat{Y}$ where the metrics are of product type as in \eqref{equation_wedge_metric_pt}, and owing to the smoothness of the variation of the these model metrics of the normal cones over the regular part $\widehat{Y}^{reg}=Y^{reg}$ of any stratum $Y$, and any other bundles describing operators $P,P^*$, the harmonic representatives of Dirac type and Laplace type operators vary smoothly along $Y$. However, the integrability of the harmonic sections change and this dictates whether or not they are in the maximal domain. It is known that for a single cone, these choices are governed by the eigenvalues of the induced operator $A$ on the link, which in turn are governed by the geometry of the link (see \cite{Albin_2017_index}). We give an example to explain this.

\begin{remark}
\label{remark_non_product_type_cohomologies_complicated}
    Consider a metric $dx^2+x^2d\theta^2|y|^2+dy^2$ on $C_x(S^1_{\theta})\times \mathbb{C}_y$, where the 1-dimensional volumes of the circular links grow larger as one moves away from $0$, and one can check that at $|y|$ increases and the eigenvalues of the Laplacian on the circle become smaller, there are more meromorphic sections $(xe^{ik\theta})^{-k}$ for $k>0$ that are $L^2$ bounded.
    Thus for $<0|y|<2$, the cohomology quotient is trivial, but for $2\leq|y|< 3$, it is one dimensional, and for $3 \leq |y|<4$ it is four dimensional, and so on. 
    Now the meromorphic sections corresponding to the $\overline{\partial}$ operator do not change with the metric if we use the usual complex structure on $\mathbb{C}^2$, but the Dolbeault-Dirac operator and more general cohomology will change accordingly.

    Here we see that if we restrict to fibrations over the annulus $3\leq |y| < 5$, there is a further stratification of the space by annuli, $3 \leq |y|<4, 4\leq |y| <5$ where the cohomology quotients are rank 3 and 4 bundles. It is easy to see that the rank 3 bundle extends as a sub-bundle over the entire annulus.
\end{remark}

More generally, the eigenvalues of the operators $A_y$ over each link varies continuously and the eigenvalues of the normal operators change continuously. Thus the dimension of the cohomology quotients are lower semi-continuous piecewise constant functions. We can find a minimal rank bundle of sections which are integrable along the entire stratum and span the cohomology quotients along each fiber. These are locally trivial bundle over any chart, and thus a finite rank vector bundle $H_Y(P)$ on the boundary hypersurface $Y$. Since the sections are in the kernel of the Dirac operator on the fibers, it is even a holomorphically flat Hermitian bundle with the metric induced from the $L^2$ inner product.

Now there are (at least) two ways to proceed. One way, along the lines of \cite{Albin_hodge_theory_cheeger_spaces,Albin_hodge_theory_on_stratified_spaces}, is to consider sections of this bundle that are flat with respect to the operator that $P+P^*$ induces on the base $Y$. This suits the broader considerations in the work in \textit{loc. cit.}, where all self-adjoint extensions of the operators are characterized (albeit for metrics satisfying certain conditions on the volumes of the links). {\textbf{The approach we take does not characterize all choices of closed extensions of $P$, instead characterizing them after identifying those that have the same cohomology}}.

We consider the sections $S$ of the bundle $H_Y(P)$ over $Y$ that satisfy $Q_Y(S)=0$. Here $Q_Y$ is the induced tangential operator on the bundle by $Q$, by restricting the derivations in the first order operator using a frame tangential to $Y$, and where the frame in the normal directions are considered to be $0$. We consider the action of this operator in the distributional sense over $Y$.
This is standard, and can be made explicit using the techniques in Subsection \ref{subsection_polyhomogenous_estimates}.
More simply, we can consider a local frame $\{\omega_k\}_k$ for the bundle, and then consider sections $S=\sum_k s_k\omega_k$ where $s_k$ are smooth functions with compact support on any local chart of $Y$, where $Q_Y(S):=Q(S)$, for $S$ considered as a section of the Hilbert space $L^2(X;\mathcal{E})$.
The space of such sections is what we call $\mathcal{H}_Y(P)$.

We next consider suitable subspaces,
\begin{equation}
    W_Y \subseteq \mathcal{H}_Y(P)
\end{equation}
which are \textbf{flat} in the sense that it is spanned over $\mathbb{C}$ by a (Schauder) basis $\{\omega_k\}_{k \in \mathbb{N}}$ where each $\omega_k$ is in the kernel of $P_Y$.

These subspaces correspond to the linear subspaces $W \subseteq \mathcal{GR}(P)$, which are the ideal boundary conditions for the Hilbert complex.

We now consider a depth two stratification $Y^2 \subset Y^1 \subset \widehat{X}$ 
and let $Z^i - H^i \overset{\phi_i}{\longrightarrow} Y^i$ denote the boundary fibration where we label the boundary hypersurfaces according to the strata.
At depth one, we pick an ideal boundary condition $W_1 \subseteq \mathcal{H}_{Y^{1}}(P)$ by repeating the construction described earlier over $Y^1 \setminus Y^2 \subset \widehat{X}$ and define 
\begin{multline}
    \mathcal{D}_{\max,W_1}(P) = \{ u \in L^2(X;\mathcal{E}) : Pu \in L^2(X;\mathcal{E}),\\
    \; \chi u \in \mathcal{D}_{W^1}(P_{\widehat{X} \setminus Y^2}), \; \forall \chi \in C^\infty_{c}(\widehat{X} \setminus Y^2) \}
\end{multline}
where the domain $\mathcal{D}_{W_1}(P_{\widehat{X} \setminus Y^2})$ is well defined for all sections supported away from $Y^2$, as the graph closure of $\mathcal{D}_{\min}(P_{\widehat{X} \setminus Y^2})\oplus W_1$. 
In particular, over any cone $C(Z^2)$ over $Y^2$, this determines a space of sections that is dense in $L^2$ since the singularities of $\widehat{Z^2}$ are those of $Y^1$ away from the cone point.
We denote this core as $\mathcal{D}_{W^1}(P_{C_{x}(Z^2)})$.
Using this we can define two domains, the first of which is a maximal extension
\begin{multline}
    \mathcal{D}_{\max,W_1}(P)_{C_{x}(Z^2)} = \{ u \in L^2(C_{x}(Z^2);\mathcal{E}|_{C_{x}(Z^2)}) : Pu \in  L^2(C_{x}(Z^2);\mathcal{E}|_{C_{x}(Z^2)}),\\
    \; \chi u \in \mathcal{D}_{W^1}(P_{C_{x}(Z^2)}), \; \forall \chi \in C^\infty_{\Phi,c}(C_x(Z^2)) \}
\end{multline}    
where $C^\infty_{\Phi,c}(C_x(Z))$ denotes the smooth sections with compact support away from the cone point at $x=0$. The second is a minimal extension that we call $\mathcal{D}_{\min,W^1}(P)$, uniquely characterized as the graph closure of the core $\mathcal{D}_{W^1}(P_{C_{x}(Z^2)})$ with respect to the graph norm of $P$.
Using these we define the quotient space
\begin{equation}
    \mathcal{GR}_{C_x(Z^2)}(P;W_1):=\frac{\mathcal{D}_{\max,W^1}(P)}{\mathcal{D}_{\min,W^1}(P)}
\end{equation}
and a quotient complex $\mathcal{GR}^{\bullet}_{C_x(Z^2)}(P;W_1)$ and similar to \ref{equation_GR_cohomology} we take
\begin{equation}
     \mathcal{H}^{\bullet}(P)_{C_x(Z^2)}:=\mathcal{H}^{\bullet} (\mathcal{GR}_{C_x(Z^2)}(P;W_1))
\end{equation}
and 
\begin{equation}
        H_{Y^2}(P) = \bigsqcup_{y\in (Y^2)^{reg}} \mathcal{H}^{\bullet}(P)_{C_x(Z^2)}
\end{equation}
We then take the flat sections of this bundle and define $\mathcal{H}_{Y^2}(P)$ and flat subspaces $W_{Y^2}$ of it similar to the case of depth 1. In this way we get an ideal boundary condition $(W_1,W_2)=:W$, and a corresponding domain.
For pseudomanifolds of higher depth, we can iterate this procedure.
In particular, the boundary data at a given stratum in general depends on the choices made at lower-depth strata, and we shall see this explicitly in Example \ref{Example_singular_toric_variety_domains}.

\begin{remark}[Iterated Gelfand-Robbin type quotients at strata]
    When the metric is $\delta>1/2$ asymptotic, the leading normal operator agrees with that of the corresponding product-type model (we postpone these computations until Subsection \ref{subsection_polyhomogenous_estimates}). By the edge calculus developed in \cite{Mazzeo_Edge_Elliptic_1,Mazzeo_Edge_Elliptic_2}, elements $u$ of the maximal domain admit a partial asymptotic expansion of the form
    \begin{equation}
    u\sim\sum_{\zeta_j}u_{\zeta_j}(y,z)x^{\zeta_j}+\widetilde u,
    \end{equation}
    where the sum is taken over the relevant indicial roots, and $\widetilde u$ is a remainder of higher asymptotic order. The coefficients $u_{\zeta_j}$ are, in general, distributional sections along the stratum, with
    \begin{equation}
    u_{\zeta_j}\in H^{-\operatorname{Re}\zeta_j}
    \bigl(Y;\ker I(xL;\zeta_j)\bigr).
    \end{equation}
    In particular, these coefficients depend continuously on $u$ with respect to the graph norm. We refer to \cite[Lemma 2.2]{Albin_hodge_theory_cheeger_spaces} for a precise formulation; see also \cite{Albin_hodge_theory_on_stratified_spaces}.

    The coefficients in this expansion constitute the asymptotic Cauchy data, and the choice of a closed domain is determined by imposing appropriate conditions on these data. For our purposes, we consider the coefficients corresponding to the $0$-eigenspace of the relevant cone operator, which requires some analysis to characterize properly.
\end{remark}

\begin{remark}[Comparison with mezzo-perversities]
    One can compare this procedure against that in \cite{Albin_hodge_theory_cheeger_spaces,Albin_hodge_theory_on_stratified_spaces} (c.f.  \cite{jayasinghe2024witten}). The cohomology of the Gelfand-Robbin complex for $P=d$, for an isolated conic singularity over a stratum, corresponds to the cohomology of the link in the degree $l/2$, where $l$ is the dimension of the link (thus trivial when the $l$ is odd, as in the case of K\"ahler varieties). The cohomology of $\mathcal{GR}(P^*=d^*)$ for the same is given by the forms $dx \wedge x^{l-1} \omega$ for $\omega$ the harmonic forms in degree $d-1$ of the cone.
    To make explicit an equivalence between the choices of domains at the level of cohomology obtained from both methods is out of the scope of this article.
\end{remark}

We make a few remarks before we study an example.

\begin{remark}[For Hamiltonians of 1-st order type]
    While it seems possible to extend our discussion for Hamiltonians of 1-st order type and the corresponding operators, for the purposes of the main results of this article it suffices to first pick equivariant Hilbert complexes of $0$-th order type corresponding to the maximal and minimal domains.
    This will depend on the choice of action. For torus actions, and those of more general compact connected Lie groups, we can find correspondences between domains for choices of circle actions. We do not explain this in detail in this article.    
\end{remark}

\begin{question}    
\label{question_cohomological_domains_dictated_at_critical_points}
When there are K\"ahler Hamiltonian actions or stratified Morse-Bott functions as in \cite{jayasinghe2024witten}, one has $\mathcal{N}=2$ SUSY complexes, and it seems plausible that the choices of extensions (at the level of cohomology) that we make are completely determined by choices of domains near the vanishing set $Z(A)$. Is this true?
\end{question}

\subsection{An example of a singular variety}

\begin{figure}[h]
    \centering
    \includegraphics[scale=0.25]{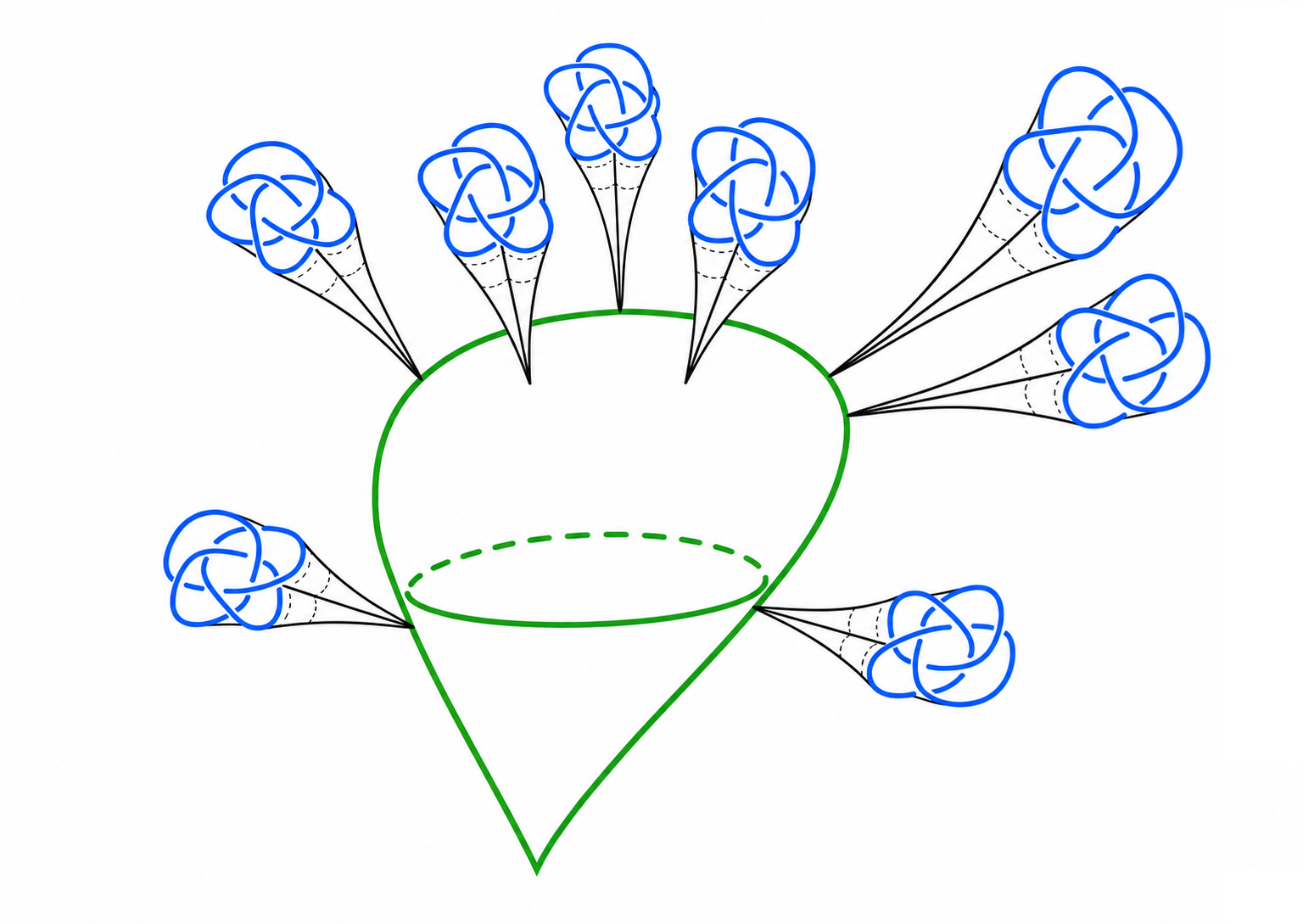}
    \caption{A schematic for a neighbourhood of the fibration over the regular part of $Y^1$, for the variety in Example \ref{Example_singular_toric_variety_domains}, also studied in Subsection \ref{subsection_non_isolated_singularity}.}
    \label{figure_singular_variety}
\end{figure}

\begin{example}
\label{Example_singular_toric_variety_domains}
    Consider the projective algebraic variety 
    \begin{equation}
    \widehat{V} = \{[U:X:Y:Z] \in \mathbb{CP}^3: Z^4 - X^3Y = 0 \},
    \end{equation}
    given in Figure \ref{figure_singular_variety}.
    Its singular locus is $\operatorname{Sing}(\widehat{V}) = \{[U:0:Y:0]\}$. On the affine chart $U=1$, with coordinates $x = \frac{X}{U}$, $y = \frac{Y}{U}$ and $z = \frac{Z}{U}$, the variety is $V = \{(x,y,z) \in \mathbb{C}^3: z^4 - x^3y = 0 \}$. The strata of the singular locus is:
    \begin{equation}
    \text{the depth 1 stratum } Y^1 = \{(0,y,0):y\neq 0\} \text{ and } \text{the depth 2 stratum } Y^2 = \{(0,0,0)\}.
    \end{equation}
    Locally near $q = \{(0,y_0,0)\} \in Y^1 \setminus Y^0$, the singularity is a product of the $Y^1$ direction with the transverse cusp $\Sigma = \{(u,v) \in \mathbb{C}^2: u^4 - v^3 = 0\}$, parametrized by $t \mapsto (u,v) = (t^3,t^4)$. 

    We use the section $\frac{u}{v} = t^{-1}$ to provide a closed domain lying strictly between the minimal and maximal domains. We first verify that $t^{-1} \in \mathcal{D}_{\max}(\overline{\partial}) \backslash \mathcal{D}_{\min}(\overline{\partial})$. The section $t^{-1}$ is unbounded as $t \to 0$ but $L^2$-integrable with respect to the pullback of the ambient Euclidean metric. Indeed,
    \begin{equation}
    g_\Sigma = (9|t|^4+ 16|t|^6) dt\otimes d\overline{t} = 9r^4 (1+\frac{16}{9}r^2)(dr^2 + r^2 d\alpha^2),
    \end{equation}
    where $t = r e^{i\alpha}, \alpha \in [0,2\pi]$. An analysis similar to that in Example 2.19 of \cite{jayasinghe2024holomorphic} shows that $|u|$ is a boundary defining function at the resolved singularity and that the link at $|u|=1$ in the model wedge will have angle $6\pi$, since $u=t^3$ and $\alpha \in [0, 2\pi]$.     
    We can check easily that
    \begin{equation}
    \Vert t^{-1} \Vert_{L^2}^2 = \int_0^{2\pi} \int_0^\varepsilon |re^{i\alpha}|^{-2} 9r^4 (1+\frac{16}{9}r^2) r dr d\alpha = 2\pi \int_0^\varepsilon r^3(9+16r^2) dr < \infty
    \end{equation}
    for finite values of $\varepsilon$.
    On $\Sigma^{reg}$, $t^{-1}$ is holomorphic, so $\overline{\partial} t^{-1} \in L^2$ distributionally. Hence $t^{-1} \in \mathcal{D}_{\max}(\overline{\partial})$.     
           
    For $u=\rho e^{i\alpha_1}=t^3$,  $g_\Sigma=(d\rho^2+\rho^2d\alpha_1^2)(1+16/9|\rho|^{2/3})$ and the model metric has angle $6\pi$, which shows that the section $t^{-1}$ has eigenvalue $-1/3$ with respect to the model operator on the circular link associated to the Dolbeault Dirac operator. Then the characterizations of the Gelfand-Robbin quotient in terms of the small eigenvalues on links (see for instance \cite{Albin_2017_index}) shows that $t^{-1}$ is not in the minimal domain.    

    In this same way we can check that $[t^{-1}], [t^{-2}]$ are distinct non-trivial elements in $\mathcal{GR}(\overline{\partial})$ for each fiber $\Sigma_{y_0}$ on $Y^1 \setminus Y^0$. We may choose a non-zero subspace $W_1$ that contains sections $s=s_1t^{-1}+s_2t^{-2}$ where $s_1,s_2$ are sections on $Y^1 \setminus Y^0$ that are in the kernel of the Dirac operator on $Y^0$. These are infinite dimensional since $Y^1 \setminus Y^0$ is not compact. 
    For instance the section $t^{-1}=u/v$ admits an extension $s=s_1t^{-1}+0t^{-2}=z/x$ where the section $s_1=1$. It also admits extensions $s^k_1t^{-1}=(1/y)^k(z/x)$, where $s^k_1=(1/y)^k$ is a holomorphic section away from $y=0$. Holomorphicity is precisely the flatness with respect to the Dolbeault Dirac operator on the base imposed in the procedure in Subsection \ref{subsubsec:ideal-bc-stratified}, when we restrict to sections of degree $0$.
    
    On the space $\widehat{V}$, the stratum $Y^1$ admits a smooth extension, and the extra choices are only to be made near $Y^0$, which we focus on now. 
    The first such choice we will make will be with $W_2$ the minimal extension, and we denote this as $\mathcal{D}_{W_{Z/X}}(P)$ where since $W_{Z/X}=(W_1,\min)$. The integral completion of the local ring at $Y^0$ has a presentation $R_1=R[1,z,xy/z,x^2y/z^2]$ as a module over the ring $R=\mathbb{C}[x,y]$, and the local cohomology at this singular point for the minimal domain is the Hilbert space completion of this ring.  
    
    It is clear that the section $(z/x)=Z/X=X^{-1/4}Y^{1/4}=t^{-1}Y^{1/4}$ on the regular part of $V$ admits a meromorphic extension near the entire singular stratum on $V$, and that this section is in the domain for $W=(W_1,\min)$ where we pick the meromorphic section $z/x$ to be in $W^1$.Moreover it is clear that it will remain $L^2$ bounded near $Y^0$.
    
    More precisely we can pick the section $Z/X$ to be flat sub-bundle in the Gelfand-Robbin quotient to form a pre-domain, and take the adjoint of the adjoint to get a closed extension. Thus it is in the graph closure of sections with support away from $Y^0$ which are in $\mathcal{D}_{W_1}$, and hence in $W_{Z/X}=(W_1,\min)$.

    We refer the reader to Subsection \ref{subsection_non_isolated_singularity} for computations of local cohomology groups and Morse polynomials for this domain, and the minimal domain, for some choices of circle actions with various fixed point sets. 
\end{example}

\subsection{Localizing domains for infinite localizing neighbourhoods} 
\label{subsec:localizing-domains}

Now it is straightforward to make the choice of domain in Definition \ref{definition_deformed_local_complexes_intro} for the Laplace-type operators precise.
We proceed by assuming that choices of self-adjoint extensions have already been made for the operators $Q_{j,\varepsilon}$. While we have given a systematic procedure for in the previous subsection for $\mathcal{N}=2$ SUSY complexes for any choice of $\varepsilon$, we shall initially work with families of $\mathcal{N}=N$ SUSY Hilbert complexes on $X$ which have Hamiltonians of $0$-th order type and potentials which have vanishing sets that have localizing neighbourhoods as we study in this article.
However we shall assume that the ideal boundary conditions for these complexes are \textit{local} in the sense that the choices are only relevent near the strata at which they are picked. This is true for the $\mathcal{N}=2$ complexes by out preceding constructions but is also know for symmetric operators $Q_{\varepsilon}$ (see, e.g. \cite[\S 3]{booss1998maslov}).
There is literature which presents conditions for Dirac-type operators to be essentially self-adjoint (see, e.g., \cite{Albin_2017_index}) and this analysis can be extended to the study of Witten deformations that we study for $\mathcal{N}=1$ complexes.

Given a domain $\mathcal{D}_W(Q_{j,\varepsilon})$ it determines a domain for $\Delta_{\varepsilon}$ as in Subsection \ref{subsec:abstract-complexes} which we also denote by $\mathcal{D}_W(\Delta_{\varepsilon})$. Then given an infinite localizing neighbourhood $U_{\infty}(F)$ for some $F \in Z(A)$ that we study, we equip it with the obvious extension of the product type metric in \eqref{equation_wedge_metric_pt}. Thus the ideal boundary condition near $F$ extends to $U_{\infty}(F)$ and we have an induced family of $\mathcal{N}=N$ SUSY Hilbert complexes 
\begin{equation}\label{eq:localized-complex}
    \mathcal{Q}_{W,\varepsilon}(U_{\infty}(F)) = \Big( L^2(U_{\infty}(F),\mathcal{E}), \{(\mathcal{D}_W(Q_{j,\varepsilon}),Q_{j,\varepsilon})\}_{j=1}^N \Big)
\end{equation}

where the $L^2$ boundedness of sections corresponds to what is known as a $L^2$ boundary condition at infinity. We call this the \textbf{localized complex} at $F \in Z(A)$.
The spectral theory results in Subsection \ref{Appendix_spectral_theory} then supply conditions which imply that these localized complexes have discrete spectrum and are Fredholm, at least for the cases where $\varepsilon>0$.
For $\varepsilon=0$, this is certainly not the case, as is apparent when looking at the spectral theory for the Hodge Laplacian on Euclidean space, which is a cone over a sphere.

\begin{remark}
    One point we now highlight is that because all of our complexes of interest are defined by beginning with an infinite localizing neighborhood $U_\infty(F)$ (c.f. Definition \ref{definition_resolved_fundamental_neighbourhood}), which is not a subset of $X$, but serves as a local model comparable to the normal bundle over $F$. As noted above, we equip $U_\infty(F)$ with the extension of the product-type metric $g_{w,pt}$ defined in \eqref{equation_wedge_metric_pt}; thus all $L^2$-based spaces are defined with respect to the (Riemannian) volume of $g_{w,pt}$. In order to compare $\mathcal{E}$-valued sections defined over \emph{actual} neighborhoods of $F$ in $X$, with these sections we try to model them by, we must then necessarily have some identification between the Clifford modules induced by the actual metric $g_w$ and our model metric $g_{w,pt}$; in keeping with our running notation we denote these two Clifford modules by $\mathcal{E}$ and ${}^{w,pt}\mathcal{E}$. We shall denote this isomorphism $\tau$, and discuss its role for our generalized Dirac operator in \eqref{eq:tau-isom-def}, and the accompanying subsection.

    In the smooth case, the corresponding canonical isomorphism has been employed in previous works \cite[eqn. (1.56)]{bismut1991complex} \cite[Prop 3.2]{wu1998equivariant}, and denoted explicitly in the latter case. In our case, because this identification is already used in the definition of our complexes, we shall typically omit $\tau$, as well conflating $\mathcal{E}$ and ${}^{w,pt}\mathcal{E}$, except when strictly necessary (e.g. when we are comparing sections of the two bundles as in \S \ref{subsection_polyhomogenous_estimates}.)
\end{remark}

\begin{remark}
If the volumes of the links of the normal conic fibers to each stratum is large enough, the eigenvalues of induced Dirac-type operators on the links are large and the Laplace type and many Dirac-type operators are essentially self-adjoint by the results in \cite{Albin_2017_index}, and Theorem 3.8 of \cite{booss1998maslov} can be used to see that the families of Dirac operators of the type $D+i\varepsilon cl(V)$ that we study are essentially self-adjoint, and the pre-Hilbert complexes extend to unique Hilbert complexes.    
\end{remark}

Consider $F \in Z(A)$ and let
\begin{equation}\label{eq:resolved-localizing-nbhd}
    \varphi_{F}: U_t(F)\longrightarrow F^{\mathrm{reg}}
\end{equation}
be its resolved localizing neighbourhood, described in definition \ref{definition_resolved_fundamental_neighbourhood}. For $y\in F^{\mathrm{reg}}$, write
\begin{equation}
    U_{a,y;t} = \varphi_{F}^{-1}(y)
\end{equation}
for the normal fiber. For $t<\infty$ its metric boundary will be denoted by $\partial U_{a,y;t}$. 
For the infinite cones, we consider the fiberwise cohomology groups. Given that the metric is of product type, we can pick an orthonormal frame $e_{i}$ of the wedge tangent bundle of the links $Z$ of the fibers $U_{a,y}=C_x(Z_y)$ over any suitable chart on $F^{reg}$ that vary smoothly in $y$, and we have a frame $\{\partial_x, \frac{1}{x}e_{i}\}$ using which we can define the \textit{normal Dirac operator}
\begin{equation}
    D_{C(Z)}:=cl(\partial_x) \nabla^\mathcal{E}_{\partial_x}+\sum cl\left(\frac{1}{x}e_{i}\right) \nabla_{\frac{1}{x}e_{i}}
\end{equation}
which also determines a \textbf{horizontal Dirac operator} $D_{F}:=D-D_{C(Z)}$ over $F^{reg}$, which can equivalently written in terms of a frame on $F^{reg}$.
This corresponds to the Dirac operators of \textbf{horizontal and vertical Clifford modules} $\mathcal{E_H}, \mathcal{E_V}$ as described in subsection \eqref{eq:clifford-mod-factorization}, when there is a decomposition of the wedge tangent bundle given by a fibration as we have here.

When there are Witten deformations of the types we study, we can write the deformed operators on localizing neighbourhoods as
\begin{equation}
    D_{\varepsilon}=D+i\varepsilon cl(V)=D_{F}+D_{C(Z),\varepsilon}, \quad D_{C(Z),\varepsilon}=D_{C(Z)}+i\varepsilon cl(V_a)
\end{equation}
where $V_a$ is the model vector field in \eqref{Morse_Bott_Kirwan_modify_hamiltonian_vf}.
For the case of $\mathcal{N}=2$ families arising from conjugation, we have
\begin{equation}
    \quad D_{C(Z),\varepsilon}= P_{C(Z),\varepsilon}+P^*_{C(Z),\varepsilon}
\end{equation}
where the operator $P_{C(Z),\varepsilon}$ is the corresponding operator on the cone, conjugated by $e^{\varepsilon h_a}$ for the function $h_a$ in Definition \ref{definition_non_degeneracy_critical_point_sets_Kahler}, restricted to the fibers.
It is easy to see that the Horizontal Dirac operators are not perturbed since $F$ is a critical point set and $D_{F}(h)=0$.
We direct the reader to Subsection \ref{subsection_discreteness_of_spectrum} for a more detailed study of the model horizontal and vertical operators.

Given we have ideal boundary conditions for the Dirac operators along the fibers, we have self-adjoint extensions for the operators $D_{C(Z)}$ and we have a family of $\mathcal{N}=1$ SUSY Hilbert complexes over each fiber
\begin{equation}\label{eq:fiberwise-cmpx}
    \mathcal{P}_{\varepsilon}(U_{a,y,\infty}) := \Big( L^2(U_{a,y,\infty};\mathcal{E_V}), \mathcal{D}_{W_y}(D_{C(Z_y),\varepsilon}), D_{C(Z_y),\varepsilon} \Big),
\end{equation}
and we call these the (family of) \textbf{fiberwise complexes}.
We can define the cohomology of such a complex in the usual way, as $\mathcal{H}(\mathcal{P}_{\varepsilon}(U_{a,y,\infty}))= \text{ker} \Delta_\varepsilon^V$ where $\Delta_{\varepsilon}^V=D_{C(Z),\varepsilon}^2$ with the induced domain. 
The case of general $\mathcal{N}=N$ complexes can be extended similarly, extending this description to all operators $Q_{j,\varepsilon}$, where the compatibility imposed for the domains over $X$ will correspond to compatibility for ideal boundary conditions.

Given such $\mathcal{N}=N$ complexes, the cohomology of the fiberwise complexes $\mathcal{H}(\mathcal{P}_{\varepsilon}(U_{a,y,\infty}))$ varies smoothly in the parameter $y \in F^{reg}$ and form a hermitian bundle with a holomorphically flat metric (see the discussion around Remark \ref{remark_non_product_type_cohomologies_complicated}), 
\begin{equation}\label{eq:normal-cohomology-bundle}
    \mathcal{H}^{\bullet}(\mathcal{P}_{\varepsilon}(U_{a,\cdot,\infty}))-\mathcal{N}_{F,\varepsilon}\rightarrow F^{reg}
\end{equation}
which we call the \textbf{normal cohomology bundle}.
The horizontal operator $D_{F}$ admits a unique self-adjoint extension as an operator over sections of the horizontal clifford module twisted by this bundle, determined by the ideal boundary conditions for $D$ on the localizing neighbourhood, which are in turn determined by those on $X$. We notate this domain as $\mathcal{D}_W(D_{F})$, and this the following complex.

\begin{definition}[Normal cohomology complex] \label{def_normal_cohomology}
    The normal cohomology complex over $F$ is the complex
    \begin{equation}
        \mathcal{P}_{W,\varepsilon}(\mathcal{N}_{F,\varepsilon}) = \Big(L^2(F^{reg}; \mathcal{E_H} \otimes \mathcal{N}_{F,\varepsilon}), \mathcal{D}_W (D_{F}), D_{F} \Big)
    \end{equation}
    where the domains and operators are as described above, and $\mathcal{E_H}$ is the horizontal Clifford module as described in \eqref{eq:clifford-mod-factorization}. Given a $\mathcal{N}=2$ complex, equivalently we have a Hilbert complex
    \begin{equation}
        \mathcal{P}_{W,\varepsilon}(\mathcal{N}_{F,\varepsilon}) = \Big(L^2(F^{reg}; \mathcal{E_H} \otimes \mathcal{N}_{F,\varepsilon}), \mathcal{D}_W (P_{F}), P_{F} \Big)
    \end{equation}
    where we abuse notation.
\end{definition}

Before formalizing the statement that the cohomology of the normal cohomology complex and the localized complex are isomorphic, we present versions of the objects we have defined here on truncated cones.

\subsection{Characterizations of cohomology from domains on truncated neighbourhoods} \label{subsec:cohomology-truncated-neighbourhoods}

We now extend the localized complexes on the infinite (non-compact) localizing neighbourhoods to those on truncated compact localizing neighbourhoods, where we need to impose boundary conditions. The point of this description is that the ideal boundary conditions inherited from the global complex and the boundary conditions imposed at the metric boundary play different roles. The former encode the choices of extension across the singular strata, whereas the latter replace the $L^2$ condition at infinity of the model problem. This is not new, and has been studied in the case where $F \in Z(A)$ are points in \cite{jayasinghe2023l2,jayasinghe2024holomorphic,jayasinghe2025supersymmetricharmonicoscillatorssingular}, and in the de Rham setting for non-isolated $F \in Z(A)$ in \cite{jayasinghe2024witten}, with and without Witten deformation.
The descriptions of the boundary conditions given in this sections work only for $\mathcal{N}=2$ SUSY complexes, where boundary conditions have been studied in \textit{loc. cit.}. There maybe extensions for $\mathcal{N}=1$ complexes with suitable boundary conditions but we do not study this, in this article.

Given a $\mathcal{N}=2$ SUSY complex on a resolved stratified space $X$, with a fixed ideal boundary condition $W$ we have, equivalently a Hilbert complex $\mathcal{P}(X) = (H=L^2(X;\mathcal{E}), \mathcal{D}_W(P), P)$. We use the same notation $W$ for its restriction to the singular boundary hypersurfaces of $U_{a,y;t}$. At $\partial U_{a,y;t}$ we impose an additional boundary condition $B$. In the applications below $B$ will be chosen from the generalized Neumann and Dirichlet boundary conditions associated with the Hilbert complex. For twisted Dolbeault complexes these correspond to the $\overline{\partial}$-Neumann boundary conditions and their adjoint conditions, and for de Rham complexes these are the absolute and relative conditions. We refer to Subsection 5.2.3 of \cite{jayasinghe2023l2}, Subsection 3.2 of \cite{jayasinghe2024holomorphic} and Subsection 3.3 of \cite{jayasinghe2024witten} for more details on these boundary conditions, including situations where they meet the strata.

We write
\begin{equation}
    \mathcal{D}_{\mathcal{B}}(P), \text{ where we use the notation } \mathcal{B} = (W,B),
\end{equation}
for the resulting domain, where it should be understood that the operator is the vertical operator, (as described in the previous subsection) on the neighbourhood ${U_{a,y;t}}$ and use the notation
\begin{equation}
    \mathcal{P}(U_{a,y;t}) := \Big( L^2(U_{a,y;t};\mathcal{E_V}), \mathcal{D}_{\mathcal{B}}(P_{U_{a,y;t}}), P_{U_{a,y;t}} \Big)
\end{equation}
for the corresponding \textbf{fiber complex}. \textbf{Now we drop the $t$ from the notation, with the understanding that our descriptions work for all $t \in (0,\infty)$}.
The normal fiber decomposes according to the \textit{attracting and expanding} directions of the Hamiltonian flow,
\begin{equation}
    U_{a,y}=U^s_{a,y}\times U^u_{a,y}.
\end{equation}
This is the assumption discussed in Remark \ref{Remark_clarify_2}, where the leading order term of the function restricted to the fiber is ``$h_{a}=(1/2) \rho_F^2 f=(1/2) \rho_1^2 f_1 + (1/2) \rho_2^2f_2$ where $f_1$ and $-f_2$ are both strictly positive functions on resolved links $Z_1, Z_2$ respectively, where $C(Z)=C(Z_1)\times C(Z_2)$ where $Z_1$ and $Z_2$ have an induced contact structures.'' Here $U^s_{a,y}=C(Z_1)$ and $U^u_{a,y}=C(Z_2)$.

Stratified Morse-Bott functions are assumed to have a similar structure (see \cite{jayasinghe2024witten}). Here the terms \textbf{attracting and expanding} correspond to the gradient flow of $h_a$ with respect to the metric gradient as opposed to the symplectic gradient, and correspond to the \textbf{stable and unstable} pseudomanifolds of this flow on the tangent cone. 

The boundary condition used in the localized complex is determined by this polarization: on the attracting factor we take the \textbf{generalized Neumann condition} associated with $P$, while on the expanding factor we take the complementary boundary condition associated with the adjoint complex, which we also refer to as the \textbf{generalized Dirichlet condition}. These domains induce domains for the corresponding SUSY operators $Q_1, Q_2$ and $\Delta$ as described abstractly earlier in this section, and the boundary conditions for the domain of the Laplace-type operator for the Dolbeault complex correspond to the $\overline{\partial}$-Neumaan conditions in the smooth setting. The references we gave earlier in this subsection have detailed expositions on this.

Given two Hilbert complexes $\mathcal{P}_1=(H_1,\mathcal{D}(P_1),P_1,\gamma), \mathcal{P}_2=(H_2,\mathcal{D}(P_2),P_2,\gamma)$ we define
\begin{equation}
    \mathcal{P}_1\widehat{\otimes}\mathcal{P}_2:=(H_1 \otimes H_2, \overline{\mathcal{D}(P_1)\otimes \mathcal{D}(P_2)}, P=P_1 \otimes 1+1 \otimes P_2, \gamma_1 \otimes 1+1 \otimes \gamma_2)
\end{equation}
where we take the graph closure of ${\mathcal{D}(P_1)\otimes \mathcal{D}(P_2)}$ with respect to $P$ to get the domain. 

Using this we denote the resulting boundary condition collectively by $\mathcal{B}$ and we have
\begin{equation} \label{eq_polarized_normal_complex}
    \mathcal{P}_{\mathcal{B}}(U_{a,y}) = \mathcal{P}_{W,N}(U^s_{a,y}) \widehat{\otimes} \mathcal{P}^{*}_{W^\perp,D}(U^u_{a,y}),
\end{equation}
Here $W^\perp$ denotes the adjoint boundary condition determined by $W$, and $N,D$ refers to the generalized Neumann and Dirichlet conditions referred to earlier. Here we assume that the Clifford module admits a factorization corresponding to the splitting of the wedge tangent bundle corresponding to the product factorization of $U_{a,y}$ (see Subsection \ref{subsection_Dirac_Operators}), and this defines the operators $Q_{j,\varepsilon}$ and $P$ on each fiber.

\textit{With this in hand, the \textbf{fiberwise complexes} for truncated fibers is complete, and the \textbf{normal cohomology bundle} and the \textbf{normal cohomology complexes} described for the infinite localizing neighbourhoods can be extended in a straightforward manner.} 

We collect the definitions we have given in this subsection and the proceeding one, notating them with uniform decorations.

\begin{definition}[localized, fiberwise, and normal cohomology complexes]
\label{Definition_all_localized_domains}
Consider a family of $\mathcal{N}=N$ SUSY Hilbert complexes $\mathcal{Q}_{W,\varepsilon}(X)= \Big(L^2(X,\mathcal{E}), \{(\mathcal{D}_W(Q_{j,\varepsilon}), Q_{j,\varepsilon})\}_{j=1}^N \Big)$, which has a Hamiltonian $\Delta_\varepsilon$ of $0$-th order type and a potential $A$.
If a given $F \in Z(A)$ has a localizing neighbourhood $U_{t}(F)$, then we have  \textbf{localized complexes}
\begin{equation}
    \mathcal{Q}_{\mathcal{B}, \varepsilon}(U_{t}(F)) = \Big(L^2(U_{t}(F),\mathcal{E}), \{ (\mathcal{D}_{\mathcal{B}}(Q_{j,\varepsilon}), Q_{j,\varepsilon}) \}_{j=1}^N \Big),
\end{equation}
defined in \eqref{eq:localized-complex}; we have \textbf{fiberwise complexes}
\begin{equation}
    \mathcal{P}_{\mathcal{B},\varepsilon}( U_{a,y,t}):= \Big(L^2(U_{a,y,t},\mathcal{E}_\mathcal{V}), \mathcal{D}_{\mathcal{B}}(D_{C(Z_y),\varepsilon}), D_{C(Z_y),\varepsilon} \Big),
\end{equation}
defined in \eqref{eq:fiberwise-cmpx}; and we have \textbf{normal cohomology complexes}
\begin{equation}\label{eq:normal-coh-cmpx}
    \mathcal{P}_{\mathcal{B},\varepsilon}\big(\mathcal{N}_{F,\varepsilon}) = \Big(L^2(F^{reg}; \mathcal{E}_{\mathcal{H}} \otimes \mathcal{N}_{F,\varepsilon}), \mathcal{D}_{\mathcal{B}} (D_{F}), D_{F} \Big),
\end{equation}
first defined in Definition \ref{def_normal_cohomology}. For $\mathcal{N}=1$, we define these for $t=\infty$ with $\mathcal{B}=(W,L^2)$ where $L^2$ denotes the $L^2$ boundary condition.
If locally there is $\mathcal{N}=2$ description, then for $t \in (0,\infty)$ this is defined with $\mathcal{B}=(W,B)$, where the boundary conditions $B$ are as described above.
\end{definition}

Let $\Delta_{\varepsilon}$ be the associated Hamiltonian, with the induced domain on the fiber. Whenever the boundary value problem is Fredholm (which is the case for twisted de Rham and Dolbeault complexes both with and without Witten deformation, the latter after passing to the $0$-th order equivariant complexes as we do), a simple consequence of the SUSY is that 
\begin{equation}
\label{eq_truncated_fiber_hodge}
    \frac{\ker P_{\varepsilon}}{\text{Im}P_{\varepsilon}} \cong \ker\!\left( \Delta_{\varepsilon} \right) = \ker\!\left( Q_{j,\varepsilon} \right)
\end{equation}
for any of the $\mathcal{N}=2$ SUSY complexes, where the latter equality holds for $\mathcal{N}=N$ complexes. This holds for any $\varepsilon \geq 0$.
Thus the cohomology determined by the local domain may equivalently be regarded as the space of harmonic sections satisfying the ideal boundary condition $W$ at the singular faces and the boundary condition $\mathcal{B}$ at the metric boundary, and we work with these harmonic sections.

\subsection{Spectral theory and cohomology of local and normal cohomology complexes}
\label{subsection_spec_cohom_normal}

We first show the correspondences between the cohomology groups for the complexes defined earlier, and then study the spectral properties of their Hamiltonians.

\begin{proposition}
\label{Proposition_cohomology_level_equivalences}
Given a family of $\mathcal{N}=2$ SUSY Hilbert complexes $\mathcal{P}_{\varepsilon}(X)$ arising from conjugation by an invertible bounded operator $M_{\varepsilon}$, then 
\begin{enumerate}
    \item $\mathcal{H}^k(\mathcal{P}_{w,\varepsilon}(X))=\mathcal{H}^k(\mathcal{P}_{W,\varepsilon'}(X))$ for all $\varepsilon, \varepsilon' \geq 0$,
    \item $\mathcal{H}^k\Big(\mathcal{P}_{\mathcal{B},\varepsilon}(U_{t}(F)) \Big) = \mathcal{H}^k\Big(\mathcal{P}_{\mathcal{B},\varepsilon'}(U_{t}(F)) \Big)$ for all $\varepsilon, \varepsilon' \geq 0$ for $t \in (0,\infty)$, and for all $\varepsilon, \varepsilon' \gneq0$ for $t =\infty$,
    \item $\mathcal{H}^k(\mathcal{P}_{\mathcal{B},\varepsilon}(U_{a,y,t}))=\mathcal{H}^k(\mathcal{P}_{\mathcal{B},\varepsilon'}(U_{a,y,t}))$ for all $\varepsilon, \varepsilon' \geq 0$ for $t \in (0,\infty)$, and for all $\varepsilon, \varepsilon' \gneq0$ for $t =\infty$,
     \item $\mathcal{H}^k(\mathcal{P}_{\mathcal{B},\varepsilon}(\mathcal{N}_{F,\varepsilon}))=\mathcal{H}^k(\mathcal{P}_{\mathcal{B},\varepsilon'}(\mathcal{N}_{F,\varepsilon'}))$ for all $\varepsilon, \varepsilon' \geq 0$ for $t \in (0,\infty)$, and for all $\varepsilon, \varepsilon' \gneq0$ for $t =\infty$.
\end{enumerate}
\end{proposition}

\begin{remark}
\label{Remark_21_on_dualities}
    We note that the cohomology of the Hilbert complexes for negative values of $\varepsilon$ correspond to dual cohomologies, corresponding to the opposite action chamber of the action, as studied in Subsection \ref{Subsection_Toric_actions}. This is also related to the time reversal symmetry symmetry that we discussed in Remark \ref{Remark_on_conjugation}.
\end{remark}

\begin{proof}
The key to understanding the first four numbered statements is that for each of the complexes for a fixed value of $t$, of which the cohomology groups are taken, there is an invertible morphism of the Hilbert complexes,
 \begin{equation}
    \begin{tikzcd}[column sep=large, row sep=large]
    \mathcal{D} (P_\varepsilon)
        \arrow[r, "P_\varepsilon"]
        \arrow[d, "M'_\varepsilon"', "\cong"]
    &
    \mathcal{D} (P_\varepsilon)
        \arrow[d, "M'_\varepsilon", "\cong"']
    \\
    \mathcal{D} (P_{\varepsilon'})
        \arrow[r, "P"']
    &
   \mathcal{D} (P_{\varepsilon'})
    \end{tikzcd}
\end{equation}
where we have expressed the domains of the operators without additional decorations which apply in each case. Here $M'_\varepsilon=M_\varepsilon M^{-1}_{\varepsilon'}$ is invertible.  
The identitities
    \begin{equation}
        P M_\varepsilon = M_\varepsilon P_\varepsilon, \quad P_{\varepsilon'} M'_\varepsilon = M'_\varepsilon P_\varepsilon
    \end{equation}
together with the definition of the conjugated domain shows that $M'_{\varepsilon}$ is an isomorphism of Hilbert complexes that is in fact a quasi-isomorphism. The induced map on cohomology is then an isomorphism. The complexes are defined for the conditions stated, as discussed in the previous subsections. 
For the complexes other than those in the first numbered statement, for different values of $t$ and for fixed $\varepsilon$, the results follows by unique continuation of harmonic sections and by restricting to smaller neighbourhoods where one can see that the boundary values are still satisfied. This is a straightforward exercise with the explicit descriptions of the harmonic sections on the fibers given in the literature cited earlier in the section on the boundary value problems.

To see how the particular choices of boundary condition enter because of the dynamics, we refer to \cite{jayasinghe2025supersymmetricharmonicoscillatorssingular}.
\end{proof}

\begin{remark}
In the case of $\mathcal{N}=2$ SUSY complexes arising by conjugation, the domain $\mathcal{D}_{\mathcal{B}}(Q_{\varepsilon})$  for the infinite neighbourhoods are not defined by conjugating the domain for $Q_{\varepsilon=0}$ on $U_\infty(F)$ by $M_\varepsilon=e^{\varepsilon h_a}$ since the Laplace-type operator for ${\varepsilon=0}$ does not have discrete spectrum on the infinite localizing neighbourhood. 

Instead, one can first take the conjugated domain for all of $X$, and then proceed as we have described above.
There are studies of continuity of Cauchy data (which applies to ideal boundary conditions as we study here) for continuous deformations (see, e.g.,  \cite[\S 3.3]{booss1998maslov}), but in our situation it is explicit.

It is similar for the normal cohomology complexes.
\end{remark}

We add a few brief notes here on the spectral theory, referring the reader to Subsection \ref{subsection-spectral-background} for more details. There we show in Theorem \ref{Theorem_spectral_discreteness_Schrodinger} that the localized complexes with infinite neighbourhoods have discrete spectrum, and in particular are Fredholm (passing to $0$-th order complexes as usual), when there is an ideal boundary condition $W$ for which the global complex on any compact stratified space $X$ has discrete spectrum. We note that the sections of the normal cohomology complexes are included in the local cohomology complexes. This is straightforward to see from the definitions of the complexes in the setting where the lengths of the fibers are infinite ($t=\infty$).

Then the results above for the fiber cohomologies and normal cohomology complexes can be used to see that the normal cohomology complexes where the fibers are truncated (and with boundary conditions $N/D$ imposed according to the dynamics of the flow) are Fredholm. In fact the conjugations give isomorphisms of Hilbert complexes as studied in \cite{bru1992hilbert}, and Lemma 2.17 of \cite{bru1992hilbert} then shows that there is discrete spectrum. This allows us to use the techniques we use in the next section in the extraction of the Witten subcomplexes as $\varepsilon$ grows.
We also refer the reader again to Subsection \ref{subsubection_Fredholm_spectral} for relevant results and remarks, and to \cite{jayasinghe2023l2,jayasinghe2024holomorphic} for more thorough expositions for the cohomology of the fibers.

\section{Constructing Witten complexes}
\label{section_Witten_complexes}

\subsection{Polyhomogenous estimates}
\label{subsection_polyhomogenous_estimates}

Let stratum $Y \in \mathcal{S}(X)$ of depth $k$ with a collar neighborhood $\mathscr{C}(\mathfrak{B}_Y) \cong [0,1)_x \times \mathfrak{B}_Y$ of $\mathfrak{B}_Y$ in $X$. Recall the form of the asymptotically $\delta$ wedge metric,
\begin{equation}
    g_w = dx^2 +x^2 g_{\mathfrak{B}_Y/Y} +\phi^*_Yg_Y +x^{2\delta}{\mathrm E}.
\end{equation}
where  $g_{w,pt} = dx^2 +x^2 g_{\mathfrak{B}_Y/Y} +\phi^*_Yg_Y$ is the model metric of the collar, the term $g_{\mathfrak{B}_Y/Y}\in \mathcal{C}^{\infty} (\mathfrak{B}_Y; S^2(^{w}T^*\mathfrak{B}_Y))$ restricts to an asymptotically $\delta$ wedge metric on the fibers of $\phi_Y$. More specifically the family of fiberwise wedge-diffeomorphisms $Z\to \phi_Y^{-1}(\{y\})\subset \mathfrak{B}_Y$ of $Z$ to the fibers of $\phi_Y$ are isometries of wedge-metrics for every $y\in Y^{reg}$ and moreover, the association $(X,Y)\mapsto g_{\mathfrak{B}_Y}(X,Y)$ for smooth wedge vector fields $X,Y\in \mathcal{C}^{\infty}(\mathscr{C}(\mathfrak{B}_Y/Y); ^{w}\ker (\phi_Y)_*)$, is nondegenerate. Also  
\begin{equation}
    x^{2\delta}\mathrm{E} = g_w - g_{w,pt} \in x^{2\delta}\mathcal{C}^{\infty} (\mathscr{C}(\mathfrak{B}_Y); S^2(^{w}T^*X))
\end{equation}
is the remainder term. Recall that off of these collar neighborhoods of the strata of depth $k$, the form of the metric is fixed by induction as an asymptotically $\delta$ wedge metric of lower depth. We denote by $\nabla^M$ (abbreviating \textit{Model}) the Levi-Civita connection on the collar with respect to the model metric $g_{w,pt}$ and by $\nabla$ the Levi-Civita connection with respect to the asymptotic $\delta$ wedge metric $g_w$. 

\begin{remark}[Notational convention]
   The index notation of the frames are as follows; we reserve Latin letters for vertical coordinates and Greek letters for horizontal coordinates. We refer to both indices with capital Greek letters.
\end{remark}
We choose submersion coordinates $\{w_A=(z_i,y_\alpha)\}_{A=(i,\alpha)}$ for the submersion $\phi_Y:\mathfrak{B}_Y\to Y$ with $\{z_i \}_i$ the fiber coordinates on the model fiber $(Z^{reg}, g_Z)$ and $\{y_\alpha\}_\alpha$ the coordinates of $Y^{reg}$. The local wedge vector fields $\{\partial_{y_\alpha}\}_\alpha$ on $\mathfrak{B}_Y$ are horizontal lifts of the corresponding coordinate vector fields of $Y$ and the $\{\frac{1}{x} \partial_{z_i}\}_i$ are vertical vector fields with respect to the submersion metric $x^2g_{\mathfrak{B}_Y/Y} +\phi^*_Y g_Y$ on the slice $\{x\}\times \mathfrak{B}_Y \xrightarrow{\phi_Y} Y$ and we refer to them as horizontal and {\bf link wedge Fermi fields}. We define the $g_{w,pt}$-orthogonal projections $p^{\mathcal V}:^{w}TX \to \ker(\phi_Y)_*$ and $p^{\mathcal H}: ^{w}TX\to (\ker(\phi_Y)_*)^{\perp_{g_{w,pt}}}$ onto vertical and horizontal wedge-Fermi fields respectively. The symbols $O_w(x^k\partial^{\mathcal H})$ and $O_w(x^k\partial^{\mathcal V})$ will denote expressions of horizontal and link model wedge Fermi fields accordingly with coefficient components vanishing up to order $x^k$ when $x\to 0^+$. 

\begin{remark}
\begin{itemize}
    \item[(a)] 
When $\delta>1/2$, the rate $k$ at which the symbols $O_w(x^k\partial^{\mathcal H})$ and $O_w(x^k\partial^{\mathcal V})$ vanish, is independent of the choice of coordinates and their associated Fermi wedge vector fields.

\item[(b)] The case $\delta=1/2$ includes the smooth case where the fields in $\{\frac{1}{x} \partial_{z_i}\}_i$ and $\{\partial_{y_\alpha}\}_\alpha$ would correspond to the spherical and tangential Fermi fields respectively.
\end{itemize}
\end{remark}

\begin{remark}
\label{rem:smooth_case_analogy}
The case of a stratum $Y$ where the space $\widehat{X}$ has smooth wedge metric $g_w$ is also covered since in Fermi coordinates, the wedge metric $g_w$ has expansion in $x$ along $Y$ of the preceding form with $\delta=1/2$ so that the model metric is $g_{w,pt} = g_N\oplus \pi^* g_F = dx^2 + x^2(\exp^N)^*g_{S(N)} + \pi^*g_F$ where $\pi^N:N\to F$ is the normal bundle map, $S(N)\subset N\to F$ is the sphere bundle, $\exp^N:N_t(F)\to U_t(F)$ is the exponential map of the ambient metric and $\pi=\pi^N\circ (\exp^N)^{-1}: U_t(F) \to F$ is the nearest point projection bundle map. Using the Gauss lemma for Fermi coordinates, the error term $E$ satisfies $E(\partial_x, \partial_x)= E(\partial_x, \frac{1}{x}\partial_{z_i}) = E(\partial_x, \partial_{y_\alpha})=0$ while $E(\frac{1}{x}\partial_{z_i}, \partial_{y_\alpha})= O(1)$ as $x\to 0^+$ in $\mathscr{C}(\mathfrak{B}_Y)$ for every $i,\alpha$, that is precisely our required assumptions stated for $g_w$ in \eqref{wedge_metric_assumption}. The two metrics coincide when restricted on $TX|_F$.
\end{remark}

Note that the radial flow corresponding to $\partial_x$ gives a geodesic ray corresponding to $\partial_x$ for the model metric $g_{w,pt}$ and the Christoffel symbols of $\nabla^M_{\partial_x}$ for the model frame vanish by \cite{Albin_2017_index} (see Section 1.2, specifically around Remark 1.4 of \textit{loc. cit.}) resulting in these frames being $\nabla^M$-parallel along the radial flow of $\partial_x$. An easy calculation also shows that $x[W_i, W_j]$ is a wedge vector field whenever $W_i, i=0,1,2$ are wedge vector fields. By the Koszul formula for the Levi-Civita connections $\nabla$ and $\nabla^M$, we obtain 
\begin{align*}
\begin{split}
g_w(\nabla_{W_0} W_1, W_2) - g_{w,pt}(\nabla^M_{W_0} W_1, W_2) & = W_0( x^{2\delta} E(W_1, W_2)) + W_1( x^{2\delta} E(W_0, W_2)) 
\\
&\quad - W_2( x^{2\delta} E(W_0, W_1)) + x^{2\delta-1} [E(x[W_0, W_1], W_2) 
\\
&\quad - E(x[W_0, W_2], W_1) - E(W_0, x[W_1, W_2])]
\end{split}
\end{align*}
and we conclude that, 
\begin{equation}
\label{eq:actual_Christofell_symbols_minus_model_Christoffel_symbols_rates}
g_w(\nabla_{W_0} W_1, W_2) - g_{w,pt}(\nabla^M_{W_0} W_1, W_2) = x^{2\delta-1}E_1(W_0,W_1, W_2), \quad \text{as $x\to 0^+$,} 
\end{equation}
for some $E_1(W_0,W_1, W_2)\in \mathcal{C}^{\infty} (\mathscr{C}(\mathfrak{B}_Y); \mathbb{R})$.

\begin{remark}
\label{rem:Parallel_transport_of_vectors_by_the_Levi_Civita_connection_of_a_wedge_metric_is_not_in_general_a_wedge_vector_field}
    Note that $\nabla_{W_0} W_1$ is not in general a wedge vector field. Likewise parallel transport of vectors by the Levi-Civita connection of a wedge metric is not in general a wedge vector field.
\end{remark}

We would like to compare $g_w$-parallel frames along the radial flow lines of $\partial_x$ and their Christoffel symbols to $g_{w,pt}$-parallel frames along the same flow lines and their Christoffel symbols as $x\to 0^+$. By Remark~\ref{rem:Parallel_transport_of_vectors_by_the_Levi_Civita_connection_of_a_wedge_metric_is_not_in_general_a_wedge_vector_field} the $g_w$-frames parallel are not in general wedge vector fields rendering formula \eqref{eq:actual_Christofell_symbols_minus_model_Christoffel_symbols_rates} useless. However the $g_{w,pt}$-parallel frames $[\partial_x, \frac{1}{x}\partial_{z_i}, \partial_{y_\alpha}]_{i,\alpha}$ along the flow lines of $\partial_x$ are wedge vector fields and the following lemma computes the asymptotics of the $g_w$-parallel frames relative to the $g_{w,pt}$-parallel frames.   

\begin{lemma}
\label{lemma_on_the_frame_expansions}
 Set $g_0 = g_{w,pt},\ g_1=g_w$ with $\nabla^0$ and $\nabla^1$ denote the associated Levi-Civita connections and let $g_0$-orthonormal wedge frame $[e_A^0]_A = [e_j^0, e^0_\alpha]_{j,\alpha}$, where $e_1^0=\partial_x$ and $[f_j^0:=x e_j^0]_{j\geq 2}$ is a $g_{\mathfrak{B}_Y/Y}$-orthonormal frame and $[e^0_\alpha]_\alpha$ is the horizontal lift of a $g_Y$-orthonormal frame. Let $[e^1_A]_A$ be a $\nabla^1$-parallel along the flow of $\partial_x$, $g_1$-orthonormal frame. The two frames are related by $[e^0_B]_B \cdot T = [e^1_A]_A$, where $T=T(x,w)\in \mathcal{C}^\infty(\mathscr{C}(\mathfrak{B}_Y);\mathrm{GL}(2n,\mathbb{R}))\footnote{The set $\mathrm{GL}(2n,\mathbb{R})$ is equipped with the trace norm $|\cdot|$.}$ is a transition matrix. 
\begin{itemize}
\item[(a)] If $\delta \geq 1/2$, there exist a uniform limit $T_0=T_0(w) \in \mathcal{C}^\infty(\mathfrak{B}_Y; \mathrm{O}(2n,\mathbb{R}))$ of $x\mapsto T(x,w)$ as $x\to 0^+$ and moreover a matrix $T_1 = T_1(w)\in \mathcal{C}^\infty(\mathfrak{B}_Y; \mathrm{GL}(2n,\mathbb{R}))$, when $\delta=1/2$ that fit into expansion,
\begin{equation}
\label{eq:transition_matrix_expansion}
T=\begin{cases} 
    T_0 + O(x^{2\delta})&, \quad \text{when $\delta>1/2$} \\
    T_0 + T_1 x +O(x^2)&,\quad \text{when $\delta=1/2$} 
\end{cases}
\end{equation}
when $x>0$. Also $T_0 \cdot T_1^t + T_1 \cdot T_0^t=0$.

\item[(b)] The $\mathrm{O}(2n,\mathbb{R})$-equivariant correspondence between the set of $\nabla^1_{\partial_x}$-parallel, $g_1$-ortho-normal frames $[e^1_A]_A$ and the set of $g_0$-orthonormal frames of the form $[e^0_j, e^0_\alpha]_{j,\alpha}\cdot (T_0+xT_1)$ with $T_0 \in\mathcal{C}^\infty(\mathfrak{B}_Y; \mathrm{O}(2n,\mathbb{R}))$ and $T_1 = T_1(w)$ as in (a). When $\delta>1/2$, the later set can be thought with $T_1=0$  
\end{itemize}
\end{lemma}

\begin{proof} 
   We can write the vectors $e^1_A$ in terms of the basis $e^0_B$ as $e^1_A = g_0(e^0_B, e^1_A) e^0_B$ where we treat $g_0(e^1_A, e^0_B)$ as a matrix, which shows that we have $T=[g_0(e^0_A, e^1_B)]_{A, B}=:[T_{AB}]_{A,B}$. We have to compute the asymptotic behavior of $x\mapsto T(x,w)$ as $x\to 0^+$. Let $S=[S_{AB}]_{A,B} :=[g_1( e_A^1, e_B^0)]_{A,B}$ so that $[e_A^0]_A = [e_B^1]_B \cdot S$. Then $S$ is the inverse of $T$ and
   \[
   \partial_x S_{AB} = g_1(\nabla_{\partial_x}^1 e_B^0, e_A^1) = T_{\Gamma A} g_1(\nabla_{\partial_x}^1 e_B^0, e_\Gamma^0) = x^{2\delta-1} T_{\Gamma A} E_1(\partial_x, e_B^0, e_\Gamma^0),
   \]
   for all $A,B$ that is,
   \begin{equation}
   \label{eq:inverse_transition_matrix_equation_1}
    \partial_x S = x^{2\delta-1} T^t\cdot E_1^t 
    \end{equation}
    where $E_1:= [E_1(\partial_x, e_E^0, e_\Delta^0)]_{E,\Delta}$ is a matrix with smooth entries. Also note that
    \[
    T-S^t = [g_0(e_A^0, e_B^1) - g_1(e_A^0, e_B^1)]_{A,B} = - x^{2\delta} E\cdot T,
    \]
    and by rearranging,
    \begin{equation}
    \label{eq:inverse_transition_matrix_equation_2}
    S^t =(I_{2n} + x^{2\delta} E) \cdot T,
    \end{equation}
    where $E:= [E(e_A^0, e_B^0)]_{A,B}$ also has smooth entries. Multiplying \eqref{eq:inverse_transition_matrix_equation_1} by $S^t$ and taking traces we obtain, 
    \[
    \frac{1}{2} \partial_x|S|^2 = x^{2\delta-1} \mathrm{Tr} E_1,
    \]
    from where we conclude that $\sup_{(x,w)\in (0,1)\times K}|S(x,w)|$ is finite and hence, by making use of \eqref{eq:inverse_transition_matrix_equation_2}, $\sup_{(x,w)\in (0,1)\times K }|T(x,w)|$ is finite in any relatively compact open subset $K\subset \mathfrak{B}_Y$. Multiplying equation \eqref{eq:inverse_transition_matrix_equation_1}  by its transpose equation, taking traces and using Cauchy-Schwartz, we obtain
    \[
    |\partial_x S|^2 = x^{4\delta-2}|E_1 T|^2\leq x^{4\delta-2}C|T|^2,  
    \]
    and
    \[
    |\partial_x T| \leq |\partial_x(I_{2n} +x^{2\delta}E)^{-1} \cdot S^t|+ |(I_{2n} +x^{2\delta}E)^{-1} \cdot \partial_xS^t| \leq C(1+ x^{2\delta-1})|T|     
    \]
    which show that the uniform limits $S_0(w):=S(0^+,w)$ and $T_0(w)=T(0^+,w)$ both exists and are orthogonal matrices, inverse to each other. By induction on the number of derivatives in $w$, we conclude that $|\partial_x \partial_w^\ell S|$ and $|\partial_x \partial^\ell_w T|$ are also bounded in $(x,w)\in (0,1)\times K$ for any $\ell$, hence $T_0, S_0$ are smooth in $w$. Also, for small $x$    
    \[
    \partial_x S = x^{2\delta-1} S\cdot (E_1\cdot(I_{2n} + x^{2\delta} E)^{-1})^t, 
    \]
    which shows that both $S_1(w):=(\partial_x S)(0,w)$ and $T_1(w):=(\partial_x T)(0,w)$ exist when $\delta \geq 1/2$ and vanish when $\delta>1/2$. This proves the claimed expansions. 
    
    When $\delta =1/2$, by using \eqref{eq:inverse_transition_matrix_equation_1} and \eqref{eq:inverse_transition_matrix_equation_2}, it follows that $T_1 = (E_1-E)|_{x=0}\cdot T_0$. From the equation \eqref{eq:actual_Christofell_symbols_minus_model_Christoffel_symbols_rates}, the matrix $E_1-E$ is skew-symmetric proving that $T_0\cdot T_1^t + T_1 \cdot T_0^t=0$.   
    
    For part (b) let $[e^1_A]_A$ and $[e^0_j, e^0_\alpha]_{j,\alpha}$ be frames as in part (a), $T$ be their transition matrix with its expansion \eqref{eq:transition_matrix_expansion}. Define the correspondence $[e^1_A]_A\mapsto [e^0_j, e^0_\alpha]_{j,\alpha}\cdot (T_0+x T_1)$. This is $\mathrm{O}(2n,\mathbb{R})$-equivariant with the right action. If $[e^1_A]_A$ and $[f^1_A]_A$ are both $g_1$-orthonormal and $\nabla^1_{\partial_x}$-parallel, wedge frames, both mapped to $[e^0_j, e^0_\alpha]_{j,\alpha}\cdot (T_0+xT_1)$ and with transition matrices $T$ and $F$ respectively, then $[e^1_A]_A\cdot T^{-1}F = [f^1_A]_A$ and since both frames are $\nabla^1_{\partial_x}$-parallel, we obtain $\partial_x(T^{-1}F)=0$. Since $T(0^+,w)=F(0^+,w)=T_0(w)$, we conclude that $T=F$. Hence the correspondence is 1-1.
\end{proof}

\begin{remark}
    \label{rem:remark_on_frame_expansions}
\begin{itemize}
\item[(a)] 
If $\mathfrak{B}_Y$ is orientable, then we choose $T\in \mathcal{C}^\infty(\mathscr{C}(\mathfrak{B}_Y); \mathrm{GL}^+(2n,\mathbb{R}))$ so that $T_0\in \mathcal{C}^\infty(\mathfrak{B}_Y; \mathrm{SO}(2n,\mathbb{R}))$.

\item[(b)] 
Multiplying by appropriate element of $\mathcal{C}^\infty(\mathfrak{B}_Y; \mathrm{SO}(2n,\mathbb{R}))$ from the right, we can always choose frame $[e^1_A]_A$ with $T_0 = I_{2n}$ and expansions of the form
\[
[e^1_A]_A =  [e^0_j, e^0_\alpha]_{j,\alpha} + O_w(x^{2\delta}\partial^{\mathcal V} + x^{2\delta}\partial^{\mathcal H})
\] 
in the case when $\delta>1/2$ or of the form
\[
[e^1_A]_A =  [e^0_j, e^0_\alpha]_{j,\alpha}\cdot(I_{2n} + T_1 x) + O_w(x^2\partial^{\mathcal V} + x^2\partial^{\mathcal H})
\]
in the case when $\delta=1/2$. In both cases, the transition matrix has expansion $T=I_{2n} + x^{2\delta}T_1  + O(x^{2\delta+1})$ and we set $[e_A^0]_A = [e_j^0, e^0_\alpha]_{j,\alpha} = :\tau [e^1_A]$. In the latter case, recall the vertical and horizontal $g_0$-orthogonal projection maps $p^{\mathcal V}$ and $p^{\mathcal H}$ onto vertical and horizontal wedge-Fermi fields. By assumption \eqref{wedge_metric_assumption} on the wedge metric $g_w$, the expression $T_1=(E_1-E)|_{x=0} = \begin{pmatrix}(T_1)_{ij} & (T_1)_{i\beta}\\ (T_1)_{\alpha i}& (T_1)_{\alpha \beta} \end{pmatrix}$ reveals that $T_{1ij}=T_{1\beta 1}=T_{11\beta}=0$, for all $i,j,\beta$. Set 
\begin{equation}
\label{eq:horozontal_lift}
H: T^{w}\mathscr{C}(\mathfrak{B}_Y) \to T^{w}\mathscr{C}(\mathfrak{B}_Y): H(e^0_A) := \begin{cases} e_i^0,& \quad \text{if $A=i$}\\ e_\alpha^0 + x T_{1j\alpha}e_j^0,& \quad \text{if $A=i$}.
\end{cases}
\end{equation}
The map $H$ is well defined and independent of the frame $[e_j^0, e_\alpha^0]_{j,\alpha}$ induced by the fibration $\phi_Y: \mathfrak{B}_Y \to Y$. Then 
\begin{align*}
e^1_A &= e^0_A + xT_{1j A} e^0_j + O_w(x^2 \partial^{\mathcal V} + x\partial^{\mathcal H})   
\\
&= \begin{cases} e^0_i &+ O_w(x^2 \partial^{\mathcal V} + x\partial^{\mathcal H}), \quad \text{if $A=i$, since  $T_{1ji}=0$}
\\
e^0_\alpha + xT_{1j\alpha} e^0_j &+ O_w(x^2 \partial^{\mathcal V} + x\partial^{\mathcal H}), \quad \text{if $A=\alpha$}
\end{cases}
\end{align*}
that is, there is expansion, 
\[
e^1_A = H(e^0_A) + O_w(x^2 \partial^{\mathcal V} + x\partial^{\mathcal H}).
\]
The dual frames $[e^{1A}]_A$ and $[e^{0B}]_B$ are transformed as $[e^{1A}]_A = [e^{0B}]_B \cdot T^{-1}$.  

\item[(c)] We can estimate the rates as $x\to 0^+$ of the Christoffel symbols of the parallel frames as follows:  
\begin{align*}
\nabla_{e^1_A}^1 e^1_B &= e^1_A(T_{Z B}) e_Z^0 + T_{\Delta A} T_{Z B} \nabla_{e_\Delta^0}^1 e_Z^0, 
\\
g_1 (\nabla_{e^1_A}^1 & e^1_B, e^1_\Gamma) = e^1_A(T_{Z B}) g_1(e_Z^0, e_\Theta^0) T_{\Theta \Gamma} +T_{\Delta A} T_{Z B} \, g_1\left(\nabla^1_{e_\Delta^0} e_Z^0, e_\Theta^0\right) T_{\Theta\Gamma} 
\\
&= O(x^{2\delta}) + e^0_A(T_{\Gamma B}) +  T_{\Delta A} T_{Z B}T_{\Theta\Gamma}\, \left[x^{2\delta-1} E_1(e^0_\Delta, e^0_Z, e^0_\Theta) + g_0\left(\nabla^0_{e_\Delta^0} e_Z^0, e_\Theta^0\right) \right]       
\\
&=\begin{cases}
    O(x^{2\delta-1}) + g_0\left(\nabla^0_{e_A^0} e_B^0, e_\Gamma^0\right), &  \text{if $\delta>1/2$,}
    \\
    O(x) + (p^{\mathcal V} e_A^0)(xT_{1\Gamma B}) + E_1(e_A^0, e_B^0, e_\Gamma^0)\big|_{x=0} + &  \text{if $\delta=1/2$,}
    \\
    \quad x\left(T_{1\Delta A}g_0\left(\nabla^0_{e_\Delta^0} e_B^0, e_\Gamma^0\right) + T_{1ZB} g_0\left(\nabla^0_{e_A^0} e_Z^0, e_\Gamma^0\right) + T_{1\Theta\Gamma} g_0\left(\nabla^0_{e_A^0} e_B^0, e_\Theta^0\right) \right)
    \\
    \quad + g_0\left(\nabla^0_{e_A^0} e_B^0, e_\Gamma^0\right), 
\end{cases}
\end{align*}
where by the skew-symmetry of the Christoffel symbol in $B$ and $\Gamma$, we assume $B<\Gamma$. Notice that by \cite{Albin_2017_index} (see Section 1.2, specifically around Remark 1.4), the model Christoffel symbols of the vertical wedge-Fermi fields $g_0\left(\nabla^0_{e_i^0} e_j^0, e_k^0\right)$ are of order $O(1/x)$ or they vanish. Therefore
\[
g_1(\nabla_{e^1_A}^1 e^1_B, e^1_\Gamma) - g_0\left(\nabla^0_{e_A^0} e_B^0, e_\Gamma^0\right) =x^{2\delta-1}G_{AB\Gamma}(x),  
\]
for some smooth function $G_{AB\Gamma}=G_{AB\Gamma}(x),\ x\in [0,\infty)$.
\end{itemize}
\end{remark}

The association of Lemma~\ref{lemma_on_the_frame_expansions}(c) defines a principal bundle isomorphism $\tau$ between the principal $\mathrm{SO}(2n,\mathbb{R})$-sub-bundles $^{w}\mathcal{O}$ and $^{w,pt}\mathcal{O}$ of $g_w$ and $g_{w,pt}$-orthogonal frames of $\prescript{w}{}{T\mathscr{C}}(\mathfrak{B}_Y)\to \mathscr{C}(\mathfrak{B}_Y)$ of the principal bundle of oriented frames respectively as follows: given a $g_w$-orthonormal frame $[e_A]_A\in \prescript{w}{}{\mathcal{O}}_{(p,s)}$ above $(p,s)\in \mathfrak{B}_Y \times (0,1)$, we $\nabla$-parallel transport along the flow-line from $(p,s)$ of $\partial_x$ and expand the resulting parallel frame in $\{p\}\times (0,1)$. We then define $\tau[e_A]_A := [\partial_x, e^0_j, e^0_\alpha]_{j,\alpha}\cdot T_0 \in \prescript{w,pt}{}{\mathcal{O}}_{(p,s)}$. Similarly the map $H:\prescript{w}{}{\mathcal{O}}\to \mathrm{GL}^+({}\prescript{w}{}{T\mathscr{C}}(\mathfrak{B}_Y))$ introduced in \eqref{eq:horozontal_lift}, is a principal bundle morphism. Clifford modules $\mathcal{E}\to \mathfrak{B}_Y\times(0,1)$ and model Clifford modules $\prescript{w,pt}{}{\mathcal{E}} \to \mathfrak{B}_Y\times(0,1)$ have frame bundles that either are induced by the $O(n)$-principal bundles $^{w}\mathcal{O}$ and $^{w,pt}\mathcal{O}$ of $g_w$-orthonormal frames and of $g_{w,pt}$-orthonormal frames directly, and hence the $\tau$ map gives an induced $\tau$-map that is an isomorphism on their corresponding frame bundles, or have frame bundles that are induced by the spin-principal bundles $\prescript{w}{}{\mathcal{P}} \to \mathfrak{B}_Y\times(0,1)$ and $\prescript{w,pt}{}{\mathcal{P}} \to \mathfrak{B}_Y\times(0,1)$ that are lifts of $g_w$-orthonormal frames and $g_{w,pt}$-orthonormal frames respectively and again the $\tau$-isomorphism lifts to $\tau$ map that is an isomorphism of the principal bundles $\tau:{}^{w}\mathcal{P}\to {}^{w,pt}\mathcal{P}$ and induce isomorphisms to the induced frame bundles of $\mathcal{E}$ and $\prescript{w,pt}{}{\mathcal{E}}$ and their Clifford actions as $cl_w$ and $cl_{w,pt}$ respectively. We continue using the same symbol $\tau: \mathcal{E} \to \prescript{w,pt}{}{\mathcal{E}}$ to denote their isometric bundle isomorphism. These isomorphisms extend over the boundary of the collar to give isomorhisms of the corresponding wedge frame bundles. Note that we have an identity,
\begin{equation}\label{eq:tau-isom-def}
   \tau\circ cl_w(e^A) = cl_{w,pt}(\tau e^A)\circ \tau, \quad \forall A. 
\end{equation}

\begin{remark}
Recall the special case where $\delta=1/2$ and the wedge metric $g_w$ is smooth as in Remark~\ref{rem:smooth_case_analogy}. The map $\tau : {}^{w}{\mathcal O} \to {}^{w,pt}{\mathcal O}$ sends an orthonormal frame of the smooth metric $g_w$ to the unique othonormal frame of $g_{w,pt}$ so that their corresponding parallel transports with their Levi-Civita connections along the radial geodesic give the same orthonormal frame of $g_N\oplus g_F$ on $TX|_F=N\oplus TF|_F$. By construction, the two frame-bundles are also isomorphic to the pullback bundle $\pi^*({\mathcal O}|_F)$ and is the one used as the model bundle (see \cite{wu1998equivariant} par. 3.2 pp.9] or \cite[\S 8g]{bismut1991complex}). In the non-smooth  wedge case the pullback bundle is absent but the $\tau$ isomorphism is still valid: Instead of defining the frames to parallel transport on the same restriction frame, we demand from their transition matrix to converge to the identity as the distance $x$ from the singular stratum $F$ vanishes.     
\end{remark}

We choose frame $[e_A]_A$ with $\tau [e_A]_A = [e_A^0]_A = [\partial_x, e^0_j, e^0_\alpha]_{j,\alpha}$. Using the construction of the $\mathrm{spin}^c$-connections $\nabla^{\mathcal{E}}$ and $\nabla^{\prescript{w,pt}{}{\mathcal{E}}}$ from the Levi-Civita connections $\nabla$ and $\nabla^M$ respectively, we obtain that on a lift frame $[f_i]_i\in \prescript{w}{}{\mathcal{P}}$ of $[e_A]_A \in \prescript{w}{}{\mathcal{O}}$,   
\begin{align*}
\tau \nabla^{\mathcal{E}}_{e_A} f_i - \nabla^{\prescript{w,pt}{}{\mathcal{E}}}_{\tau e_A} (\tau f_i) &= -\frac{1}{2}\sum_{B<\Gamma}(^{w}\omega_{AB}^\Gamma \tau (cl_w(e^B) cl_w(e^\Gamma)f_i - ^{w,pt}\omega_{A B}^\Gamma cl_{w,pt}(\tau e^B) cl_{w,pt}(\tau e^\Gamma) \tau f_i) 
\\
&= -\frac{1}{2}\sum_{B<\Gamma} [g_w(\nabla_{e_A} e_B, e_\Gamma) - g_{w,pt}(\nabla^M_{\tau e_A} \tau e_B, \tau e_\Gamma)]  \tau (cl_w(e^B) cl_w(e^\Gamma) f_i)
\end{align*}
and by Remark~\eqref{rem:remark_on_frame_expansions}(c), we obtain
\begin{equation}
\label{eq:actual_Christofell_spinor_symbols_minus_model_Christoffel_spinor_symbols_rates}
\tau \nabla^{\mathcal{E}}_{e_A} f_i - \nabla^{\prescript{w,pt}{}{\mathcal{E}}}_{\tau e_A} (\tau f_i) = \begin{cases}
    [\tau f_j]_j \cdot O(x^{2\delta-1}),& \quad \text{if $\delta>1/2$,}
    \\
    {\mathcal B}_{\tau e_A}(\tau f_i) + [\tau f_j]_j\cdot O(x),&\quad \text{if $\delta=1/2$,}
\end{cases} 
\end{equation}
where 
\[
{\mathcal B}_{e_A^0} := -\frac{1}{2} \sum_{B<\Gamma} G_{AB\Gamma}(0)  cl_{w,pt}( e^{0B}) cl_{w, pt}( e^{0\Gamma}): ^{w,pt}{\mathcal E} \to {}^{w,pt}{\mathcal E}, 
\]
is skew-symmetric. By construction the $\mathrm{spin}^c$ as the pullback connection through the principal bundle morphisms ${}\prescript{w}{}{\mathcal{P}}\to \prescript{w}{}{\mathcal{O}}$, the parallelism properties of the frame $[e_A]_A\in \prescript{w}{}{\mathcal{O}}_{(p,s)}$ transfer to the lifted frame $[f_i]_i\in \prescript{w}{}{\mathcal{P}}$. In particular, the frame $[f_i]_i\in  \mathcal{C}^\infty\left(\mathfrak{B}_Y\times(0,1); {}\prescript{w}{}{\mathcal{P}}\big|_{\mathfrak{B}_Y\times(0,1)}\right)$ is $\nabla^{\mathcal{E}}$-parallel along the flow line of $\partial_x$. An analogue statement is true for the model frames $[\tau e_A]_A$ and $[\tau f_i]_i$.

Following Definition \ref{definition_resolved_fundamental_neighbourhood} , let $U_t:=U_t(F)$ be a resolved localizing neighborhood of this component $F$ of the vanishing locus $Z(|V|^2)$ of the potential. We consider  two spin Dirac operators, the first being $D_X$ restricted to $U_t$ which we will continue to denote as $D_X$. The second is the \textbf{model operator $D_M$}, which is the composition of the connection and the Clifford action for $U_t$ with the model metric $g_{w,pt}$. The two operators are written explicitly in the preceding frames as
\[
D_X = cl_w(e^A) \nabla^{\mathcal{E}}_{e_A} \quad \text{and} \quad D_M= D^{\mathcal V} +  D^{\mathcal H},
\]
where the vertical and horizontal operators are defined locally as follows: given a $g_{w,pt}$-orthonormal frame $[e_i^0, e_\alpha^0]_{i,\alpha}\in \prescript{w,pt}{}{\mathcal{O}}$ with lift $[f_j]_j\in {}^{w,pt}{\mathcal E}$ and on a basic section $\xi = f f_j$ we have the expressions
\begin{equation}\label{eq:vert-ops}
D^{\mathcal V}\xi= cl_{w,pt}(e^{0i}) \nabla^{\prescript{w,pt}{}{\mathcal{E}}}_{e_i^0}\xi
\end{equation}
and 
\begin{equation}\label{eq:hor-ops}
D^{\mathcal H}\xi = \begin{cases} cl_{w,pt}(e^{0\alpha}) \nabla^{\prescript{w,pt}{}{\mathcal{E}}}_{e_\alpha^0} \xi&, \quad  \text{if $\delta>1/2$} \\  cl_{w,pt}(e^{0\alpha})[ H (e_\alpha^0)(f) f_j + f \nabla^{\prescript{w,pt}{}{\mathcal{E}}}_{e^0_\alpha}f_j] &, \quad \text{if $\delta=1/2$}\end{cases} .
\end{equation}
The operators extend to all local sections by linearity, they are independent of the choice of $g_{w,pt}$-orthonormal frames $[e^0_i, e^0_\alpha]_{i,\alpha}$ associated with the fibration $\phi_Y:\mathfrak{B}_Y\to Y$ and therefore are defined on sections of a tubular neighborhood of the $Y$ stratum.

We have the following proposition, which gives the polyhomogenous estimates for the (undeformed) Dirac operator.

\begin{proposition}\label{prop:poly-est-delta-geq-1/2}
Let $\delta \geq 1/2$. Then for every $\xi\in\mathcal{C}^\infty(\mathscr{C}(\mathfrak{B}_Y); {}^{w,pt}{\mathcal{E}})$     \begin{equation}
      (\tau D_X \tau^{-1}) \xi - D_M\xi = \begin{cases} O_w(x^{2\delta}\partial^{\mathcal V} + x^{2\delta}\partial^{\mathcal H} + x^{2\delta-1})\xi,& \quad \text{if $\delta>1/2$,}
      \\
       {\mathcal B}\xi + O_w(x^2\partial^{\mathcal V} + x\partial^{\mathcal H} + x)\xi,& \quad \text{if $\delta=1/2$,}
      \end{cases}
    \end{equation}
    where $\mathcal{B} := cl_{w,pt} (\tau e^A)\mathcal{B}_{\tau e_A}$.
\end{proposition}

\begin{proof}
    We choose a local frame $[e_A]_A\in \mathcal{C}^\infty(U_t; \prescript{w}{}{\mathcal{O}}\big|_{U_t})$ with $\tau [e_A]_A = [\partial_x, e^0_\alpha, e^0_j]_{\alpha,j}$ so that $p^{\mathcal V}(\tau e_A) = \tau e_A$ or $p^{\mathcal H}(\tau e_A)=0$ depending on whether $\tau e_A$ is vertical or horizontal Fermi field. Let $[f_i]_i \in \mathcal{C}^\infty(U_t; {}\prescript{w}{}{\mathcal{P}}\big|_{U_t})$ be a lift of $[e_A]_A$. By linearity, it is enough to prove the expansion for local sections $\xi$ with $\tau^{-1} \xi = f f_i$ for some $f\in \mathcal{C}^\infty(U_t)$. When $\delta>1/2$ we apply \eqref{eq:actual_Christofell_spinor_symbols_minus_model_Christoffel_spinor_symbols_rates} to obtain, 
    \begin{align*}
        (\tau D_X \tau^{-1})\xi &=  \tau cl_w (e^A) \nabla^{\mathcal{E}}_{e_A}(f f_i) = cl_{w,pt} ( \tau e^A)[ e_A(f) \tau f_i + f \tau\nabla^{\mathcal{E}}_{e_A} f_i]   
        \\
        &= cl_{w,pt} ( \tau e^A)[ (\tau e_A)(f) \tau f_i + f \nabla^{\prescript{w,pt}{}{\mathcal{E}}}_{\tau e_A} \tau f_i] + O_w(x^{2\delta}\partial^{\mathcal V} + x^{2\delta}\partial^{\mathcal H} + x^{2\delta-1})\xi  
        \\
        &= D_M\xi + O_w(x^{2\delta}\partial^{\mathcal V} + x^{2\delta}\partial^{\mathcal H} + x^{2\delta-1})\xi. 
    \end{align*}

When $\delta = 1/2$ we recall Remark~\ref{rem:remark_on_frame_expansions}(b) and we estimate
\begin{align*}
e_A(f) &= H(\tau e_A)(f) + O(x^2\partial^{\mathcal V}+ x \partial^{\mathcal H}) f
\end{align*}
Then, by applying \eqref{eq:actual_Christofell_spinor_symbols_minus_model_Christoffel_spinor_symbols_rates}, the same calculation as before gives,
\begin{align*}
        (\tau D_X \tau^{-1}) \xi &= cl_{w,pt} (\tau e^A)[ e_A(f)  + f\tau\nabla^{\mathcal{E}}_{e_A}f_i] 
\\
&= cl_{w,pt} (\tau e^A)[ H(\tau e_A)(f)  + f\nabla^{\prescript{w,pt}{}{\mathcal{E}}}_{\tau e_A}f_i + \mathcal{B}_{\tau e_A}] +  O_w(x^2\partial^{\mathcal V} + x\partial^{\mathcal H} + x)\xi
\\
&= cl_{w,pt}(p^{\mathcal V}\tau e^A)[p^{\mathcal V} (\tau e_A)(f) \tau f_i + f \nabla^{\prescript{w,pt}{}{\mathcal{E}}}_{p^{\mathcal V}(\tau e_A)} \tau f_i]  
\\
&\quad +cl_{w, pt} (p^{\mathcal H}\tau e^A)[(H(p^{\mathcal H} (\tau e_A)))(f) \tau f_i  + f \nabla^{\prescript{w,pt}{}{\mathcal{E}}}_{p^{\mathcal H}(\tau e_A)} \tau f_i] 
\\
&\quad + \mathcal{B}\xi + O_w(x^2\partial^{\mathcal V}+ x \partial^{\mathcal H}+x) \xi 
\\
&= D_M\xi+ \mathcal{B} \xi + O_w(x^2\partial^{\mathcal V} + x\partial^{\mathcal H} + x)\xi. 
    \end{align*}
    The result now follows by linearity.
\end{proof}

The ambient Sobolev norm and the model Sobolev norm we subsequently use, become comparable by introducing the ratio of the ambient density with the model density:

\begin{definition}
\label{defn:jacobian_ratio}
We define the Jacobian factor $k_F= \frac{d\mathrm{vol}_w}{d\mathrm{vol}_{w,pt}}$ as the ratio of the product-type $L^2$-density on $U_{t_0}(F)$ used to define the local model, and the density obtained by taking the restriction of the ambient $L^2(X)$-volume to $\overline{\beta^{-1}\left(\{\rho_{F}<t\}\right)} \subset X$ and pushing it forward  along $\beta^{-1}: \overline{\{\rho_F<t\}} \to U_t(F)$. 
\end{definition}

In the preceding local frames, the two volume forms are represented as $d\mathrm{vol}_w = e^1\wedge \dots \wedge e^{2n}$ and $d\mathrm{vol}_{w,pt} = \tau e^1 \wedge \dots \wedge \tau e^{2n}$. We have that $d\mathrm{vol}_w = \det T(x,w)\, d\mathrm{vol}_{w,pt}$ so that 
\begin{equation}
\label{eq:jacobian_rates}
k(x,w)=\det T(x,w)= 1+O(x^{2\delta})\quad \text{as $x\to 0^+$.}
\end{equation}

We fix a component $F$ of the zero set of $A$ with fibration ${U} - {U(F)} \xrightarrow{{\phi_F}} F^{reg}$ of the resolved localizing neighbourhood and the fiber decomposition $U_y = U_y^s \times U_y^u$ into the attracting and expanding factors for the gradient vector field $V|_{U(F)}$ and distance that decompose to the stable and unstable cones as $x^2 = \sum_j x_j^2$. Together with expansion \eqref{Morse_Bott_Kirwan_modify_hamiltonian_vf} we obtain 

\begin{corollary}
Let $\delta \geq 1/2$. Then for every $\xi\in\mathcal{C}^\infty(U_{t}(F); {}^{w,pt}\mathcal{E})$     
\begin{multline}
      \tau (D_X + \sqrt{-1}\varepsilon cl_w(V))\tau^{-1}\xi - (D_M + \sqrt{-1}\varepsilon\gamma_j x_j cl_{w,pt}(V_j)) \xi
      \\
      = \begin{cases} O_w(x^{2\delta}\partial^{\mathcal V} + x^{2\delta}\partial^{\mathcal H} + x^{2\delta-1} + s x^{2\delta+1})\xi,& \quad \text{if $\delta>1/2$,}
      \\
        {\mathcal B}\xi + \sqrt{-1}\varepsilon x^2 cl_{w,pt}(\tau V_E)\xi   
       \\
       \quad \; + O_w(x^2\partial^{\mathcal V} + x\partial^{\mathcal H} + x+\varepsilon x^3)\xi,& \quad \text{if $\delta=1/2$,}
      \end{cases}
    \end{multline} 
\end{corollary}

We define the following scaling operator.
\begin{definition}[Scaling operator]\label{def:scaling_operator}
    Let $\varepsilon>0$. The principal bundle $\mathcal{G^+}( \prescript{w}{}{T\mathscr{C}}\mathfrak{B}_Y)$ of positively oriented frames can be trivialized along rays of the radial field $\partial_x$ and we define in a trivialization chart $(x,z)\in U_t(F)$ of the bundle the following rescaling operator: given frame $g_w$-orthonorormal frame  $[e_A]_A \in \mathcal{C}^\infty(U_t(F); \prescript{w}{}{\mathcal{O}})$, we define
    \begin{equation}
    (S_{\varepsilon} [e_A]_A)(x,z) = \sqrt{\varepsilon} p^{\mathcal V}([\tau e_A]_A)(x/\sqrt{\varepsilon}, z)+ p^{\mathcal H}([\tau e_A]_A)(x/\sqrt{\varepsilon}, z) \in \mathcal{C}^\infty(U_{t\sqrt{\varepsilon}}(F); \mathcal{G^+}( {}\prescript{w}{}{T\mathscr{C}}\mathfrak{B}_Y)).
    \end{equation}
    Likewise, if $u \in L^2(U_{t}(F); {\mathcal E})$, let
    \begin{equation}
    S_{\varepsilon}u(x,z):=\tau \circ u(x/{\sqrt{\varepsilon}},z)\in L^2\big(U_{t\sqrt{\varepsilon}}(F); \mathcal{E}\big).
    \end{equation}
    
\end{definition}

The following is the Witten deformed version of the polyhomogenous estimates.

\begin{proposition}
\label{proposition_rescaled_estimates-0}
Let $\delta \geq 1/2$. Consider a Hamiltonian function $h$, and the rescaled setting. Recall the expansion \eqref{Morse_Bott_Kirwan_modify_hamiltonian_vf}. Then for every $\xi\in\mathcal{C}^\infty(U_{t\sqrt{\varepsilon}}(F); ^{w,pt}{\mathcal{E}})$     
\begin{multline}\label{eq:rescaling-grad-term}
    S_{\varepsilon}(D_X+ \sqrt{-1}\varepsilon cl_w(V)) S^{-1}_{\varepsilon}\xi - (D^{\mathcal{H}}+\varepsilon^{1/2}(D^{\mathcal{V}}+ \sqrt{-1}\gamma_j x_j cl_{w,pt}(V_j)) \xi 
    \\
=    \begin{cases}O_w(\varepsilon^{1/2-\delta}x^{2\delta}\partial^{\mathcal V} + \varepsilon^{-\delta}x^{2\delta}\partial^{\mathcal H} + \varepsilon^{1/2-\delta}x^{2\delta-1})\xi, & \quad \text{if $\delta>1/2$,}
      \\
      \mathcal{B}\xi+ \sqrt{-1}x^2 cl_{w,pt}(\tau V_E)\xi +\varepsilon^{-1/2}O_w(x^2\partial^{\mathcal V} + x\partial^{\mathcal H} + x)\xi, & \quad \text{if $\delta=1/2$}
      \end{cases}.
    \end{multline}
\end{proposition}

We will use this to prove the polyhomogenous estimates for the Witten deformed Dirac operator in Proposition \ref{prop:block-estimates}. We end with analogues of \cite[Thm 8.21]{bismut1991complex}, \cite[Prop 3.4]{wu1998equivariant}:
\begin{lemma}\label{lem:hor-op-conj}
    Using the $\tau$ isomorphism defined above, and $\phi_F:U_t(F)\to F$ the resolved localizing neighborhood \eqref{eq:resolved-localizing-nbhd}, and $D_F$ the connection of the normal cohomology complex $\mathcal{P}_{\mathcal{B},\varepsilon}(\mathcal{N}_{F,\varepsilon})$, we have for all $\alpha\in L^2(F\,;\mathcal{E}_\mathcal{H})$ that
    \begin{equation}
        D^\mathcal{H}\phi_F^*(\alpha) = \phi_F^*(D_F \alpha) .
    \end{equation}    
\end{lemma}
\begin{proof}
    From the definition of $D^\mathcal{H}$ \eqref{eq:hor-ops}, and the identity \eqref{eq:tau-isom-def}, this is immediate after observing that our construction of the normal cohomology complex $\mathcal{N}_{F,\varepsilon}$ already incorporates the $\tau$ isomorphism on its fibers, and the definitions of $D^\mathcal{H}$ and $D_F$ otherwise only differ based on the source of their sections ($U_t(F)$ and $F$ respectively). 
\end{proof}

\subsection{Estimates of the operator and resolvent as $\varepsilon\to \infty$}\label{subsec:resolvent-ests}

In this section we will construct the \textit{quasimode approximation} which will furnish the remaining estimates to be proven. Beginning with a family of $\mathcal{N}=1$ or $\mathcal{N}=2$ SUSY Hilbert complexes $\mathcal{P}_{W,\varepsilon}(X)=(L^2(X,\mathcal{E}), \mathcal{D}_W(Q_{j,\varepsilon}), Q_{j,\varepsilon})_{j=1}^\mathcal{N}$, associated to a Hamiltonian $\Delta_\varepsilon$ of $0$-th order type and a potential $A$. Recall that if $\mathcal{N}=1$ then $Q_1=Q=D$, and if $\mathcal{N}=2$ then $Q_1=D=P+P^*$ and $\mathcal{D}(D)=\mathcal{D}(P)\cap\mathcal{D}(P^*)$; the proceeding arguments in this subsection will hold for the cases of $\mathcal{N}=1$ and 2 simultaneously and will be presented as such.

Throughout this subsection we will denote 
\begin{equation}\label{eq:kth-domain}
H_{\mathcal{B}}^k (X;\mathcal{E}) := \mathcal{D}_{\mathcal{B}}(D^k):=\{u \in \mathcal{D}_{\mathcal{B}}(D^{k-1}) \, | \, D^ku \in  \mathcal{D}_{\mathcal{B}}(D^{k-1}) \}\subset L^2(X;\mathcal{E})   
\end{equation}
to be the iteratively defined domain of the $k$-th power of the Dirac operator $D_X$ with boundary conditions $\mathcal{B}$ described in subsection \ref{subsec:cohomology-truncated-neighbourhoods}. The choice of notation is suggestive, as in the smooth case this collapses to the corresponding $k$th Sobolev space of sections of our (generalized) spinor bundle. For example, for an $\mathcal{N}=1$ SUSY Hilbert complex arising from Witten deformation by a K\"ahler-Hamiltonian vector field $V$, as in definition \ref{def:witten-def-cmpx-vfs}, if $X$ is smooth and K\"ahler then \eqref{eq:kth-domain} is exactly the $k$-th Sobolev space of $\Lambda^{(0,q)}T^*X\otimes \mathcal{E}$ (with the corresponding boundary conditions). Notice also that because our action acts on these sections we can restrict $D^k$ to the $\mu$-eigenspace of $\mathcal{L}_V$; we denote the corresponding domain of sections by $H_{\mathcal{B},\mu}^k(X;E)$. Finally, we note that (for each fixed $\mu$), by Proposition \ref{Proposition_equivariant_complexes} we can reduce our 1st order complexes to zeroth order complexes. For this reason, the operator obtained by restricting $D_X$ to the $\ker(iL_V-\mu)$ eigenspace can be obtained as a zeroth order operator, and we shall otherwise \emph{omit $\mu$ in what follows}, in keeping with our presenting the proceeding arguments for the $\mathcal{N}=1$ and 2 cases. \\

Now, for each connected component $F\in Z(A)$ of the vanishing locus, fix \(t_0\in (0,\infty]\) and let
\begin{equation}
U_{t_0}(F)\subset U(F)
\end{equation}
denote the corresponding resolved localizing neighbourhood defined in \ref{definition_resolved_fundamental_neighbourhood}. By construction there is
a coordinate identification of the dense smooth part of $U_{t_0}(F)$ with a neighbourhood $\mathcal{T}_{F^{\text{reg}}}\subset X$, which allows us to identify sections over $U_{t_0}(F)$ with sections over $X$, via pullback. This should be viewed as a singular replacement for a tubular coordinate system.

Recall the existence of a $t_0>0$ such that corresponding truncated localizing neighborhoods $U_{t_0}(F)$ are all disjoint. Then for each such $F$ we introduce a smooth cutoff $\gamma_{t_0}^F: \mathbb{R}\to [0,1]$ such that
\begin{equation}
     \gamma_{t_0}^F(t) = \begin{cases}
    1 & t< t_0/2 \\
    0 & t>1
\end{cases}  
\end{equation}

Now, we recall the construction given in definition \ref{Definition_all_localized_domains} (see also the discussion surrounding \ref{def_normal_cohomology}) which introduced the Witten-deformed normal cohomology complex; this will be crucial to the estimates that follow and is used to define the \emph{quasimode-type approximation} applied in this section. We refer to this as such because the cohomology of the normal cohomology complex $\mathcal{P}_{\mathcal{B},\varepsilon}\big(\mathcal{N}_{F,\varepsilon})$ will define $\varepsilon^f$-scale quasi-modes of $\Delta_\varepsilon$. However, rather than approximating the Witten complex by the cohomology of $\mathcal{P}_{\mathcal{B},\varepsilon}\big(\mathcal{N}_{F,\varepsilon})$, we approximate global sections by such $\mathcal{N}_{F,\varepsilon}$-twisted sections, and show that these converge to the Witten complex. This construction should be viewed as the singular analogue of ``small eigenvalue complex" approximation first considered in \cite{bismut1991complex, wu1998equivariant}, which themselves were inspired by the ``Born-Oppenheimer" approximation considered in \cite{witten1984holomorphic}.
\begin{definition}\label{def:adapted-normal-cohom-complex}
We recall the normal cohomology complex is given by
\begin{equation}
\mathcal{P}_{\mathcal{B},\varepsilon}\big(\mathcal{N}_{F,\varepsilon})
=
\left( L^2\left(F^{\text{reg}}; \mathcal{E}_{\mathcal{H}}\otimes \mathcal{N}_{F,\varepsilon}\right),
\mathcal{D}_{\mathcal{B}}(D_\varepsilon^F),
D_{\varepsilon}^F
\right) , 
\end{equation}

which are the sections over $F^{\text{reg}}$ valued in the normal cohomology bundle $\mathcal{N}_{F,\varepsilon}$ \eqref{eq:normal-cohomology-bundle}, (whose fibers come from the vertical Hodge cohomology $\mathcal H^{\bullet}_{\mathcal{B},\varepsilon}(U_{F,\varepsilon})$ of $\Delta_{\varepsilon}^V|_{F}$), and twisted by the horizontal Clifford module $\mathcal{E}_{\mathcal{H}}$. 
\end{definition}
Following \eqref{eq:kth-domain} we then define the corresponding higher regularity sections as
\begin{equation}
H_{\mathcal{B},\varepsilon}^k(\mathcal{N}(F); E):=\mathcal{D}_\mathcal{B}((D_{F,\varepsilon})^k), \quad \text{where} \quad H_{\mathcal{B},\varepsilon}^0(\mathcal{N}(F); E):= L^2\left(F; \mathcal{E}_{\mathcal{H}}\otimes \mathcal{N}_{F,\varepsilon}\right)
\end{equation}
as the notation suggests these should be viewed 

Now introduce the embedding 
\begin{equation}
    J_{\varepsilon}: H_{\mathcal{B},\varepsilon}^1(F\, ; \mathcal{E}_{\mathcal{H}}\otimes \mathcal{N}_{F,\varepsilon}) \hookrightarrow H_\mathcal{B}^1(X;\mathcal{E}), 
\end{equation}
and we denote the image of $J_{\varepsilon}$ by $H_{\mathcal{B},\varepsilon}^1(X;E)$. Explicitly, we consider the pushforward on sections induced by the inclusion $\iota_a: U_{t_0}(\widehat{F})\hookrightarrow X$, and $\beta_F: U_{t_0}(F) \to U_{t_0}(\widehat{F})$ the resolution of the localizing neighborhoods, then we define for $\sigma\in H_{\mathcal{B},\varepsilon}^1(F\, ; \mathcal{E}_{\mathcal{H}}\otimes \mathcal{N}_{F,\varepsilon})$
\begin{equation}
\label{equation_local_expression_quasimode}
J_{\varepsilon,F}(\sigma)
=
(\iota_{F}\circ \beta_F)_\ast \!\left(
k_a^{-1/2}\gamma_{t_0}^F(t)\,
\sum_{j=1}^{n_F}
\phi_{F}^*(\alpha_{j})\wedge \omega_{j,\varepsilon}
\right).
\end{equation}
where $k_F$ is the density defined in definition \ref{defn:jacobian_ratio}; the $\phi_{F}^*(\alpha_{a,j})$ are lifts of horizontal forms $\alpha_{j}\in H_{\mathcal{B},\varepsilon}^1(F^{reg}; \mathcal{E}_\mathcal{H})$ over the vanishing set $F$; and where the forms
\begin{equation}
\{\omega_{j,\varepsilon}(y,\cdot)\}_{j=1}^{n_F} \subset L^2(F^{\text{reg}};\mathcal{N}_{F,\varepsilon})
\end{equation}
comprise a local orthonormal frame for the vertical Hodge cohomology $\mathcal H^{\bullet}_{\mathcal{B},\varepsilon} = \ker \Delta_\varepsilon^V$ over \(y\in F^{\text{reg}}\). In this case each
\(\omega_{j,\varepsilon}(y,\cdot)\) can be  realized as a fibrewise harmonic representative on the normal model over \(y\).

Also define the orthogonal complement of $H_{\mathcal{B},\varepsilon}^k(X;\mathcal{E})$ in $H_{\mathcal{B}}^k(X;\mathcal{E})$ by 
\begin{equation}
    H_{\mathcal{B},\varepsilon}^{k,\perp}(X;\mathcal{E}) := H_{\mathcal{B}}^k(X;\mathcal{E}) \cap H_{\mathcal{B},\varepsilon}^{0,\perp}(X;\mathcal{E}).
\end{equation}
Then given this embedding we can define the orthogonal projectors $\Pi_{\varepsilon}, \Pi_{\varepsilon}^\perp$ from $H_{\mathcal{B}}^k(X;\mathcal{E})$ onto $H_{\mathcal{B},\varepsilon}^{k}(X;E)$ and $H_{\mathcal{B},\varepsilon}^{k,\perp}(X;E)$ respectively. These projectors induce a block decomposition on $D_X$, with respect to the splitting $H_{\mathcal{B}}^k(X;\mathcal{E}) = H_{\mathcal{B},\varepsilon}^{k}(X;\mathcal{E}) \oplus H_{\mathcal{B},\varepsilon}^{k,\perp}(X;\mathcal{E})$ we have
\[ D_{X,\varepsilon} = \begin{pmatrix}
    D_{\varepsilon}^{(1)} & D_{\varepsilon}^{(2)} \\ D_{\varepsilon}^{(3)} & D_{\varepsilon}^{(4)}
\end{pmatrix} . \] 
We record the following block estimates for the induced operators below:

\begin{proposition}[Block estimates]\label{prop:block-estimates}

For the $\delta>0$ defined in Proposition \ref{proposition_rescaled_estimates-0}, we define
\begin{equation}\label{eqn:f_delta-defn}
    f=f(\delta)=\begin{cases}
    \tfrac{1}{2} -\delta & \delta > \tfrac{1}{2} \\
    -\tfrac{1}{2}    & \delta = \tfrac{1}{2}
    \end{cases} ,
\end{equation}
based on the form of our asymptotically $\delta$ wedge metric. Having fixed this parameter, in our setting, the following estimates hold.
\begin{enumerate}
\item There exist constants \(\varepsilon_0>0\) and \(C_0>0\) such that for every
\(\varepsilon>\varepsilon_0\),
\begin{equation}
J_{\varepsilon}^{-1}D_{\varepsilon}^{(1)}J_{\varepsilon}
=
D_F+\varepsilon^f R_\varepsilon,
\end{equation}
where $R_\varepsilon$ is a first-order differential operator on
\( \mathcal N_{\mathcal{B},\varepsilon}(F) \) whose coefficients are uniformly bounded independently of
\(\varepsilon\). Equivalently, for every
\begin{equation}
\sigma\in H_{\mathcal{B},\varepsilon}^1(F\, ; \mathcal{E}_{\mathcal{H}}\otimes \mathcal{N}_{F,\varepsilon}),
\end{equation}
one has
\begin{equation}
\bigl\|
J_{\varepsilon}^{-1}D_{\varepsilon}^{(1)}J_{\varepsilon}\sigma
-
D_{F,\varepsilon}\sigma
\bigr\|_{L^2_{\varepsilon}}
\le
C_0\,\varepsilon^f\,\|\sigma\|_{H^{1}_{\varepsilon}}.
\end{equation}

\item There exist constants \(\varepsilon_1>0\) and \(C_1>0\) such that for every
\(\varepsilon>\varepsilon_1\), every
\begin{equation}
s'\in H_{\mathcal{B},\varepsilon}^{1,\perp}(X;\mathcal{E}),
\qquad
s\in H_{\mathcal{B},\varepsilon}^{1}(X;\mathcal{E}),
\end{equation}
one has
\begin{equation}
\|D_{\varepsilon}^{(2)}s'\|_{L^2_{\varepsilon}}
\le
C_1(\|s'\|_{L^2_{\varepsilon}} + \varepsilon^f\|s'\|_{H^{1}_{\varepsilon}} ),
\qquad
\|D_{\varepsilon}^{(3)}s\|_{L^2_{\varepsilon}}
\le
C_1 (\|s\|_{L^2_{\varepsilon}} + \varepsilon^f \|s\|_{H^{1}_{\varepsilon}}) .
\end{equation}

\item There exist constants \(\varepsilon_2>0\) and \(C_2>0\) such that for every
\(\varepsilon>\varepsilon_2\) and every
\begin{equation}
s\in H_{\mathcal{B},\varepsilon}^{1,\perp}(X;\mathcal{E}),
\end{equation}
one has
\begin{equation}\label{eq:rescaled-Garding-est}
\|D_{\varepsilon}^{(4)}s\|_{L^2_{\varepsilon}}
\ge
C\Bigl(\|s\|_{H_{\varepsilon}^1}+\sqrt{\varepsilon}\,\|s\|_{L^2_{\varepsilon}}\Bigr).
\end{equation}
\end{enumerate}
\end{proposition}

\begin{remark}
    Before we prove the block estimates, we observe that (similarly as in the proof of \cite[Thm 8.21]{bismut1991complex}) the zeroth order terms $\mathcal{B}$ and $x^2 cl_{w,pt}(\tau V_{E})$ arising in \eqref{eq:rescaling-grad-term} will vanish in the expressions of the operator  $D_\varepsilon^{(1)}$ for very strong (i.e. pointwise) reasons, related to degree shifting. Both of these zeroth-order terms arise from Clifford multiplication, and are in particular degree-shifting. On the other hand, the $D_\varepsilon^{(1)}$ term is projected onto a section of the same Clifford degree (as a vertical Gaussian multiple of the volume form). Hence such quantities do not contribute to the local model, providing a statement of the form $J^{-1}D_\varepsilon^{(1)}J = D_F + \varepsilon^{f}\mathcal{E}$, rather than $J^{-1}D_\varepsilon^{(1)}J = D_\varepsilon^F + \mathcal{B} + x^2 cl_{w,pt}(\tau V_E)+\varepsilon^{f}\mathcal{E}$
\end{remark} 

\begin{proof}
\textbf{Proof of (1).} We consider the operator $J_{\varepsilon}^{-1} D_{\varepsilon}^{(1)}J_{\varepsilon}$. Since the sections in the image of $J_{\varepsilon}$ are the elements $H_{\mathcal{B},\varepsilon}^1(X;\mathcal{E})$, and since the projection $\Pi_{\varepsilon}$ restricts to the identity map on $H_{\mathcal{B},\varepsilon}^1(X;\mathcal{E})$, the operator reduces to $\Pi_{\varepsilon} D_{X,\varepsilon
}$ on $H_{\mathcal{B},\varepsilon}^1(F\, ; \mathcal{E}_{\mathcal{H}}\otimes \mathcal{N}_{F,\varepsilon})$. 

Given a section $\sigma \in H_{\mathcal{B},\varepsilon}^1(F\, ; \mathcal{E}_{\mathcal{H}}\otimes \mathcal{N}_{F,\varepsilon})$ it restricts to a local neighbourhood $U(F)$ as an expression of the form
\begin{equation}
\sigma_F=\!\left(k^{-1/2}\gamma_{t_0}^F(t)\,
\sum_{j=1}^{n_F}
\phi_{F}^*(\alpha_{j})\wedge \omega_{j,\varepsilon}
\right)
\end{equation}
where we recall that $\alpha_{j}\in H_{\mathcal{B},\varepsilon}^1(F^{reg}; \mathcal{E}_\mathcal{H})$ and the $\{\omega_{j,\varepsilon}(y,\cdot)\}_{j=1}^{n_F}$ are an orthonormal frame for the vertical Hodge cohomology \(\mathcal H^{\bullet}_{\mathcal{B},\varepsilon} \) over \(y\in F^{reg}\). 

By Proposition \ref{proposition_rescaled_estimates-0}, the difference between $J_{\varepsilon}^{-1}D_X J_{\varepsilon}$ and the model operator is contained in the error operator $\varepsilon^f R_\varepsilon$. Using this fact, and the Leibniz rule applied to \eqref{equation_local_expression_quasimode} we compute
\begin{equation}\label{equation_First_estimate_computation_cero}
    J_{\varepsilon}^{-1}\Pi_{\varepsilon}D_{X,\varepsilon} J_{\varepsilon}\sigma_a = \Pi_{\varepsilon}(D_{\mathcal{H}}+\varepsilon^{1/2}D_{\mathcal{V}}+\varepsilon^{1/2}cl(V_a))\Pi_{\varepsilon}\sigma_F + \varepsilon^fR_\varepsilon\sigma_F \\
    =: \Pi_{\varepsilon} D_{M,\varepsilon}\Pi_{\varepsilon}\sigma_F + \varepsilon^f R_\varepsilon\sigma_F
\end{equation}

Now we will use that the model operator $D_M$ on $H_\mathcal{B,\varepsilon}^1(X;E)$ splits into the vertical \eqref{eq:vert-ops} and horizontal \eqref{eq:hor-ops} operators $D_{\mathcal{H}}+D_{\mathcal{V}}$, that this latter operator acts trivially on the local sections $\omega_{j,\varepsilon}(x(y),z(y))$ of the frame, and that $\Pi_{\varepsilon} D_{\mathcal{H},\varepsilon} = D_F$. These claims are summarized in the following lemma.

\begin{lemma}\label{lem:1st-est}
Given a section $\sigma_F \in H_{\mathcal{B},\varepsilon}^1(F\, ; \mathcal{E}_{\mathcal{H}}\otimes \mathcal{N}_{F,\varepsilon})$
of the form
\begin{equation}
\sigma_F=\!\left(
u\,
\phi_{F}^*(\alpha_{j})\wedge \omega_{j,\varepsilon,\mu}
\right)
\end{equation}
we have that
\begin{equation}
    D_{X,\varepsilon} \sigma_F = cl(du)(\phi_{F}^*(\alpha_{j})\wedge \omega_{j,\varepsilon})+u\,(\phi_{F}^*(D_F\alpha_{j})\wedge \omega_{j,\varepsilon}) +\varepsilon^fR_\varepsilon\sigma_F
\end{equation}
\end{lemma}

\begin{proof}
We have that the vertical model operator acts on $\sigma$ by
\begin{multline} 
D_{\mathcal{V},\varepsilon}(u\,\omega_{j,\varepsilon})=D_{\mathcal{V}}(u\,\omega_{j,\varepsilon})+\varepsilon^{1/2}\sqrt{-1}\gamma_j x_j cl_{w,pt}(V_j)(u\,\omega_{j,\varepsilon})=cl(du)\omega_{j,\varepsilon} +u\, D_{\mathcal{V},\varepsilon}(\omega_{j,\varepsilon})=cl(du)\omega_{j,\varepsilon}
\end{multline}
where the last equality follows since $\omega_{j,\varepsilon}$ is in the null space of the model vertical operator.
With this we can see that
\begin{equation}\label{equation_First_estimate_computation_tres}
D_{M,\varepsilon}(u\,
\phi_{F}^*(\alpha_{j})\wedge \omega_{j,\varepsilon})=cl(du)(
\phi_{F}^*(\alpha_{j})\wedge \omega_{j,\varepsilon})+u\, D_{M,\varepsilon}(
\phi_{F}^*(\alpha_{j})\wedge \omega_{j,\varepsilon})
\end{equation}
where 
\begin{equation}
D_{M,\varepsilon}(
\phi_{F}^*(\alpha_{j})\wedge \omega_{j,\varepsilon})=
D_{\mathcal{H}}(\phi_{F}^*(\alpha_{j}))\wedge \omega_{j,\varepsilon}=(\phi_{F}^*(D_F\alpha_{j}))\wedge \omega_{j,\varepsilon}
\end{equation}
where in the last equality we use lemma \ref{lem:hor-op-conj}. The error operator acts boundedly by a 1st order differential operator, scaled by $\varepsilon^{f}$ where $f<0$. This proves the result.
\end{proof} 

From this we can now compute that \eqref{equation_First_estimate_computation_cero} reduces to 
\begin{multline}
\label{equation_First_estimate_computation_uno}
D_{X,\varepsilon} \sigma_F=D_{M,\varepsilon}(k_F^{-1/2}\gamma_{t_0}^F(t)) \!\left(\sum_{j=1}^{n_F}
\phi_{F}^*(\alpha_{j})\wedge \omega_{j,\varepsilon}
\right)\\+k_F^{-1/2}\gamma_{t_0}^F(t)  \left(\sum_{j=1}^{n_F}
\phi_{F}^*(D_F\alpha_{j})\wedge \omega_{j,\varepsilon}
\right)+\varepsilon^f R_\varepsilon\sigma_F.
\end{multline}
Now observe that
\begin{equation}
    D_{M,\varepsilon}(k_F^{-1/2}\gamma_{t_0}^F(t))\!\left(\sum_{j=1}^{n_F}
\phi_{F}^*(\alpha_{j})\wedge \omega_{j,\varepsilon}
\right) \leq C \exp(-\varepsilon C')
\end{equation}
for some constants $C, C'>0$ as the cutoff function $\gamma_{t_0}^F$ is identically $1$ on a smaller neighbourhood $t<t_0/5$ where the expression on the left hand side vanishes, and since away from this set there is an exponential damping given by the forms $\omega_{j,\varepsilon}$. This shows we can absorb the terms $D_{M,\varepsilon}(k_F^{-1/2}\gamma_{t_0}^F(t)) \!\left(\sum_{j=1}^{n_F}
\phi_{F}^*(\alpha_{j})\wedge \omega_{j,\varepsilon}
\right)$ in \eqref{equation_First_estimate_computation_uno} into $\varepsilon^{f} R_\varepsilon\sigma_a$. Further, arguing as in the lemma shows that $D_F \left( k_F^{-1/2}\gamma_{t_0}^F(t) \sum_{j=1}^{n_F}
\phi_{F}^*(\alpha_{j})\wedge \omega_{j,\varepsilon}
\right)$ differs from the expression $k_F^{-1/2}\gamma_{t_0}^F(t) D_F \left(\sum_{j=1}^{n_F}
\phi_{F}^*(\alpha_{j})\wedge \omega_{j,\varepsilon}
\right)$ up to a term which can be absorbed into $\varepsilon^{f} R_\varepsilon\sigma_F$.

Thus we have shown that 
\begin{equation}
\label{equation_First_estimate_computation_two}
   J^{-1} D_{\varepsilon}^{(1)} J \sigma_F=D_F \sigma_F + \varepsilon^{f} R_\varepsilon\sigma_F.
\end{equation}
We can easily take sums over sections $\sigma_F$ supported near different vanishing loci $F$ to extend this estimate for all sections in $H_{\mathcal{B},\varepsilon}^1(X;\mathcal{E})$. Since the projection $\Pi_{\varepsilon}$ is a bounded operator which projects onto the image of $J_{\varepsilon}$, on which $J_{\varepsilon}^{-1}$ is also a bounded operator, we can compose the expressions on both sides of \eqref{equation_First_estimate_computation_two} by $J_{\varepsilon}^{-1} \Pi_{\varepsilon}$ to get the required estimate in the statement.

\textbf{Proof of (2).}
It suffices to prove the claim for $D^{(3)}$, as the the other estimate will follow from self-adjointness.

First, if $s\in H_{\mathcal{B},\varepsilon}^{1}(X;\mathcal{E})$ then $s$ can be taken to be in the image of $J_{\varepsilon}$, so
\begin{equation}
   \Pi_\varepsilon^\perp D_{X,\varepsilon} \Pi_\varepsilon s = \Pi_\varepsilon^\perp D_{X,\varepsilon} J_{\varepsilon}(\tilde s) = \Pi_\varepsilon^\perp (D_F + \varepsilon^{f} R_\varepsilon) \tilde s,
\end{equation}
where we have in the last equality used \eqref{equation_First_estimate_computation_cero} and the fact that $D_{\varepsilon}^{\mathcal{V}}$ annihilates elements in the image of $J_{\varepsilon}$. Now because $\Pi_{\varepsilon}^\perp$ is bounded, it suffice to prove the desired estimate for $\Pi_\varepsilon^\perp (D_F + \varepsilon^f R_\varepsilon) \tilde s$.

The error polynomial term will contribute the $\varepsilon^{f}$ term in our estimate, so we focus now on the $D_F$ term. Applying lemma \ref{lem:1st-est} and \eqref{equation_First_estimate_computation_uno} we seek to bound
\begin{multline}
    \Pi_{\varepsilon}^\perp D_F\tilde{s} = \Pi_{\varepsilon}^\perp D_F\left( k_F^{-1/2}\gamma_{t_0}^F(t) \sum_{j=1}^{n_F}
\phi_{F}^*(\alpha_{j})\wedge \omega_{j,\varepsilon}
\right) \\
= \Pi_{\varepsilon}^\perp\left( D_{M,\varepsilon}(k_F^{-1/2}\gamma_{t_0}^F(t))\!\left(\sum_{j=1}^{n_F}
\phi_{F}^*(\alpha_{j})\wedge \omega_{j,\varepsilon}
\right)  \right) + 0
\end{multline}
as all other terms in the right hand side of \eqref{equation_First_estimate_computation_uno} lie in $\text{Im}(\Pi_{\varepsilon,\mu})$. However this remaining term can be estimated by \eqref{equation_First_estimate_computation_tres}, which combined with boundedness of $\Pi_{\varepsilon}^\perp$ completes the proof.

\textbf{Proof of (3)}
    To prove this final result we need to first compute how $D_{X,\varepsilon}$ acts on elements in the image of $\Pi_{\varepsilon}^\perp$. This reduces to considering three sub-cases: if $s\in H_{\mathcal{B},\varepsilon}^{1,\perp}(X;\mathcal{E})$ (i.e. in the image of $\Pi_{\varepsilon}^\perp$) then such a section $s$ must satisfy at least one of the following
    \begin{enumerate}
        \item $s$ is supported away from the critical point set $F$, i.e. $\text{supp}(s)\cap \text{supp}(\gamma)=\emptyset$. In this case invertibility of $cl(V)$ away from $Z(|V|^2)$ gives claim
        \item $s\in H_{\mathcal{B},\varepsilon}^{1,\perp}(X;\mathcal{E})$ is supported near $F$, and is acting on the subbundle of $H_{\mathcal{B},\varepsilon}^{0,\perp}(X;\mathcal{E})$ on which $\varepsilon^2 |V|^2$ is a strictly positive operator. In this case we can show $D_{\varepsilon}^M$ is invertible as a generalized Harmonic oscillator would be.
        \item $s\in H_{\mathcal{B},\varepsilon}^{1,\perp}(X;\mathcal{E})$ is supported near $F$, and is acting on sections living in a subbundle of $H_{\mathcal{B},\varepsilon}^{0,\perp}(X;\mathcal{E})$ in which $\varepsilon^2|V|^2$ fails to be uniformly bounded below; in this case $D_{\varepsilon}^a$ is no longer uniformly positive for all $\varepsilon$. But for each fixed $\mu$, we can show there exists a $\varepsilon>0$ sufficiently small (possibly depending on $\mu$) such that the sign-changing term $\sqrt{\varepsilon}cl(V_a)$ can be absorbed into the error term $\varepsilon^{f}R_\varepsilon$, for every $\delta\geq \tfrac{1}{2}$.
    \end{enumerate}
Proof of (1): first we note that the expansion 
\[  (D_{X,\varepsilon})^2 = (D_X + \sqrt{-1}\varepsilon cl(V))^2|_{H_{\mathcal{B}}^1} = (D_X)^2 + \varepsilon^2 |V|^2 + \varepsilon[D_X,cl(V)] \]
holds globally, and the commutator term is bounded on $H_{\mathcal{B}}^k$ for all $k$ and all $\varepsilon>0$. By a partition of unity argument we may assume that $s$ is supported away from the entire vanishing loci $Z(|V|^2)$; for our fixed $t_0>0$ we have that $|V|^2$ is invertible on $X\setminus \bigcup_{F\in Z(|V|^2)} U_{t_0}(F)$, and is in fact uniformly bounded below on sections supported in this set. From boundedness of the commutator term and ellipticity of $D_X$ (acting on sections with boundary conditions specified by $\mathcal{B}$) we can conclude, as in \cite[Prop 9.13]{bismut1991complex} that there exists $c, C>0$ such that
\begin{equation}
    \| D_{X,\varepsilon} s \|_{H_{\mathcal{B},\varepsilon}^0}^2   \geq C\left( \| s \|_{H_{\mathcal{B},\varepsilon}^1}^2 + (\varepsilon - c)\| s \|_{H_{\mathcal{B},\varepsilon}^0}^2 \right) .
\end{equation} 

Proof of (2): if $s$ is supported near some $F$ then $|V|^2$ may no longer be uniformly positive as in case (1). However, depending on our choice of $t_0>0$ defining $U_{t_0}(F)$ we can find a subbundle of $E$, denoted $E_+$, such that over $\bigcup_{F\in Z(|V|^2)} U_{t_0}(F)$ we have that $|V|^2$ is positive on $E_+$.\footnote{This can be achieved by noting that as an endomorphism on $E_+$, the eigenvalues of $|V|^2$ are continuous functions of $x\in X$ (recall that our bundles are defined over the regular part $F^{reg}$ and then extended by continuity). Then we may choose $E_+$ to be the subbundle consisting of eigenspaces with eigenvalue larger than the smallest non-zero eigenvalue of $|V|^2$ over $\bigcup_{F\in Z(|V|^2)} U_{t_0}(F)$.} According to the splitting $E=E_+\oplus E_-$, for sections in $E_+$ we are in the case of (1) and argue as before, so we now consider sections which land in $E_-$. 

Proof of (3): In this final step we will no longer have positivity to give the coercive lower bounds and will need to consider instead the elliptic estimates that arise from $D_{\varepsilon}^a$ acting on $E_-$. Now from proposition \ref{proposition_rescaled_estimates-0}, we have that if $s\in H_{\mathcal{B},\varepsilon}^{1,\perp}(X;E_-)$ then, letting the right hand side of \eqref{eq:rescaling-grad-term} be denoted $\varepsilon^{f}R_\varepsilon$ we have 
\begin{equation}
    (D_{X,\varepsilon} -\varepsilon^{f}R_\varepsilon) s = [D_\mathcal{H} + \varepsilon^{1/2}(D_{\mathcal{V}} + \sqrt{-1}\gamma_\alpha x_j cl_{w,pt}(V_j))]s ,  
\end{equation}
(From our polyhomogeneous expansion \eqref{Morse_Bott_Kirwan_modify_hamiltonian_vf}, the quantity $\sqrt{-1}\gamma_\alpha x_j cl_{w,pt}(V_j)$ should be compared with the quantity $\widetilde\nabla_Z^{\xi} V(y)$ arising in e.g. \cite[eqn (9.56)]{bismut1991complex}). Now because $\varepsilon^{f}R_\varepsilon$ satisfies the conclusion of the claim, we only need to verify the estimate \eqref{eq:rescaled-Garding-est} for the right hand side above. Restricted to sections landing in $E_-$ we have that $[D_{\mathcal{H}}, \sqrt{-1}\gamma_\alpha x_j cl_{w,pt}(V_j))]=0$ and $[D_{\mathcal{H}}, D_{\mathcal{V}}]=0$, and thus
\begin{equation}
    [D_\mathcal{H} + \varepsilon^{1/2}(D_{\mathcal{V}} + \sqrt{-1}\gamma_\alpha x_j cl_{w,pt}(V_j))]^2 = (D_\mathcal{H})^2 + \varepsilon \Big(D_{\mathcal{V}} + \sqrt{-1}\gamma_\alpha x_j cl_{w,pt}(V_j)) \Big)^2, 
\end{equation}
and hence 
\begin{equation}
   \| [D_\mathcal{H} + \varepsilon^{1/2}(D_{\mathcal{V}} + \sqrt{-1}\gamma_\alpha x_j cl_{w,pt}(V_j))]s \|_{L_{\varepsilon}^2}^2 = \| D_{\mathcal{H}}s \|_{L_{\varepsilon}^2}^2 + \varepsilon \| (D_{\mathcal{V}} + \sqrt{-1}\gamma_\alpha x_j cl_{w,pt}(V_j))s \|_{L_{\varepsilon}^2}^2.
\end{equation}
From here we observe that by the local ellipticity of $D_{X,\varepsilon}$ restricted to the local neighborhood, the operator $D^{\mathcal{V}} + \sqrt{-1}\gamma_\alpha x_j cl_{w,pt}(V_j)$ has discrete spectrum (see our Assumption \ref{assumption_on_metrics_and_domains}). From here the argument proceeds exactly as in the arguments following  \cite[eqn (9.75)]{bismut1991complex} to conclude this final case. This completes the proof.
\end{proof}

Now that we have established our blockwise estimates on $D_{X,\varepsilon}$ we can consider the construction of the spectral projectors. In order to do so we must first establish a comparison between the resolvents of the operators $D_{X,\varepsilon}$ and $D_F$. 
\begin{remark}
    The following demonstrates one situation where at minimum $D_F$ is forced to have discrete spectrum; in the proof of the proceeding proposition, we consider the trace norm of $(\lambda-D_F)^{-1}$, whose finiteness implies $D_F$ has purely discrete spectrum.
\end{remark}

For two bounded operators $A\in \mathcal{L}(H_{\mathcal{B}}^0(X;E))$ and $B\in \mathcal{L}(H_{\mathcal{B}}^0(\mathcal{N}(F);E))$ we define 
\begin{equation}
    d(A,B) = \| J_{\varepsilon}^{-1} A^{(1)}J_{\varepsilon} - B \|_1 + \sum_{j=2}^4 \| A^{(j)} \|_1 .
\end{equation}
where $\|-\|_1$ denotes the trace norm. Now fixing a constant $\kappa>0$ such that
\begin{equation}
    \text{Spec}(D_F)\cap [-2\kappa, 2\kappa] \subset \{0\},
\end{equation}
we can show

\begin{proposition}\label{prop:spec-comparison}
There exists $\varepsilon_0\geq 1$ such that for any $\varepsilon\geq \varepsilon_0$, $\lambda\in \mathbb{C}$ with $|\lambda|=\kappa$, then $\lambda-D_{X,\varepsilon}$ is invertible on $H_{\mathcal{B}}^0(X;E)$. Moreover for any integer $p\geq 2\,\emph{dim}X+2$, we have
\begin{equation}
    d\left((\lambda-D_{X,\varepsilon})^{-p},(\lambda-D_F)^{-p}\right) \leq C(1+|\lambda|^{p+1})\varepsilon^{f}
\end{equation} 
\end{proposition}
\begin{proof}
    Using proposition \ref{prop:block-estimates} this result can proved almost identically to the purely functional analytic arguments of \cite[\S 9.c-9.e]{bismut1991complex}, with one small modification; in our case the estimates involving $D_{\varepsilon}^{\ell}$ for $\ell=1,2,3$ involve $\varepsilon^{f}$ instead of $1/\sqrt{\varepsilon}$. Because $f<0$ for $\delta\geq \tfrac{1}{2}$, the proofs proceed exactly as written with the modification of the sets therein to be 
    \begin{equation}
        U_\varepsilon :=\left\{\lambda\in \mathbb{C} \, : \, |\lambda|\leq \frac{\kappa}{2} \varepsilon^{-f} , \; \inf_{\lambda_0\in \text{Spec}(D_F)} |\lambda-\lambda_0| \geq \frac{\kappa}{4} \right\} ,
    \end{equation} 
    at which point the claim follows in the same manner.
\end{proof}

\subsection{Witten subcomplexes and equivariant index theorems}

Now we prove the existence of Witten subcomplexes for families of $\mathcal{N}=1,2$ SUSY complexes. In the case of $\mathcal{N}=2$, for families arising from conjugation there are subcomplexes of the global complexes quasi-isomorphic to the Witten subcomplexes for all $\varepsilon \geq 0$, in particular one can prove Morse inequalities for the original complex. 
For $\mathcal{N}=1$, the Witten subcomplexes for large values of $\varepsilon \geq 0$ do not correspond to those for $\varepsilon=0$ in general, but the indices of the complexes are the same for all $\varepsilon \geq 0$, so we have a generalized Atiyah-Bott-Lefschetz fixed point formula (i.e.- an equivariant index formula). We formalize this as follows.

The following is the restatement of Theorem \ref{Theorem_Witten_complex_N_2_abstract_intro} that we stated in the introduction.

\begin{theorem}
\label{Theorem_Witten_complex_N_2_abstract}
Consider a family $\mathcal{P}_{\varepsilon}(X)$ of $\mathcal{N}=2$ SUSY complexes arising from a conjugation, with a Hamiltonian of zero-th order type where each element $F \in Z(W)$ has a localizing neighbourhood $U(F)$ which admits a local $\mathcal{N}=2$ refinement. Then for every $\varepsilon \geq 0$ there exists a Witten subcomplex which is finite dimensional.
\end{theorem}

\begin{proof}[Proof of theorem \ref{Theorem_Witten_complex_N_2_abstract}]
Proposition \ref{prop:spec-comparison} shows that the small eigenvalues of the Dirac operator on the entire space $X$ are well approximated by the spectrum of the Dirac operator $D_F$ of the normal cohomology complexes
$\mathcal{P}_{\mathcal{B}_a,\varepsilon}(\mathcal{N}_{F_a,\varepsilon})$ over all $F \in Z(A)$, for metrics which asymptotically $2\delta$ for $\delta>1/2$.

Thus the harmonic sections are (the cohomology) of the normal cohomology complex correspond are quasimodes for the small eigenvalue eigensections of the entire operator. We take the smallest non-trivial eigenvalues $\lambda'_F$ of $D_F$ on the horizontal Clifford modules on the vanishing sets (without twisting by the normal cohomology bundles), and minimize them over all the vanishing sets $F$ to get $\lambda'$, then for all $\varepsilon >\varepsilon_0$ where $\varepsilon_0$ is large enough, then the eigensections for the eigenvalues in $[0,\lambda'-\varepsilon_0^f)$ gives a Witten subcomplex $\mathcal{W}(\mathcal{P}_{W,\varepsilon}(X))$ of the complex $\mathcal{P}_{W,\varepsilon}(X)$.
The finiteness of the cohomology follows from the fact that each localized complex has finite dimensional cohomology, and that the number of connected components of $Z(A)$ are finite.

This proves the result for large $\varepsilon$.
The conjugation gives a quasi-isomorphism of complexes for $\mathcal{P}_{\varepsilon}(X)$ for each $\varepsilon \geq 0$ (see Proposition \ref{Proposition_cohomology_level_equivalences}). Thus taking the Witten subcomplexes $\mathcal{W}(\mathcal{P}_{W,\varepsilon}(X))$ and applying the quasi-isomorphism we get a subcomplex for all $\varepsilon \geq 0$.
\end{proof}

We now consider $\mathcal{N}=1$ SUSY complexes with operators $Q_{\varepsilon}=D_{\varepsilon}=D+i\varepsilon cl(V)$.

\begin{theorem}
\label{Theorem_N_1_Witten_complex}
Consider a $\mathcal{N}=1$ SUSY complex $\mathcal{P}_{\varepsilon}=(L^2(X;\mathcal{E}),\mathcal{D}(D_{\varepsilon}),D_{\varepsilon},\gamma)$ where $\gamma$ is a $\mathbb{Z}_2$ grading, where the family of Hamiltonians is of $0$-th order type with a potential $A$, and where $Z(A) \neq \emptyset$.
Then there exists $\varepsilon_0>0$ such that for each $\varepsilon>\varepsilon_0$, there is a Witten subcomplex.

In particular the index of $D_{\varepsilon}$ restricted to the Witten subcomplex is equal to the index of the operator $D_{\varepsilon}$ for any $\varepsilon >\varepsilon_0$ and for any $\varepsilon' \geq 0$, and
\begin{equation}
    \operatorname{Ind} _g(D)=\emph{Str}\,T_g|_{\mathcal{H}^k(\mathcal{P}_{0}(X))}=\sum_{F \in Z(A)}\emph{Str}\, T_g|_{\mathcal{H}^k(\mathcal{P}_{\varepsilon'}(U(F)))}
\end{equation}
\end{theorem}

\begin{proof}[Proof of Theorem \ref{Theorem_N_1_Witten_complex}]
The arguments in the proof we gave for Theorem \ref{Theorem_Witten_complex_N_2_abstract} shows that there exists a Witten subcomplex in this setting, for $\varepsilon>\varepsilon_0$ for some $\varepsilon_0>0$ which includes the harmonic sections of $D_{\varepsilon}$. However we no longer have an isomorphism of Hilbert complexes given by conjugation now, and thus cannot guarantee the existence of a Witten subcomplex for the undeformed complex.

Instead we have that for all $\varepsilon \geq 0$ $D_{\varepsilon}-D$ is a bounded 0-th order operator, in particular a compact operator from $\mathcal{D}(D_{\varepsilon})$ to $L^2(X;\mathcal{E})$ (here we can use the Kato-Rellich theorem \ref{Theorem_Kato_Rellich} to see that $\mathcal{D}(D_{\varepsilon'})=\mathcal{D}(D_0)$). This shows that the operators $D_{\varepsilon'}$ have the same Fredholm index for all $\varepsilon' \geq 0$.

The last equality follows, and we note that the localized complexes $\mathcal{H}^k \Big(\mathcal{P}_{\varepsilon'}(U(F)) \Big)$ have cohomologies
isomorphic for all $\varepsilon'$ since we assume the actions are K\"ahler Hamiltonian near each $F \in Z(A)$.
\end{proof}

\begin{remark}
Proposition 4.10 of \cite{Zhanglectures} generalizes in this context. That is, there is a homotopy of Fredholm operators
\begin{equation}
    D_{\varepsilon}(u)=D_{1,\varepsilon}+D_{4,\varepsilon}+u(D_{2,\varepsilon}+D_{3,\varepsilon}): \mathcal{D}(D_{\varepsilon}) \rightarrow L^2(X;\mathcal{E})
\end{equation}
for $\varepsilon \geq \varepsilon_0$ and $0 \leq u \leq 1$, and that $D_{4,\varepsilon}$ is invertible
\end{remark}

If $V$ has no zeros, then the index of the SUSY complex vanishes. This follows in our singular setting, for instance from the Lefschetz fixed point theorem in \cite{jayasinghe2023l2}. 
In fact, since in this case $Z(A=|V|^2)$ is empty, we have a general result for any $\mathcal{N}=N$. We formulate the following Theorem for completion.

\begin{theorem}
\label{Theorem_no_vanishing_sets}
Given a family $\mathcal{P}_{\varepsilon}$ of $\mathcal{N}=N$ SUSY complexes with a Hamiltonian of $0$-th order or 1-st order type with a potential $A$ which is positive ($Z(A)=\emptyset$), then the there is some $\varepsilon_0 \geq 0$ such that for all $\varepsilon >\varepsilon_0$, the Hamiltonian $\Delta_{\varepsilon}$ is a positive operator, and in particular is invertible.
\end{theorem}

\begin{proof}
For complexes of first order type we can use Proposition \ref{Proposition_equivariant_complexes} and reduce it to $0$-th order complexes, for which the result follows from Proposition \ref{Propostion_growth_estimate_witten_deformed}.
\end{proof}

\begin{remark}[SUSY and equivariant cohomology]
In equivariant cohomology on smooth manifolds $M$, if there is a vector field $V$ that has does not vanish, one can construct a form $\theta=V^{\flat}/|V|^2$, so that $\iota_V \theta=1$. Then
then given form $\alpha$ that is equivariantly closed (i.e closed with respect to $d_V=d-\iota_V$), we see that it is equivariantly exact, and realized as $\alpha=d_V\beta$ for $\beta=\theta \wedge (d_V \theta)^{-1} \wedge \alpha$, where $(d_V\theta)^{-1}=\sum_{k=0}^N (d\theta)^k$ where $N=\lfloor \dim M/2 \rfloor$.
Thus $\alpha,\beta$ are superpartners in this context.

When there are fixed points, the contributions come from fixed point sets, and the local contributions there correspond to computations in local equivariant cohomology groups.
Similar theories are yet to be studied in detail on stratified spaces, where the different choices of objects near strata (such as ideal boundary conditions) play roles in the analysis.
\end{remark}

\section{Morse inequalities, equivariant index theorems and applications}
\label{section_Morse_inequalities}

We proved the holomorphic Morse inequalities using the Witten complexes in Theorem \ref{Theorem_introduction_holomorphic_Morse_inequal_not_abstract}, using Theorem \ref{Theorem_main_Witten_complex_Kahler_case}, for circle actions.
In Theorem 5.13 of \cite{jayasinghe2024holomorphic}, the holomorphic Morse inequalities for the Serre dual complexes associated to a Dolbeault complex were made explicit, generalizing the version in the smooth setting that was first formulated in \cite{witten1984holomorphic}. This extends easily to the case of non-isolated critical points that we study in this article, and we leave it to the interested reader. 

In the next subsection we will study the case of Hamiltonian torus actions where we formulate Morse inequalities corresponding to each action chamber of the Lie algebra, as well as the walls that separate them. 
The Holomorphic Morse inequality above and the dual Morse inequality corresponds to two action chambers $\mathbb{R}^{\pm}$ for a given circle action (see Remark \ref{Remark_dual_Morse_polynomial_and_witten_complex}).

\subsection{K\"ahler Hamiltonian Torus actions}
\label{Subsection_Toric_actions}

Given a Hamiltonian action of a connected compact Lie group $G$ with moment map $M$, we can pick a maximal torus $T$ with Lie algebra $\textbf{t}=\sqrt{-1} Lie(T)$, with $\mathcal{L}$ the integral lattice in $\textbf{t}$. The dual lattice $\mathcal{L}^*$ in $\textbf{t}^*$ is the weight lattice, where $\textbf{t}_\mathbb{R}$ is the real Lie algebra. The following is Definition 2.1 of \cite{wu1998equivariant}.
\begin{definition}
    Let $\mathbb{Z}[\mathcal{L}^*]$ be the ring of formal characters of $T$, consisting of the elements $q=\sum_{\xi \in \mathcal{L}^*} q_{\xi} e^{\xi}$, for $q_{\xi} \in \mathbb{Z}$. We say $q \geq 0$ if $q_{\xi} \geq 0$ for all $\xi \in \mathcal{L}^*$. 
\end{definition}

While the ring of formal characters is sometimes denoted $\mathbb{Z}[[\mathcal{L}^*]]$ (to emphasize a completion), we will use the notation in \cite{wu1998equivariant}. This ring, expressed in terms of the dual integral lattice is an efficient way of formulating certain results, especially related to geometric quantization. In the smooth setting, the holomorphic Morse inequalities are polynomials with values in the character variety, up to some factors corresponding to the twists by vector bundles $E$. We incorporate this into our notation since these twists appear in the local cohomology groups.

\begin{example}
    Consider the Dolbeault complex on $\mathbb{CP}^1$ with coefficients in the prequantum line bundle $E_2$, which in this case is a square root of the canonical line bundle (the Dolbeault Dirac operator is in this case isomorphic to the spin Dirac operator). Consider the action $(\lambda)[Z_0:Z_1]=[Z_0:\lambda Z_1]$.
    The Morse polynomial in this case is $(1+b) \sum_{j=0}^\infty \lambda^{1/2+j}$, which is both a polynomial in $b$ with coefficients in $\mathbb{Z}[\Lambda]$ up to a factor $\lambda^{1/2}$. If we consider $\lambda=e^{i\theta}$, this series converges for complex values of $\theta$ with positive imaginary part. Here $\lambda=e^{\xi}$ for $\xi \in \sqrt{-1}\mathbb{R}$, and where $\Lambda$ is the integral lattice in the dual Lie algebra $(\sqrt{-1}\mathbb{R})^*$
    It is also an element of $\mathbb{Z}_{E_2}[b][\Lambda]=\mathbb{Z}[b][\Lambda+1/2]$.

    Similarly, if one considers twists by $E_N$ which are $N$-th roots of the canonical bundle, the Morse polynomials are elements of $\mathbb{Z}_{E_N}[b][\Lambda]=\mathbb{Z}[b][\Lambda+1/N]$.

    For twists by $E=E_{2}\oplus E_{3}$, the Morse polynomial is the sum of the Morse polynomials of $E_{N}$ for $N=2,3$. This is an element of $\mathbb{Z}_{E}[b][\Lambda]=\mathbb{Z}_{E_2}[b][\Lambda]\oplus\mathbb{Z}_{E_3}[b][\Lambda]$. This is an example where the lattice is a direct sum of shifted integral lattices, which is not itself an integral lattice.    
\end{example}

The splitting principle for holomorphic bundles over the regular part allows us to define $\mathbb{Z}_E[b][\mathcal{L}^*]$ for general Dolbeault complexes. 
\begin{definition}
Given an action of a torus on a vector bundle $E$ on a space $X^{reg}$ as we consider in this article, we define
\begin{equation}
    \mathbb{Z}_E[b][\mathcal{L}^*]:=\bigoplus_{j} \mathbb{Z}_{L_j}[b][\mathcal{L}^*]
\end{equation}
where $L_j$ are the summands for $E$ given by the splitting principle on $X^{reg}$, where $\mathbb{Z}_{L_j}[b][\mathcal{L}^*]$ are the rings of formal power series with a shifted integral lattice.

\noindent We use a similar notation for the ring of formal power series with the complexified integral lattice $\mathbb{Z}_{E}[b][\mathcal{L}^* \otimes \mathbb{C}]$.
\end{definition}

\begin{remark}
The local cohomology groups differ as the ideal boundary conditions change. In most cases, it appears that the choice $\mathbb{Z}_E[b][\mathcal{L}^*]$ is independent of the choice of ideal boundary conditions, including in all examples studied in this article, but we do not attempt to prove this.

Hence we note that this ring may $E$ depend on the ideal boundary conditions in general, even though we do not denote this explicitly since the notation is already complicated.
\end{remark}

\begin{definition}[Non-negativity and positivity]
\label{definition_non_negativity_and_positivity}
    We say $Q(b) \geq 0$ if $q_k \geq 0$ in $\mathbb{Z}_E[\mathcal{L}^*]$ for all $k$. For two such polynomials $Q_1(b), Q_2(b)$, we define that $Q_1(b) \leq Q_2(b)$ if $Q_1(b)-Q_2(b) \geq 0$, and $Q_1(b) > Q_2(b)$ if $Q_1(b)-Q_2(b) > 0$.
\end{definition}

\begin{definition}[Intersection of polynomials]
\label{definition_intersection_of_polynomials}
Given finitely many polynomials $M_1,M_2,... M_l$ of $\mathbb{Z}_E[b][\mathcal{L}^*]$, each of which can be expressed as 
\begin{equation}
\label{equation_polynomial_temp_121}
    M_k=\sum_q b^q \sum_{\xi \in \mathcal{L}^*+\sigma} q^{k}_{\xi} e^{\xi}
\end{equation}
where $\sigma$ corresponds to shifts of the lattice corresponding to the twisting by $E$. We define $M_0=\bigcap_{k=1}^l M_k$ to be the polynomial where in the same notation as \eqref{equation_polynomial_temp_121} the coefficients $q^{0}_{\xi}:= \min \{ q^{1}_{\xi}, q^{2}_{\xi},...,q^l_{\xi} \}$.    
\end{definition}

\begin{example}
Consider the polynomials 
\begin{equation}
    M_1=1+(1+b)\sum_{\xi \in \mathbb{Z}^+} e^{\xi}, \quad M_2=1+(1+b) \sum_{\xi \in \mathbb{Z}^-} e^{\xi},
\end{equation}
where the lattice is $\mathcal{L}^*$. This is simply the Morse polynomials for the two action chambers for the standard toric action on $\mathbb{CP}^1$. If we pick a vector field in $i\theta V \in \sqrt{-1}Lie(S^1)=\sqrt{-1}\mathbb{R}^1$, we have that $\xi(i\theta V)=i\theta$, we can compare this with Example 1.6 of \cite{jayasinghe2024holomorphic}.
Here it is easy to see that the only $q^0_{0}=1$ is non-trivial, and thus the polynomial $M_1\cap M_2=1$ is the Poincar\'e polynomial.
\end{example}

We briefly describe what an \textbf{action chamber} is and how it relates to the Hamiltonians on the tangent cones. 
Let $G$ be a compact connected Lie group acting on $X$, and
choose a maximal torus $T\subset G$, with Lie algebra
$\mathfrak t$. For $\xi\in\mathfrak t$, denote by $V_\xi$
the corresponding vector field on $X$. Every element
of $\mathfrak g$ is conjugate to an element of $\mathfrak t$, and
the remaining ambiguity is the action of the Weyl group
$\mathcal{W}=N_G(T)/T$.

We denote the generically fixed points by $X^T=\bigsqcup_a F_a$ for all $F_a \in Z(A)$  and that the normal fibers at each $F_a$ admits a $T$-invariant decomposition into conic factors $C_{x_{a,j}}(Z_{a,j})$, with fixed
normalized Reeb vector fields $R_{a,j}$, such that
\begin{equation}
    V_\xi\big|_{C_{x_{a,j}}(Z_{a,j})}
      =\alpha_{a,j}(\xi)R_{a,j},
    \qquad \alpha_{a,j}\in\mathfrak t^*\setminus\{0\}.
\end{equation}   
Here the linear functionals $\alpha_{a,j}$ incorporate
the constant rescalings in our assumption on the action.
With the convention $\iota_{V_\xi}\omega=-dh_\xi$,
the Hamiltonian on the normal tangent cone is
\begin{equation}
    h_\xi-h_\xi|_{F_a}
      =\frac12\sum_j\alpha_{a,j}(\xi)x_{a,j}^2,
\end{equation}   
matching the definitions in Subsection \ref{subsection_fixed_point_sets}.

\begin{definition}
Assume that the collection of hyperplanes
$\ker\alpha_{a,j}$ is finite. An \emph{action chamber}
for the $G$-action, relative to $T$, is a connected
component of
\begin{equation}
    \mathfrak t\setminus
      \bigcup_{a,j}\ker\alpha_{a,j}.
\end{equation}   
Equivalently, two elements belong to the same action
chamber precisely when all the coefficients
$\alpha_{a,j}(\xi)$ are nonzero and have the same signs.
Thus the action chambers are open convex polyhedral
cones, obtained by taking the common refinement of
the sign decompositions at all components of $X^T$.
\end{definition}

Assume moreover that the normal conic data are compatible
with the action of $N_G(T)$. Then the Weyl group permutes
the wall hyperplanes and hence the action chambers.
Changing the maximal torus transports this decomposition
by conjugation and perturbations within a chamber preserve
the signs of the normal Hamiltonians at every fixed
component.
 
Given an Hamiltonian $G$ action with moment map $\mathcal{M}$ on a space with an extension of the Dolbeault complex $W$, with generic fixed point sets indexed by $crit(\mathcal{M})$, we define the \textbf{local equivariant Morse polynomials}
\begin{equation}
    M_F(g)(b)=\sum_{q=0}^n b^q \operatorname{Tr} T_g|_{\mathcal{H}^q(\mathcal{P}_{W,B(g)}(U(F)))},
\end{equation}
and so $M_F \in \mathbb{Z}_E[b][\mathcal{L}^* \otimes \mathbb{C}]$, and $M_a$ converges for positive imaginary valued elements of $\xi \in \mathcal{L}^* \otimes \mathbb{C}$.

Here the boundary conditions $B(g)$ depend only on the choice of action chamber $C$ containing the generator of $g$. Since the local Morse polynomials are analytic in the dual Lie algebra for a fixed action, it is easy to see that the coefficients 
$q_{\xi}$ in \eqref{equation_polynomial_temp_121} also depend only on the action chamber. Thus for a given action chamber $C$, we define
\begin{equation}
    M_{F,C}(g)(b)=\sum_{q=0}^n b^q \operatorname{Tr} T_g|_{\mathcal{H}^q(\mathcal{P}_{W,B(C)}(U(F)))}.
\end{equation}

For a given action, since the Lie algebra of the maximal torus $T$ is also that for $G$ we see that the \textbf{(global) equivariant Morse polynomial} $M_C=\sum_{F \in crit(\mathcal{\mathcal{M}(\xi^*)})} M_{F,C}$ (where $\xi^*$ is a generic element of the Lie algebra corresponding to generic fixed points) is an element of $\mathbb{Z}_E[b][\mathcal{L}^* \otimes \mathbb{C}]$ as well.

Similarly we define the \textbf{equivariant Poincar\'e polynomial}
\begin{equation}
    P(g)(b)=\sum_{q=0}^n b^q \operatorname{Tr} T_g|_{\mathcal{H}^q(\mathcal{P}_{W}(X))}
\end{equation}
which is a finite multi-variable polynomial in $\mathbb{Z}[\mathcal{L}^* \otimes \mathbb{C}]$ for each fixed $b$, defined independent of the action.
Then for each element $g$ which lies in a circle action generated by an element in an action chamber $C$, we have a Morse inequality
\begin{equation}
    M_C(g)(b)=P(g)(b)+(1+b)\sum_{q=0}^{n-1}b^qQ_q(g)
\end{equation}
where $Q_q \in \mathbb{Z}_E[b][\mathcal{L}^* \otimes \mathbb{C}]$ are non-negative in the sense of \ref{definition_non_negativity_and_positivity}.

We formulate the following strengthening of the Morse inequalities for general actions.

\begin{theorem}[Combined upper bound]
\label{Theorem_strong_Morse_Dolbeault_version_groups}
Given a stratified pseudomanifold with a holomorphic bundle $E$, and a twisted Dolbeault complex with a choice of self-adjoint extension indexed by $W$, and a K\"ahler Hamiltonian moment map $\mathcal{M}$, we define the \textbf{combined upper bound}
\begin{equation}
\label{equation_total_Morse_polynomial}
M_{\cap}(b):=\bigcap_C M_C(b),
\end{equation}
where the intersection takes the minimum coefficient in each degree
and weight. Then
\begin{equation}
P(b)\leq M_{\cap}(b),
\qquad
M_{\cap}(b)=P(b)+R(b),
\qquad R(b)\geq 0.
\end{equation}
\end{theorem}

\begin{proof}
Since \(Q_C(b)\geq 0\), every coefficient of $M_C(b)-P(b)=(1+b)Q_C(b)$ is nonnegative. Taking the minimum over chambers in each degree and weight gives the result.
\end{proof}

The action of the Weyl group permutes the action chambers, which can be used to understand the different choices of action chambers for the action. We refer the reader to Section \ref{section_examples} for some examples where Morse inequalities for different action chambers are discussed, and also to Example 6.1 of \cite{jayasinghe2024holomorphic}.

\begin{remark}
This result can be used to generalize Proposition 2.6 of \cite{wu1998equivariant} for the singular spaces we study. 
The result is also the analog of the refined Morse polynomial for the de Rham complex in Corollary 1.6 of \cite{jayasinghe2024witten}.
\end{remark}

\begin{proof}[Proof of Theorem \ref{Theorem_strong_Morse_Dolbeault_version_groups}]
    Since each Morse polynomial corresponding to an action chamber participates in a Morse inequality, taking the intersections of them in the sense of Definition \ref{definition_intersection_of_polynomials} will still participate in a Morse inequality. In particular, the terms $Q^q$ of the intersected error polynomial will be given by intersecting the corresponding terms of the error polynomials for the action chambers, while the Poincar\'e polynomial does not depend on the action chamber chosen. 
\end{proof}

\begin{remark}[Dual Morse polynomials and dual Witten complexes]
\label{Remark_dual_Morse_polynomial_and_witten_complex}
In the case of isolated critical points in \cite{jayasinghe2024holomorphic}, dual Morse polynomials were defined in terms of Serre dual and adjoint Hilbert complexes for given complexes. It is easy to check that the correct dual Witten complex (corresponding to the dual Morse inequalities presented in \cite{witten1984holomorphic}) corresponding to an action $C$ is given by the Witten complex for $-C$.
We leave it to the reader to compare this with the dual Morse complexes in \cite{jayasinghe2024holomorphic}.
\end{remark}

\begin{remark}
    Note that one cannot construct instanton complexes that categorify these stronger Morse polynomial where the categorification is by local cohomology complexes that are localized near a given fixed point. 
    This is already clear in the example of spinning the sphere, here for two different actions, the non-trivial global harmonic section localizes near the two different fixed points for the rotations in different directions.

    Concretely, for $\mathbb{CP}^1$ with the action $(\lambda)[Z_0: Z_1]=[Z_0:\lambda Z_1]$, for $\lambda=e^{i\theta}$ with positive $\theta$, the Witten deformed constant section localizes near $[1:0]$ in the deformed limit, while it localizes near $[0:1]$ when $\theta$ is negative.

    Of course, the quantum mechanical state space can be treated as an abstract Hilbert space and an abstract finite dimensional complex corresponding to these stronger Morse inequalities can be formulated, but at present it lacks motivation.
\end{remark}

\subsubsection{Comparison with other equivariant Riemann-Roch formulae}
\label{Subsec_comparison_Lefschetz_smooth}

Since it is well known that the (holomorphic) Morse inequalities reduce to the (holomorphic) Lefschetz fixed point theorem for $b=-1$, it is instructive to compare our formulas with other known local index formulas. 
Even in the case of isolated critical points, objects used to describe the local Lefschetz numbers as equivariant eta invariants, linear algebraic formulas (for orbifolds) are different from the objects describing instanton complexes that categorify the (holomorphic) Morse inequalities. 
This was investigated in \cite{jayasinghe2026spectralasymmetrysupersymmetryequivariant} where a companion paper is being compiled where the relation to $\mathcal{N}=2$ SUSY is explained (see Appendix \ref{Appendix_SUSY_and_conjectures}). 
A quick description of the results there, in the complex setting, is that while the eta invariants are alternating traces over the spectrum of an associated Dirac-type operator on the link at a singularity, not all the eigensections contribute to the trace. There are additional supersymmetric cancellations in the trace (corresponding to the second supersymmetry), and only the sections corresponding to local cohomology groups on the cones contribute non-trivially to the trace.

The special formulas for orbifold equivariant Lefschetz numbers \cite{vergne1994quantification} can be reformulated as supertraces over the local cohomology group of a smooth local cover at the fixed points, up to extra factors corresponding to the covering degree. However this point of view is not intrinsic to the singular space, and does readily relate to descriptions of other singular analytic objects arising in the analytic descriptions of equivariant index defects (such as eta and xi invariants \cite{jayasinghe2026spectralasymmetrysupersymmetryequivariant}) or Poincar\'e residues arising in algebraic versions. In fact the work of \cite{Degeratuthesis} that derives orbifold equivariant indices were given categorical explanations at the level of harmonic section in \cite{jayasinghe2026spectralasymmetrysupersymmetryequivariant}. We explain the role of SUSY cancellations related to that work in Appendix \ref{Appendix_SUSY_and_conjectures}.

In Remark 2.3 of \cite{wu1998equivariant}, it is observed that the local holomorphic Lefschetz numbers over a smooth fixed set $F$ are given by
\begin{equation}
\label{equation_smooth_equivariant_index}
    \sum_{q=0}^n (-1)^q \operatorname{Tr} T|_{\mathcal{H}^q(U(F);E)}=\int_{F} ch_T(K^C(N_r) \otimes E|_{F}) td(F)=\int_{F} ch_T \Bigg(\frac{E|_{F}}{det(1-N^{*(1,0)})} \Bigg)td(F)
\end{equation}
and we refer to \textit{loc. cit.} for more details, where we have rephrased the cohomology group in the notation and formulations of our article. In the case of a Spin Dirac operator, where $E|_{F}$ would correspond to a square-root bundle $L$ of the canonical bundle (i.e. $L^{\otimes 2} =\mathcal{K}$), it is well known that the product of the Chern class of $L|_{F}$ with the Todd class $td(F)$ gives the $\widehat{A}$ class on $F$. We work out an example on a smooth toric manifold in subsection \ref{subsection_smooth_toric_action} where the expression on the right hand side of \eqref{equation_smooth_equivariant_index} are also explained following Duistermaat \cite{duistermaat2013heat}.

Now consider the index formula for stratified pseudomanifolds in \cite{Albin_2017_index}, where the contribution from the regular part of a singular stratum $Y$ of any depth (where the local fibration is as in Definition \ref{iterated_fibration_structure}) is given by the integral 
\begin{equation}
\label{singular_index_contrib}
    \int_{Y} \widehat{A}(Y)\mathcal{J}(\mathcal{B}_Y/Y)
\end{equation}
where we refer the reader to \textit{loc. cit.} for the definition of the Bismut-Cheeger $\mathcal{J}$ form over $Y$. Comparing with equation \ref{equation_smooth_equivariant_index} with a \textbf{conjectural singular equivariant formula} generalizing equation \ref{singular_index_contrib} (at least at a strata of maximal depth), it is easy to see that the equivariant Chern class $ch_T(K^C(N_r)|_{F})$ admits a description as an equivariant Bismut-Cheeger $\mathcal{J}$ form for the normal disc bundle over $Y$, viewed as a boundary fibration. For the case of isolated fixed points this follows from the results in \cite{jayasinghe2023l2}, and from that it is easy to extend to the case of flat normal bundles in the sense of \cite{jayasinghe2024witten} as well. 

In an upcoming revision of \cite{jayasinghe2026spectralasymmetrysupersymmetryequivariant}, we show that just as the $\xi_T$ invariant (which records the equivariant index defect at isolated fixed points) corresponds to a supertrace over the local cohomology group, as appearing in the left hand side of \ref{equation_smooth_equivariant_index}, the non-equivariant $\xi$ invariant at an isolated singularity  can be phrased entirely in terms of local \textit{Dirac cohomology groups}.

\textbf{We conjecture that this perspective extends even to non-isolated singularities, giving a reformulation of the index theorem of \cite{Albin_2017_index} in suitable settings (see Section \ref{Appendix_SUSY_and_conjectures}).}

We mention that there are links to other formulations of equivariant Riemann-Roch theorems, particularly the work of Edidin-Graham in \cite{equivariant_intersection_edidin,localization_algebraic_Graham_Edidin} (see Remark 7.32 of \cite{jayasinghe2023l2}) that generalize the work of Baum-Fulton-Macpherson-Quart \cite{baum1979lefschetz,baumfultonmacKtheoryRiemannRoch,baumfultonmacpheresonRiemannRoch,Baumformula81}, which were discussed in the precursor work of \cite{jayasinghe2024holomorphic}, \textbf{in particular formulating conjectural Morse inequalities for these algebraic equivariant Riemann-Roch results} (see Section 6.2 of \cite{jayasinghe2024holomorphic}).

In \cite{wu1998equivariant}, Wu constructs examples which show that the holomorphic Morse inequalities fail for certain toric actions on complex manifolds which are not K\"ahler Hamiltonian actions. It is also identified that compact, connected Lie group actions are the correct general setting for the result, where in \cite{wu2003instanton} he shows that a meromorphic group action on a complex manifold with a filterable Bialynicki-Birula decomposition (weaker than the K\"ahler assumption in this paper) suffices.
In section 2 of \cite{wu1998equivariant}, the author shows how more information can be drawn in the case of a Lie group action. This was already used to prove vanishing results for harmonic spinors in certain algebraic settings (see Section 7.3 of \cite{jayasinghe2024holomorphic}), but we extend this formulation more precisely in Theorem \ref{Theorem_strong_Morse_Dolbeault_version_groups} in our geometric setting.

In Section 5 of \cite{wu1996equivariant}, Wu also studies the case of non-Abelian group actions including relationships to the Weyl character formula. These are beautiful results in the smooth setting, but potential generalizations are more complicated in our setting as can be checked in examples and a more careful study is not in the scope of this article.

\subsection{Rigidity and Poincar\'e Hodge numbers}

The rigidity of equivariant invariants is a topic of great interest, and one can trace its roots to at least the rigidity for Spin Dirac operators by Atiyah-Hirzebruch in \cite{atiyah1970spin}. Witten's questions in \cite{witten1983fermion} gave the impetus to prove rigidity of Rarita-Schwinger operators and more general elliptic genera, and we refer to Section 7 of \cite{jayasinghe2024holomorphic} for a more thorough review of the literature.

\begin{definition}[Rigid equivariant indices and rigid complexes]
Given the equivariant index $\operatorname{Ind} _{g}(D)$ of a Fredholm operator $D$, where $g \in G$ for some group $G$, we say the \textbf{equivariant index is rigid} if $\operatorname{Ind} _{g}(D)$ is independent of $g$ for all $g \in G$, in which case it is equal to the index $\operatorname{Ind} _e(D)$ where $e$ is the identity element.

\noindent Given a family of Hilbert complexes as in Proposition \ref{Proposition_equivariant_complexes} with a Hamiltonian of 1st order type, we refer to the equivariant Hilbert complex $\mathcal{P}_{\varepsilon}^{\mu=0}$ as the \textbf{rigid subcomplex} of $\mathcal{P}_{\varepsilon}$.
\end{definition}

In this subsection, we study the Hilbert complexes 
\begin{equation}
\label{equation_Hilbert_complexes_rigidity_section}
    \mathcal{P}^d_{\min}(X)=(L^2\Omega^{p,q}(X;E),\mathcal{D}_{\min}(d_E),d_E), \quad \prescript{p}{}{\mathcal{P}}_{W}(X)=(L^2\Omega^{p,\bullet}(X;E),\mathcal{D}_{W}(P),P)
\end{equation}
where $E$ is a bundle that is \text{both flat and Hermitian}, $P=\overline{\partial}_{\Lambda^{p,0} \otimes E}$ and $\min$. We also study versions for fundamental neighbourhoods $U(F)$ with various choices of the boundary conditions $\mathcal{B}$ that arise in our study.
We shall denote by $\prescript{\bullet}{}{\mathcal{P}}_{W}$ the direct sum over all $p$ of the complexes over $\prescript{p}{}{\mathcal{P}}_{W}$, and denote their cohomology with the $(p,q)$ bi-degree by $\mathcal{H}^{p,q}(\prescript{\bullet}{}{\mathcal{P}}_{W})=\mathcal{H}^{q}(\prescript{p}{}{\mathcal{P}}_{W})$.

\subsubsection{Rigidity of de Rham and Dolbeault cohomology}

We recall that the de Rham complexes on complex spaces with only even codimensional strata have $L^2$ cohomology that is independent of the choice of wedge metric (while varying as one changes more broadly), indeed isomorphic to intersection cohomology with middle perversity (both lower middle and upper middle perversities give the same cohomology in this setting (see \cite{jayasinghe2023l2} for a survey of the literature). 
We also recall that the Dolbeault cohomology of a smooth manifold is independent of metrics, and for any Hermitian metric has a unique harmonic representative, which in the K\"ahler setting are the same representatives as those for the de Rham cohomology.

The following result considers the generalization of this to the singular setting and formulates a related statement in local cohomology.

\begin{theorem}
\label{theorem_rigidity_preliminary}
    Given a resolved stratified space $X$ with a wedge Hermitian metric, then the $L^2$ cohomology groups $\mathcal{H}^{q}(\prescript{p}{}{\mathcal{P}}_{\min}(X))$ are independent of the choice of wedge Hermitian metric, and in the K\"ahler setting are isomorphic to the cohomology groups $\mathcal{H}^{p,q}(\mathcal{P}^d_{\min}(X))$.

   \noindent Moreover, given a fundamental neighbourhood $U(F)$ of a critical point set of a local K\"ahler isometric action, 
    \begin{equation}
    \mathcal{H}^{p,q} \Big(\mathcal{P}^d_{\min,B}(U(F)) \Big) \subset \mathcal{H}^{q} \Big( \prescript{p}{}{\mathcal{P}}^{\mu=0}_{W,B}(U(F)) \Big)
    \end{equation}
    where the choice of generalized Neumann/Dirichlet conditions for $d_E$ and $P$ are the same for the respective complexes on the factors $C(Z_j)$ appearing in Definition \ref{definition_non_degeneracy_critical_point_sets_Kahler}.
\end{theorem}

\begin{remark}
    We emphasize that the local rigidity result holds for locally K\"ahler actions near fixed point sets of the type we study.
\end{remark}

\begin{proof}
    Given any wedge Hermitian metric compatible with the complex structure $J$, the minimal domain for $P$ gives rise to a rolled up extension for the Dirac operator $P+P^*$ and the harmonic sections of this operator are isomorphic to $ker(P)/im(P)$ (with the domain). Since the operator $P$ is given by the projection to the anti-holomorphic form parts after operating by $d$, the formal operator is independent of the metric and this persists for the case of twists by holomorphic vector bundles that are flat. 
    
    In the K\"ahler setting, these groups are isomorphic to the de Rham cohomology groups with the $(p,q)$ decomposition. This is clear if the metrics of the links at the strata have small volumes so that the spectral Witt condition is satisfied and the Dirac operators are essentially self-adjoint. Then the Laplace-type operators are equal for the de Rham and Dolbeault complexes (up to a constant factor), and the isomorphism holds at the level of harmonic sections. 

    For general metrics, while the Laplace-type operator for the de Rham complex is not essentially self adjoint, it is well known that the cohomology groups are independent of the choice of metric. The smooth sections with compact support in the interior give a smooth core and the cohomology of the core gives the cohomology of the complex. Since the formal operators are independent of the metric, we see that the cohomology of the minimal domains are independent of the choice of Hermitian metric, and since the core is the same for both de Rham and Dolbeault complexes and the Laplacians are formally equal upto a constant, we have the isomorphism of cohomology groups as stated in the theorem statement.
        
    To understand the inclusion of the local de Rham cohomology into the rigid Dolbeault cohomology with the minimal ideal boundary conditions, we first consider the case of a cone $U=C_x(Z)$ where the local K\"ahler isometry is generated by the Reeb vector field $V$ on the Sasaki link.

    The de Rham cohomology of a cone $\mathcal{H}^{p,q}(\mathcal{P}^d_{\min,N}(U))$ was computed in Lemma 6.16 of \cite{jayasinghe2023l2} where the harmonic representatives $\phi$ were shown to be given by the pullback under $\pi_Z: C(Z) \rightarrow Z=\{x=1\}$ of the harmonic representatives of the cohomology of $Z$, of degree $k \leq [\frac{l-1}{2}]$ where $l$ is the dimension of the link $Z$.
    In particular that ${\partial_x} \phi=0$ (short for $\nabla_{\partial_x} \phi=0$).
    We observed in the proof of Proposition 7.1 of \cite{jayasinghe2023l2} that these harmonic representatives, when restricted to the link $Z$ are in the basic cohomology corresponding to the Reeb foliation, in particular $L_V \phi=0$. This shows that they are in the rigid Dolbeault cohomology, and it is easy to check that they satisfy the generalized Neumann conditions for $P$.

    The case for $-V$ follows by duality, and we have that $\mathcal{H}^{p,q}(\mathcal{P}^d_{\min,D}(U))$ also satisfies the inclusion.
    If the fixed point set $F$ is not isolated, then we have the inclusions above on the normal fibers. The cohomology groups for the localized complexes are isomorphic to the cohomology groups of the normal cohomology complexes, and the above arguments for a cone $U$ shows that the inclusion holds for the fiber cohomology groups. Then since the Laplace operators restricted to $F$ differ by a constant factor, the inclusion holds as in the statement of the Theorem.
\end{proof}

The following example shows that Theorem \ref{theorem_rigidity_preliminary} fails in general for domains that are not minimal, even with just isolated singularities and isolated fixed points.

\begin{example}[Spherical suspension of a circle of length $4\pi$]
\label{example_spherical_suspension_circle_4pi}
A spherical suspension of a stratified pseudomanifold $Z$ with a wedge metric $g_Z$ is the stratified pseudomanifold given by the metric closure of $[0,\pi]_{\phi} \times Z$ with respect to the metric $d\phi^2+\sin^2(\phi)g_Z$.

We consider the spherical suspension of the circle of length $4\pi$, which has the K\"ahler wedge metric $d\phi^2+\sin^2(\phi)d\theta^2$ or $\phi \in [0,\pi]$. The circle action generated by $\partial_{\theta}$ is Hamiltonian with moment map $h=\cos(\phi)$. We pick the maximal domain for $\overline{\partial}$. The two singular points have tangent cones that are both cones $C_r(S^1)$ over circles with length $4\pi$. The cohomology of the maximal domain for this operator has a basis $\{z^{k}\}_{k=-1}^\infty$, where $z=r^{1/2}e^{i\theta/2}$. The cohomology of the adjoint complex has basis $\{\overline{z}^kd\overline{z}\}_{k=1}^{\infty}$.

Since the two holomorphic coordinates at the two singularities are inverses of each other (akin to the Riemann Sphere), for an action that take $z$ to $\lambda z$, where $\lambda=e^{i\alpha}$ for $\alpha>0, \alpha \neq 4m\pi$ for $m \in \mathbb{Z}$, we have that the Morse polynomial is $\lambda^{-1}+1+\lambda +(1+b)(\lambda^2/(1-\lambda))$. For the opposite action chamber, with $\lambda^{-1}$ the Morse polynomial is $\lambda+1+\lambda^{-1} +(1+b)(\lambda^{-2}/(1-\lambda^{-1}))$, and we see 
that the total Morse polynomial introduced in Theorem \ref{Theorem_strong_Morse_Dolbeault_version_groups} is $\lambda+1+\lambda^{-1}$. Since there are no $1+b$ factors this is exactly the Poincar\'e polynomial.

When considering the branch of the Hodge diamond for $p=1$, the local cohomology groups on the tangent cones are $\{z^{k}dz\}_{k=-1}^\infty$ and $\{\overline{z}^kd\overline{z}\wedge dz\}_{k=1}^{\infty}$. In particular this shows that the volume form is not in the local cohomology group for the minimal domain of $\overline{\partial}^*$. The Morse polynomials can be computed as earlier to see that the Poincar\'e polynomial is $1$.

In particular, the Dolbeault cohomology with this domain is not contained in the de Rham cohomology, and the opposite inclusion does not hold either.
\end{example}

We now use Theorem \ref{theorem_rigidity_preliminary} to prove the following rigidity result.

\begin{theorem}
\label{Theorem_instanton_correspondence_deRham_Dolbeault_non_intro}
Let $X$ be a resolved stratified pseudomanifold with a wedge K\"ahler metric and a K\"ahler Hamiltonian action. 
Then the Witten deformed de Rham complex in Theorem \ref{theorem_de_Rham_Morse_Bott} is isomorphic to a Hilbert sub-complex of the Witten deformed complex for $\prescript{\bullet}{}{\mathcal{P}}_{\min}$. 

In particular a global K\"ahler Hamiltonian action on $X$, the global de Rham and Dolbeault cohomology groups, in addition to being isomorphic, are rigid.
\end{theorem}

\begin{proof}
Given a global K\"ahler Hamiltonian action on $X$, the inclusion in Theorem \ref{theorem_rigidity_preliminary} holds for the Witten complexes since it holds for all the localized complexes at the fixed points. Thus the Witten deformed de Rham complex is a subcomplex of the Witten deformed complex $\prescript{\bullet}{}{\mathcal{P}}_{W}$.

Since the cohomology of the de Rham complex is the same as its Witten deformed one, it is rigid. Since the de Rham and Dolbeault cohomology groups are isomorphic on $X$, we have that the latter is rigid as well. 
\end{proof}

\textbf{The $L^2$ cohomology of the de Rham and Dolbeault complexes (with minimal domain) in the presence of K\"ahler Hamiltonian functions is rigid, that is the representations of the Hamiltonians are trivial on the global cohomology.}

In the examples we have studied so far, it seems that the Witten complex for the de Rham and rigid Dolbeault complexes are isomorphic, but so far we can only prove one inclusion rigorously.

\subsubsection{Poincar\'e Hodge polynomials}
\label{subsubsection_Poincare_Hodge_polynomial}

Here we follow \cite{wu1998equivariant}, and \cite{jayasinghe2024holomorphic} and refer to those articles for more details and motivations for the study of Poincar\'e Hodge polynomials.

\begin{definition}[Poincar\'e Hodge polynomials]
\label{equivariant_Poincare_Hodge_definition}
Given a Hilbert complex $\mathcal{P}$ with a $(p,q)$-bidegree
 define the equivariant \textit{\textbf{Poincar\'e Hodge polynomial}} to be 
\begin{equation}
    \chi_{y,b}(\mathcal{P},T_g):=\sum_{p=0}^n y^p \sum_{q=0}^n b^q \operatorname{Tr} T_g |_{\mathcal{H}^{p,q}(\mathcal{P})}
\end{equation}
\end{definition}

The following Theorem follows naturally from the results we have proven in this article.

\begin{theorem}
\label{theorem_Morse_polynomials_for_Poincare_Hodge_polynomials}
Given a K\"ahler Hamiltonian $h$ on a resolved stratified pseudomanifold $X$ with a flat Hermitian bundle $E$, given geometric endomorphisms $T$ we have that for each of the complexes $\mathcal{P}_{\min}=\mathcal{P}^d_{\min},\prescript{\bullet}{}{\mathcal{P}}_{\min}, \prescript{\bullet}{}{\mathcal{P}}^{\mu=0}_{\min}$
\begin{equation}
\label{equation_Poincare_polynomials_Morse_inequality}
    \sum_{F \in Z(A)} \chi_{y,b}(\mathcal{P}_{\min,B}(U(F)),T)=\chi_{y,b}(\mathcal{P}_{\min}(X),T)+(1+b) \sum_{q=0}^{n-1} R_q b^q
\end{equation}
are independent of $T$, in particular equal to the case where $T=Id$. For $\mathcal{P}^d_{\min}, \prescript{\bullet}{}{\mathcal{P}}^{\mu=0}_{\min}$, $R_q$ are non-negative integers.

\noindent The polynomials $\chi_{y,b}(\mathcal{P}_{\min}(X)):=\chi_{y,b}(\mathcal{P}_{\min}(X),Id)$ are equal for all three complexes.
\end{theorem}

\begin{proof}
    The results follow easily by taking the Witten complexes of the complexes, and applying Theorem \ref{Theorem_instanton_correspondence_deRham_Dolbeault_non_intro}.
\end{proof}

\begin{proposition}
In the same setting as Theorem \ref{theorem_Morse_polynomials_for_Poincare_Hodge_polynomials}
\begin{equation}
    \sum_{F \in Z(A)} \chi_{y,b}(\mathcal{P}_{\min,B}(U(F)),Id)  =\chi_{y,b}(X)\chi_{y,b}(\mathcal{P}_{\min}(X)) +(1+b)(1+y)\sum_{q}\sum_{p} y^p b^q R^{q,p}
\end{equation}    
where $R^{q,p} \in \mathbb{Z}, R^{p,q} \geq 0$.
\end{proposition}

\begin{proof}
This is a simple consequence of the symmetries of the Hodge diamond in the K\"ahler setting, which implies that $\chi_{y,b}(X)=\chi_{b,y}(X)=(yb)^n\chi_{y^{-1},b^{-1}}$ and similarly for the local cohomology for all choices of boundary conditions.
\end{proof}

The equivariant Hirzebruch invariant for a geometric endomorphism $T$ is given by 
\begin{equation}
    \chi_{y}(\mathcal{P}^d_{\min}(X),T):=\chi_{y,-1}(\mathcal{P}^d_{\min}(X),T),
\end{equation}
and we shall denote these by $\chi_y(X,T)$ for short,
which for the actions we study are rigid by the rigidity results of this section. The Euler characteristic is given by $y=-1$, and the signature invariant in the K\"ahler setting is given by $y=+1$. For $X$ of real dimension 4, the equivariant self-dual and anti-self dual complexes have indices given by $(\chi_{-1}(X,T)+\chi_1(X,T))/2$. 
In \cite{jayasinghe2023l2} the first author showed that these results continue to hold in the wedge setting, and formalizing results such as Self-duality under Serre duality (see Remark 6.8 of \cite{jayasinghe2024holomorphic}) which reduce to trivial observations in the smooth setting in terms of polynomials but are more nuanced in the singular setting. Indeed, symmetries of these invariants are related to symmetries of elliptic genera coming from modular invariance (see for instance \cite{hirzebruch1992manifolds}) and symmetries of Nekrasov partition functions (see Section 7.2.8 of \cite{jayasinghe2023l2}) arising in numerous studies of QFT. 
It is easy to see that many of the arguments there extend to the case of non-isolated critical point sets that we study, and we leave this to the interested reader.

\subsubsection{A vanishing theorem}
\label{subsubsection_vanishing_results}

In \cite{witten1984holomorphic}, Witten \textit{partly retrieves} a result of Carell and Liberman, that given a holomorphic vector field with only isolated zeros on a K\"ahler manifold, then the cohomology groups $H^{k,l}$ vanish for $k\neq l$. This is a simple corollary of the Holomorphic Morse inequalities in the smooth case, but more complicated in the singular setting (c.f. Section 7.2.9 of \cite{jayasinghe2023l2}).
Theorem 4.1 of \cite{wu1998equivariant} generalizes this the Poincar\'e Hodge polynomials to the case of non-isolated critical points, and is used to draw Corollary 4.5 of that paper, which shows that the cohomology groups $H^{k,l}$ vanish for $|k-l| > N$ where $N$ is the maximum of the complex dimension over the components of the critical point set. 
We introduce the following notion, which will allow us to formulate a general version of this result in a succinct form.

\begin{definition}[Hodge deviation (local and global)]
If the twisted cohomology of the de Rham complex $\mathcal{P}_{\min}$ (local or global) is supported in Hodge bi-degrees $(k,l)$, we say that $\max\{|k-l|\}$ over all such pairs $(k,l)$ is the \textbf{Hodge deviation} $\mathcal{HD}(\mathcal{P}^d_{\min})$ of the complex.
\end{definition}

Now we have the following vanishing theorem.

\begin{proposition}
Given a twisted de Rham complex on a stratified pseudomanifold with a K\"ahler wedge metric and a hamiltonian action as we study in this article,
the Hodge deviation of the space satisfies
\begin{multline}
    \mathcal{HD}(\mathcal{P}^d_{W}(X))= \max_{F \in Z(A)}\mathcal{HD} \Big(\mathcal{P}^d_{W,B(g)}(U(F)) \Big)\\ \leq \max_{F \in Z(A)}\mathcal{HD}(\mathcal{P}^d_{W}(F) )+\max_{F \in Z(A)}\mathcal{HD}(\mathcal{P}^d_{W,B(g)}(U_a)).
\end{multline}
Moreover, the Hodge deviation $\mathcal{HD}(\mathcal{P}^d_{W,B(g)}(U_a)$ of a fiber is independent of the choices of action/boundary conditions $\mathcal{B}(g)$.
\end{proposition}

We observe that $\mathcal{HD}(\mathcal{P}^d_{W}(F)) \leq ({\dim F})$, and the vanishing of local cohomology on fibers give many restrictions, and this can be used to recover the results in \cite{wu1998equivariant} in the smooth setting described at the beginning of this Subsection \ref{subsubsection_vanishing_results}.

\begin{proof}
It suffices to prove this for a circle action. The Witten deformed complex for the de Rham cohomology can be used to see that the Hodge deviation has to be equal to the maximum of the Hodge deviations of the localized complexes.

Any harmonic section in the local cohomology for the de Rham complex is locally a product of a harmonic section over $F$ and a harmonic section on a fiber, and we can see that the sums of the Hodge deviations of the factors give an upper bound as in the statement of the proposition.

Serre duality on the fibers given by the Hodge star operator shows that the Hodge deviation is the same for both the generalized Neumann and Dirichlet conditions for a fiber. Since this holds for any factor $C(Z_{a,j})$ of the fibers (see Definition \ref{definition_non_degeneracy_critical_point_sets_Kahler}), and the sums of the Hodge deviations on factors give the Hodge deviation on the product, we see that the Hodge deviation is independent of the choice of domain.
\end{proof}

\begin{remark}
In the smooth setting, the Hodge deviation of the normal fibers is $0$. Conditions on the cohomology of normal fibers can be used to get more precise vanishing results for the cohomology.
For instance the pure and odd vanishing conditions in the case of isolated fixed points, there is related work in Section 7.2.9 of \cite{jayasinghe2023l2}.
\end{remark}

\subsubsection{Bolt charges and Bott residues}
\label{subsection_bolt_charges}

The case of NUT charges arising in the study of supergravity, and the study of gravitational instantons was investigated in \cite{jayasinghe2024holomorphic}, highlighting how the ideas in Section 5 of \cite{bott1967vector} could be used to compute more general invariants using localization.

We now review the equivariant signature and NUT and Bolt charge conventions of
\cite[\S 7.2]{jayasinghe2024holomorphic}, extending the local contributions to
positive-dimensional fixed components (also known as bolts).

Let $F$ be a connected fixed component, and let
$\tau_F(\theta)$ denote its local contribution to the equivariant
signature (i.e. $\chi_{y=+1,b=-1}(\mathcal{P}(U(F),T_{2\theta})$, with circle parameter $g=e^{2i\theta} \in S^1$. In real dimension four, when this contribution has the expansion
\[
\tau_F(\theta)
=
\frac{N_F}{\theta^2}+T_F+O(\theta^2),
\]
we call $N_F$ its charge and $T_F$ its constant signature contribution.
For an isolated fixed point this is the NUT charge, for a
two-dimensional fixed component we call it the bolt charge.
This definition applies also when the fixed component meets the
singular locus, using the local equivariant signature contribution of the chosen signature complex.

Consider the singular surface in Example \ref{Example_singular_toric_variety_domains}, and in Subsection \ref{subsection_non_isolated_singularity} with the toric action considered there. The equivariant signatures at the 3 generic fixed points are 
\begin{equation}
    A_1=\frac{
1+\lambda^3\mu+\lambda^2\mu^2+\lambda\mu^3
+y(\lambda+\mu)^2(\lambda^2+\mu^2)
+y^2(\lambda^4\mu^4+\lambda^3\mu+\lambda^2\mu^2+\lambda\mu^3)
}{
(1-\lambda^4)(1-\mu^4)
}
\end{equation}
at $a_1$, and 
\begin{equation}
\frac{(1+y\lambda^{-4})(1+y\mu/\lambda)}
     {(1-\lambda^{-4})(1-\mu/\lambda)},
\qquad
\frac{(1+y\lambda/\mu)(1+y\mu^{-4})}
     {(1-\lambda/\mu)(1-\mu^{-4})}
\end{equation}
at $a_2$ and $a_3$ respectively.
Writing
\[
\lambda^4=e^{2ip\theta},\qquad
\mu^4=e^{2iq\theta},
\qquad pq(p-q)\ne0,
\]
with the corresponding continuous choices
$\lambda=e^{ip\theta/2}$ and $\mu=e^{iq\theta/2}$.
The local signature contributions at $a_2$ and $a_3$ are
\begin{equation}
    \tau_{a_2}(\theta)=\cot(p\theta)\cot (({q-p})\theta/{4}),
\qquad
\tau_{a_3}(\theta)=\cot(q\theta)\cot (({p-q})\theta/{4}).
\end{equation}
Expanding these expressions and $A_1$ gives
\[
\begin{array}{c|c|c}
F & N_F & T_F\\ \hline
a_1
&
-\dfrac{4}{pq}
&
\dfrac{1}{12}\left(\dfrac pq+\dfrac qp\right)-\dfrac12
\\[6pt]
a_2
&
\dfrac{4}{p(q-p)}
&
-\dfrac13\left(\dfrac{4p}{q-p}+\dfrac{q-p}{4p}\right)
\\[6pt]
a_3
&
\dfrac{4}{q(p-q)}
&
-\dfrac13\left(\dfrac{4q}{p-q}+\dfrac{p-q}{4q}\right).
\end{array}
\]
In particular, the charge at $a_1$ is $-4/(pq)$, correcting
the coefficient $-3/(pq)$ in \cite[\S 7.2]{jayasinghe2024holomorphic}.
These values satisfy
\[
\sum_{j=1}^3N_{a_j}=0,
\qquad
\sum_{j=1}^3T_{a_j}=1.
\]

For the subgroup $\mu=1$, the fixed set consists of the bolt $B_0=\{X=Z=0\}$
and the isolated point $a_2$.
It is easy to reverse engineer and compute the equivariant $\chi_1$ invariant at the Bolt, which can be used to see that $N_{B_0}=1/4$ $T_{B_0}=-5/12.$
The remaining isolated contribution is
\[
\tau_{a_2}(\theta)
=
-\cot(4\theta)\cot\theta
=
-\frac{1}{4\theta^2}+\frac{17}{12}+O(\theta^2),
\]
and their sum is identically $1$.
The vanishing of the higher order terms is due to rigidity that we have proven earlier in this section.

In the smooth setting, the NUT and Bolt charges are computed in \cite{gibbons1979classification} using the equivariant signature formula and the influence of Atiyah in that paper is clear, not just by the explicit acknowledgement by name. The principle behind the computations is described by Bott in the last section of \cite{bott1967vector}. 
The Bott residue theorem shows that given a global integral of an equivariant polynomial $\Phi$ of the Chern classes of a smooth manifold can be computed in terms of local quantities near fixed points.
Corollary 3 of \cite{bott1967vector} states that if the weight of $\Phi$ is weight less than $m$ (half the dimension of the manifold), the global integral vanishes on a smooth manifold of dimension $m$.
The integral  global polynomial $1$ corresponds precisely to the NUT and Bolt charges, hence why the sum of all NUT and Bolt charges over all fixed points of a space vanish. 
In fact it can be shown that for a general almost complex manifold of dimension $2m$, the equivariant signature has a power series expansion in powers of $\theta$ with negative terms that all add up to $0$, that have topological significance.

In the case of a 2 dimensional Bolt in 4 dimensions, the Bolt charge corresponds to the self intersection number of the Bolt (see \cite{gibbons1979classification}). On singular spaces this was generalized in the non-equivariant case by the index formulas in \cite{atiyah2013curvature}, also studied in \cite{Albin_2016_Index} (see Theorem 6.2 of that article). Our work here is a step towards understanding this equivariantly, and the non-equivariant work described above is related to Conjecture \ref{Conjecture_brave_2}.

The entropy functional in equation (6.9) of \cite{gibbons1979classification} is in terms of the NUT charges at fixed points (of both NUTs and Bolts), and in terms of the area of Bolts. In the case where the Bolts are arising from the actions at walls of toric actions, the Duistermaat-Heckman theorem for the Bolt with a refined action on the Bolt can be used to compute the equivariant area, and in particular the area of the Bolt. This can be used to compute a generalized entropy functional on singular spaces. There are other, more sophisticated equivariant quantities computed in more recent studies of supergravity, including on spaces with orbifold singularities in \cite{genolini_sparks2019localization,Equivariant_loc_Sparks_2024}, where NUT and Bolt computations are done. This requires a better understanding of equivariant localization, in a manner that can capture the behaviour at singularities, including choices of domains and the geometry and topology that they see at singularities.
Remark 7.12 of \cite{jayasinghe2024holomorphic} can be used to see that the NUT and Bolt charges can be computed using the equivariant signature for any choice of domain, but other equvariant quantities that can be read off of power series expansions of the equivariant index formulae may depend on these choices in general.

\subsection{The Kervaire semi-characteristic}
\label{subsection_Kervaire}

In the smooth setting, Atiyah proved a vanishing theorem for the Kervaire semi-characteristic on a $4q+1$ dimensional space admitting two linearly independent non-vanishing vector fields using a perturbation of a Dirac type operator.
This was generalized other mod 2 indices and to a counting formula by Zhang in \cite{zhangkervaire,zhangkervaire99} (see also \cite[\S 7]{Zhanglectures}). These methods are also related to the work in \cite{zhang2000n,chen2018poincar} where the Poincar\'e Hopf index theorem is derived using spectral flow.

\begin{definition}[Kervaire semi-characteristic]
\label{Definition_Kervaire_SC}
Given a stratified pseudomanifold $X$ of dimension $4q+1$, the \textbf{Kervaire semi-characteristic} is 
\begin{equation}
    k_W(X):= \sum_{i=0}^{2q} \dim \mathcal{H}^{2i}(\mathcal{P}_W^d(X)) (\text{mod } 2\mathbb{Z})
\end{equation}
where we denote the de Rham cohomology by $\mathcal{H}^*(\mathcal{P}_W^d(X))$ consistent with the notation in this section.    
\end{definition}
In the smooth case the Hodge Dirac operator is essentially self-adjoint and we can drop $W$. However the above definition extends to resolved stratified pseudomanifolds that we study.

Consider a resolved stratified pseudomanifold $X$ of dimension $4q+1$ where there is one non-vanishing vector field $V_1$ and another $V_2$ that is transverse to it, with potential vanishing loci around circles fixed by the flow of $V_1$. In this context, the Euler characteristic vanishes since $V_1$ has no fixed points, and $k_W$ can be thought of as a secondary invariant.

Note that there are such examples where the second vector field is the gradient flow of a stratified Morse-Bott function $h$ which is \textit{perfect} (i.e. the Morse polynomial and the Poincar\'e polynomial are equal), for instance taking a stratified Morse function on a stratified pseudomanifold $X'$ of dimension $4q$ and extending it in the obvious manner to $X=X' \times S^1$, where the critical point sets will be circles. 
Some examples of perfect stratified Morse functions are the height functions $\cos(\phi)$ of spherical suspensions (see Example \ref{example_spherical_suspension_circle_4pi}). 
In this setting, the Witten deformed de Rham complex in \cite{jayasinghe2024witten} can be used, even without the improvements in this article since the critical point sets are products of circles with cones.

More generally one can consider perfect stratified Morse-Bott functions $h$ on a space $X$ where the critical point sets are circles, and carry circle actions. 
We present a formula for $k(X)$ which generalizes that of Zhang to singular spaces, albeit in the context of stratified Morse-Bott functions arising as above. 

\begin{theorem}
Given a stratified pseudomanifold $X$ of dimension $4q+1$ and a perfect stratified Morse-Bott function $h$ where $Z(|dh|^2)$ consists of circles, then 
\begin{equation}
    k_W(X)= \sum_{F \in Z(|dh|^2)} \sum_{i=0}^{2q} \dim \mathcal{H}^{2i} \Big(\mathcal{P}_{W,B}^d(U(F)) \Big) (\text{mod } 2\mathbb{Z})
\end{equation}
where the de Rham complex on $X$ and where $\mathcal{P}_W^d$ is a de Rham complex on $X$ with a choice of mezzo-perversity $W$. 
\end{theorem}

\begin{proof}
This is a corollary of the Witten subcomplexes introduced Theorem \ref{theorem_de_Rham_Morse_Bott}), where the local cohomology groups at the fixed points are $\mathcal{H}^i \Big(\mathcal{P}_{W,B}^d(U(F)) \Big)$. Since the function $h$ is perfect, the global cohomology is isomorphic to the vector spaces in the Witten subcomplex, and the result follows from using this to replace the cohomology groups in Definition \ref{Definition_Kervaire_SC}.    
\end{proof}

That the local cohomology groups are more general can be seen, for instance from the work in \cite[\S 6]{jayasinghe2024witten}, where the cohomologies of the fibers over the circles appear in the formulas. It is also easy to check that this result generalizes that of Zhang (for instance in \cite[\S 7]{Zhanglectures}, and we leave this to the reader).

\textbf{This shows how Zhang's formulas for mod 2 indices can be reformulated in terms of local cohomology groups, in ways suitable for generalizations.}

\section{Examples and guide for computations}
\label{section_examples}

Since we did not find any worked out examples of the holomorphic Morse polynomials in the literature in the smooth setting, we will work out a smooth example in detail in Subsection \ref{subsection_smooth_toric_action}, comparing and contrasting with the formulas of Atiyah-Segal-Singer. 

\subsection{Toric action on $\mathbb{CP}^2$}
\label{subsection_smooth_toric_action}

We study the following example, relating it to a study by Duistermaat \cite[\S 13.4]{duistermaat2013heat}. 
We consider the action $(\lambda,\mu)[ Z_0:  Z_1: Z_2]=[\lambda Z_0:  \mu Z_1: Z_2]$ on $\mathbb{CP}^2$, which for generic values of $(\lambda,\mu)=(e^{i\theta_1},e^{i\theta_2}) \in (\mathbb{C}^*)^2$ has three fixed points on the origins of the standard affine charts. 
For this algebraic two torus action, the holomorphic Lefschetz fixed point theorem yields
\begin{equation}
    \frac{1}{(1-\lambda)(1-\mu)}+\frac{1}{(1-\lambda/\mu)(1-1/\mu)}+\frac{1}{(1-1/\lambda)(1-\mu/\lambda)}=1
\end{equation}
If one picks an action chamber for the action, we can upgrade this to a Morse polynomial. However at the walls of the action chambers, the denominators appearing in the local Lefschetz numbers go to zero. Consider the case where $\mu \rightarrow 1$, for which the first two terms have vanishing denominators. This is because at $\mu=1$, the action fixes the complex line $Z_0=0$, and the isolated fixed points on $Z_0=0$ for generic actions are connected by a fixed complex sphere.

We see that sum of these first two terms can be simplified as 
\begin{align}
\label{equation_pole_zero_cancellations_smooth}
    &=\frac{1}{(1-\lambda)(1-\mu)}+\frac{-\mu}{(1-\lambda/\mu)(1-\mu)}\\
    &=\frac{1}{(1-\mu)} \bigg(\frac{1-\lambda/\mu-\mu+\lambda\mu}{(1-\lambda)(1-\lambda/\mu)} \bigg)
    =\frac{1}{(1-\mu)} \bigg(\frac{1-\mu +\lambda(\mu-1/\mu)}{(1-\lambda)(1-\lambda/\mu)} \bigg)\\
    &=\frac{1}{(1-\mu)} \bigg(\frac{(1-\mu) +\lambda/\mu(\mu-1)(\mu+1)}{(1-\lambda)(1-\lambda/\mu)} \bigg)=\frac{1 -\lambda/\mu(\mu+1)}{(1-\lambda)(1-\lambda/\mu)}
\end{align}
where we can take the limit as $\mu$ goes to $1$ to get $(1-2\lambda)/(1-\lambda)^2$. \textbf{This form of \textit{reverse engineering} can be used to compute holomorphic Lefschetz numbers for non-isolated fixed point sets corresponding to the circle actions at walls of torus actions which generically have isolated fixed point sets, even in the singular setting.}

\textbf{In fact even the Morse polynomial can be obtained using this}, with limited information of the cohomology of the non-isolated critical point set. We consider the case where $\lambda=e^{i\theta}$ where $\theta>0$. Then the local Lefschetz number and the local Morse polynomial at the isolated fixed point can be easily checked to be the \textbf{unique convergent taylor series expansions for $|\lambda|<1$} of
\begin{equation}
    \frac{1}{(1-1/\lambda)^2}=\frac{\lambda^2}{(1-\lambda)^2}, \quad \text{ and } \quad b^2\frac{\lambda^2}{(1-\lambda)^2}
\end{equation}
respectively.    
The global Lefschetz number (the global Poincar\'e polynomial) is
\begin{equation}
    1=P(\lambda,\mu=1)(b)=\frac{\lambda^2}{(1-\lambda)^2}+\frac{-2\lambda+1}{(1-\lambda)^2}.
\end{equation}
The second term on the right is the index defect which together with the Lefschetz number at the isolated fixed point adds up to $1$.

The contribution from the degree $0$ cohomology of the fixed complex line $Z_0=0$, with values in the local Dolbeault cohomology on the normal fibers is ${1}/{(1-\lambda)}$. This is because the holomorphic functions have the basis given by the Hilbert space completion of $\mathbb{C}[[Z_0]]$, which is the local cohomology in degree $0$ of the fundamental neighbourhood. 

Since we know the contributions from degrees $0,2$ for the local Morse polynomials, and since we know the Poincar\'e polynomial is $1$, the Morse polynomial completely determines the contribution of degree $1$, and it is
\begin{equation}
    b \Bigg(\frac{\lambda}{(1-\lambda)}+\frac{\lambda^2}{(1-\lambda)^2} \Bigg)=b \Bigg(\frac{\lambda(1-\lambda)+\lambda^2}{(1-\lambda)^2} \Bigg)=b \frac{\lambda}{(1-\lambda)^2}.
\end{equation}
The Morse polynomial is 
\begin{equation}
   M_C(\lambda=e^{i\theta},\mu=1)(b)=b^2\frac{\lambda^2}{(1-\lambda)^2}+b \Bigg(\frac{\lambda^2}{(1-\lambda)^2}+\frac{\lambda}{(1-\lambda)}\Bigg)+b^0\frac{1}{(1-\lambda)}
\end{equation}
where $C$ is the action chamber where $\theta>0$ for $\lambda=e^{i\theta}$. This is actually a wall on the action chamber for the well known $SU(2)$ action on $\mathbb{CP}^2$ for general $\lambda,\mu$.
The dual Morse polynomial corresponding to $-C$ (i.e. $\theta<0$) is
\begin{equation}
   M_{-C}(\lambda=e^{-i\theta},\mu=1)(b)=\frac{1}{(1-1/\lambda)^2}+b \Bigg(\frac{1}{(1-1/\lambda)^2}+\frac{1}{(1-1/\lambda)}\Bigg)+b^2\frac{1/\lambda}{(1-1/\lambda)}.
\end{equation}
These expressions can be easily used to see that the error polynomial has a factor of $(1+b)$, and we can recover the Lefschetz numbers by setting $b=-1$ in either of these Morse polynomials.

\textbf{We now compute the local cohomology group at the fixed line $\{Z_0=0\}$}.
The local cohomology group with the $\overline{\partial}$-Neumann boundary conditions on the boundary of the fundamental neighbourhood can be expressed in terms of their harmonic representatives as follows.
On the chart where $Z_2=1$, the cohomology for fixed eigenvalue $k$ of the operator $\sqrt{-1}L_{V}$ (where $V$ is the Hamiltonian vector field) is given by
\begin{equation}
    \mathcal{H}^{q}(\mathbb{CP}^1;\mathcal{H}^*(\mathbb{D}^2))=\begin{cases}
      (Z_0/Z_2)^k & \text{if  $q=0$} \\
      \frac{\overline{w}^jd\overline{w}}{(1+|w|^2)^{k+1}} \otimes  \frac{Z_0^k}{Z_2^{k+1}} \text{   for $j=0,..., k-1$} &\text{ if $q=1$}
    \end{cases}
\end{equation}
where $w=Z_1/Z_2$, and corresponds to the holomorphic functions on each fiber of the normal disc bundle.

On the chart where $Z_1=1$, we have
\begin{equation}
    \mathcal{H}^{q}(\mathbb{CP}^1;\mathcal{H}^*(\mathbb{D}^2))=\begin{cases}
      (Z_0/Z_1)^k & \text{if  $q=0$} \\
      -\frac{\overline{v}^{k-1-j}d\overline{v}}{(1+|v|^2)^{k+1}} \otimes  \frac{Z_0^k}{Z_1^{k+1}} \text{   for $j=0,..., k-1$} &\text{ if $q=1$}
    \end{cases}
\end{equation}
where $v=Z_2/Z_1$.

\begin{remark}
    In \cite{wu1998equivariant}, these cohomology groups are expressed as the cohomology of the bundle of small eigenvalue eigensections which are well understood in the smooth setting (see equation (2.1) of \textit{loc. cit.}.

    The Lichnerowicz-Schr\"odinger identity in equation \eqref{equation_operator_identities_3} (see Remark \ref{Remark_Lichnerowicz_Schrodinger}) has a term $\varepsilon iL_V=\varepsilon \mu$ which grows as we take large values of $\mu=k$ in this example, and can be thought of as a curvature term. This explains the denominators $(1+|v|^2)^{k+1}$ arising in the cohomology for fixed $k$, which for the Fubini-Study metric on the fixed spheres correspond to curvature of higher and higher tensor powers of the line bundle associated to the normal bundle.

    This is related to similar Lichnerowicz-Schr\"odinger identities arising in the study of asymptotic holomorphic Morse inequalities pioneered by Demailly in \cite{demailly1991holomorphic} where the \textit{deformation} is due to taking tensor powers of line bundles. 
\end{remark}

It is easy to see that these sections glue on the two charts on the fixed $\mathbb{CP}^1$, and are in fact the cohomology of a line bundle corresponding to the sheaf $(Z_0/Z_2)^k$ and $(Z_0/Z_1)^k$ on the two charts (with the corresponding transition map).
Then it is easy to see that the trace over the local cohomology in degree $1$ is given by 
\begin{equation}
    \sum_{k \geq 1} k\lambda^k=\lambda \sum k\lambda^{k-1}=\lambda \frac{d}{d\lambda}\bigg(\frac{1}{1-\lambda} \bigg)=\frac{\lambda}{(1-\lambda)^2}
\end{equation}
where we appeal to the analyticity of the formulas to justify the calculus for the simplication we have used.

We briefly review for the benefit of the reader the Atiyah-Segal formula from Duistermaat \cite[\S 13.4]{duistermaat2013heat} for the non-isolated fixed point set $F_1$. Denoting the Todd class of $F_1$ and the Chern class of the normal bundle of $F_1$ in the space as $t,c \in H^2(F_1;\mathbb{C})$, we get the contribution to be 
\begin{equation}
\label{equation_duistermaat_example_contribution}
    \int_{F_1} \text{area} \Big( \frac{1+t}{1-\lambda e^{-c}} \Big)
\end{equation}
Using that $c^2=0$ on $F_1$ we get
\begin{equation}
    (1-\lambda e^{-c})^{-1}=(1-\lambda+\lambda c)^{-1}=(1-\lambda)^{-1}(1 +  (1-\lambda)^{-1} {\lambda c} )^{-1}
\end{equation}
where again using that $c^2=0$, we see that this is equal to
\begin{equation}
    (1-\lambda)^{-1}\bigg(1-\frac{\lambda c}{(1-\lambda)} \bigg)=(1-\lambda)^{-2}(1-\lambda-\lambda c )
\end{equation}
from which we see that the $2$ form part of the integrand of \eqref{equation_duistermaat_example_contribution} is 
\begin{equation}
     \text{area} \bigg( \frac{(1+t)(1-\lambda-\lambda c)}{(1-\lambda)^2}\bigg)= \frac{(-(t+c)\lambda+t)}{(1-\lambda)^2}  
\end{equation}
and since the integrals of $t$ and $c$ over $F_1$ are both $1$, we see that the expression \eqref{equation_duistermaat_example_contribution} simplifies to
${(-2\lambda+1)}/{(1-\lambda)^2}$.
Duistermaat infers the fact that $T+C=2$ and $T=1$ by reverse engineering, using that fact that $\chi(\lambda)=1$.

\subsection{Singular example}
\label{subsection_non_isolated_singularity}

We consider the example of the singular toric variety in Example \ref{Example_singular_toric_variety_domains}, and study various invariants that we have defined in this article explicitly. In section 7.3.6 of \cite{jayasinghe2023l2}, the algebraic Lefschetz numbers and the $L^2$ holomorphic Lefschetz numbers for the minimal domain were studied, in particular showing how to compute the equivariant $\chi_y$ invariant at the singular point of depth 2 for generic actions.
In section 6.2 of \cite{jayasinghe2024holomorphic}, holomorphic Morse inequalities generalizing the algebraic equivariant Riemann-Roch results of \cite{baum1979lefschetz,Baumformula81} were conjectured. The ideas in this article can be used to conjecture versions for non-isolated critical point sets, and in amalgamation with the ideas in \cite{jayasinghe2026spectralasymmetrysupersymmetryequivariant} can be used to understand even non-equivariant extensions of Riemann-Roch defects, corresponding to \cite{baumfultonmacKtheoryRiemannRoch,baumfultonmacpheresonRiemannRoch}.
We leave this to the inquisitive reader, and limit our scope for brevity.

The space $\widehat{V}$ in Example \ref{Example_singular_toric_variety_domains}, that is the vanishing locus of $Z^4-X^3Y=0$ in $\mathbb{CP}^3$ with coordinates $[U:X:Y:Z]$, admits the algebraic torus action $(\lambda,\mu)\cdot [U:X:Y:Z]=[U:\lambda^4X:\mu^4Y:\lambda^3\mu Z]$.

We compute the equivariant index for the $\mathcal{P}=(\Omega^{0,q}(X), \mathcal{D}_{\min,\min}(\overline{\partial}),\overline{\partial})$.
Generically there are three isolated fixed points at $a_1=[1:0:0:0]$, $a_2=[0:1:0:0]$ and $a_3=[0:0:1:0]$, where $a_2$ is a smooth fixed point while $a_1$ and $a_3$ are singular. The singular set is the locus $[U:0:Y:0]$.

This space is the compactification of the affine cone over the singular algebraic curve $Z^4-X^3Y=0$ in $\mathbb{CP}^2$ and $a_1$ is a stratum of depth $2$ within the singular set.

The equivariant Riemann-Roch numbers in $L^2$ cohomology for $p=0$, corresponding to the minimal domain, were computed in section 7.3.6 of \cite{jayasinghe2023l2} and we refer the reader to that article for more details on local cohomology groups at the isolated fixed points that we discuss here.
The smooth point gives the contribution $1/[(1-1/\lambda^4)(1-\mu/\lambda)]$. The integral completion of the local ring at the singular fixed point $a_3$ is $Q=\mathbb{C}[X/Z,U/Y]$.
The local cohomology in degree $0$ at this point is $\overline{Q}$ and the local Lefschetz number is $1/[(1-\lambda/\mu)(1-1/\mu^4)]$.
On the chart where $U=1$, we consider coordinates $x=X/U$, $y=Y/U$ and $z=Z/U$. The local cohomology group at the singular fixed point $a_1$ is $\overline{R_1}$, where $R_1=R[1,z,xy/z,x^2y/z^2]$, where $R=\mathbb{C}[x,y]$, which can be used to compute the local Lefschetz number $[1+(\lambda^3\mu)+(\lambda^2\mu^2)+(\lambda\mu^3)]/[(1-\lambda^4)(1-\mu^4)]$. Summing up, the global $L^2$ Lefschetz Riemann-Roch number is
\begin{equation}
\label{equation_Lefschtez_number_some_blah_blah}
    \frac{1+(\lambda^3\mu)+(\lambda^2\mu^2)+(\lambda\mu^3)}{(1-\lambda^4)(1-\mu^4)}+\frac{1}{(1-1/\lambda^4)(1-\mu/\lambda)}+\frac{1}{(1-\lambda/\mu)(1-1/\mu^4)}=1.
\end{equation}

The Morse polynomial for the chamber $C$, given by $0< Arg(\lambda) < Arg(\mu)$ for $\lambda, \mu$ in the first quadrant can be checked to be 
\begin{equation}
    \frac{1+(\lambda^3\mu)+(\lambda^2\mu^2)+(\lambda\mu^3)}{(1-\lambda^4)(1-\mu^4)}+b\frac{\lambda^4}{(1-\lambda^4)(1-\mu/\lambda)}+b^2\frac{\mu^5/\lambda}{(1-\mu/\lambda)(1-\mu^4)}=1+(1+b) \sum_{q=0}^1 Q_C^q(\lambda,\mu)
\end{equation}
We consider the second domain $\mathcal{D}_{W_{Z/X}}(P)$ in Example \ref{Example_singular_toric_variety_domains}, and the presentation there shows that this domain has one additional meromorphic section $Z/X$, on which the group acts by the endomorphism $\mu/\lambda$.

The local cohomology groups for the $\overline{\partial}$ operator at the critical points $a_1, a_3$ corresponding to this domain can be presented as the Hilbert space completion of the modules
\begin{equation}
    R_1\oplus \mathbb{C}[y](z/x]) \quad \mathbb{C}[X/Z,U/Y]\oplus\mathbb{C}[U/Y](Z/X)
\end{equation}
at the points $a_1,a_3$ respectively,
where $R_1=R(1,z,xy/z,x^2y/z^2)$ for $R=\mathbb{C}[x,y]$, as in the case of the minimal domain. Thus the local Lefschetz numbers at these points are
\begin{equation}
     \frac{1+(\lambda^3\mu)+(\lambda^2\mu^2)+(\lambda\mu^3)}{(1-\lambda^4)(1-\mu^4)} +\frac{\mu/\lambda}{(1-\mu^4)},\quad\frac{1}{(1-\lambda/\mu)(1-1/\mu^4)}+\frac{\mu/\lambda}{(1-1/\mu^{4})}
\end{equation}
respectively, and the local Morse polynomials are
\begin{equation}
     \frac{1+(\lambda^3\mu)+(\lambda^2\mu^2)+(\lambda\mu^3)}{(1-\lambda^4)(1-\mu^4)} +\frac{\mu/\lambda}{(1-\mu^4)}, \quad b^2\bigg(\frac{\mu^5/\lambda}{(1-\mu/\lambda)(1-\mu^4)}+\frac{\mu^5/\lambda}{(1-\mu^{4})}\bigg)
\end{equation}
Since the third point is smooth, we have the same Lefschetz numbers and Morse polynomials. We see that the $L^2$ Lefschetz Riemann-Roch number for this domain is 
\begin{equation}
    1+\frac{\mu/\lambda}{(1-\mu^4)}+\frac{\mu/\lambda}{(1-1/\mu^{4})}=1+\mu/\lambda
\end{equation}
where we have used \eqref{equation_Lefschtez_number_some_blah_blah} to simplify the other terms in the local Lefschetz numbers to give the number $1$ on the left hand side. The local and global Morse polynomials are easy enough to work out and we leave it to the reader.

We now consider non-isolated critical point sets.
The circle actions for the non-generic actions at various walls of action chambers are of different types.
\begin{enumerate}
    \item $\mu = \lambda$ fixes an algebraic curve given by the intersection of $U=0$ with $\widehat{V}$ in $\mathbb{CP}^3$ that intersects the singular curve $X=Z=0$ at a single point away from the singularity of depth $2$. 
    \item $\mu = 1$ fixes the entire singular set $X=Z=0$. 
    \item $\lambda = 1$ fixes an algebraic curve given by the intersection of $Y=Z=0$ with $\widehat{V}$ in $\mathbb{CP}^3$ that intersects the singular curve $X=Z=0$ at a single point away from the singularity of depth $2$. 
\end{enumerate}
Each of these are critical point sets of type II, and do not satisfy the normally non-singular condition, for instance used in \cite{banagl2025gsignaturetheoremwittspaces}.

In each of these cases, \textbf{the reverse engineering techniques introduced in the subsection above can be used to obtain the Morse polynomials}. 
We demonstrate this for the 2nd case, where $\mu=1$.

We first consider the minimal domain, and use the same method used to derive the local Lefschetz numbers for the non-isolated fixed point set in equation \ref{equation_pole_zero_cancellations_smooth}, for the first and last terms in equation \ref{equation_Lefschtez_number_some_blah_blah} to get 
\begin{align}
    &=\frac{1+(\lambda^3\mu)+(\lambda^2\mu^2)+(\lambda\mu^3)}{(1-\lambda^4)(1-\mu^4)}+\frac{-\mu^4}{(1-\mu^4)(1-\lambda/\mu)}\\&=\frac{(1+\lambda^3\mu+\lambda^2\mu^2+\lambda\mu^3-\lambda/\mu-\lambda^4-\lambda^3\mu-\lambda^2\mu^2-\mu^4+\lambda^4\mu^4)}{(1-\lambda^4)(1-\mu^4)(1-\lambda/\mu)}\\&=\frac{(1-\mu^4)(1-\lambda^4-\lambda/\mu)}{(1-\lambda^4)(1-\mu^4)(1-\lambda/\mu)}=\frac{(1-\lambda^4-\lambda/\mu)}{(1-\lambda^4)(1-\lambda/\mu)}
\end{align}
which for $\mu=1$ gives
\begin{equation}
\label{equation_blah_what}
    \frac{1-\lambda^4-\lambda}{(1-\lambda^4)(1-\lambda)}.
\end{equation}

The local cohomology in degree $0$ at the critical point set is given by sections which on the chart $Y=1$ corresponds to
\begin{equation}
\label{equation_singular_cohomology_11}
    =\begin{cases}
      (X/Y)^k & \text{if  $q=0$} \\
      \text{span}[1,Z/Y,(X/Z), (X/Z)^2] \frac{\overline{w}^jd\overline{w}}{(1+|w|^2)^{k+1}} \otimes  \frac{X^k}{Y^{k+1}} \text{   for $j=0,..., k-1$} &\text{ if $q=1$}
    \end{cases}
\end{equation}
where $w=U/Y$.
Here we note that the sections $(X/Z)d\overline{w}, (X/Z)^2d\overline{w}$ are $L^2$ bounded. 

On the chart $U=1$, we have
\begin{equation}
\label{equation_singular_cohomology_12}
    =\begin{cases}
      (X/U)^k & \text{if  $q=0$} \\
      \text{span}[1, Z/U,(X/Z), (X/Z)^2]\frac{-\overline{v}^{k-1-j}d\overline{v}}{(1+|v|^2)^{k+1}} \otimes  \frac{X^k}{U^{k+1}} \text{   for $j=0,..., k-1$} &\text{ if $q=1$}
    \end{cases}
\end{equation}
where $v=Y/U$. Recall that the integral completion of the local ring at the depth two fixed point has a presentation given by $R[1,Z/U, (XY)/(ZU), (X^2Y)/(Z^2U)]$, where $R=\mathbb{C}[X/U,Y/U]$. One can check that $(X/Z)d\overline{w},(X/Z)^2d\overline{w}$ are $L^2$ bounded.

It is easy to check that the trace of the induced map on the cohomology in degree $0$ is given by ${1}/{(1-\lambda^4)}$. 
The trace over the cohomology in degree $1$ is given by
\begin{equation}
    \frac{\lambda^4}{(1-\lambda)(1-\lambda^4)}=\frac{\lambda^4(1+\lambda+\lambda^2+\lambda^4)}{(1-\lambda^4)^2}.
\end{equation}
This can be checked by taking the trace over the sections with eigenvalue of $\sqrt{-1}L_V$ equal to $k$, which can be checked to be 
\begin{equation}
    k(1+\lambda^3+\lambda^2+\lambda)\lambda^{4k}
\end{equation}
for positive integers $k$. Summing over all $k$, we get 
\begin{equation}
    (1+\lambda+\lambda^2+\lambda^4)\lambda^4 \sum_{k=1}^{\infty}  k \lambda^{4(k-1)} =\frac{\lambda^4}{(1-\lambda)(1-\lambda^4)}
\end{equation}
where we use the fact that $(1+\lambda+\lambda^2+\lambda^4)=(1-\lambda^4)/(1-\lambda)$, and that for $x=\lambda^4$ we have
\begin{equation}
    \sum_{k=1}^{\infty}   kx^{k-1}=\sum_{k=0}^{\infty} \frac{dx^k}{dx}  =\frac{d}{dx}\frac{1}{(1-x)}=\frac{1}{(1-x)^2}.
\end{equation}
Summing these two contributions, with a minus sign for the contributions from degree $1$, we easily recover the expression in equation \ref{equation_blah_what}. 
The global Morse polynomial can also be checked to be 
\begin{equation}
    \frac{1-\lambda+b\lambda^4}{(1-\lambda^4)(1-\lambda)}+\frac{b^2\lambda^5}{(1-\lambda^4)(1-\lambda)}.
\end{equation}
For the case of the domain $\mathcal{D}_{W_{Z/X}}$, the local cohomology group has a Schauder basis that is the union of the basis in equations \eqref{equation_singular_cohomology_11} and \eqref{equation_singular_cohomology_12}, together with the section
$Z/X$ in degree $0$. Thus it is easy to see that the global Morse polynomial is 
\begin{equation}
    \frac{1-\lambda+b\lambda^4}{(1-\lambda^4)(1-\lambda)}+\frac{b^2\lambda^5}{(1-\lambda^4)(1-\lambda)}+1/\lambda.
\end{equation}

\section{Appendix A: Spectral theory for Schr\"odinger operators}
\label{Appendix_spectral_theory}

Here we present the spectral theory for Schr\"odinger operators corresponding to the Hamiltonians that arise in Witten deformation, for fundamental neighbourhoods of the type we consider.

\subsection{Discreteness of spectrum on infinite neighbourhoods}
\label{subsection_discreteness_of_spectrum}

In this subsection we establish that the spectrum of model Schr\"odinger operators on infinite fundamental neighbourhoods $U_{\infty}(F)$ is discrete, equivalently that they have compact resolvent. This follows from standard arguments, modified for the stratified setting, and its utility in the paper is discussed in Subsection \ref{subsubection_Fredholm_spectral}.

\subsubsection{Relative boundedness}

We begin by reviewing the notions of relative boundedness and relative compactness of operators. We refer the reader to \cite[\S 13.1]{hislop2012introduction} for the following definition and the statement and proof of the Kato-Rellich theorem (Theorem \ref{Theorem_Kato_Rellich}).
\begin{definition}[Relatively bounded operators]
    An operator $B$ is called \textbf{$A$-bounded} if $\mathcal{D}(A) \subseteq \mathcal{D}(B)$. If $res(A) \neq \emptyset$ and $B$ is $A$ bounded, then there exist non-negative constants $a,b$ such that 
    \begin{equation}
        \label{definition_relative_bounded_ops}
        ||Bu|| \leq a||Au||+b||u||
    \end{equation}
    for all $u \in \mathcal{D}(A)$.
    The smallest non-negative constant $a$ such that \eqref{definition_relative_bounded_ops} holds for all $u \in \mathcal{D}(A)$ is called the \textbf{bound of $B$ relative to $A$}, or the \textbf{$A$ bound of $B$}.
\end{definition}

\begin{theorem}[Kato-Rellich theorem]
\label{Theorem_Kato_Rellich}
    Let $A$ be a self-adjoint operator in a Hilbert space $H$ and let $B$ be a symmetric and $A$-bounded with relative $A$-bound less than 1, then the operator $A+B$ with $\mathcal{D}(A+B)=\mathcal{D}(A)$ is self-adjoint.
\end{theorem}

Loosely speaking, this shows that perturbations of operators by relatively bounded operators as in the statement of the theorem preserve self-adjoint domains of operators. We will use this implicitly in many of the arguments of this section.

\subsubsection{Compactifications of truncated model spaces}

Let $\widehat{U_t(F)}$ be the localizing neighbourhood of a vanishing set in $\widehat{X}$, which resolves to the localizing neighbourhood $U_t(F)$ where we freeze coefficients of the metric in the fiber directions to obtain a model metric. Note that this model metric does not necessarily have isometric fibers $C(Z)$ as assumed in \cite{jayasinghe2024holomorphic}. 

Since $\widehat{U_t(F)}$ is the metric closure of the regular part of the set, which admits a fibration over $F^{reg}$ where the radial distance from the base is $\rho_{F}$, we can construct a topological suspension (or simply the resolved cylinder) of the boundary fibration which we will label as $\Sigma(\mathfrak{B}F)$. Given the model metric 
\begin{equation}    g_{U(F)}=d\rho^2_{F}+\rho_{F}^2g_Z+\phi_{F}^*g_{F^{reg}}
\end{equation}
we can construct a new metric on $\Sigma(\mathfrak{B}F)=[0,3R]\times \mathfrak{B}F$ where $Z- \mathfrak{B}F \rightarrow F$
\begin{equation}    g_{\Sigma(\mathfrak{B}F)}=d\rho^2_{F}+f(\rho_{F})^2g_Z+\phi_{F}^*g_{F^{reg}}
\end{equation}
where $f$ is a smooth function that satisfies $f(x)=x$ on $x<R$, $f(x)=3R-x$ on $x>2R$ and is positive elsewhere. A Dirac-type operator on $U_R(F)$ can be realized as the restriction of a Dirac-type operator on $\Sigma(\mathfrak{B}F)$ (when Clifford modules, and even spin-c structures can be extended in a suitable way) and the ideal boundary conditions on $\mathfrak{B}F$ and $U_R(F)$ match those on $\Sigma(\mathfrak{B}F)$ when one identifies $U_R(F)$ as with the subset $\rho_{F}< R$ in $\Sigma(\mathfrak{B}F)$.

\subsubsection{Discrete spectrum for model Schr\"odinger operators}

\begin{theorem}[Theorem]
\label{Theorem_spectral_discreteness_Schrodinger}
Let $\mathscr{H}' = \Delta + V$ on $L^2(U_{\infty}(F);\mathcal{E})$ be equipped with a self-adjoint domain corresponding to ideal boundary conditions $W$ for which the domain $\mathcal{D}_W(\mathscr{H}')$ is compactly contained in $L^2(\Sigma(\mathfrak{B}F;\mathcal{E})$.
Assume further that the potential $V \in C_{\Phi}(U_{\infty}(F); \mathbb{R})$ is bounded from below and $V(x) \to +\infty$ as $|x| \to \infty$, where $x$ is the radial distance from $F$. Then $\operatorname{Spec}(\mathscr{H}')$ is discrete.
\end{theorem}

\begin{proof}
It suffices to consider the shifted operator $\mathscr{H} = \Delta + V + C$ where $C = 1 - \inf_{\mathbb{R}} V(x)$, and hence $\tilde{V} = V + C \geq 1$. Since the shift is by a constant, any eigensection of $\mathscr{H}$ is one for $\mathscr{H}'$ with the eigenvalue shifted by $C$, and an orthonormal basis of eigensections for $\mathscr{H}$ is also one for $\mathscr{H}'$.      
Since $\mathscr{H} \geq 1$, $\mathscr{H}$ is an invertible operator.

\begin{remark}[details on invertibility]
The operator $\mathscr{H}$ is self-adjoint, and the operator $\Delta + \tilde{V} + i$ is bijective from $\mathcal{D}(\mathscr{H})$ to $L^2(U_{\infty}(F);\mathcal{E})$.
Furthermore, we have the resolvent bound
\begin{equation}
    \|(\Delta + \tilde{V} + i)^{-1}\|_{L^2 \to L^2} \leq 1.
\end{equation}
As the resolvent $(\Delta + \tilde{V} - z)^{-1}$ is holomorphic in the parameter $z$ for $\operatorname{Im} z > 0$, we can choose a path from $i$ to $0$ contained in the upper half-plane. As $z \to 0$, the inequality $(\mathscr{H}f, f) \geq \|f\|^2$ for $f \in \mathcal{D}(\mathscr{H})$ shows that the resolvent cannot have a pole at $0$, and hence $\mathscr{H}$ is invertible.   
\end{remark}

The statement of the theorem is then equivalent to the claim that $\mathscr{H}^{-1}$ is a compact operator on $L^2(U_{\infty}(F);\mathcal{E})$. In fact it suffices to show that the set
\begin{equation}
    \mathcal{M} = \{f \in \mathcal{D}(\mathscr{H}) : |(\mathscr{H}f, f)| \leq 1\}
\end{equation}
is relatively compact in $L^2(U_{\infty}(F);\mathcal{E})$.
If it is, then given any bounded sequence $\|f_n\| \leq C$ in $L^2(U_{\infty}(F);\mathcal{E})$, we have that $|(\mathscr{H}\mathscr{H}^{-1} f_n, \mathscr{H}^{-1} f_n)| \leq \|f_n\|^2 \|\mathscr{H}^{-1}\|_{L^2 \to L^2}$. Hence, by dividing each $f_n$ by a fixed constant $c_1$ we get an associated sequence $c_1^{-1}\mathscr{H}^{-1}f_n$ that is in $\mathcal{M}$, and thus admits a convergent subsequence, showing that $\mathscr{H}^{-1}f_n$ has a convergent subsequence since $\mathcal{M}$ is relatively compact, and thus $\mathscr{H}^{-1}$ is compact.

\begin{lemma}
In the setting described above, for ideal boundary conditions $W$ such that $\mathcal{D}_W(\Delta)$ on $\Sigma(\mathfrak{B}F)$ which are compactly contained in $L^2(\Sigma(\mathfrak{B}F);\mathcal{E})$, $\mathcal{M}$ is relatively compact in $L^2(U_{\infty}(F);\mathcal{E})$.
\end{lemma}

\begin{proof}
Define the set $\mathcal{M}_N = \{f \in \mathcal{M} : f = 0 \text{ for } |x| \geq N\}$. Then, for $f \in \mathcal{M}_N$, we have
\begin{equation}
\label{equation_square_pointer_23}
(\mathscr{H}f, f) = \int_{U_{\infty}(F)} |D f|^2 + \tilde{V} |f|^2 \, dvol.
\end{equation}
where the pointwise norms $|\cdot|$ are those on the fibers $\mathcal{E}$.

Now let $\chi \in C_{\Phi}^\infty(U_{\infty}(F);\mathbb{R})$ be equal to one on the truncated neighbourhood $U_{1/2}(F)$ and vanishes outside of $U_{1}(F)$, and define $B = \|D \chi\|_{L^\infty}$. We now show that for $f \in \mathcal{M}$ and $N > 2B$, $\chi(\cdot/N)f \in 4\mathcal{M}_N$.

We define $\chi_N = \chi(\cdot/N)$ and $g_N = \chi_N f$ and we get
\begin{align}
(Hg_N, g_N) &= \int_{U_{\infty}(F)} |D \chi_N f + \chi_N D f|^2 + \tilde{V} |\chi_N f|^2 \, dvol \\
&= \int_{U_{\infty}(F)} \left( |D \chi_N|^2 + \tilde{V} \chi^2_N \right) |f|^2 + \chi_N^2 |D f|^2 + 2\chi_N \langle f D \chi_N , D f \rangle_\mathcal{E} \, dvol
\end{align}

We estimate the terms to get,
\begin{align}
\left(|D \chi_N|^2 + \tilde{V}\right) |f|^2 &\leq \left(\frac{B^2}{N^2} + \tilde{V}\right) |f|^2 \leq \frac{5}{4} \tilde{V} |f|^2, \,\, \text{since  } \widetilde{V} \geq 1 \,\, \text{and  }  B^2/N^2 <4,\\
\chi_N^2 |D f|^2 &\leq |D f|^2 \\
2  \langle \chi_N f D \chi_N , D f \rangle_\mathcal{E} &\leq |D \chi_N| (|\chi_N f|^2 + |D f|^2) \leq \frac{B}{N} (|f|^2 + |D f|^2) \leq \frac{1}{2} (\tilde{V} |f|^2 + |D f|^2)
\end{align}
where in the last estimate we used the fact that $B/N <1/2 <\widetilde{V}/2$. 
Integrating these bounds we get
\begin{equation}
    (\mathscr{H} g_N, g_N) \leq \int_{U_{\infty}(F)} \frac{3}{2} |D f|^2 + \frac{7}{4} \tilde{V} |f|^2 \leq \frac{13}{4}<4,
\end{equation}
since $f \in \mathcal{M}$. We also have for $f \in \mathcal{M}$ that
\begin{equation}
    \int_{U_{\infty}(F)} |f - \chi_Nf|^2 \, dvol \leq \int_{|x| \geq N/2} |f|^2 dvol \leq \left( \min_{|x| \geq N/2} V(x) \right)^{-1},
\end{equation}
which tends to zero as $N \to \infty$ by our assumption on $V$.

Since any element of $\mathcal{M}$ can be made arbitrarily close to an element in $4\mathcal{M}_N$, it suffices to show $\mathcal{M}_N$ is compact (as this will imply $\mathcal{M}$ is totally bounded).

We first observe that the operator $\widetilde{V} \widetilde{\chi_N}$ (where $\widetilde{\chi_N}$ is the characteristic function of the subset where $|x|<N$ on $U_{\infty}(F)$) is a relatively $\Delta$-compact operator (as a multiplier operator) and has relative bound $0$ (since it has compact support) and by the Kato-Rellich theorem we have that $\mathcal{D}(\mathscr{H})=\mathcal{D}(\Delta)$ for sections in $\mathcal{M}_N$.
Equation \eqref{equation_square_pointer_23} implies that $\mathcal{M}_N$ is bounded in $\mathcal{D}(\Delta)=\mathcal{D}(\mathscr{H})$.
Thus if $f \in \mathcal{M}_N$, then $f \in \mathcal{D}(\Delta)$, which is compactly contained in $L^2(\Sigma(\mathfrak{B}F);E)$.
Since the sections with compact support (in $\mathcal{M}_N$ in $L^2(\Sigma(\mathfrak{B}F);\mathcal{E})$ and $L^2(U_{\infty}(F);\mathcal{E})$ can be identified, $\mathcal{M}_N \subset L^2(U_{\infty}(F);\mathcal{E})$ is compact.
\end{proof}

The claim about eigenfunctions and eigenvalues then follows from the spectral theorem for compact self-adjoint operators.
\end{proof}

\subsection{Spectral gaps and Fredholm complexes in various geometries}\label{subsection-spectral-background}

Having discussed some of the background spectral theory definitions and results we review the various conclusions we can make for different geometric contexts considered in our results.

\begin{theorem}[Spectral gap and Fredholm property for isolated vanishing sets]
\label{Theorem_spectral_properties_of_isolated_vanishing} $ $\\
Given a resolved stratified space $X$ where the singular stratum has depth $1$, a K\"ahler wedge metric and a holomorphic bundle $E$, that admits a stratified K\"ahler Hamiltonian Morse function $h$ (i.e. $h$ has only isolated critical points) and where the group action lifts to one on $E$, then for any choice of ideal boundary condition $W$ that gives a rolled up extension \eqref{equation_rolled_up_extensions} for the Dolbeault Dirac operator, the corresponding Hamiltonian operator $\mathcal{D}_{W,\varepsilon}(\Delta)$ on $X$ is Fredholm, and has a spectral gap.

This result holds for K\"ahler wedge metrics that are asymptotically $\delta$ for $\delta>1/2$.
\end{theorem}

\begin{remark}
This was proven for for K\"ahler wedge metrics which are locally conformally asymptotically $\delta=1/2$ (called locally conformally totally geodesic in \cite{jayasinghe2024holomorphic}).

The locally conformally totally geodesic condition is used in \cite{jayasinghe2024holomorphic} since one can use the conformal covariance of the Dolbeault operators to find the local cohomology groups of the Dolbeault operators for a metric conformal to the model metric at isolated fixed points, and then use them to construct the quasi-modes, which in this article correspond to the cohomology of the normal cohomology complex at the isolated vanishing points of $|dh|^2$. The cohomology for the $P^*$ complex can also be obtained using dualities. This observation can be used to generalize these results for even more general situations.

Once good candidates for the quasi-modes near critical points are found the semi-classical methods used in those articles go through. We shall comment on the $\delta=1/2$ condition in the sketch of the proof that follows.    
\end{remark}

\begin{proof}[Sketch of proof of Theorem \ref{Theorem_spectral_properties_of_isolated_vanishing}]
The result follows from observing that the proof of the spectral gap results for the Witten deformed complex in \cite{jayasinghe2024holomorphic} (a generalization of the proofs in \cite{Zhanglectures}, which itself builds on a lot of literature including \cite{bismut1991complex}) goes through for any choice of ideal boundary condition, and does not assume that there is discrete spectrum to begin with. We explain this now.

In Section 5.2 of \cite{jayasinghe2024holomorphic}, the analog of the block estimates in Proposition \ref{prop:block-estimates} are proven for the setting in the Theorem statement. There the estimates can be used to show that for metrics that are locally asymptotically $\delta$, where the near the critical points the model operators are equal to the actual operators, there exists $C>0$ and $\varepsilon_0>0$ such that for all $\varepsilon>\varepsilon_0$ there exists a spectral gap 
\begin{equation}
    \text{Spec}(\Delta_{\varepsilon})=[0, C\varepsilon^{{(1-\delta)/2}}) \cup (C\varepsilon^{1/2}, \infty)
\end{equation}
which for $\delta \in (0,1/2)$ is one where the upper bound of the lower interval (corresponding to the Witten subcomplex) may grow, but at a rate smaller than the rate of growth of the lower bound for the complement of the sections of the Witten subcomplex.

\begin{figure}[h]
    \centering
    \includegraphics[scale=0.60]{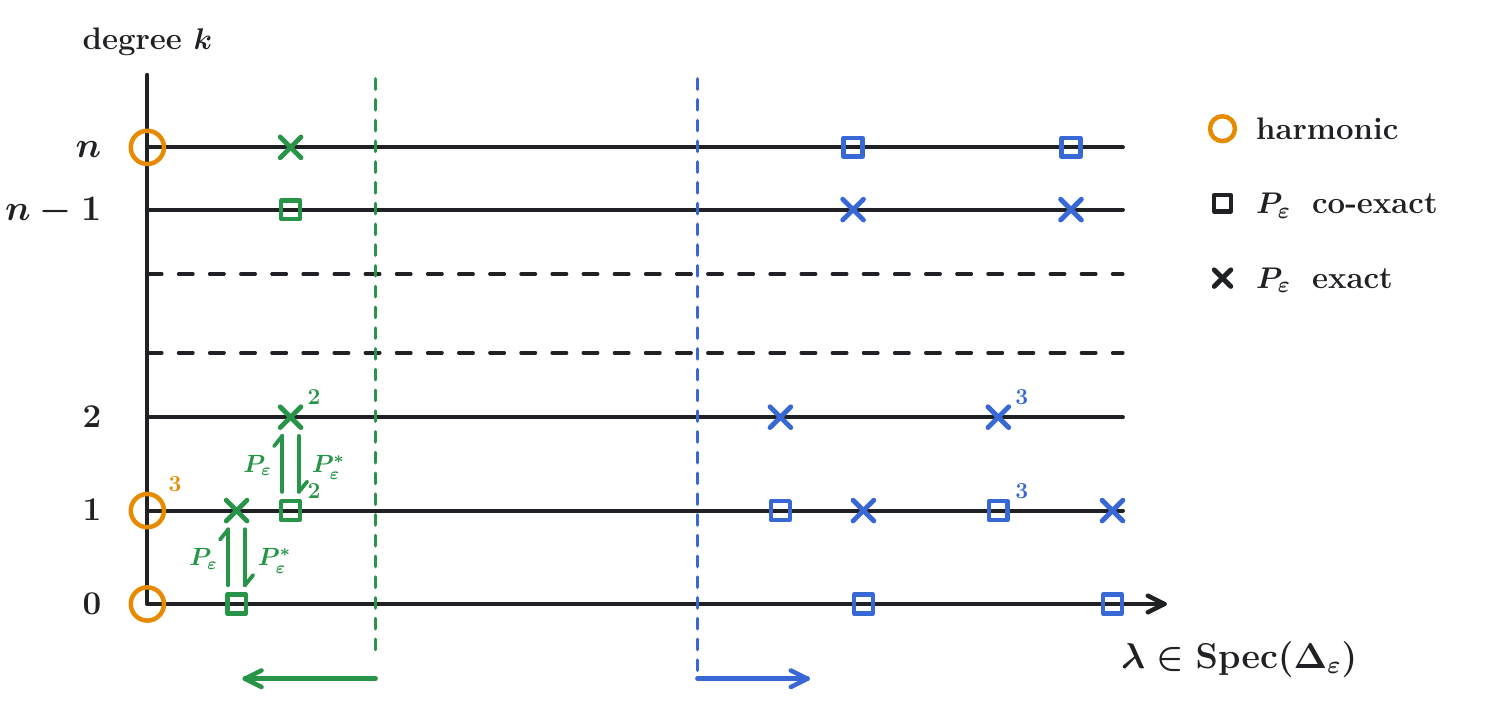}
    \caption{A spectral abacus for a $\mathbb{Z}$ graded supersymmetric complex corresponding to isolated critical points. The eigen sections with eigenvalues $\lambda$ are marked (with multiplicities in superscripts). The sections marked in green and orange constitute the Witten complex.}
    \label{figure_spectral_abacus_isolated}
\end{figure}

The methods of the present article show that this can be improved for all metrics that are asymptotically $\delta$ for $\delta>1/2$.

The proof of the result in \cite{jayasinghe2024holomorphic}, while stated assuming there is discrete spectrum for the operator $\mathcal{D}_W(\Delta)$ on $X$ (there for certain choices of domains), goes through without this assumption, as long as there is discrete spectrum for the model harmonic oscillators on the localizing neighbourhoods of critical point sets. We now help the reader to trace through the proof to see this.

In the case of isolated vanishing sets, the method of the Cheeger ansatz with generalized Dirichlet/Neumann boundary conditions on truncated cones and with $L^2$ conditions on the infinite cones, can be used to show that the local equivariant Hilbert complexes have discrete spectrum (see \cite{jayasinghe2023l2,jayasinghe2025supersymmetricharmonicoscillatorssingular} for more details), as long as there is discrete spectrum for the corresponding operators on the links. This is proven by giving explicit orthonormal bases for the Hilbert complexes, using separation of variables and Sturm-Liouville theory. Here we use the assumption in the theorem statement that the singular strata are of depth 1, which implies that the links are smooth, and so the relevant operators on the links have discrete spectrum.

Then the quasimodes notated as $G_{\varepsilon,\mu}$ in Subsection 5.2 of \cite{jayasinghe2024holomorphic} are finite for each $\mu$, and so give a closed subspace of the domain of the Hilbert complex, and the Hilbert space itself. Then the proofs of Propositions 5.2 and 5.9 of that article go through without the assumption of discrete spectrum for the operators on $X$, since they only use projection valued measures, which are amenable for the self-adjoint operators that are studied, and show that there is a spectral gap, and that the spaces of eigensections below that are finite. This proves the statement.
\end{proof}

More generally, there are certain scenarios with non-isolated vanishing loci in which one can make the same conclusions as theorem \ref{Theorem_spectral_properties_of_isolated_vanishing}: one notable class of examples are those considered in \cite{jayasinghe2024witten}, in which hypotheses on the Morse-Bott function $h$ on a stratified pseudomanifold satisfies a certain flatness condition. Namely, if for every connected component $F\in Z(|dh|^2):=\text{Crit}(h)$ of the vanishing loci, its corresponding the singular localizing neighborhood $U(\widehat{F})$ is a flat (fiber) bundle, which implies we can find a wedge metric $g'$ over the resolved localizing neighborhood $U(F)$ which admits a strong product decomposition; this in turn implies a certain splitting of the spectrum of the Witten-deformed Laplacian $\Delta_\varepsilon = \Delta^H+\Delta_\varepsilon^V$. We state this implication in terms of this spectral condition and then discuss other geometric hypotheses which imply this condition. We recall the following from \cite[Prop 4.8]{jayasinghe2024witten}

\begin{proposition}[Spectral gap and Fredholm property for certain non-isolated vanishing sets.]\label{prop:specgap-flatness}
    Let $X$ be a resolved stratified space $X$ with a wedge metric, and Witten-deformed Laplace-type operator $\Delta_\varepsilon$ such that there exist constants $C_1, C_2>0$ such that the restriction of $\Delta_\varepsilon$ to the localized complex $\mathcal{Q}_{\mathcal{B},\varepsilon}(U_t(F))$ is self-adjoint and has spectrum satisfying
    \begin{equation}
        \emph{Spec}(\Delta_\varepsilon) \subset \{0\}\cup [C_1,\infty) \cup [C_2\varepsilon,\infty) ,
    \end{equation}
    and thus the operator participates in a Fredholm complex and has a spectral gap. In particular for large $\varepsilon>0$, $[C_2\varepsilon,\infty) \subseteq [C_1,\infty)$, and we have a fixed spectral gap.
\end{proposition}
Note that while the second and third intervals appear to have some large intersection, the upshot of this result is that while the operator $\Delta_\varepsilon$ may fail to have a \emph{global} spectral gap on the full Hilbert complex $L^2(X;\mathcal{E})$, it nevertheless possesses such a gap on any of the localized complexes $\mathcal{Q}_{\mathcal{B},\varepsilon}(U_t(F))$ for every $F\in Z(A)$. 
In particular whenever the conclusions of the Proposition \ref{Theorem_spectral_properties_of_isolated_vanishing} or \ref{prop:specgap-flatness}, we can prove (stronger) versions of the block estimates of proposition \ref{prop:block-estimates}, just as is done in \cite[Prop 5.2]{jayasinghe2024witten}. Thus the Witten inequalities of this paper follow from the arguments presented there. \\

\emph{However}, for an arbitrary stratified pseudomanifold $X$ equipped with a $\mathcal{N}=N$ SUSY Hilbert complex, such that the associated Hamiltonian satisfies $\Delta_\varepsilon=\Delta+\varepsilon K+\varepsilon^2 A$, it is possible for the vanishing loci to fail to satisfy the flatness condition used in the above result. For example even the fixed points of the toric action on $\mathbb{CP}^2$ considered in Subsection \ref{subsection_smooth_toric_action} (see also Example \ref{Example_singular_toric_variety_domains}) does not enjoy the flatness condition required in \cite{jayasinghe2024witten} and described above. In such cases, we must assume that our domain for operator $\mathcal{D}_W(D_\varepsilon)$ on $X$ is such that $D_\varepsilon$ has discrete spectrum. One family of examples of such domains are considered for example in \cite{Albin_2017_index} who prove that certain perturbations of Dirac operators with the Vertical APS domain form a self-adjoint complex with discrete spectrum. Other examples include the minimal domain for the Dolbeault complex on spaces with wedge K\"ahler metrics (see Proposition 3.6 of \cite{AnalyticToddBeiPiazza}).

In all such cases, Theorem \ref{Theorem_spectral_discreteness_Schrodinger} can be applied to see that the normal cohomology complexes are Fredholm and have a spectral gap. Once this property of the normal cohomology complexes is established, we can proceed as in Subsection \ref{subsec:resolvent-ests} and prove the conclusions of proposition \ref{prop:block-estimates}.

\section{Appendix B: Perspectives from SUSY in Index theory}
\label{Appendix_SUSY_and_conjectures}

Jointly written with Alex R. Taylor.

\begin{center}
\centering
\itshape
``And the point of my lectures is to particularly motivate young people, when they're still young, to learn both the requisite mathematics from experts in the field and the requisite physics from the experts in the physics field so that they not only can speak both languages, but they can begin to think the way these two different groups think.''
    
\vspace{0.7em}

\hfill --- Isadore Singer, Simons lectures, Oct 1996 \cite{Singer1996}
\end{center}

\begin{center}
\textbf{To Pierre Albin, with admiration.}
\end{center}

In this appendix we explore some new perspectives on the index theory of $\mathcal{N}=1$ and $\mathcal{N}=2$ SUSY complexes. We begin by reminding the reader of some of the deep and mysterious interactions of ideas in mathematical physics and index theory, placing SUSY at the very heart of index theory.
Their investigations have proven fruitful and have shed light on previously known phenomena, as well as certain hitherto unknown phenomena which we explain in this article.

\begin{figure}[h]
    \centering
    \includegraphics[scale=0.9]{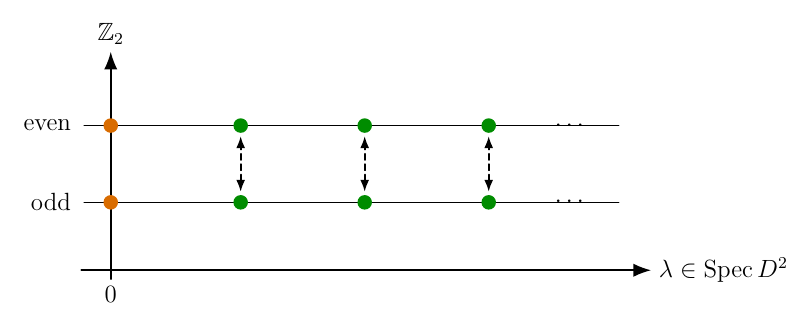}
    \caption{Supersymmetry of Dirac operators corresponding to McKean Singer cancellations}
    \label{figure_McKean_Singer}
\end{figure}

In their article ``Curvature and the eigen\textit{values} of the Laplacian'' \cite{mckean1967curvature}, McKean and Singer explained how the non-zero eigenvalues of the Laplace type operators (squares of Dirac-type operators) appear as pairs, and thus taking spectral expansions of alternating sums of traces of the heat kernels, we have that
\begin{equation}
    \sum_{k=0}^{1} (-1)^k \operatorname{Tr} \mathcal{K}(t,x,y)=\sum_{k=0}^{1} (-1)^k \sum_{j} e^{-t\lambda_j^2}\langle \phi_j^k,\phi_j^k \rangle= \sum_{k=0}^{1} (-1)^k \sum_{j|\lambda_j=0} e^{-t0}\langle \phi_j^k,\phi_j^k \rangle
\end{equation}
where $\{\phi_j^k\}_j$ form an orthonormal basis of eigensections of the operator $D^2$ in even/odd degree ($k=\pm 1$) with eigenvalues $\lambda_j$. 
Here the non-zero eigenvalues cancel off in the alternating trace since for each positive eigenvalue of $D^2$ of odd degree, there is one of even degree, and we get an alternating trace on the harmonic sections on the right hand side.
This is what was later termed SUSY cancellations, with the positive eigensections arising as SUSY pairs. This gives a model for how in SUSY field theories, for every bosonic state of a given energy, the SUSY yields a fermionic state of the same energy, for positive energy. The alternating trace is called a supertrace.

The question of whether the integrands of the heat kernel could be expressed in purely topological terms using such cancellations was raised by them, and was addressed by Patodi in ``Curvature and the eigen\textit{forms} of the Laplace operator'' \cite{Patodi1971CurvatureAT}. This was later approached more systematically by Getzler in ``Pseudodifferential operators on supermanifolds and the Atiyah-Singer index theorem'' \cite{getzler1983pseudodifferential}, where he developed a more rigorous version of the work of several physicists including ``Supersymmetry and the Atiyah-Singer index theorem'' \cite{alvarez1983supersymmetry} by Alvarez-Gaum\'e, who studied both the index and equivariant index theorems using perspectives from Supersymmetric Quantum field theory (the notion of $\mathcal{N}=N$ SUSY we use in this article corresponds to $\mathcal{N}=N/2$ in the aforementioned article).

Getzler considers what we call a $\mathcal{N}=1$ SUSY complex ($\mathcal{N}=1/2$ in the sense used in \cite{alvarez1983supersymmetry}), and models a proof of the local index theorem on smooth manifolds on a proof of Weyl's theorem on a compact Riemannian manifold. Weyl's theorem states 
\begin{equation}
\label{equation_Tr_Weyl}
    \operatorname{Tr} e^{-t\Delta}
    = \left(\frac{\pi}{t}\right)^{n/2}
      \operatorname{vol}(M) + O(t^{-n/2+1}),
\end{equation}
and using the calculus of pseudodifferential operators, one can shows that the symbol $\sigma$ of the heat kernel satisfies,
\begin{equation}
    \sigma(e^{-t\Delta})
    = e^{-t|\xi|^2} + \text{small error},
\end{equation}
where $|\xi|^2$ is in fact the symbol of $\Delta$.
Combining this with a formula for the trace of a pseudodifferential
operator
\begin{equation}
    \operatorname{Tr} P
    = \int_{T^*M} \sigma(P)\, dx\, d\xi,
\end{equation}
the theorem follows, by evaluating the Gaussian integrals along each fibre. The principal symbol of $P$ gives an integrand corresponding to a Gaussian cutoff of the volume measure of the cotangent bundle, and the integral may give extra terms with higher order powers of $t$.

Getzler considers a supermanifold over the classical cotangent bundle, and builds a pseudodifferential calculus, in which the principal symbol is precisely the local index formula; the integral of the principal symbol gives the index, where there are no higher order error terms as in equation \eqref{equation_Tr_Weyl} due to the supersymmetric cancellations.

The equivariant index theorem on smooth manifolds for actions with non-isolated singularities also has a proof using Getzler's method (see \cite{berline2003heat}), and this allows one to prove the Duistermaat-Heckman theorem, which has been phrased as the fact that \text{stationary phase is exact} when there are sufficient symmetries (see \cite{atiyah1985circular}). Here the reason for the vanishing of the higher order contributions to the stationary phase expansion of Hamiltonian functions at their critical points can be attributed to SUSY cancellations.

In fact, in \cite{atiyah1985circular}, Atiyah describes Witten's conjecture with a heuristic ``proof'' explaining how the Duistermaat-Heckman theorem on the loop space of a spin manifold yields the local index for the spin Dirac operator on the manifold, where the differential form corresponding to the local index at a point is obtained in terms of geometric objects on the tangent space to the constant loop corresponding to the said point.

In their series of papers titled ``Spectral asymmetry and Riemannian Geometry'' starting with \cite{atiyah1975spectral}, Atiyah-Patodi-Singer (APS) proved their index theorem on a manifold with boundary $M$, where the contribution from the boundary corresponded to the \textit{reduced eta invariant}, also called the $\xi$ invariant. This invariant is typically expressed in terms of the spectral data of a Dirac operator $A=cl(dr)D|_{\partial M}$ on the boundary, induced by the Dirac operator $D$ on $M$ (where in a collar neighbourhood of the boundary one considers a cylindrical metric), where only the eigensections of the even degree sub-bundle are considered. Later this was extended to an equivariant index theorem by Donnelly \cite{donnelly1978eta}, where the information of the group action on the eigenspaces provided the formula at the boundary, and this boundary defect is called the equivariant $\xi$ invariant (also known as the reduced eta invariant) and we denote it by $\xi_T$ for geometric endomorphisms corresponding to group actions ($\xi_{Id}=\xi$).

Proposition 3.11 of \cite{atiyah1975spectral} can be used to express this invariant as a supertrace over the $L^2$ bounded harmonic sections on an infinite cylindrical extension of the boundary (this has also been studied using the b-calculus by Melrose in \cite{melrose1992calculus}). 
This suggests a natural categorification here which has not been formalized in the literature, which is to consider the above mentioned $L^2$-bounded harmonic sections on the cylinder as a \textbf{Dirac cohomology group $\mathcal{H}^k(\mathcal{D}_\mathcal{B}(D))$} where $\mathcal{B}$ denotes the APS boundary conditions and the cohomology are considered not just in the even degree bundle, but also in the odd degree bundle. In fact the negative eigenvalue eigensections of the operator $A$ correspond to $L^2$ bounded harmonic sections in odd degree, while the positive ones correspond to those in even degree; one can show that with the correct definition of the Dirac cohomology, the $\xi_T$ is the supertrace of $T$ over the Dirac cohomology divided by $2$.
For now, this seems like a reformulation of the $\xi_T$ invariant, but we shall soon study how this simplifies at the category of Dirac cohomology when there are extra supersymmetries.

Before we study this, let us consider a few other results where the contributions to the equivariant and non-equivariant indices of Dirac type operators have been studied using spectral theoretic methods, with spectral formulae.

In \cite{cheeger1983spectral} (c.f. \cite{cheeger1980hodge,cheeger1982l2}), Cheeger studied how for the Hodge Dirac and signature operators on spaces with isolated conic singularities, the index contribution at a cone point came from an eta invariant. Moreover in the case of the Hodge Dirac operator, the equivariant index contribution matches that in the intersection cohomology of Goresky and MacPherson \cite{Goresky_1985_Lefschetz,Goresky_1993_local_Lefschetz}, which is given by a supertrace of the geometric endomorphism on the local intersection cohomology of a neighbourhood of the singular fixed point.
In the equivariant Riemann-Roch formula of Baum-Fulton-Quart in \cite{baum1979lefschetz,Baumformula81}, the formula at a singular isolated fixed point of a locally complete intersection is given by a trace of the geometric endomorphism over the local ring. There are no such explicit formulas given for non-equivariant Riemann-Roch numbers at singularities in the work of Baum-Fulton-MacPherson (BFM) \cite{baumfultonmacKtheoryRiemannRoch,baumfultonmacpheresonRiemannRoch}. We note that the formula of Wu and Zhang in \cite{wu1998equivariant} (see also \cite{wu2003instanton}) is different from that of Atiyah-Segal-Singer \cite{atiyah1965index_2} (c.f. \cite{berline2003heat,duistermaat2013heat}) for spin-c Dirac operators on K\"ahler manifolds with Hamiltonian actions, where the contributions at the fixed points are supertraces over local cohomology groups, and we note that no such formula is even hinted at in the work of BFM.
There have been generalizations of these, for instance by Edidin and Graham \cite{equivariant_intersection_edidin,localization_algebraic_Graham_Edidin, MaximJorgCharacteristicsingulartoric2015}, which are slightly different in nature, using more algebraic tools.

In \cite{weiping1990note}, Zhang expresses the equivariant index of a general Dirac-type operator (under some mild conditions) at an isolated conic fixed point in terms of a $\xi_T$ invariant.
In \cite{nazaikinskii2005elliptic} spectral zeta functions are used to explain contributions to the equivariant index at fixed points, and this is extended by Bei for Hilbert complexes in \cite{Bei_2012_L2atiyahbottlefschetz}. In the special case of the Hodge de Rham operator, the equivariant index Bei matches it with that of Goresky and MacPherson in intersection cohomology. 

In addition, there are many Morse inequalities for the de Rham case, for instance in \cite{Jesus2018Wittens,Jesus2018Wittensgeneral,ludwig2013analytic,ludwig2017index} where the formulas in particular give Lefschetz numbers for the gradient flow of stratified Morse functions on spaces with conic and other singular warped product metrics in terms of local cohomology groups for de Rham operators.

Indeed, since the equivariant index formula of \cite{weiping1990note} as a $\xi_T$ invariant applies to certain smooth settings as well, where the tangent space of a smooth isolated fixed point can be considered as a cone over a sphere, there is question on how the Atiyah-Bott formulas \cite{AtiyahBott1,AtiyahBott2} correspond to each of these cases.

All these cases have been studied disparately, and there had been no attempts to unite these perspectives in the literature as far as we know, until recent work by the authors. 
We have now come to realize that it is a lack of proper understanding of SUSY that has held back insight and perspective as to the relationships between these invariants which are closely related to each other.

In \cite{jayasinghe2023l2}, a framework to understand $L^2$ Lefschetz fixed point theorems and Morse inequalities for stratified spaces was laid out, tackling many cases of Dirac operators. In particular it was shown that in the case of both the de Rham ($P=d$) and Dolbeault ($P=\overline{\partial}$) complexes, the equivariant index formulas for the corresponding Dirac operators could be phrased in terms of local cohomology groups for the local operators $P,P^*$, where for the local Morse polynomials the particular choice of $P,P^*$ was related to the dynamics of the gradient flow of the Morse function (see also \cite{jayasinghe2024holomorphic}). This should be compared with the formulas for the equivariant index for a general Dirac operator in \cite{jayasinghe2023l2} (when the local neighbourhoods of fixed points do not admit $\mathcal{N}=2$ refinements in the language of this article), where it is formulated in terms of a equivariant Bismut-Cheeger $\mathcal{J}$ form, turning it into an equivariant version of the index theorem in \cite{Albin_2017_index} (c.f. \cite{Albin_2016_Index}).

While these insights were facilitated through viewing them through the lens of the formalism of Witten deformation, the true role of the second SUSY was only made clearer in the investigations in \cite{jayasinghe2026spectralasymmetrysupersymmetryequivariant}, where the relationship between spectral asymmetry and supersymmetry were first investigated. \textbf{We now proceed to give an intuitive explanation of the insights of these investigations, some of which will be laid out with more precision in an upcoming expansion of \cite{jayasinghe2026spectralasymmetrysupersymmetryequivariant}.}

\subsection{On the second SUSY and residual cancellations}
\label{subsection_appendix_second_susy_cancellation}

Consider a manifold with boundary, where the boundary has a neighbourhood isometric to a metric cylinder $C_x(\partial M)$. We equip it with a Dirac operator, and as in \cite{atiyah1975spectral} equip it with the APS boundary conditions. We can also consider topological cones, with conic metrics or other singular warped product metrics \cite{Jesus2018Wittens,Jesus2018Wittensgeneral,jayasinghe2025supersymmetricharmonicoscillatorssingular} where for simplicity we can consider the situations where the Dirac operators are essentially self-adjoint.
In such cases, the boundary contribution to the equivariant and non-equivariant indices is the $\xi_T$ invariant,
\begin{equation}
\label{equation_xi_T}
    \xi_T(s):=(1/2)\text{Str}_s\, T |_{\mathcal{H}^k(\mathcal{D}_{\mathcal{B}}(D))}, \quad \xi_T=\xi_T(0).
\end{equation}
This is a simple reformulation in the case of the manifold with boundary of the main result in \cite{atiyah1975spectral} as discussed earlier. Here there are various equivalent choices of regularization that have been used in different contexts (and indeed some of the issues with the regularization seem to have been distracting previous researchers of these objects from studying their symmetries), but we will focus on the symmetries at the level of Dirac cohomology.

It can be similarly seen that such a statement holds for conic singularities, where the reader can get a better idea of the statements we shall make since we study these in greater detail in this article. Here for the Dirac cohomology on a truncated cone, one can think of it as the span of the harmonic extensions of the eigensections of the link to the cone, which are $L^2$ bounded near the singularity. When the operator is not essentially self-adjoint, then additional conditions have to be stated which we do not consider here.

Given a $\mathcal{N}=2$ SUSY refinement of the Dirac operator as $D=P+P^*$ the Dirac cohomology admits a threefold splitting
\begin{equation}
\label{equation_Dirac_cohomology_decomposition_threefold}
    \mathcal{H}^k(\mathcal{D}_\mathcal{B}(D))=\mathcal{H}^k_P\oplus \mathcal{H}^k_{P^*}\oplus \mathcal{H}_3
\end{equation}
and a simple way to state what happens is that \textbf{regularized supertraces over the third partition cancel off and the first two cohomology groups are dual to each other}. In particular, another good name for the cancellations can be \textbf{cohomological SUSY cancellations}. Here $\mathcal{H}^k_P$ corresponds to the cohomology of the Hilbert complex for $P$ with generalized Neumann conditions for $P$ at the boundary of the truncated cone. Similarly $\mathcal{H}^k_{P^*}$ is the cohomology of $P^*$ with the generalized Dirichlet conditions for $P^*$ at the boundary. The third partition then is uniquely characterized as the orthogonal complement of the first two spaces in the Dirac cohomology, but we present a different description.

For any Dirac-type operator near a boundary/ an isolated conic singularity with metric $g=dx^2+x^{2a}g_Z$ (where for boundaries $a=0$, for cones $a=1$), we assume that there is an expansion
\begin{equation}
\label{equation_Dirac_112}
    D=\sigma\left(\nabla_{\partial_x}-\frac{1}{x^a}A\right)
\end{equation}
where $\sigma=\sigma_1(D)(dx)$ which switches the parity of the $\mathbb{Z}_2$ grading, and an operator $A$ on sections of a Clifford module (suitably rescaled in the wedge setting).
A wedge Dirac-type operator can be put in this form  after some rescaling of the sections on the bundle (see \cite{Albin_2017_index}, where the sign for $A$ is picked differently).
However when there is a second SUSY operator, then the operator $P$ admits a splitting $P=P_1+P_2$ where $P_1$ can be understood as follows.

The Clifford module $\mathcal{E}$ admits a decomposition
\begin{equation}
    \mathcal{E}=\mathcal{E}' \oplus \sigma_1(P)(dx) \mathcal{E'}
\end{equation}
where $\mathcal{E}'$ is determined by the condition that all the sections of this sub-bundle are annihilated by $\sigma_1(P^*)(dx)$.
This determines the splitting of the operator as $P=P_1+P_2$ by
\begin{equation}
P = \begin{pmatrix}
    P_2 & 0 \\ P_1 & P_2
\end{pmatrix}.
\end{equation}

Then the sections in $\mathcal{H}(D)$ correspond to $L^2$ bounded Dirac harmonic sections not in the minimal domain of the operator on the cone. This leads to the typical analysis of small eigenvalues of the operator $A$, and we refer to \cite{Albin_2017_index} for more details.
This yields a threefold splitting of the operator $D$
\begin{equation}
    D=P_1+P_1^*+(P_2+P_2^*)=\sigma\left(\nabla_{\partial_x}-\frac{1}{x^a}(A_1+A_2)\right)
\end{equation}
where $P_2,P_2^*$ are nilpotent (square to 0) as well. Moreover, $P_2+P_2^*$ +commutes with $D^2$ and we can find a mutually diagonalizable basis, similar to the case explained in Lemma \ref{lem:exact-coexact-eigenvectors} where $P_2+P_2^*$ gives SUSY pairs. 
Thus we see that any element in the kernel of $D$ is necessarily in the kernel of $P_2+P_2^*$. 
The operator $A$ is equal to $A_1+A_2$ where $\sigma A_2$ is equal to $(P_2+P_2^*)$, up to certain rescalings of bundles and up to multiplication by the bundle morphism $\sigma=cl(dx)$.
We have that the third space is $\mathcal{H}_3=\{\text{ker}(P_2+P_2)^*\}^{\perp}$, while $\mathcal{H}_P$ are the sections of the Dirac cohomology which take values in $\mathcal{E}'$ and are in $\text{ker}(P_2+P_2)^*$, while $\mathcal{H}_{P^*}$ are the section of the Dirac cohomology with values in $\sigma_1(P)(dx)\mathcal{E}'$ in the kernel of $P_2+P_2^*$.
\begin{remark}
\label{Remark_appendix_deRham_vs_Dolbeault_1}
    In the de Rham setting, $P_2=d_Z$, and $A_1=0$, $A_2=A$.    
    In the Dolbeault setting $P_2=\overline{\partial}_b$, which is the transversal Dolbeault operator on the contact distribution transverse to the foliation generated by $x^aJ(\partial_x)$. Thus $A_1$ corresponds to an operator where the derivatives are only with respect to the vector field $J(\partial_x)$. 
\end{remark}

\begin{figure}[h]
    \centering
    \includegraphics[scale=0.8]{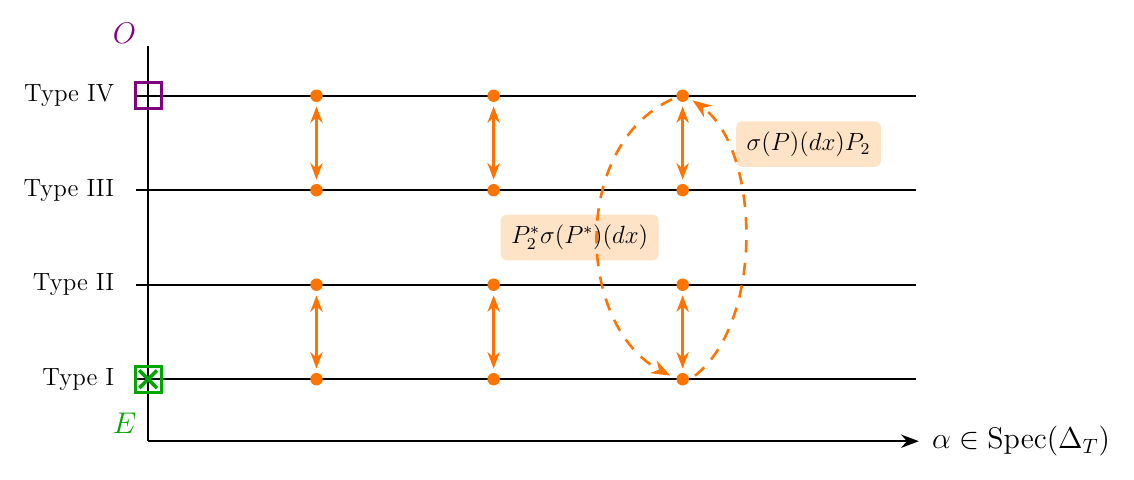}
    \caption{Sections of $\mathcal{H}^k(\mathcal{D}_{\mathcal{B}}(D))$ under the splitting of the Cheeger ansatz. Each orange dot corresponds to a harmonic section of $D^2$ not in the kernel of $A_2^2=\Delta_T$, and the solid arrows between them are the maps corresponding to the operator $A_2$. The green crosses correspond to $\mathcal{H}^k_P$ while the purple boxes correspond to $\mathcal{H}^k_{P^*}$.}
    \label{figure_Cheeger_abacus_2}
\end{figure}

In fact, these decompositions and the corresponding analysis of the spectrum in the case of the Hodge de Rham complex on cones were already hidden in the work of Cheeger in \cite{cheeger1983spectral}, where the study included even non harmonic eigensections on the cones. This was extended to the Dolbeault complex on cones, and for both de Rham and Dolbeault complexes on truncated cones, with boundary conditions in \cite{jayasinghe2023l2}. We also refer to \cite{jayasinghe2025supersymmetricharmonicoscillatorssingular} where a comprehensive study was made, including with Witten deformation, for more general metrics.

\begin{figure}[h]
    \centering
    \includegraphics[scale=0.7]{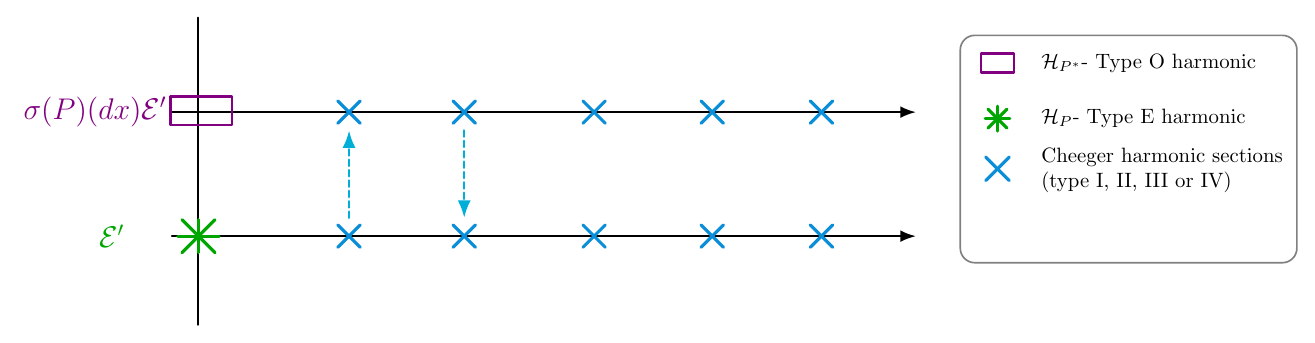}
    \caption{Sections of $\mathcal{H}^k(\mathcal{D}_{\mathcal{B}}(D))$ under the splitting of the Clifford module}
    \label{figure_Cheeger_abacus_1}
\end{figure}

We explain this briefly. The eigensections of the Dirac operator $A$ on the boundary/link come in 6 types, $E, O, I,II,III,IV$.
The eigensections $E,O$ are those that are in the kernel of $P_2+P_2^*$, and take values in $\mathcal{E}',\sigma_1(P)(dx)\mathcal{E}'$ respectively.
The eigensections of type $I$ are those that are $P_2$-coexact (in the image of $P_2^*$) and take values in $\mathcal{E}'$, and those of type $II$ are precisely their superpartners, obtained by applying $P$ to those of type $I$.
The eigensections of type $IV$ are those that are $P_2$-exact (in the image of $P_2$) and take values in $\mathcal{E}'$, and those of type $III$ are precisely their superpartners, obtained by applying $P^*$ to those of type $IV$.
We record the information of the harmonic sections of these types in Figure \ref{figure_Cheeger_abacus_2}. There, the operator $\sigma(P)(dx)P_2$ maps sections of type $I$ to those of type $IV$, and $P_2^*\sigma(P^*)(dx)$ maps those of type $IV$ to type $I$. 
Figure \ref{figure_Cheeger_abacus_1} shows the same sections, but corresponding to the splitting of the Clifford module.

\begin{figure}[h]
    \centering
    \includegraphics[scale=0.8]{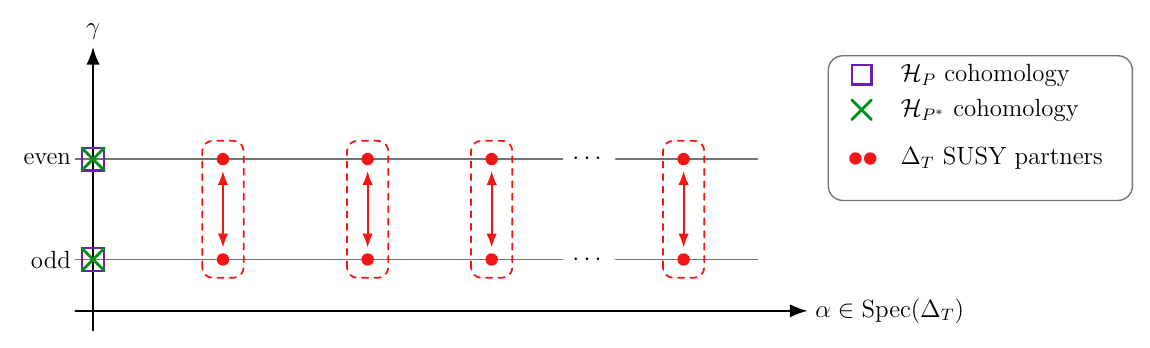}
    \caption{The sections of $\mathcal{H}^k(\mathcal{D}_{\mathcal{B}}(D))$, organized in odd and even degrees for the eigenvalues with respect to $\Delta_T=A_2^2$.}
    \label{figure_Cheeger_Sections_odd_even_delta_T}
\end{figure}

Figure \ref{figure_Cheeger_Sections_odd_even_delta_T} indicates the same sections arranged by the eigenvalues of $\Delta_T=A_2^2$, and for the odd and even degrees, where there are superpartners for all the non-zero eigenvalues. In the de Rham case, since $A_1=0$, $A_2=A$ is simply the Hodge de Rham operator on the link, which is the second SUSY operator.

\begin{figure}[h]
    \centering
    \includegraphics[scale=0.8]{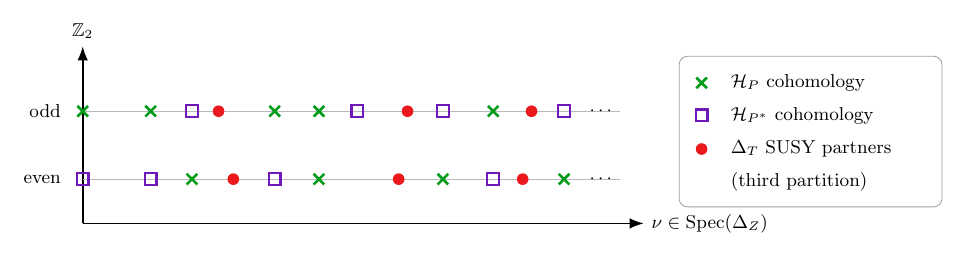}
    \caption{The sections of $\mathcal{H}^k(\mathcal{D}_{\mathcal{B}}(D))$, organized in odd and even degrees for the eigenvalues with respect to $\Delta_Z=A^2$, when $A_1 \neq 0$. The eigenvalues of the superpartners of $A_2^2$ may shift slightly since $A=A_1+A_2$.}
    \label{figure_Cheeger_Sections_odd_even_delta_Z}
\end{figure}

In the Dolbeault setting, as explained in Remark \ref{Remark_appendix_deRham_vs_Dolbeault_1}, the operator $A_1$ is non-trivial, and when arranged by the eigenvalues of $\Delta_Z=A^2$, we can no longer see the SUSY as clearly. We illustrate this in Figure \ref{figure_Cheeger_Sections_odd_even_delta_Z}. In particular this plays a role when regularizing equivariant eta invariants, $\xi_T$ invariants, and Lefschetz fixed point contributions (see \cite[\S 5.2.5]{jayasinghe2023l2}, and \cite{jayasinghe2026spectralasymmetrysupersymmetryequivariant}).
The de Rham setting is much simpler since the $P,P^*$ cohomology groups are finite dimensional.

It should be clear that the supertrace in \eqref{equation_xi_T} should be regularized suitably to respect the symmetries of the eigenvalues. These spectral symmetries are obscured if one uses the spectrum of the Dirac operator on the link $A$ instead of $A_2$ for the Dolbeault complex. This may very well have been one reason these symmetries were not spotted earlier. 

Since the third summand of the partition in \eqref{equation_Dirac_cohomology_decomposition_threefold} has SUSY, the supertraces cancel and one can reduce the $\xi_T$ invariant to (regularized if infinite) supertraces over the $P,P^*$ cohomology groups (see \cite{jayasinghe2026spectralasymmetrysupersymmetryequivariant}). In fact, if domains that respect the $\mathcal{N}=2$ SUSY structure are picked for the Dirac operators, then the duality of the Hilbert complexes will yield a duality for the supertraces over the $P,P^*$ cohomology groups, and the equivariant index contributions over both will be the same. Then the factor of 1/2 in \eqref{equation_xi_T} shows how the supertrace over the Dirac cohomology group is equal to the supertrace on the cohomology groups $(P,P^*)$ considered separately (both of which are themselves equal to one another).

In the case of the Dolbeault complex, this splitting in \eqref{equation_Dirac_cohomology_decomposition_threefold} corresponds to the splitting of the harmonic sections as 
\begin{equation}
    \text{harmonic}=\text{holomorphic}+\text{anti-holomorphic}+\text{(neither holomorphic nor anti-holomorphic)}
\end{equation}
and reflects the natural intuition that the third partition should have nothing to do with the Riemann-Roch numbers. That the $P,P^*$ cohomologies have the same information is why it suffices to look at just one, as is commonly done in algebraic geometry where anti-holomorphic sections are not studied as keenly as their holomorphic counterparts.
In the de Rham setting, the splitting corresponds to the cohomology of $d_Z$ and the wedge product of the cohomology of $d_Z$ by $x^{aN} dx$, and the third summand in the partition consists of sections that are not in the kernel of $d_Z+\frac{1}{x^{2a}}d_Z^*$.
Here $N$ is the \textit{number operator}, which is just multiplication of the form by its degree.
This can be used together with the cohomology of a sphere to obtain Poincar\'e-Hopf contributions to fixed points from the $\xi_T$ invariant on the tangent space of the fixed point considered as a cone over a sphere, which can be a nice exercise for the motivated reader. For a different but related approach see \cite{zhang2000n,chen2018poincar}.

\begin{figure}[h]
    \centering
    \includegraphics[scale=0.8]{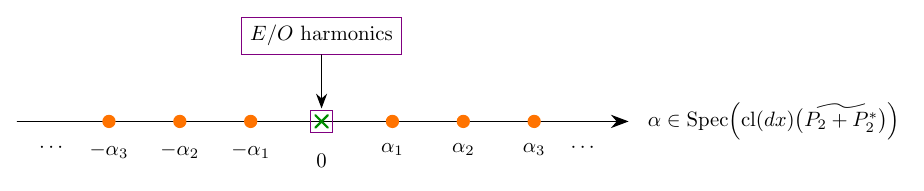}
    \caption{The restriction of the sections in $\mathcal{H}^k(\mathcal{D}_{\mathcal{B}}(D))$ to the even degree sub-bundle on the link, arranged by the eigenvalues of the residual SUSY operator on the link. The restriction to the odd degree sub-bundle is similar.}
    \label{figure_Eta_spectrum}
\end{figure}

Typically people have tended to work with the spectrum of the operator $A$ on the link, restricted to the even degree sub-bundle, when studying the $\xi_T$ invariants, so we provide Figure \ref{figure_Eta_spectrum} to show that the eigenvalues of $cl(dx)(P_2+P_2^*)$ (which after rescaling in the conic case we indicate as $cl(dx)(\widetilde{(P_2+P_2^*)})$ in the figure) have spectral symmetries on the third partition. Thus regularizing eta and $\xi_T$ invariants with this SUSY operator leads to a better understanding of such invariants.
In fact, as we noted earlier and as reflected in Figure \ref{figure_Cheeger_abacus_2}, there is a SUSY operator intertwining the eigensections of types $I$ and $IV$, given by multiplying by $cl(dx)(P_2+P_2^*)$, which when restricted to the type $I, IV$ sections simplifies to the operators denoted in that figure. Together with the other SUSY intertwining the sections of types $I,II$ and $III,IV$, for each eigenvalue of $A_2^2=\Delta_T$ in the even/odd sub-bundle, there is another in the even/odd sub-bundle if it is in the third partition.
However, these can be shown to come in pairs for the associated first-order operator $A_2$ where the eigenvalues differ by a sign. This can be used to show how regularizing the eta invariants and $\xi_T$ invariants defined in terms of the spectrum of the operator $A_2$ on the link will still yield the same cancellations as those for the supertraces on the Dirac cohomology.

We give a slightly different description of the same phenomenon, in terms of measures of spectral asymmetry. When the eta invariant of a boundary/link vanishes, it indicates that there is a symmetry of the eigensections about $0$, and for each positive eigenvalue, there is a negative eigenvalue, leading to a cancellation. This is called \textbf{spectral symmetry}, and the eta invariant (which is essentially half of the supertrace over the eigensections of $A$ on the link which are not in the kernel of $A$) is a measure of spectral asymmetry. It is clear from the discussion above (in particular relevant to the comparisons of Figure \ref{figure_Cheeger_Sections_odd_even_delta_Z} and Figure \ref{figure_Cheeger_Sections_odd_even_delta_T}) that the second SUSY in the case of twisted Dolbeault complexes does not correspond to a spectral symmetry of $cl(dx)A$, but a spectral symmetry of $cl(dx)A_2$, and the cancellations are only away from the kernel of $cl(dx)A_2$. Thus the supertrace over the third partition would be the measure of spectral asymmetry of the operator $cl(dx)A_2$, which of course vanishes.
Thus the second SUSY corresponds to the spectral symmetry of the operator $cl(dx)A_2$ in general, where we note that in the de Rham setting, $A_2=A$. This also explains the discrepancies in the eigenvalues for the spin Dirac operator on the sphere (and for orbifolds) for third partition (of what we call the Dirac cohomology), computed in \cite[\S 3.3]{Degeratuthesis} (see Corollary 3.4.2 of \textit{loc. cit.}), where representation theoretic decompositions are utilized. 

Yet another way of viewing this is to consider the correspondence of the eigensections of the operators $cl(dx)A_2, cl(dx)A$ on the boundaries/links $Z$ of the cylinders/cones, and see there is a one to one correspondence with the Dirac cohomology. Here the spectral symmetry of operators on the links are in one to one correspondence with secondary SUSY of the operators on the cylinders/cones. Indeed the reason one can restrict to either the even or the odd degree sections to compute the $\xi_T$ invariant using spectral asymmetry is because there is at least one supersymmetry. The second spectral symmetry that leads to more cancellations in $\xi_T$ invariants is due to the second supersymmetry. 
Since the proofs expressing the index contributions in terms of $\xi_T$ invariants consume one SUSY in the McKean singer cancellations, this is the second/residual SUSY that now cancels out the third partition. 
This point of view is quite useful in understanding the role of boundary conditions and ideal boundary conditions for $\mathcal{N}=2$ SUSY complexes as well.

\textit{We now enter into a fairly high level discussion on insights on picking boundary conditions, especially with spectral cuts, which is by no means self contained.}

The Green-Stokes 
pairing arising from integration by parts for formally symmetric operators $D$ is given by
\begin{equation}
    \omega_D(u,v)=\langle Du, v\rangle_{L^2(M;\mathcal{E})}-\langle u, Dv\rangle_{L^2(M;\mathcal{E})}=\int_{x=0} \langle \sigma(D)(dx)u,v\rangle_\mathcal{E} x^{al}dvol(g_Z) 
\end{equation}
where $x$ is a boundary defining functions for the boundary (resolved boundary if singular) so that $\{x=0\}=Z$ is the boundary. Here $l$ is the dimension of $Z$.
Notice that the inner products $\langle \sigma(D)(dx)u,v\rangle_\mathcal{E}$ vanish unless both $u,v$ have the same parity (are both even/odd), and that $\sigma(D)(dx)$ is a bundle isomorphism that switches the parity. 

One way to pick self-adjoint extensions of the operator $D$ is to carry out the following procedure (where we refer the reader to \cite{Albin_2017_index,bar2022boundary} for more rigorous treatments). Consider a closed subspace $S_1 \subseteq \Gamma(C_x(Z);\mathcal{E}^{even})$, where $\Gamma(C_x(Z);\mathcal{E}^{even})$ is a space of sections whose graph closure includes the minimal domain of $D$, but has more sections which are determined by a spectral condition, which we \textit{discuss separately for the boundary and conic settings}.

Consider the \textit{boundary setting}: the sections $s \in S_1$, when restricted to $x=0$, should be spanned by eigensections $\{\phi_j\}_{j \in L}$ of the induced operator $A$ on the boundary/link, where $\{\phi_j\}_{j \in J}$ are an orthonormal basis of eigensections for the operator $A$ on the even subbundle $\mathcal{E}^{even}$, and $L \subseteq J$. Then take $cl(dx)S_1 \subseteq \Gamma(C_x(Z);\mathcal{E}^{odd})$ and denote the complement of this closed subspace as $S_2'$. Taking the direct sum of $S_2'$ with the minimal domain for $D$ on the odd bundle, and taking the graph closure (with respect to the graph norm of $D$) we get the domain on the odd bundle which we denote $\mathcal{D}_L(D)$. The domain on the even bundle $\mathcal{D}_L(D)$ is given by the graph closure of $S_1$.

One can check that $\{cl(dx)\phi_j\}_{j \in J}$ are an orthonormal basis of eigensections for the odd subundle. Then one can check that the restriction of the odd degree sections of $\mathcal{D}_L(D)$ to $x=0$ corresponds to the span of the eigensections $\{cl(dx)\phi_j\}_{j \in L^{\perp}}$, where $L^{\perp}$ corresponds to the complement of $L \subseteq J$. These are sometimes called spectral cuts, when $L$ and $L^{\perp}$ are taken to be the indices corresponding to eigenvalues in disjoint semi-infinite intervals.

It is easy to see how the boundary pairing for sections which restrict to a span of eigensections at the boundary at $x=0$ vanishes with such choices of domains, since the basis of eigensections are orthonormal, and we have picked orthogonal complements on the even and odd sub-bundles. For general sections, the result follows from the density of these cores and functional analysis, similar to the results in Section \ref{sec:domains}.

Now we consider the \textit{conic setting}, where the operator $D=\sigma(\nabla_{\partial_x}-\frac{1}{x}A)$ has harmonic sections on the model cone $x^{\mu_j}\phi_j$, where $\mu_j$ are the eigenvalues of $A$ for $\phi_j$. While these are locally in $L^2$ near $x=0$ for all positive $\mu_j$, small eigenvalues of $\mu_j$ also lead to integrable harmonic sections due to the factor $x^l$. Together with the rescaling of the bundle of sections needed to put the Dirac operator in the form $D=\sigma(\nabla_{\partial_x}-\frac{1}{x}A)$, one can get a strict lower bound for the eigenvalues that may give rise to integrable harmonic sections.
There are other non-zero eigenvalue eigen sections on the cone that may also be in the Gelfand-Robbin complex, but have typically been ignored when studying index theory, since it is the harmonic sections that feature in the formulas for index theory. 

In the study of Mezzo-perversities in \cite{Albin_hodge_theory_cheeger_spaces,Albin_hodge_theory_on_stratified_spaces}, by rescaling the links of the metrics so that the eigenvalues $\mu_j$ grew, it was arranged that the choices of self-adjoint extensions at the level of the cohomology were only coming from the $0$ eigenvalues of $A$, and this is known as the \textit{spectral Witt condition}. Other methods such as that in \cite{Albin_2017_index} are closer to the case of the APS domain for cylinders with boundary. There, assuming that $A$ has no harmonic sections on the sections in $Z$, the eigenvalues corresponding to those that give rise to 
$L^2$ integrable harmonic sections are cut in 2, using a spectral cut similar to that in the case of the APS condition on cylinders. We refer the reader to \textit{loc. cit.} for more details. In particular, this domain is a rolled up extension of the Dolbeault-Dirac operator (this can be inferred from the duality properties for the domains introduced in \cite{jayasinghe2024holomorphic}).

Thus, in \textit{both the cylindrical and conic settings, the choices of domains have been picked using spectral cuts with respect to the operator $A$}. When there is a second supersymmetry, as we have studied in this article, it suffices to look at the kernels of the normal operators $A_2$ (corresponding to the model metrics on the boundary fibrations over the strata). Indeed, in the de Rham setting, since $A_2=A$, this corresponds to the method in \cite{Albin_hodge_theory_cheeger_spaces,Albin_hodge_theory_on_stratified_spaces}. However in the Dolbeault case, when $A_2 \neq A$, even after looking at the kernel of the operator $A_2$, one has to consider spectral cuts for the operator $A_1$, and some examples in \cite{jayasinghe2026spectralasymmetrysupersymmetryequivariant} already show how this can be done for rolled up extensions. In the case of a 1-dimensional link, $A=A_1$, and the problem simplifies to taking spectral cuts corresponding to Fourier eigenmodes on a circle. Now there is an analogous problem even on cones with K\"ahler metrics (equivalently when the links are Sasaki), where the problem boils down to picking eigensections of the operator $A_1$, basically an operator generated with derivatives with respect to the Reeb-vector field, up to Clifford multiplication. 

This begs the question, \textit{how does such a procedure compare with that of the procedure of splitting the eigenvalues of $A$, say in \cite{Albin_2017_index} using the VAPS domain?}
The answer is again given by SUSY. When restricted to the kernel of $A_2$, the operator $A$ is exactly equal to $A_1$, and the choices away from the third partition are the same, while the choices made for the eigensections in the third partition do not change the index formulas (as long as they are choices corresponding to rolled up extensions) due to SUSY cancellations.

A similar statement holds for spectral cuts for cylindrical problems, and even beyond that for horn metrics and other metrics with some form of radial symmetry (see for instance \cite{jayasinghe2025supersymmetricharmonicoscillatorssingular,Jesus2018Wittens,Jesus2018Wittensgeneral}. These preliminary investigations show that SUSY has a lot to contribute when studying even more general singular metrics.

\subsection{Comparison of the formulas with Lefschetz formulas in the smooth setting}
\label{Subsection_comparison_Lefschetz_smooth}

Since it is well known that the (holomorphic) Morse inequalities reduce to the (holomorphic) Lefschetz fixed point theorem for $b=-1$, it is instructive to compare our formulas with other known local index formulas. 
Even in the case of isolated critical points, objects used to describe the local Lefschetz numbers as equivariant eta invariants, linear algebraic formulas (for orbifolds) are different from the objects describing instanton complexes that categorify the (holomorphic) Morse inequalities. 
This was investigated in \cite{jayasinghe2026spectralasymmetrysupersymmetryequivariant} where a companion paper is being compiled where the relation to $\mathcal{N}=2$ SUSY is explained. 
A quick description of the results there, in the complex setting, is that while the eta invariants are alternating traces over the spectrum of an associated Dirac-type operator on the link at a singularity, not all the eigensections contribute to the trace. There are additional supersymmetric cancellations in the trace (corresponding to the second supersymmetry), and only the sections corresponding to local cohomology groups on the cones contribute non-trivially to the trace.

The special formulas for orbifold equivariant Lefschetz numbers \cite{vergne1994quantification} can be reformulated as supertraces over the local cohomology group of a smooth local cover at the fixed points, up to extra factors corresponding to the covering degree. However this point of view is not intrinsic to the singular space, and does readily relate to descriptions of other singular analytic objects arising in the analytic descriptions of equivariant index defects (such as eta and xi invariants \cite{jayasinghe2026spectralasymmetrysupersymmetryequivariant}) or Poincar\'e residues arising in algebraic versions.

In Remark 2.3 of \cite{wu1998equivariant}, it is observed that the local holomorphic Lefschetz numbers over a smooth fixed set $F_a$ are given by
\begin{equation}
\label{equation_smooth_equivariant_index-2}
    \sum_{q=0}^n (-1)^q \operatorname{Tr} T|_{\mathcal{H}^q(U(F_a);E)}=\int_{F_a} ch_T(K^C(N_r) \otimes E|_{F_a}) td(F_a)=\int_{F_a} ch_T \Bigg(\frac{E|_{F_a}}{det(1-N^{*(1,0)}_a)} \Bigg)td(F_a)
\end{equation}
and we refer to \textit{loc. cit.} for more details, where we have rephrased the cohomology group in the notation and formulations of our article. In the case of a Spin Dirac operator, where $E|_{F_a}$ would correspond to a square-root bundle $L$ of the canonical bundle (i.e. $L^{\otimes 2} =\mathcal{K}$), it is well known that the product of the Chern class of $L|_{F_a}$ with the Todd class $td(F_a)$ gives the $\widehat{A}$ class on $F_a$.

Now consider the index formula in \cite{Albin_2017_index}, where the contribution from a singular stratum $N$ of any depth (where the local fibration is as in Definition \ref{iterated_fibration_structure}) is given by the integral 
\begin{equation}
\label{singular_index_contributions}
    \int_{Y} \widehat{A}(Y)\mathcal{J}(\mathcal{B}_Y/Y)
\end{equation}
where we refer the reader to \textit{loc. cit.} for the definition of the Bismut-Cheeger $\mathcal{J}$ form over $Y$. Comparing with equation \ref{equation_smooth_equivariant_index-2} with a \textbf{conjectural singular equivariant formula} generalizing equation \eqref{singular_index_contributions} (at least at a strata of maximal depth), it is easy to see that the equivariant Chern class $ch_T(K^C(N_r)|_{F_a})$ admits a description as an equivariant Bismut-Cheeger $\mathcal{J}$ form for the normal disc bundle over $Y$, viewed as a boundary fibration. For the case of isolated fixed points this follows from the results in \cite{jayasinghe2023l2}, and from that it is easy to extend to the case of flat normal bundles in the sense of \cite{jayasinghe2024witten} as well. 

In an upcoming revision of \cite{jayasinghe2026spectralasymmetrysupersymmetryequivariant}, we show that just as the equivariant $\xi$ invariant which records the equivariant index defect at isolated fixed points corresponds to a supertrace over the local cohomology group, as appearing in the left hand side of \ref{equation_smooth_equivariant_index-2}, the non-equivariant $\xi$ invariant at an isolated singularity  can be phrased entirely in terms of local cohomology groups.

\textbf{We conjecture that this perspective extends even to non-isolated singularities, giving a reformulation of the index theorem of \cite{Albin_2017_index} in suitable settings (see Section \ref{Appendix_SUSY_and_conjectures}).}

We mention that there are links to other formulations of equivariant Riemann-Roch theorems, particularly the work of Edidin-Graham in \cite{equivariant_intersection_edidin,localization_algebraic_Graham_Edidin} (see Remark 7.32 of \cite{jayasinghe2023l2}) that generalize the work of Baum-Fulton-Macpherson-Quart \cite{baum1979lefschetz,baumfultonmacKtheoryRiemannRoch,baumfultonmacpheresonRiemannRoch,Baumformula81}, which were discussed in the precursor work of \cite{jayasinghe2024holomorphic}, in particular formulating conjectural Morse inequalities for these algebraic equivariant Riemann-Roch results (see Section 6.2 of \cite{jayasinghe2024holomorphic}).

In \cite{wu1998equivariant}, Wu constructs examples which show that the holomorphic Morse inequalities fail for certain toric actions on complex manifolds which are not K\"ahler Hamiltonian actions. It is also identified that compact, connected Lie group actions are the correct general setting for the result, where in \cite{wu2003instanton} he shows that a meromorphic group action on a complex manifold with a filterable Bialynicki-Birula decomposition (weaker than the K\"ahler assumption in this paper) suffices.
In section 2 of \cite{wu1998equivariant}, the author shows how more information can be drawn in the case of a Lie group action. This was already used to prove vanishing results for harmonic spinors in certain algebraic settings (see Section 7.3 of \cite{jayasinghe2024holomorphic}), but we extend this formulation more precisely in Subsection \ref{Subsection_Toric_actions} to our geometric setting.

In Section 5 of \cite{wu1996equivariant}, Wu also studies the case of non-Abelian group actions including relationships to the Weyl character formula. These are beautiful results in the smooth setting, but potential generalizations are more complicated in our setting as can be checked in examples and a more careful study is not in the scope of this article.

\subsection{Conjectures on index theory}
\label{subsection_appendix_index_theory}

Now we put the work in our article in a broader context. 
Recall that one framing of the equivariant index theorem on smooth manifolds is that there is a functor
from the category of equivariant K-theory to $R(G)$, which for the index is a functor from K-theory to the integers.
This factors through cohomology, via the Chern character, which for rational cohomology is invertible. Thus it suffices to prove localization at the level of K-theory.

In the context of Morse theory and Witten complexes, we prove localization at the level of Witten subcomplexes, in the category of $\mathcal{N}=2$ SUSY Hilbert complexes. We can pass to the local cohomology groups which can also be thought of as infinite rank bundles, of finite rank after passing to the equivariant complexes to representatives in K-theory with finite rank bundles.
Of course one only has a $\mathbb{Z}_2$ grading in K-theory.

Consider a resolved stratified pseudomanifold $X$ with a wedge K\"ahler metric, 
and a holomorphic bundle $E$ on which there is a Hamiltonian action of a compact connected Lie group $G$, where the connected components $F$ of the fixed point sets of $g \in G$ each have suitable localizing neighbourhoods $U(F)$. 
Then we have shown that there are (regularized/formal) supertrace functors
\begin{equation}
\label{functor_1}
    \text{Str}_s\, T_g : \mathcal{H}^k(\mathcal{P}_{\mathcal{B}}(\mathcal{N}_F)) \longrightarrow R(G)
\end{equation}
where $\mathcal{H}^k(\mathcal{P}_{\mathcal{B}}(\mathcal{N}_F))$ is the cohomology of the normal cohomology complex of the twisted Dolbeault complex, and these give the twisted equivariant Riemann-Roch numbers. Here $s$ is a parameter that comes up in the regularization, something well known in the study of equivariant eta and $\xi_T$ invariants (see \cite{jayasinghe2025supersymmetricharmonicoscillatorssingular}).
There are obvious lifts of this functor to the category of SUSY Hilbert complexes, but we shall study those at the level of cohomology to match the correspondences we pointed out in the earlier part of this Appendix.

Here there are two canonical choices of boundary conditions $B=N,D$ for $\mathcal{B}=(W,B)$, for any choice of ideal boundary condition $W$ on the space $X$, which give the same (regularized) supertrace. Thus we can form a symmetrized version of this functor
\begin{equation}
\label{functor_2}
   (1/2)\text{Str}_s\, T_g : \mathcal{H}^k(\mathcal{P}_{W,N}(\mathcal{N}_F)) \oplus \mathcal{H}^k(\mathcal{P}_{W,D}(\mathcal{N}_F)) \longrightarrow R(G)
\end{equation}
which gives the same series in the representation ring.

The results in \cite{jayasinghe2026spectralasymmetrysupersymmetryequivariant} show that in the case of isolated fixed points at isolated conic singularities, these correspond to the $\xi_T$ invariants- that is-there is a regularized supertrace functor (see equation \eqref{equation_xi_T}) from Dirac cohomology to $R(G)$
\begin{equation}
    (1/2)\text{Str}_s T_g : \mathcal{H}^k(\mathcal{D}_{\mathcal{B}}(D)) \longrightarrow R(G)
\end{equation}
which is equal to the functors \eqref{functor_1}, \eqref{functor_2}. In an upcoming revision of \cite{jayasinghe2026spectralasymmetrysupersymmetryequivariant}, we shall give the details for these functors when $g$ is the identity element, that is the non-equivariant index formulae, and for the case of compact but non-connected Lie groups $G$ (including discrete actions).

Thus we make the following conjecture.
\begin{conjecture}
\label{Conjecture_1}
Given a $\mathcal{N}=1$ complex $\mathcal{P}_{W}(X)=(L^2(X;\mathcal{E}),\mathcal{D}_W(D),D)$ on a resolved stratified space $X$ with a wedge metric, with an action of a compact Lie group $G$, where there are localizing $U(F)$ neighbourhoods for the connected components $F$ of the fixed point set, we conjecture that there are Dirac cohomology groups $\mathcal{D}_{\mathcal{B}}(D_{U(F)})$ associated to the localizing neighbourhoods, and functors
\begin{equation}
    (1/2)\emph{Str}_s\, T_g : \mathcal{H}^k(\mathcal{D}_{\mathcal{B}}(D_{U(F)})) \longrightarrow R(G)
\end{equation}
that correspond to the equivariant index functor.
Moreover when there are local $\mathcal{N}=2$ refinements these functors are equal to 
\begin{equation}
\label{functor_2-new}
   (1/2)\emph{Str}_s\, T_g : \mathcal{H}^k(\mathcal{P}_{W,N}(\mathcal{N}_F)) \oplus \mathcal{H}^k(\mathcal{P}_{W,D}(\mathcal{N}_F)) \longrightarrow R(G).
\end{equation}
\end{conjecture}

The above conjecture seems quite plausible since we have a proof for the case of isolated fixed points on cone points, and for non-isolated fixed points where the metric is a product of the fixed point set and the normal fibers (where one can do separation of variables to reduce the result to the case of isolated fixed points). This will be presented in an upcoming revision of \cite{jayasinghe2026spectralasymmetrysupersymmetryequivariant}.

We separate this next question from the earlier, though they are related.

\begin{question/conjecture}
\label{Conjecture_brave_2}
Given a $\mathcal{N}=1$ complex $\mathcal{P}_{W}(X)=(L^2(X;\mathcal{E}),\mathcal{D}_W(D),D)$ on a resolved stratified space $X$ with a wedge metric. When is the contribution to the index of the Dirac operator at singular stratum $S$ given by functor 
\begin{equation}
    (1/2)\emph{Str}_s\, T_{\emph{Id}} : \mathcal{H}^k(\mathcal{D}_{\mathcal{B}}(D_{S})) \longrightarrow R(G)
\end{equation}
for a suitable cohomology group $\mathcal{D}_{\mathcal{B}}(D_{S})$?

If there are $\mathcal{N}=2$ refinements, are there SUSY cancellations of such functors, or can we formulate such functors at least when the complex has $\mathcal{N}=2$ SUSY?
\end{question/conjecture}

Part of the question here is to come up with proper definitions for the cohomology groups. There are some examples where it seems one can formulate such groups, and the fact that we have proofs of this for isolated singular strata and those with product metrics makes it a compelling question. The discussion in Subsections \ref{Subsection_comparison_Lefschetz_smooth} and \ref{Subsec_comparison_Lefschetz_smooth} where the relationships to the work in \cite{Albin_2017_index} were discussed are highly relevant, since this conjecture would reformulate the results proven in terms of differential forms, in terms of local cohomology.
What it clear is that there is a study to be done to relate eta forms and Bismut-Cheeger $\mathcal{J}$ forms on stratified spaces (both equivariant and non-equivariant) to local cohomology groups.

\subsection{Some questions related to physics}

Index theorems are used in various contexts in physics already, including settings with manifolds with boundary and singularities. There are many instances where mathematics and physics have complemented each other fruitfully and we encourage the reader to discover this with modern tools that are far superior to anything we write here. 

In questions related to gravitational anomalies \cite{witten1985global}, Witten explained how $\xi$ invariants (also called reduced eta invariants) correspond to the holonomy of determinant bundles, which lead to several interesting investigations on smooth manifolds \cite{dai1994eta,bismut1986analysis,bismut1986analysisII} and singular spaces \cite{bismut1989eta,cheeger1987eta}. In each of these contexts, when the corresponding space has complex structures, the content described previously in this appendix applies.
Equivariant index theorems are used in SUSY localization more broadly where the equivariant indices of operators correspond to various partition functions such as the Nekrasov partition functions (we refer to \cite[\S 7.2.8]{jayasinghe2023l2} for a literature review as well as studies in the singular setting), and gravitational entropy functional studied by Gibbons-Hawking and others (we refer to \cite{jayasinghe2024holomorphic} and Subsection \ref{subsection_bolt_charges}). Since partition functions are counting the relevant states in the supersymmetric theory, we are tempted to ask some natural questions.

\begin{question}
When using equivariant index formulas for orbifolds and smooth manifolds which implicitly use $\mathcal{N}=2$ SUSY cancellations to simplify the equivariant eta invariants on spheres and other links on orbifolds (see \cite{jayasinghe2026spectralasymmetrysupersymmetryequivariant,Degeratuthesis})
is the relevant count to be understood to be the superstates of the theory before the extra supersymmetric cancellations, which at the particular fixed point sets cancel off due to the $\mathcal{N}=2$ SUSY refinement at the fixed points?

What are the qualitative aspects of local contributions to partition functions that can advance ones understanding of physics?
\end{question}

Moreover there are studies in the physics literature (see, e.g., \cite{cvetivc2025extra} for studies of anomaly inflows, cf.  \cite{witten2022anomaly}), where the view is expressed that the results they study in the case of orbifolds where the eta invariants are computed using the results in \cite{Degeratuthesis} extend to more general singularities.
Since these extra SUSY cancellations in the eta invariant are not taken into consideration in these studies, we ask the same questions there, but we also ask \textbf{do the choices of domains, which play important roles when studying these invariants on singular spaces, have a role to play in the study of physics?}

\bibliographystyle{alpha}
\bibliography{references}
\end{document}